\documentclass[11pt]{amsart}

\usepackage[T1]{fontenc}
\usepackage[english]{babel}
\usepackage[utf8]{inputenc}

\usepackage{amsmath,amsthm,amssymb,amsfonts}
\usepackage{enumerate,amscd,amsxtra,MnSymbol}
\usepackage{mathrsfs}
\usepackage{bbm}

\usepackage{geometry}
\usepackage{color}
\usepackage[cmtip,all]{xy}
\usepackage{tikz} 
\usepackage{tikz-cd} 
\usepackage{tikz-qtree} 
\usepackage{forest}
\usetikzlibrary{arrows,calc}

\usepackage{graphicx}

\usepackage{xparse}

\usepackage[normalem]{ulem}

\usepackage{hyperref}
\usepackage[capitalise]{cleveref} 

\tikzstyle{d}=[double distance=.3ex]
\tikzstyle{w}=[preaction={draw=white, -,line width=4pt}]
\NewDocumentCommand{\punctuation}{ m m O{5pt} }{\node at ($(#1.east)-(0,#3)$) {#2};}

\newcounter{diagram}
\renewcommand{\thediagram}{\thetheorem}

\tikzset{%
node distance=1.5cm, la/.style={scale=0.8}, rr/.style={xshift=1.5cm},
space/.style={xshift=.5cm}, over/.style={auto=false,fill=white,inner sep=1.5pt, minimum size=0, outer sep=0},
    symbol/.style={%
        draw=none,
        every to/.append style={%
            edge node={node [sloped, allow upside down, auto=false]{$#1$}}},
            
    }, pro/.style={postaction={decorate,decoration={
        markings,
        mark=at position .5 with {\node at (0,0) {$\bullet$};}
      }},
      inner sep=.9ex,
      },
  n/.style={double equal sign distance, -implies}, t/.style={double distance=2.5pt, -implies, postaction={draw,-}},
}

\usepackage{mathtools}

\theoremstyle{definition}

\newtheorem*{defn}{Definition} 

\newtheorem{theorem}{Theorem}[section]
\newtheorem*{theorem*}{Theorem}
\newtheorem{lemma}[theorem]{Lemma}
\newtheorem{proposition}[theorem]{Proposition}
\newtheorem*{proposition*}{Proposition}
\newtheorem{corollary}[theorem]{Corollary}
\newtheorem*{corollary*}{Corollary}

\newtheorem{construction}[theorem]{Construction}
\newtheorem{example}[theorem]{Example}

\newtheorem{introtheorem}{Theorem}

\newtheorem{definition*}{Definition}

\newtheorem{conjecture*}{Conjecture}

\newtheorem{notation}[theorem]{Notation}

\newtheorem{remark}[theorem]{Remark}
\newtheorem{remark*}{Remark}
\newtheorem{Terminology}[theorem]{Terminology}

\newtheorem{definition}[theorem]{Definition}

\newcommand{\caR}{{\mathcal R}}
\newcommand{\caC}{{\mathcal C}}
\newcommand{\caD}{{\mathcal D}}
\newcommand{\caB}{{\mathcal B}}
\newcommand{\caA}{{\mathcal A}}
\newcommand{\caE}{{\mathcal E}}

\newcommand{\caS}{{\mathcal S}}
\newcommand{\caO}{{\mathcal O}}
\newcommand{\caH}{{\mathcal H}}
\newcommand{\caZ}{{\mathcal Z}}

\newcommand{\caP}{{\mathcal P}}
\newcommand{\caK}{{\mathcal K}}
\newcommand{\caQ}{{\mathcal Q}}
\newcommand{\caL}{{\mathcal L}}
\newcommand{\caW}{{\mathcal W}}
\newcommand{\caF}{{\mathcal F}}
\newcommand{\caV}{{\mathcal V}}

\newcommand{\D}{{\mathsf D}}
\renewcommand{\L}{{\mathbb L}}

\newcommand{\fib}{\mathsf{fib}}
\newcommand{\uHom}{\underline{\mathsf{Hom}}}

\newcommand{\op}{\mathsf{op}}
\newcommand{\dual}{\vee}

\newcommand{\Hom}{\mathsf{Hom}}

\newcommand{\s}{\mathsf{s}}

\newcommand{\C}{\mathsf{C}}

\newcommand{\pr}{\mathsf{pr}}

\newcommand{\Sp}{\mathsf{Sp}}

\newcommand{\f}{\mathsf{f}}

\newcommand{\colim}{\mathrm{colim}}

\newcommand{\Mod}{{\mathsf{Mod}}}
\newcommand{\LMod}{{\mathsf{LMod}}}

\newcommand{\Alg}{{\mathsf{Alg}}}

\newcommand{\X}{\mathrm{X}}
\newcommand{\Y}{\mathrm{Y}}
\newcommand{\Pre}{\mathsf{Pre}}

\renewcommand{\H}{\mathcal{H}}

\newcommand{\id}{\mathsf{id}}

\newcommand{\Cat}{\mathsf{Cat}}

\newcommand{\Set}{\mathsf{Set}}
\newcommand{\sSet}{\mathsf{sSet}}
\newcommand{\rsSet}{\mathsf{rsSet}}

\newcommand{\map}{\mathrm{map}}

\renewcommand{\Pr}{\mathsf{Pr}}

\newcommand{\Spc}{\mathsf{Spc}}
\newcommand{\gdpreadd}{\mathsf{Pre}\mathsf{Add}^\mathsf{gd}}
\newcommand{\gdadd}{\mathsf{Add}^\mathsf{gd}}
\newcommand{\Add}{\mathsf{Add}}
\newcommand{\St}{\mathsf{St}}
\newcommand{\Mon}{\mathsf{Mon}}

\newcommand{\Fun}{\mathsf{Fun}}

\newcommand{\Cart}{\mathsf{Cart}}
\newcommand{\Grp}{\mathsf{Grp}}

\newcommand{\add}{\mathsf{add}}

\newcommand{\Fin}{\mathsf{Fin}}

\newcommand{\lax}{\mathsf{lax}}
\newcommand{\SegSpc}{\mathsf{SegSpc}}
\newcommand{\rSegSpc}{\mathsf{rSegSpc}}
\newcommand{\crSegSpc}{\mathsf{crSegSpc}}

\newcommand{\Ar}{\mathsf{Ar}}

\newcommand{\Tw}{\mathsf{Tw}}

\newcommand{\KR}{\mathsf{KR}}
\newcommand{\KH}{\mathsf{KH}}
\newcommand{\Wald}{\mathsf{Wald}}

\newcommand{\fiss}{\mathsf{fiss}}
\newcommand{\gd}{\mathsf{gd}}
\newcommand{\Exact}{\mathsf{Exact}}
\renewcommand{\d}{\mathsf{d}}
\newcommand{\Seg}{\mathsf{Seg}}

\newcommand{\Exc}{\mathsf{Exc}}

\newcommand{\rsSpc}{\mathsf{rsSpc}}

\newcommand{\brSegSpc}{\mathsf{brSegSpc}}
\newcommand{\ev}{\mathsf{ev}}
\newcommand{\rev}{\mathsf{rev}}
\newcommand{\Cmon}{\mathsf{CMon}}

\newcommand{\Calg}{\mathsf{CAlg}}
\newcommand{\exc}{\mathsf{exc}}
\newcommand{\rs}{\mathsf{rs}}

\newcommand{\cof}{\mathsf{cof}}
\newcommand{\PreCat}{\mathsf{PreCat}}
\newcommand{\Ass}{\mathsf{Ass}}
\newcommand{\Op}{\mathsf{Op}}
\newcommand{\caJ}{\mathcal{J}}

\usepackage[textwidth=3cm,
textsize=small,
colorinlistoftodos]
{todonotes} 

\title{Real $K$-theory for Waldhausen $\infty$-categories with genuine duality}

\author[H.\ Heine]{Hadrian Heine}
\address{University of Oslo, Oslo, Norway}
\email{hadriah@math.uio.no}

\author[M.\ Spitzweck]{Markus Spitzweck}
\address{Fakult\"at f\"ur Mathematik, Universit\"at Osnabr\"uck, Osnabr\"uck, Germany}
\email{markus.spitzweck@uni-osnabrueck.de}

\author[P.\ Verdugo]{Paula Verdugo}
\address{Max Planck Institute f\"ur Mathematik, Bonn, Germany}
\email{verdugo@mpim-bonn.mpg.de}

\begin{document}

\begin{abstract}
    We develop a new framework to study real $K$-theory in the context of $\infty$-categories. For this, we introduce Waldhausen $\infty$-categories with genuine duality, which will be the input for such $K$-theory. These are Waldhausen $\infty$-categories in the sense of Barwick equipped with a compatible duality and a refinement of their hermitian objects generalizing the concept of Poincar\'e $\infty$-categories of Lurie. They may also be thought of as a version of complete Segal spaces enriched in genuine $C_2$-spaces whose underlying $\infty$-category carries a compatible Waldhausen structure, since we show that their respective $\infty$-categories are equivalent. 
    We define the real $K$-theory genuine $C_2$-spaces by means of an enriched version of the  $S_\bullet$-construction, defined for Waldhausen $\infty$-categories with genuine duality. Moreover, we prove an Additivity Theorem for this $S_\bullet$-construction which leads to an Additivity Theorem for real $K$-theory. Furthermore, such real $K$-theory satisfy a universal property \textemdash analogous to that proved by Barwick for algebraic $K$-theory of Waldhausen $\infty$-categories\textemdash: we prove that every theory can be universally turned into an additive theory and identify our real K-theory with the universal additive theory associated to the functor that associates to a Waldhausen $\infty$-category with genuine duality its maximal subspace.
    
    Finally, we promote the real $K$-theory genuine $C_2$-spaces to genuine $C_2$-spectra.
\end{abstract}

\maketitle

\tableofcontents

\section{Introduction}

Algebraic $K$-theory serves as a central organizing principle linking geometric and algebraic data, and its influence can be found in a variety of fields, from number theory and algebraic geometry to manifold topology and mathematical physics. Its most basic component, the group $K_0$, was introduced by Grothendieck in \cite{grothendieck1968classes} motivated by studying a generalization of Riemann-Roch theorem. It has since proven a powerful invariant that, at its core, captures relevant information of rings by means of their category of modules.
Higher $K$-groups were introduced by Quillen \cite{quillen.higher}, and encode rich information about rings, schemes, and categories. 

For much of its history, the emphasis in the study of algebraic $K$-theory was computational: determining $K$-groups for concrete rings, schemes, or categories. In the last decade, with the advent of higher category theory, a new conceptual perspective has materialized. This perspective consists in looking at algebraic $K$-theory as a \emph{universal functor}, and not only as a its specific values, defined and studied through its formal and homotopical properties; notably with the works of Blumberg-Gepner-Tabuada \cite{bgt} and Barwick \cite{barwick.wald}. This shift of focus from computation to universality sits in a broader transformation in modern homotopy theory: the formalization that many classical constructions arise as initial or terminal objects satisfying certain excision or additivity properties (e.g.\ topological cyclic homology \cite{harpaz2024trace}, laced $K$-theory and topological Hochschild homology \cite{nikolaus2018topological}).

The study of \emph{hermitian} and \emph{real algebraic $K$-theory}\textemdash that is, algebraic $K$-theory in the presence of dualities-has a long and intricate history. Karoubi and Villamayor \cite{karoubi1971k} \cite{karoubi1980le}
developed hermitian $K$-theory for rings with involutions, which was later extended to exact categories with dualities. Later, it was shown by Schlichting to be a well-behaved invariant of schemes in \cite{schlichting.mv}.
Real algebraic $K$-theory developed by \cite{hesselholt-madsen}
unifies algebraic $K$-theory and hermitian $K$-theory which appear, respectively, as the underlying space and $C_2$-fixed points of a genuine $C_2$-space or spectrum.

In recent years, real and hermitian algebraic $K$-theory has been going through the same process of using $\infty$-categories to obtain both a more conceptual presentation and expand applications. One of the approaches to this, that has proven computationally very powerful, was developed by Calmès-Dotto-Harpaz-Hebestreit-Moi-Land-Nardin-Nikolaus-Steimle in \cite{Calmes_etal1,Calmes_etal2,Calmes_etal3} which parallels the work in \cite{bgt}, another approach is the one presented in this work which parallels the work in \cite{barwick.wald}. Before diving into the content of this work, we look at the non-real case for context.  

\subsection{Algebraic $K$-theory in the presence of higher categories}

A fundamental property of algebraic $K$-theory, that was long considered to be a core piece, is that it splits cofiber sequences, a property known as additivity. This was, in fact, the starting point for the pursuit of a universal characterization that was achieved independently by Blumberg--Gepner-Tabuada\textemdash in the context of stable $\infty$ categories \textemdash and Barwick \textemdash in the context of Waldhausen $\infty$-categories.

Blumberg--Gepner--Tabuada construct $K$-theory of stable $\infty$-categories and formulate additivity as a homological property of preserving short exact sequences.
As a main result the authors prove that algebraic $K$-theory universally extends to a corepresentable functor on a certain stable $\infty$-category. The authors call this stable $\infty$-category the $\infty$-category of non-commutative motives drawing an analogy to Grothendieck's philosophy of motives, the universal category to which homological invariants of algebraic varieties extend.

Barwick developed the \emph{analogue in the Waldhausen setting}.
He extends the concept of Waldhausen categories to the realm of $\infty$-categories. In this setting, Barwick proves that algebraic $K$-theory is the \emph{initial additive theory}, i.e.\ the initial functor splitting short exact sequences,
under the functor sending a Waldhausen $\infty$-category to its space of objects. 
As a key result Barwick constructs an analogue of Blumberg--Gepner--Tabuada's stable $\infty$-category of non-commutative motives, which he calls the fissile derived $\infty$-category of Waldhausen $\infty$-categorie. He proves that algebraic $K$-theory universally extends to an
excisive functor on this $\infty$-category, which Barwick insightful characterizes as the first Goodwillie derivative of the functor that takes the maximal subspace.

In this way, both Blumberg--Gepner--Tabuada and Barwick revealed $K$-theory as an essentially \emph{formal} construction—determined entirely by its universal property—rather than by specific models or presentations. 
This conceptual viewpoint will serve as the foundation for our real refinement.

\subsection{Real and Hermitian $K$-theory}

Similarly to the non-hermitian case, the technology of $\infty$-categories allows for a jump both, in formal statements, and in understanding. As said before, two modern frameworks for real/hermitian $K$-theory have been developed in parallel.

On one hand, Calmès-Dotto-Harpaz-Hebestreit-Moi-Land-Nardin-Nikolaus-Steimle develop in \cite{Calmes_etal1, Calmes_etal2, Calmes_etal3} a framework for hermitian $K$-theory of stable $\infty$-categories equipped with a non-degenerate \emph{quadratic functor}, a categorification of a quadratic form originally proposed by Lurie \cite{lurie.L4}. Their theory extends the stable approach of \cite{bgt} to the equivariant setting and provides a powerful link between algebraic and geometric topology: hermitian $K$-theory and $L$-theory arise respectively as the genuine and geometric fixed points of a real $K$-theory genuine $C_2$-spectrum. In recent work, Schlichting \cite{schlichting1} constructs hermitian $K$-theory for so called exact form categories with weak equivalences, and Schlichting-Marlowe \cite{schlichting2} prove an equivalence between this construction and the one of Calmès-Dotto-Harpaz-Hebestreit-Moi-Land-Nardin-Nikolaus-Steimle.

On the other hand, in this work we develop a framework for real algebraic $K$-theory that generalizes that of Barwick \cite{barwick.wald}, which is what we introduce in the rest of the introduction. We develop a general theory of \emph{Waldhausen $\infty$-categories with genuine duality}, and within this setting, we prove a \emph{real additivity theorem} and a \emph{universal characterization} of real algebraic $K$-theory as the initial real additive functor. A preliminary version of this manuscript was made publicly available on the arXiv in 2019. Since then, we have significantly improved its readability and presentation, and also have proved that the two approaches above coincide \cite{comparisonKth}.

We know start introducing the context of our framwork.

\subsubsection{Real $\infty$-categories and the real $S$-construction}

A guiding principle of this work is that \emph{real phenomena in algebraic $K$-theory are best understood through enrichment in genuine $C_2$-spaces}. 
The original constructions of algebraic $K$-theory, due to Quillen for exact categories and Waldhausen for categories with cofibrations, are rooted in simplicial methods that encode additivity and homotopy coherence. 
In the real setting, we adopt an analogous perspective: real algebraic $K$-theory is studied via \emph{real simplicial methods}, where simplicial structures are replaced by their $C_2$-equivariant analogues. 
For instance, our proof of the real additivity theorem relies crucially on the construction and careful analysis of highly intricate real simplicial homotopies.
Real simplicial methods find their natural expression in the language of \emph{real $\infty$-categories}—$\infty$-categories enriched in $\Spc^{C_2}$, the $\infty$-category of genuine $C_2$-spaces.

The philosophy of our approach is to \emph{enrich classical constructions from Quillen and Waldhausen in genuine $C_2$-spaces}.  
This endows all constructions with an intrinsic equivariant structure reflecting the presence of dualities.
Within this setting, the real algebraic $K$-theory functor arises as the universal real additive invariant of real $\infty$-categories.

For instance, we construct a real version of the $S_\bullet$-construction of Waldhausen  \cite{waldhausen1985algebraic}, which was extended to $\infty$-categories by Barwick \cite{barwick.wald}.
Waldhausen's $S_\bullet$-construction produces from a category with cofibrations a simplicial space. 
We extend Barwick's $S_\bullet$-construction to a real $S_\bullet$-construction that, in analogy, produces from a Waldhausen $\infty$-category equipped with suitable duality data, a real simplicial space.
This is a real functor $\Delta^\op \to \Spc^{C_2}$ to the real $\infty$-category of genuine $C_2$-spaces,
where we endow the simplex category $\Delta$ with a natural enrichment in discrete $C_2$-spaces arising from the natural involution on the simplex category.

\subsubsection{Waldhausen $\infty$-categories with genuine duality}

A central component of our work is the introduction and systematic study of \emph{Waldhausen $\infty$-categories with genuine duality}. 
These are Waldhausen $\infty$-categories in the sense of Barwick \cite{barwick.wald} endowed with a compatible duality and a \emph{genuine refinement} of the duality, 
playing the role of the space of genuine $C_2$-fixed points in equivariant homotopy theory.

We devote substantial effort to developing the notion of an $\infty$-category with genuine duality, which generalizes the concept of an $\infty$-category with duality introduced in \cite{hls}. 

To motivate this notion, it is instructive to recall the analogy with equivariant homotopy theory:
by Elmendorf's theorem \cite{Elmendorf}, a genuine $C_2$-space can be described by a $C_2$-space $X$ together with a map of spaces $$Y \to X^{hC_2}.$$
In this presentation the space $Y$ is the space of $C_2$-fixed points of a genuine $C_2$-space.
In direct analogy, an $\infty$-category with genuine duality consists of an $\infty$-category with duality $\caC$ equipped with a right fibration over the $\infty$-category of hermitian objects $\caH(\caC)$ (\cref{hermitian}), 
which serves as the categorical counterpart of the space of $C_2$-fixed points in the setting of dualities.

Our key conceptual framework is summarized in the following dictionary:

\begin{center}
\begin{tabular}{ | l | l | c | c}
\hline
Space & $\infty$-category \\
\hline
$C_2$-space & $\infty$-category with duality \\ \hline
Genuine $C_2$-space & $\infty$-category with genuine duality \\ \hline
Homotopy $C_2$-fixed points & Hermitian objects \\ \hline
$C_2$-fixed points & Genuine refinement (genuine hermitian objects) \\
\hline
\end{tabular}
\end{center}

\vspace{3mm}

We provide a concrete model for $\infty$-categories with genuine duality in terms of \emph{complete real Segal spaces},
which extends the classical Segal space model for $\infty$-categories to the equivariant setting.
A well-known model for $\infty$-categories are complete Segal spaces, simplicial spaces satisfying the Segal condition and a completeness condition ensuring the existence of compositions and that equivalences are well-behaved.
In direct analogy, we identify $\infty$-categories with genuine duality with \emph{complete real Segal spaces}:
these are real simplicial spaces satisfying $C_2$-equivariant analogues of both the Segal condition and the completeness condition:

\begin{introtheorem}[\cref{equiv_crSegSpc_GD}]
There is an equivalence of $\infty$-categories between that of  $\infty$-categories with genuine duality and that of complete real Segal spaces. 
\end{introtheorem}

The model for $\infty$-categories with genuine duality via complete real Segal spaces plays a central role in the conceptual understanding of the real $S$-construction. Indeed, Waldhausen’s $S$-construction associates to a category with cofibrations $\caC$ a simplicial category with cofibrations
$S_\bullet(\caC)^\simeq$,
whose simplicial space of equivalences $S_\bullet(\caC)^\simeq$ is obtained by discarding all non-invertible morphisms in each degree.
This simplicial space freely determines a complete Segal space whose classifying space identifies with the geometric realization of $S_\bullet(\caC)^\simeq$.
In particular, it provides a delooping of algebraic $K$-theory:
Waldhausen defines algebraic $K$-theory of a category with cofibrations $\caC$ as the loop space
$$\mathrm{K}(\caC)\coloneqq \Omega \vert {S_\bullet (\caC)}^\simeq \vert.$$

In analogy, we assign to any Waldhausen $\infty$-category with genuine duality $\caC$ a real $S_\bullet$-construction $S_\bullet(\caC) $, which is a real simplicial Waldhausen $\infty$-category with genuine duality.
By passing to equivalences in each simplicial degree, we obtain a real simplicial space $S_\bullet(\caC)^\simeq$,
which freely determines a complete real Segal space.
This complete real Segal space admits a classifying genuine $C_2$-space, which identifies with the real geometric realization of
$S_\bullet(\caC)^\simeq$ and thus provides a $C_2$-equivariant delooping of real algebraic $K$-theory, defined as 
$$\KR(\caC)\coloneqq \Omega^{1,1} \vert {S_\bullet (\caC)}^\simeq \vert.$$
Here loops are formed with respect to the sphere equipped with sign representation.

In other words real algebraic $K$-theory is the $C_2$-equivariant loop space of the 
classifying genuine $C_2$-space of the complete real Segal space freely generated by the real simplicial space arising from the real Waldhausen $S_\bullet$-construction.

Another important perspective on $\infty$-categories with genuine duality is the following, which establishes a bridge between our framework and that of \emph{Poincar\'e $\infty$-categories} of 
\cite{Calmes_etal1}, \cite{Calmes_etal2}, \cite{Calmes_etal3}. 

By work of \cite{Herm} there is an equivalence between dualities on an $\infty$-category $\caC$ and non-degenerate symmetric bilinear functors on $\caC$, a categorification of a bilinear form.
A symmetric bilinear functor on $\caC$ is a $C_2$-equivariant functor $ \alpha: \caC^\op \times \caC^\op \to \Spc$ from the flip action on
$\caC^\op \times \caC^\op$ to the trivial action.
It is non-degenerate if the corresponding functor
$\caC^\op \to \Fun(\caC^\op, \Spc)$ factors trough the Yoneda-embedding by an equivalence $d: \caC^\op \to \caC,$
the corresponding duality. 

An $\infty$-category with genuine duality is an $\infty$-category  $\caC$ equipped with a duality $d\colon \caC^\mathrm{op} \simeq \caC $
together with a lift of the functor
$$ \caC^\op \to \Spc[C_2], X \mapsto \map_\caC(X,d(X))$$ to the $\infty$-category $\Spc^{C_2} $ of genuine $C_2$-spaces.

We extend the correspondence between dualities and symmetric bilinear functors to an equivalence between genuine dualities and quadratic functors on a stable $\infty$-category. This way we establish a bridge between our framework and that of \emph{Poincar\'e $\infty$-categories}. 

\begin{introtheorem}[\cref{equiv_qu_gd}]
There is an equivalence of $\infty$-categories between that of stable $\infty$-categories with genuine duality and that of  
Poincar\'e $\infty$-categories.
\end{introtheorem}

This result is not only interesting form a theoretical perspective of quadratic functors, but also prepares the ground to show in \cite{comparisonKth} that the real algebraic $K$-theory spaces of \cite{Calmes_etal1} and ours coincide. An extension to the level of genuine $C_2$-spectra is current work in progress.

\subsection{The real additivity theorem}

Arguably, the most fundamental result of algebraic $K$-theory is the additivity theorem, from which several other fundamental results on $K$-theory can be derived; see \cite{grayson1987exact, staffeldt1989fundamental,mccarthy1993fundamental}.

The additivity theorem was first proven by Quillen for exact categories \cite{quillen.higher}, and later by Waldhausen \cite{waldhausen1985algebraic} using the $S$-construction; later it was generalized to the context of $\infty$-categories by \cite{bgt} and \cite{barwick.wald}, respectively. 
The additivity theorem asserts that algebraic $K$-theory splits cofiber sequences in a precise sense. 
In Barwick's language of Waldhausen $\infty$-categories, this means that, after applying the $S$-construction, the functor of Waldhausen $\infty$-categories 
$$\gamma\colon \caC^{[1]}_\cof \to \mathcal{C} \times \mathcal{C}, \ (A \hookrightarrow B) \mapsto (A,B/A),$$
induces an equivalence on geometric realizations taken in the derived $\infty$-category of Waldhausen $\infty$-categories \cite[Theorem 7.12]{barwick.wald}.
Here $ \caC^{[1]}_\cof \subset \caC^{[1]}$ is the full subcategory of the Waldhausen $\infty$-category of arrows spanned by the cofibrations in $\caC.$

The additivity theorem is not only computationally meaningful but also conceptually fundamental: it characterizes algebraic $K$-theory among all invariants of Waldhausen $\infty$-categories. 
This guiding idea had long been part of the folklore surrounding $K$-theory, but it was first made precise in the work of Barwick, who introduced the notion of \emph{additive theories}—functors satisfying the additivity theorem—and proved that algebraic $K$-theory is the initial such additive theory. 
In other words, Barwick formalized the long-suspected universal property of algebraic $K$-theory within the framework of Waldhausen $\infty$-categories.

We prove a real version of the additivity theorem.
For that we construct for every Waldhausen $\infty$-category with genuine duality $\mathcal{C}$ a canonical map of Waldhausen $\infty$-categories 
$$\gamma\colon \caC^{[2]}_\cof  \to \widetilde{\mathcal{C} \times \mathcal{C}} \times \mathcal{C},$$
that preserves genuine dualities. 
Here $\caC^{[2]}_\cof  \subset \caC^{[2]} $
is the full subcategory of composable pairs of two cofibrations,
and $\widetilde{\mathcal{C} \times \mathcal{C}}$ is the product $\caC \times \caC$ equipped with the duality $(X,Y) \mapsto (Y^\dual, X^\dual)$ and the canonical genuine refinement.

More explicitly, the functor $\gamma$ maps the sequence $(A \to B \to C) $ of two cofibrations to the triple $ (A, C/B,B/A)$. Regarding this map, we prove a real version of the additivity theorem expressed in terms of the real $S_\bullet$-construction, which leads to a real additivity theorem for real $K$-theory:

\begin{introtheorem}[\cref{thm:add}]
For every Waldhausen $\infty$-category with genuine duality, 
the map 
$$S(\gamma)\colon S(\Fun([2],\caC)^\cof) \to S(\widetilde{\mathcal{C} \times \mathcal{C}}) \times S(\mathcal{C})$$
becomes an equivalence after realization in the derived
real $\infty$-category of Waldhausen $\infty$-categories with genuine duality.
\end{introtheorem}

To prove the real additivity theorem we develop real simplicial homotopy theory, an extension of simplicial homotopy theory to the $C_2$-equivariant setting, and construct intricate real simplicial homotopies.
Our arguments are inspired by
Lectures 3 and 4 of the Felix Klein Lectures “Around topological Hochschild homology”
given by Lars Hesselholt at the HCM Bonn in fall 2016. In these lectures Lars Hesselholt
discusses proofs of the additivity theorem for the K-theory of Waldhausen categories and
for real K-theory of exact categories with weak equivalences and duality, which are based on methods of \cite{mccarthy1993fundamental}.

We use the real additivity theorem to prove a universal characterization of real algebraic $K$-theory analogous to Barwick's characterization of algebraic $K$-theory. Moreover, we combine this universal property with a universal property of the real $\infty$-category of genuine $C_2$-spectra to lift the real $K$-theory functor to genuine $C_2$-spectra:

\begin{introtheorem}[\cref{thm:lifting_Omega11_S11}]
The real $K$-theory functor $\mathrm{KR}\colon \Wald_\infty^\gd \to \Spc^{C_2}$ lifts to a real functor to connective genuine $C_2$-spectra  $$\mathrm{KR}\colon \Wald_\infty^\gd \to \Sp_{\geq0}^{C_2}.$$

\end{introtheorem}

Hermitian $K$-theory of an $\infty$-category with genuine duality arises as the connective spectrum of $C_2$-fixed points of the real algebraic $K$-theory spectrum.

\subsection{A universal characterization of real algebraic $K$-theory}

Barwick characterizes algebraic $K$-theory as the initial additive theory sending a Waldhausen $\infty$-category to its space of objects. Moreover, inspired by Grothendieck's philosophy of motives, he constructs a universal
additive invariant of Waldhausen $\infty$-categories and lifts algebraic $K$-theory
to a motivic invariant, which he further relates to Goodwillie calculus.

Barwick calls the corresponding $\infty$-category of motives
the fissile derived $\infty$-category of Waldhausen $\infty$-categories since it is 
the universal environment where cofiber sequences split.
Barwick's fissile derived 
$\infty$-category of Waldhausen $\infty$-categories is the analogue of Blumberg-Gepner-Tabuada's stable $\infty$-category of non-commutative additive motives, to which every additive invariant of stable $\infty$-categories descends.

As a key result Barwick proves that every additive invariant of Waldhausen $\infty$-categories descends to an excisive functor in the sense of Goodwillie calculus on the fissile derived $\infty$-category of Waldhausen $\infty$-categories.
This mirrors the corresponding result of
Blumberg-Gepner-Tabuada \cite{bgt}, who show that every additive invariant of stable $\infty$-categories universally factors through the stable $\infty$-category of non-commutative additive motives by a left adjoint corepresentable functor.

We extend Barwick's results to the real setting and establish a real version of Barwick's universal characterization of algebraic $K$-theory. We characterize real algebraic $K$-theory by a universal property among real theories, the $C_2$-equivariant version of Barwick's notion of theories, that are real additive, i.e.\ satisfy the real additivity theorem:

\begin{introtheorem}[\cref{univ}]
Real algebraic $K$-theory $\KR\colon \Wald_\infty^\gd \to \Spc^{C_2}$ is the universal real additive theory on the real functor $\Wald_\infty^\gd \to \Spc^{C_2}$ 
sending a Waldhausen $\infty$-category with genuine duality to its space of objects.
\end{introtheorem}

Analogously to Barwick's fissile derived $\infty$-category of Waldhausen $\infty$-categories, which plays the role of the $\infty$-category of motives, we define a real version of this construction and prove a similar result to \cite[Theorem 7.8 and Corollary 6.9.1]{barwick.wald}. Explicitly, we prove that real algebraic $K$-theory $\KR$ extends to the $\infty$-category of fissile Walhausen $\infty$-categories with genuine duality, denoted by $\D\Wald_\infty^{\gd,\fiss} $, which is the respective real $\infty$-category of motives in our setting. We identify such extension with the genuine excisive approximation\textemdash which is the excisive approximation in the sense of $C_2$-equivariant Goodwillie calculus with respect to the representation sphere $S^{1,1}$\textemdash of the natural real functor
$\iota\colon \D\Wald_\infty^{\gd,\fiss} \subset \D\Wald_\infty^{\gd} \to \Spc^{C_2}$ sending a fissile Waldhausen $\infty$-category with genuine duality to its space of objects.

\begin{introtheorem}(\cref{prop:Kth_Goodwillie}, \cref{prop:Kexcisive_approximation}) Real algebraic $K$-theory $\KR\colon\Wald_\infty^\gd \to \Spc^{C_2}$ factors as 
\[
\begin{tikzcd}
\Wald_\infty^\gd\ar[rr, "\KR"]\ar[dr]        &&\Spc^{C_2}\\
&\D\Wald_\infty^{\gd,\fiss}\ar[ur, "\KR'"']  &
\end{tikzcd}
\]
where $\KR'$ is the genuine excisive approximation of the real functor $\iota\colon \D\Wald_\infty^{\gd, \fiss} \to \Spc^{C_2}$
that sends $\caC\in \D\Wald_\infty^{\gd, \fiss}$ to its space of objects. Moreover, $\KR'(\caC)$ can be computed as the colimit of the following sequential diagram of genuine $C_2$-spaces
$$\iota(\caC) \to \Omega^{S^{1,1}}(\iota(\Sigma^{S^{1,1}}(\caC))) \to \Omega^{S^{1,1}} (\Omega^{S^{1,1}} (\iota(\Sigma^{S^{1,1}}(\Sigma^{S^{1,1}}(\caC))))) \to\cdots$$ 
where $S^{1,1}$ is the representation sphere.
\end{introtheorem}

\vspace{5mm}

\subsection*{Relation to other work} We have mentioned these relations throughout the introduction, but we amalgamate them here for easy reference. 

In the $1$-categorical setting, Hesselholdt-Madsen \cite{hesselholt-madsen} define real algebraic $K$-theory for exact categories with weak equivalences and duality via a hermitian $S$-construction and prove a real additivity theorem.

Our work is built on the treatment of algebraic $K$-theory for Waldhausen $\infty$-categories developed by Barwick \cite{barwick.wald}, incorporating and building up on work of the second author \cite{spitzweck.gw} who constructs Grothendieck-Witt spaces for stable $\infty$-categories with duality. 

In \cite{Calmes_etal1,Calmes_etal2,Calmes_etal3} 
Calmes-Dotto-Harpaz-Hebestreit-Moi-Land-Nardin-Nikolaus-Steimle
build an extensive framework for real algebraic $K$-theory of Poincar\'e $\infty$-categories, stable $\infty$-categories equipped with a non-degenerate quadratic functor. The latter is a categorification of a quadratic form suggested by Lurie \cite{lurie.L4}.
Their framework extends the framework of \cite{bgt} and \cite{spitzweck.gw} and
supports applications in algebraic surgery theory connecting 
hermitian $K$-theory to $L$-theory, both appearing respectively as the genuine $C_2$-fixed points and the geometric $C_2$-fixed points of a real $K$-theory genuine $C_2$-spectrum.

The following table illustrate the parallelisms established above. 

\begin{center}
\begin{tabular}{ | l | l | c | c}
\hline
 & Non-real & Real \\
\hline
Stable & Blumberg-Gepner-Tabuada & Calm\`es at al. \\ \hline
Non-stable & Barwick & This work \\
\hline
\end{tabular}
\end{center}

Moreover, in \cite{comparisonKth} we prove that the genuine $C_2$-spaces of real algebraic $K$-theory of Calm\`es et al. and our framework  are canonically equivalent.
An extension to the level of genuine $C_2$-spectra is current work in progress.

In recent work, Schlichting \cite{schlichting1} constructs hermitian $K$-theory for exact form categories with weak equivalences, exact categories with weak equivalences and duality equipped with quadratic data, and proves an additivity theorem for hermitian $K$-theory via the hermitian $\caQ$-construction.
Moreover Schlichting-Marlowe \cite{schlichting2} prove an equivalence between this construction of hermitian K-theory and the one of \cite{Calmes_etal1}.

\subsection{Outline}
This work is divided in two parts plus three appendices. In \textbf{\cref{part:1}} we set the stage by presenting two equivalent ways of thinking about our base objects: $\infty$-categories with duality. \textbf{\cref{part:2}} build on that work in order to define and motivate our input categories, namely, Waldhausen $\infty$-categories with genuine duality. In \textbf{\cref{part:3}} we finally dive into real $K$-theory. Each of the first three parts is comprised by three sections itself.

\textbf{\cref{part:1}} includes \cref{sec:GD,sec:real_structures,sec:GD_complete_rSegalspaces}. \textbf{\cref{sec:GD}} is devoted to introduce the notion of $\infty$-category with genuine duality, which we do in \cref{subsec:GD}. For this we work first on $\infty$-categories with duality (\cref{subsec:infty_w_duality}) and define hermitian objects, the means by which we refine dualities, in \cref{subsec:lax_hermitian_objects}.

In \textbf{\cref{sec:real_structures}} we focus on presenting a real notion of complete Segal spaces \textemdash by which we mean a version enriched in genuine $C_2$-spaces. We dedicate \cref{subsec:real_infty_categories} to define real $\infty$-categories and study their main features. In \cref{subsec:real_Segal_objects} we define (complete) real Segal objects. 

In the last section of \cref{part:1}, \textbf{\cref{sec:GD_complete_rSegalspaces}}, we identify complete real Segal spaces as a model of $\infty$-categories with genuine duality.

\textbf{\cref{part:2}} exists to introduce, justify, and prove the essential features of our input categories: Waldhausen $\infty$-categories with genuine duality. We start this part with \textbf{\cref{sec:genuine_add_preadd_stable}} were we present the theory of (pre)additive and stable $\infty$-categories with genuine duality in \cref{subsec:preadd_with_gd,subsec:additive_with_gd,subsec:stable_with_gd} respectively. Remarkably, we prove that the $\infty$-categories of (pre)additive and stable $\infty$-categories with genuine duality are presentable, genuine preadditive \textemdash a $\Spc^{C_2}$-enriched and stronger version of preadditivity\textemdash and carry closed symmetric monoidal structures; these results can be found in \cref{subsec:structural_results_of_additive_infty_categores_w_gd}.

\textbf{\cref{sec:quadratic_functors}} is devoted to show an equivalence between quadratic functors and stable $\infty$-categories with duality, \cref{subsec:stable_wgd_are_quadratic_functors}. After a brief recapitulation of quadratic functors in \cref{subsec:recap_qu}, we give an interpretation of dualities in terms of symmetric monoidal right fibrations in \cref{subsec:dualities_as_symmetric_monoidal_righ_fibs}. Those subsections together allow us to show the aforementioned equivalence, \cref{fghkklll}.

Such a result plays the role of a hinge; it justifies our definition of stable $\infty$-categories and in doing so motivates the generalization to \emph{Waldhausen} $\infty$-categories with genuine duality, which are finally introduced in \textbf{\cref{sec:Wald_gd}}, closing \cref{part:2}.

Real $K$-theory enters as the protagonist of \textbf{\cref{part:3}}, which comprises \cref{sec:additivity_thm,sec:universal_prop,sec:spectrum}. In \textbf{\cref{sec:additivity_thm}} we prove an additivity theorem for an enriched version of the hermitian $S_\bullet$-construction for Waldhausen $\infty$-categories with genuine duality. Since the proof of the main theorem (Additivity Theorem, \cref{thm:add}) is intrincated, we dedicate \cref{subsec:add_statement_and_directions} to present the statement and directions to read the proof. Each of the subsequent three subsections, \cref{subsec:S_construction,subsec:DWald,subsec:real_s_htpies} take care of important ingredients for the proof. These are, respectively, a real version of the $S_\bullet$-construction in our setting \textemdash for which we strongly use that $\Wald_\infty^\gd$ is cotensored over the $\infty$-categories of $\infty$-categories with genuine duality, result that is deduced from  \cref{subsec:structural_results_of_additive_infty_categores_w_gd} \textemdash; the non-abelian derived category of $\Wald_\infty^\gd$, which will correct completeness failures of such; and real simplicial homotopies. In \cref{subsec:add_thm} we present the core of the proof of the Additivity Theorem, and in \cref{subsec:lemmas_add} can be found several lemmas needed in the proof of this result that were deferred to simplify the read. 

\textbf{\cref{sec:universal_prop}} is where it first appears the definition of the real $K$-theory genuine $C_2$-space associated to a Waldhausen $\infty$-category with genuine duality. The goal of this section is to show a universal property of such real $K$-theory. In \cref{subsec:pre_add_theoroes} we define preadditive and additive theories, where additive theories are roughly preadditive theories satisfying the additivity theorem. \cref{subsec:additivization} shows that there is a canonical way to associate an additivity theory to every preadditive theory $\phi$, we call this its additivization and denote it by $\add(\phi)$. To show that such an association is universal we dedicate \cref{subsec:universal_add_th}; precisely, we show that $\add(\phi)$ is the initial additive theory
equipped with a map of theories $\lambda\colon \phi \to \add(\phi)$ (\cref{univ}). Since the real $K$-theory functor is no more than the additivization of the real functor that associates the maximal subspace to any Waldhausen $\infty$-category, we obtain a universal property for real $K$-theory. 

The last piece of \cref{part:3} is \textbf{\cref{sec:spectrum}}. The aim of this section is to promote the real $K$-theory genuine $C_2$-space to a genuine $C_2$-spectrum, notion that is defined therein. We begin by recalling the definition of the real $K$-theory genuine $C_2$-space and describe both its $C_2$-fixed points and underlying space in \cref{subsec:KR_space}. In \cref{subsec:enriched_spectra} we introduce generalities about enriched spectra, walking the path to present the definition and main features of genuine $C_2$-spectra. Finally, we promote out reak $K$-theory functor to a functor whose target is the $\infty$-category of genuine $C_2$-spectra in \cref{subsec:KR_spectrum}.

Last, we have the appendices. In \textbf{\cref{Appx:enriched}} we write a sort of consultation guide for the results we make use of about the category theory of enriched $\infty$-categories. In \textbf{\cref{Appx:real_spine_inc}} we provide further evidence that our definition of real spine inclusions, used to define complete real Segal objects, compares with its usual counterpart in the expected manner.

\subsection*{Acknowledgements}
We would like to thank Lars Hesselholt for giving the inspiring lecture series at HCM Bonn on topological Hochschild homology. Part of these notes are based on Lectures 3 and 4.
We would like to thank Clark Barwick, Hongyi Chu, David Gepner, Kristian Moi, Thomas Nikolaus, Oliver R\"ondigs and Manfred Stelzer for very helpful discussions and suggestions on the subject.

This material is partially based upon work supported by the National Science Foundation under Grant No. 1440140, while the second author was in residence at the Mathematical Sciences Research Institute in Berkeley, California, during the spring semester 2019. The third named authour gratefully acknowledges the support of the DFG-GK 1916 while at Osnabr\"uck, an international Macquarie University Research Excellence Scholarship, and of the National Science Foundation grant DMS-1652600 while at Johns Hopkins.

\part{\texorpdfstring{$\infty$-categories with genuine duality}{Infinity-categories with genuine duality}}\label{part:1}

\section{\texorpdfstring{$\infty$-categories with genuine duality}{Infinity-categories with genuine duality}}\label{sec:GD}

The aim of this section is to introduce and study the notion  of $\infty$-category with genuine duality (\cref{def:GD}), a refinement of that of $\infty$-category with duality. We begin by studying in \cref{subsec:infty_w_duality} $\infty$-categories with duality as presented in \cite{hls}. In particular, we devote \cref{subsubsec:infty_w_duality_as_Segal_Spc} to lay the groundwork for \cref{sec:GD_complete_rSegalspaces} by showing that $\infty$-categories with duality can be seen as an $\Spc[C_2]$-enriched version of complete Segal spaces.

Coming back to the main point of this section, the reader should have in the back of their minds that the notion of $\infty$-category with genuine duality that we present generalizes that of genuine $C_2$-space, an idea that we make precise in \cref{subsubsec:genuine_C2_spaces}. Moreover, an analogy can be drawn between \emph{the sense} in which the concept of genuine $C_2$-space is a refinement of that of $C_2$-spaces and the concept of genuine $\infty$-category with duality is a refinement of that of $\infty$-category with duality; we explain this in the next paragraphs.

Recall that genuine $C_2$-spaces are often considered in topology, and consist of $C_2$-spaces equipped with a notion of $C_2$-fixed points that refine the homotopy $C_2$-fixed points in the following sense. A genuine $C_2$-space is a $C_2$-space $X$ together with a map of spaces $Y\to X^{hC_2}$, and $Y$ is called the space of $C_2$-fixed points of $X$.

The homotopy $C_2$-fixed points of a genuine $C_2$-space $X$ are, roughly, given by an object $x\in X$ and an equivalence $x\simeq \tau(x)$ where $\tau$ is the involution on $X$. To take their place, given an $\infty$-category with duality $\caC$ we define its hermitian objects (see \cref{subsec:lax_hermitian_objects}). These are roughly given by an object $C\in \caC$ and a morphism $\alpha\colon C \to C^\dual$ for $(-)^\dual$ the duality on $\caC$. We denote by $\caH(\caC)$ the $\infty$-category of hermitian objects on $\caC$.

We then define $\infty$-categories with genuine duality as
$\infty$-categories with duality $\caC$ that are equipped with a functor 
$\phi\colon H \to \caH(\caC)$, where we ask for the fibers of $\phi$ to vary contravariantly \textemdash formally, we ask that $\phi$ be a right fibration. So $\infty$-categories with genuine duality $\caC$ are $\infty$-categories with duality equipped with a genuine refinement $H$ of the hermitian objects of $\caC$.  As expected, when $\caC$ is a space, the notion of $\infty$-category with duality reduces to the notion of genuine $C_2$-space.

The $\infty$-category $\Cat_\infty^\gd$
of small $\infty$-categories with genuine duality fits into a commutative square 
\[
\begin{tikzcd}
\Cat_\infty^\gd \ar[r] \ar[d]      & \Cat_\infty^{hC_2} \ar[d] \\
\Spc^{C_2} \ar[r]      & \Spc[C_2],
\end{tikzcd}
\]

where the right vertical functor takes the maximal subspace and the left vertical functor takes the maximal subspace of the genuine refinement. We conclude the section by showing that the $\infty$-category of small $\infty$-categories with genuine duality $\Cat_\infty^\gd$ is cartesian closed, as well as the similarly defined $\infty$-category of small simplicial spaces with genuine duality. This is the content of \cref{subsubsec:GD_is_cartesian_closed}.

\subsection{Complete Segal spaces} This model of $\infty$-categories will have a relevant role in this section. For this reason, despite the fact that the reader is surely familiar with complete Segal spaces, we make a brief summary to fix notation. 

\begin{definition}
We call a simplicial space $X\colon  \Delta^\op \to\Spc$ a Segal space if for any $[n] \in \Delta$ 
the canonical map $$X_n \to X_1 \times_{X_0} ... \times_{X_0} \times X_1$$
is an equivalence.
\end{definition}

\begin{definition}
We say that a Segal space $X$ is complete if it is local
with respect to the unique map $\mathcal{J} \to \ast$,
where $\mathcal{J}$ is the nerve of the groupoid with two elements and unique isomorphism between them.
\end{definition}

\begin{proposition}[{\cite[Theorem 14.6]{barwick-schommer-pries.uni}}]\label{prop:localization_cSegSpc}
The composite
$$\Cat_\infty\xrightarrow{y}\caP(\Cat_\infty)\to\caP(\Delta)$$
where the first map is the Yoneda embedding and the second one restricts along the embedding $\Delta\subset\Cat_\infty$,
admits a left adjoint, and induces an equivalence to the full subcategory of complete Segal spaces.
\end{proposition}

\subsection{$\infty$-categories with duality}\label{subsec:infty_w_duality}

We start by recalling the definition of $\infty$-categories with duality, originally presented in \cite{hls}. Roughly, we equip the $\infty$-category $\Cat_\infty$ of small $\infty$-categories with the canonical $C_2$-action that sends a small $\infty$-category to its opposite $\infty$-category, and define small $\infty$-categories with duality as $C_2$-homotopy fixed points with respect to this $C_2$-action.

\begin{construction}[{$C_2$-action on the simplex category}] We consider the $C_2$-action in $\Delta$ that in objects sends $[n]\mapsto [n]$, and sends a morphism $[n]\xrightarrow{f}[m]$ in $\Delta$ to the morphism
$$[n]\simeq[n]^\op\xrightarrow{f^\op}[m]^\op\simeq [m]$$
ehere the equivalence $ [n] \simeq [n]^\op$ sends $i\mapsto n-i$.
\end{construction}

\begin{construction}[{$C_2$-action on $\Cat_\infty$}] Observe that the $C_2$-action in $\Delta$ above, induces a $C_2$-action on presheaves $\caP(\Delta)$ by precomposition: 
\[
\begin{tikzcd}[row sep=tiny, column sep=small]
    BC_2 \ar[r]       &\Cat_\infty\ar[r]    & \widehat{\Cat_\infty}\\
    \ast  \ar[r,mapsto]      &\Delta \ar[r,mapsto]        &\caP(\Delta).
\end{tikzcd}
\]
Now, we consider the composite
$$\Cat_\infty\xrightarrow{y}\caP(\Cat_\infty)\to\caP(\Delta)$$
where the first map is the Yoneda embedding and the second one comes from the embedding $\Delta\subset\Cat_\infty$.

By \cref{prop:localization_cSegSpc} the composite induces an equivalence to the full subcategory of complete Segal spaces. Now observe that the full subcategory of complete Segal spaces is closed under the $C_2$-action since for every complete Segal space $X$ the presheaf $X \circ (-)^\op$ is evidently a complete Segal space.
This implies that the $C_2$-action on $\caP(\Delta)$
induces a $C_2$-action on $\Cat_\infty$ and the embedding
$\Cat_\infty \to \caP(\Delta)$ is $C_2$-equivariant.
\end{construction}

\begin{remark}
Note that on objects, this action is simply taking the opposite category. 
\end{remark}

\begin{definition}
We define the $\infty$-category of small $\infty$-categories with duality as the homotopy fixed points of the non-trivial $C_2$-action on $\Cat_\infty$ presented above. We denote this $\infty$-category by $\Cat_\infty^{hC_2}$.
\end{definition}

\begin{definition}\label{def:sspc_w_duality}We define the $\infty$-category of simplicial spaces with duality as the homotopy fixed points of the non-trivial $C_2$-action on $\caP(\Delta)$. We denote this $\infty$-category by $\caP(\Delta)^{hC_2}$.
\end{definition}

One could ask why do we choose this action to define duality. The answer is unambiguous: it's the only non-trivial one. 

\begin{lemma}\label{unicity}
The $C_2$-actions on $\Delta$ and $\Cat_\infty$ constructed above are the only non-trivial ones.
\end{lemma}
\begin{proof}
By the same argument found in \cite[Proposition 6.2]{toen.vers-ax}, every autoequivalence of $\Cat_\infty$ restricts to an autoequivalence of $\Delta \subset \Cat_\infty$. This restriction defines a map $\phi\colon \mathrm{Aut}({\Cat_\infty}) \to \mathrm{Aut}({\Delta})$ of grouplike $A_\infty$-spaces. Since $\Cat_\infty$ is a localization of $\caP(\Delta)$, the map $\phi$ is an embedding. One can see that there is an isomorphism of groups $\mathrm{Aut}({\Delta}) \cong C_2 $, which then implies that $\phi$ is an equivalence.

The spaces of $C_2$-actions on $\Cat_\infty$ and on $\Delta$ are, respectively, the spaces of $A_\infty$-maps $C_2 \to \mathrm{Aut}({\Cat_\infty})$ and $C_2 \to \mathrm{Aut}(\Delta).$

Consequently, the spaces of $C_2$-actions on $\Cat_\infty$ and on $\Delta$ are both equivalent to the set of group endomorphisms of $C_2$, which consists of two elements.
\end{proof}

In what remains of this subsection we recall properties of $\infty$-categories with duality, and their relation with simplicial spaces with duality.

\begin{remark}
The embedding $ \Cat_\infty^{hC_2} \subset \caP(\Delta)^{hC_2}$ is an accessible localization of presentable $\infty$-categories
since the embedding $\Cat_\infty \subset \caP(\Delta)$ is an accessible localization of presentable $\infty$-categories and the inclusion of the subcategory $\Pr^R$ of presentable $\infty$-categories and right adjoint functors into all large $\infty$-categories preserves small limits.
\end{remark}

\begin{remark}We can see any space as an $\infty$-category via the canonical embedding $\Spc \to \caP(\Delta) \to \Cat_\infty$, where the first map assigns the constant presheaf and the second is the left adjoint to the embedding $\Cat_\infty\subset\caP(\Delta)$.

Now the non-trivial $C_2$-action on $\Cat_\infty$
restricts to the trivial $C_2$-action on $\Spc$,
as it is the only one that it admits, for $\mathrm{Aut}(\Spc)$ is contractible. Therefore, the embedding $\Spc\subset\Cat_\infty$ is $C_2$-equivariant.

Since the homotopy fixed points of the trivial $C_2$-action on $\Spc$ are the spaces with a $C_2$-action, this is $\Spc^{hC_2}=\Spc[C_2]$, we obtain the induced embedding $$\Spc[C_2] \subset \Cat_\infty^{hC_2},$$

which admits left and right adjoints, the latter given by taking maximal subspace.

\end{remark}

We proceed now to show that the $\infty$-categories $\Cat_\infty^{hC_2}$ and $\caP(\Delta)^{hC_2}$ are cartesian closed. 

\begin{lemma}\label{forgetful_hasadjoints}
Let $\caC$ be an $\infty$-category with $C_2$-action that admits finite coproducts. Then the forgetful functor $\caC^{hC_2} \to \caC$ admits a left adjoint sending $X \in \caC$ to $X \coprod \tau(X)$, where $\tau$ denotes the involution on $\caC$ given by its $C_2$-action. 
\end{lemma}
\begin{proof}
    A proof of this (stated in the case that $\caC$ is $\caC^\op$ is given in \cite[Lemma 2.10]{hls}).
\end{proof}

\begin{corollary}\label{cor:hC2_functor_cat}Let $\caA,\caB$ be $\infty$-categories with $C_2$-action such that $\caB$ admits finite products.
The forgetful functor $$\Fun(\caA,\caB)^{hC_2} \to \Fun(\caA,\caB)$$ admits a right adjoint $R$ that sends a functor $F\colon \caA \to \caB$ to $F \prod (\sigma \circ F \circ \tau),$
where $\sigma$ is the involution of $\caB$ and $\tau$ is the involution of $\caA$.
\end{corollary}

\begin{proof}
For any $\infty$-categories $\caA,\caB$ with $C_2$-action the $\infty$-category $\Fun(\caA,\caB)$ carries an induced $C_2$-action, which sends a functor $F\colon\caA \to \caB$ to $\sigma \circ F \circ \tau,$
where $\sigma$ is the involution of $\caB$ and $\tau$ is the involution of $\caA$. Hence the statement follows directly from \cref{forgetful_hasadjoints}.
\end{proof}

\begin{proposition}\label{hC2_inherits_monoidal}Let $\caC$ be a closed monoidal $\infty$-category with a (compatible) $C_2$-action, that admits finite coproducts. Then the $\infty$-category $\caC^{hC_2}$ is also closed monoidal, and the forgetful functor $\caC^{hC_2} \to \caC$ preserves internal homs.
\end{proposition}

\begin{proof}
We first note we may reduce to the case that $\caC$
is presentable: the monoidal $C_2$-action on $\caC$ yields a closed monoidal $C_2$-action on the $\infty$-category $\caP(\caC)$ of presheaves on $\caC$ endowed with Day-convolution such that the Yoneda-embedding $\caC \subset \caP(\caC) $ is a $C_2$-equivariant monoidal functor.
With $\caP(\caC)$ also $\caP(\caC)^{hC_2}$
is presentable and so monoidal closed as the monoidal forgetful functor $\caP(\caC)^{hC_2} \to \caP(\caC)$ preserves small colimits.
As the Yoneda-embedding $\caC \subset \caP(\caC)$ preserves internal homs, the internal hom in $\caP(\caC)^{hC_2}$ of two objects of $\caC^{hC_2} \subset \caP(\caC)^{hC_2}$ belongs to $\caC^{hC_2}$
if we have proved that the forgetful functor $\caP(\caC)^{hC_2} \to \caP(\caC)$ preserves internal homs.

By \cref{forgetful_hasadjoints} the conservative forgetful functor $\caC^{hC_2} \to \caC$ admits a left adjoint sending $X \in \caC$ to $X \coprod \tau(X) $.
Thus the forgetful functor $\caP(\caC)^{hC_2} \to \caP(\caC)$ preserves internal homs since for every $X \in \caC^{hC_2}, Y \in \caC$ the canonical morphism $$ (X \otimes Y) \coprod \tau(X \otimes Y) \to (X \coprod \tau(X)) \otimes (Y \coprod \tau(Y)) \to X \otimes (Y \coprod \tau(Y))$$ is an equivalence.
\end{proof}

\begin{corollary}\label{Cat_infty_SpcC2_enrichment}The $\infty$-categories $\Cat_\infty^{hC_2}$ and $\caP(\Delta)^{hC_2}$ are cartesian closed and so carry closed actions of $\Spc[C_2]$.
\end{corollary}

\begin{proof}
From the proof of \cref{hC2_inherits_monoidal} we can see that the inherited closed monoidal structure is given by the cartesian product, since both $\Cat_\infty$ and $\caP(\Delta)$ are cartesian closed. This implies that they admit closed actions over themselves and so over $\Spc[C_2]$ by restriction along the canonical embedding, which is left adjoint.
\end{proof}

\subsubsection{Interpretation via complete Segal spaces}\label{subsubsec:infty_w_duality_as_Segal_Spc}
As next we will present a different way to think about $\infty$-categories with duality by adding structure to an existing model of $\infty$-categories, that of complete Segal spaces. Indeed, we will show that the $\infty$-category of small $\infty$-categories with duality is equivalent to the $\infty$-category of an $\Spc[C_2]$-enriched version of complete Segal spaces.

The first step towards this new characterization is to enrich the simplex category $\Delta$ over $\Spc[C_2]$. For this, we rely on the theory of enriched $\infty$-category theory, developed for example in \cite{GEPNER2015575,heine-enriched1,heine-enriched3}; for ease of reading, we include the results we need with their proofs or references in \cref{Appx:enriched}.

The $C_2$-equivariant embedding $\Delta \subset \Cat_\infty$ yields an embedding $\Delta^{hC_2} \subset \Cat_\infty^{hC_2}$. We then can restrict the $\Spc[C_2]$-enrichment on $\Cat_\infty^{hC_2}$ of \cref{Cat_infty_SpcC2_enrichment} to the $\infty$-category $\Delta^{hC_2}$. For this, we will use the following notation.

\begin{notation} We denote by $\underline{\Delta}$ the $\Spc[C_2]$-enriched $\infty$-category $\Delta^{hC_2}$. 
\end{notation}

\begin{remark}\label{rmk:Delta_as_Spc[C2]_enriched}
Taking homotopy $C_2$-fixed points on all hom $C_2$-spaces of $\underline{\Delta} $ we get $\Delta^{hC_2},$ applying the forgetful functor $\Spc[C_2] \to \Spc$ on all hom $C_2$-spaces of $\underline{\Delta} $ we get a category canonically equivalent to $\Delta$ using that the forgetful functor $\Delta^{hC_2} \to \Delta$ is essentially surjective. In particular, $\underline{\Delta}$ is a $\Set[C_2]$-enriched category.

Explicitely, for $[n], [m] \in \Delta$ the $C_2$-action on the hom-set $\map_\Delta([n],[m])$ sends $[n] \xrightarrow{f} [m] $ to $$ [n] \simeq [n]^\op \xrightarrow{f^\op} [m]^\op \simeq [m], $$ where the equivalence $ [n] \simeq [n]^\op$ sends $i$ to $ n-i$ and similar for $[m].$

\end{remark}
Finally, we are in condition of proving the announced interpretation of $\infty$-categories with duality via complete Segal spaces.
\begin{proposition}\label{dhfvghjjk}
There is a canonical equivalence
$$ \caP(\Delta)^{hC_2} \simeq \Fun_{\Spc[C_2]}(\underline{\Delta}^\op, \Spc[C_2])$$
of $\Spc[C_2]$-enriched $\infty$-categories. Moreover this equivalence restricts to an equivalence 
$$ \Cat_\infty^{hC_2} \simeq \Fun_{\Spc[C_2]}(\underline{\Delta}^\op, \Spc[C_2])_{Seg}$$
where $\Fun_{\Spc[C_2]}(\underline{\Delta}^\op, \Spc[C_2])_{Seg}$ is the full subcategory of $\Spc[C_2]$-functors whose underlying functors are complete Segal spaces.

\end{proposition}

\begin{proof}

By \cref{gghhjjn,fghjkfb} for $\caC\coloneqq B (C_2)$ and $ \caV\coloneqq\Spc$ we have an adjunction
$$L\colon \Cat_\infty^{\Spc[C_2]} \rightleftarrows \Cat_\infty[C_2]\colon R. $$

The right adjoint $R$ sends an $\infty$-category with $C_2$-action $\caD$ to the canonical $\Spc[C_2]$-enrichment on $\caD^{hC_2}.$ Especially it sends an $\infty$-category $\caB$ equipped with trivial $C_2$-action to a canonical $\Spc[C_2]$-enrichment on $\caB^{hC_2} \simeq \caB[C_2].$
The left adjoint $L$ sends a $\Spc[C_2]$-enriched $\infty$-category $\caC$ to a canonical $C_2$-action on the $\infty$-category arising from $\caC$ by forgetting the $C_2$-actions on the mapping spaces.
The left adjoint $L$ sends $\underline{\Delta} $ to $\Delta$ with its unique non-trivial $C_2$-action.
So we get a canonical equivalence
$$ \caP(\Delta)^{hC_2}= \Fun(\Delta^\op, \Spc)^{hC_2} \simeq \Fun_{\Spc[C_2]}(\underline{\Delta}^\op, \Spc[C_2])$$
of $\Spc[C_2]$-enriched $\infty$-categories that seen as an equivalence of $\infty$-categories fits into a commutative triangle

\begin{equation*}
\begin{tikzcd}[column sep=small]
 \caP(\Delta)^{hC_2} \ar[rr, "\simeq"]\ar[rd]       &       &\Fun_{\Spc[C_2]}(\underline{\Delta}^\op, \Spc[C_2]) \ar[ld]\\
        &\caP(\Delta),      &
\end{tikzcd}
\end{equation*}
where the right hand functor forgets the $C_2$-actions on the mapping spaces.

Consequently under this equivalence the full subcategory $ \Cat_\infty^{hC_2} \subset  \caP(\Delta)^{hC_2}$ corresponds to those $\Spc[C_2]$-enriched functors $\underline{\Delta}^\op \to \Spc[C_2]$
that induce a complete Segal space after forgetting the $C_2$-actions on the mapping spaces.
\end{proof}

\subsection{Hermitian objects}\label{subsec:lax_hermitian_objects}

Generalizing the notion of hermitian form on a vector space, we define a hermitian structure on an object $X$ in an $\infty$-category $\caC$ with duality. For us, such thing will be a homotopy $C_2$-fixed point of the $C_2$-action on $\map_\caC(X, X^\dual)$, where $X^\dual$ is the dual of $X$, sending a morphism $f\colon X \to X^\dual $ to $$f^\dual\colon X \simeq (X^\dual)^\dual \to X^\dual.$$ If such a homotopy $C_2$-fixed point of $\map_\caC(X, X^\dual)$ lies over an equivalence $X \simeq X^\dual,$ we call the hermitian structure on $X$ a non-degenerated structure.

\subsubsection{Twisted arrow $\infty$-category}
We now present an $\infty$-categorical version of the twisted arrow category of a category; see also \cite[Sec.\ 2]{barwick.mack1} for a similar definition in the context of quasi-categories and \cite[Definition 3.1]{haugseng2020coends} for this one.

\begin{definition} The edgewise subdivision functor is the functor $e\colon\Delta \to \Delta$, given by $[n]\mapsto [n]*[n]^\op \simeq [2n+1]$. 
\end{definition}

\begin{definition}For a simplicial space $\caC$ we define its simplicial space of twisted arrows by the composition $\Tw(\caC)\coloneqq\caC\circ e^\op.$
\end{definition}

\begin{remark}\label{rmk:functor_Tw}
The above definition yields a functor $\Tw\colon\caP(\Delta)\to \caP(\Delta)$.
\end{remark}

\begin{remark}
One can check that if $\caC$ is a (conventional) category, $\Tw(N\caC)$ is an $\infty$-category that coincides with the nerve of the (ordinary) twisted arrow category of $\caC.$
\end{remark}

\begin{remark}
Observe that there are two natural transformations $\id\to e$ and $(-)^\op\to e$, given respectively by $[n] = [n] * \emptyset \to [n] * [n]^\op$ and $[n]^\op = [n]^\op* \emptyset \to [n] * [n]^\op$. These natural transformations allow us to define a map of simplicial spaces $\Tw(\caC)\to\caC\times\caC^\op$ as follows:
$$\Tw(\caC)=\caC \circ e^\op \to (\caC \circ id) \times (\caC \circ (-)^\op) = \caC \times \caC^\op. $$

\end{remark}
We conclude the introduction of the twisted arrow simplicial space by proving some result that will later be useful to explore the behavior of hermitian objects of an $\infty$-category with duality; for this we first introduce the notion of right fibration between simplicial spaces.

\begin{definition}
We call a map of simplicial spaces $X \to Y$
a right fibration if for any $n \geq 1$
the map $\{n\} \subset [n] $ induces an equivalence 
$$ X_n \to X_0 \times_{Y_0} Y_n.$$
\end{definition}

\begin{lemma}\label{rightfib}
Let $\caC$ be a Segal space, then the map $\Tw(\caC)\to\caC\times\caC^\op$ is a right fibration of simplicial spaces.
\end{lemma}
\begin{proof}
For any $n \geq 1$ the map $\{n\} \subset [n] $ induces the map 
$$\Tw(\caC)_n \to \Tw(\caC)_0 \times_{(\caC_0 \times \caC^\op_0)} (\caC_n \times \caC^\op_n)$$
that identifies with the equivalence
$$\caC_{2n+1} \to \caC_1 \times_{(\caC_0 \times \caC_0)} (\caC_n \times \caC_n) \simeq\caC_n\times_{\caC_0}\caC_1\times_{\caC_0} \caC_n.$$
\end{proof}

\begin{lemma}\label{completeright} Let $X$ and $Y$ be simplicial spaces, and $f\colon X \to Y $ a right fibration between them. If $Y$ is a (complete) Segal space, so is $X$.
\end{lemma}

\begin{proof}
Using that equivalences are stable under pullbacks, to show that $Y$ is a Segal space it is enough to verify that for every $n \geq 2 $ the evident commutative square depicted below is a pullback square.
\begin{equation}
\begin{tikzcd}\label{sghk}
X_n\ar[r]\ar[d]     &X_{n-1} \times_{X_0} X_1\ar[d]\\
Y_n\ar[r]           &Y_{n-1} \times_{Y_0} Y_1
\end{tikzcd}
\end{equation}

Since $f\colon X\to Y$ is a right fibration, for every $k\geq 0$ the commutative square below, induced by the map $\{k\} \subset [k]$, is a pullback square.

\begin{equation}
\begin{tikzcd}\label{pullb}
X_k\ar[r]\ar[d]     &X_0\ar[d]\\
Y_k\ar[r]           &Y_0
\end{tikzcd}
\end{equation}

Now, observe that when we take $k=n$, the pullback square (\ref{pullb}) factors as follows
\[
\begin{tikzcd}
X_k\ar[r]\ar[d]     &X_0\ar[d]\\
Y_k\ar[r]           &Y_0
\end{tikzcd}
=
\begin{tikzcd}
X_n\ar[r]\ar[d]     &X_{n-1} \times_{X_0} X_1\ar[d]\ar[r]       &X_1\ar[r]\ar[d]        &X_0\ar[d]\\
Y_n\ar[r]           &Y_{n-1} \times_{Y_0} Y_1\ar[r]     &Y_1\ar[r]  &Y_0,
\end{tikzcd}
\]
where the leftmost square of the factorization is diagram (\ref{sghk}) and the rightmost is diagram (\ref{pullb}) with $k=1$. Therefore, by pasting of pullbacks, it is enough to see that the middle square in the factorization, that is
\begin{equation}
\begin{tikzcd}\label{gghhjj}
X_{n-1} \times_{X_0} X_1\ar[r]\ar[d]        &X_1\ar[d]\\
Y_{n-1} \times_{Y_0} Y_1\ar[r]              &Y_1,
\end{tikzcd}
\end{equation}
is a pullback square.

That square (\ref{gghhjj}) is in fact a pullback follows by observing that it arises by applying the pullback preserving endofunctor $(-) \times_{X_0} X_1$ of the $\infty$-category of spaces to the pullback square \ref{pullb} for $k= n-1.$ This concludes the proof that $X$ is a Segal space.

For the completeness axiom, we consider the commutative diagram below, induced by the map $[1]\to[0]$ in $\Delta$.

\begin{equation}
\begin{tikzcd}\label{fgd}
X_0\ar[r]\ar[d]     &X_1\ar[d]\\
Y_0\ar[r]           &Y_1.
\end{tikzcd}
\end{equation}  
Since the pasting of diagram (\ref{fgd}) with the pullback square (\ref{pullb}) for $k=1$ is the identity, we deduce that diagram (\ref{fgd}) must also be a pullback square. Hence the map $X_0 \to X_1$ is an embedding with image the morphisms of $X$ that go under $f$ to invertible morphisms in $Y$.

To complete the proof, we need to prove that a right fibration $X \to Y$ of simplicial spaces reflects invertible morphisms, i.e.\ a morphism of $X$ admits an inverse if its image in $Y$ does.

To check this, it is enough to see that a morphism of $X$, whose image in $Y$ admits a right inverse, also admits a right inverse in $X$ that lies over the given right inverse in $Y$. Let $\alpha\colon A \to B$ be a morphism of $X$ and $\alpha'\colon A' \to B'$ its image in $Y$, whose right inverse is $\beta'\colon B' \to A'\in Y$. Since $X \to Y$ is a right fibration of simplicial spaces, there is a morphism $\beta\colon C \to A$ in $X$ lying over $\beta'.$ The composition $\alpha \circ \beta\colon C \to B $ of $X$ lies over the morphism  $\alpha' \circ \beta'= \id_{B'}$ of $Y$ and by the uniqueness of lifts, it is the identity of $B$. Hence $\alpha$ admits a right inverse.
\end{proof}

\begin{remark}\label{rem:Tw_cat_is_inftycat} It follows directly from \cref{rightfib} and \cref{completeright} that when $\caC$ is an $\infty$-category, also $\Tw(\caC)$ is one.
\end{remark}

\begin{notation}
In view of the remark above, when $\caC$ is an $\infty$-category, we will refer to $\Tw(\caC)$ as its twisted arrow $\infty$-category.
\end{notation}

\subsubsection{Hermitian objects}\label{subsubsec:lax_hermitian_objects} We are now ready to define hermitian objects. We recall the edgewise subdivision functor $e\colon  \Delta \to \Delta$, given by $[n]\mapsto [n]*[n]^\op \simeq [2n+1]$. Note that, if we endow the $\Delta$ in the source of the edgewise subdivision with the trivial action and the $\Delta$ in its target with the unique non-trivial action, then the functor $e$ is $C_2$-equivariant. This, together with the functoriality of taking presheaves, yields a $C_2$-equivariant functor $\caP(\Delta)^{hC_2} \xrightarrow{\Tw } \caP(\Delta)[C_2]$, which we will also call $\Tw$.

\begin{definition}
We define the hermitian object functor, $\H$, as the composition
$$\caP(\Delta)^{hC_2} \xrightarrow{\Tw } \caP(\Delta)[C_2]\xrightarrow{(-)^{hC_2} } \caP(\Delta).$$
\end{definition}

\begin{remark}
 By \cref{rem:Tw_cat_is_inftycat}
the functor $\H$ restricts to a functor $\Cat_\infty^{hC_2} \to \Cat_\infty.$ We will often refer to this restriction when we write $\H$.
\end{remark}

\begin{definition}\label{hermitian}
We call an object of $ \H(\caC)$ a hermitian object of $\caC.$
\end{definition}
\begin{definition}
A hermitian object of $\caC$ is a non-degenerated object if its image under the map $\H(\caC) \to \Tw(\caC)$ corresponds to an equivalence in $\caC.$
\end{definition}

In what follows we construct a forgetful functor $\caH(\caC)\to\caC$, for every $\infty$-category with duality $\caC$. Further, we show that it is a right fibration, which should be interpreted as an internal compatibility of taking hermitian structures over different objects of $\caC$. For that, we start with a technical lemma \textemdash that will also allow us to show a different perspective on hermitian objects.

\begin{lemma}\label{gghjjnv}
The canonical natural transformation $\Tw \to \id \times (-)^\op$ of functors $\caP(\Delta) \to \caP(\Delta)$ is  $C_2$-equivariant if we endow the $\Delta$ in the source with the unique non-trivial $C_2$-action and the $\Delta$ in the target with the trivial $C_2$-action.
\end{lemma}

\begin{proof} 
By \cref{cor:hC2_functor_cat} the forgetful functor 
$$U\colon \Fun(\caP(\Delta),\caP(\Delta))^{hC_2} \to \Fun(\caP(\Delta),\caP(\Delta))$$ admits a right adjoint $R$ that sends a functor $F\colon \caP(\Delta) \to \caP(\Delta)$ to $F \prod ((-)^\op \circ F ),$ where we view $\caP(\Delta)$
in the source as equipped with the induced $C_2$-action 
and $\caP(\Delta)$ in the target as equipped with the trivial $C_2$-action.

The functor $e^*\colon \caP(\Delta) \to \caP(\Delta)$ is
$C_2$-equivariant when we endow the source with the induced $C_2$-action and the target with the trivial $C_2$-action,
and therefore is an object of $\Fun(\caP(\Delta),\caP(\Delta))^{hC_2}.$
By functoriality of taking presheaves, the canonical natural transformation $\id \to e$
induces a natural transformation $U(e^*) \to \id^* = \id$
of functors $\caP(\Delta) \to \caP(\Delta),$ 
which corresponds to a $C_2$-equivariant natural transformation 
$e^* \to R(\id) =\id \prod (-)^\op$,
which refines the canonical natural transformation.
\end{proof}

To construct the forgetful functor, observe that by \cref{gghjjnv}, given any simplicial space with duality $\caC$, the canonical map $\Tw(\caC) \to \caC \times \caC^\op$ is $C_2$-equivariant and so yields a map of simplicial spaces $$\caH(\caC)= \Tw(\caC)^{hC_2} \to (\caC \times \caC^\op)^{hC_2}= \caC.$$ This functor behaves well in the following sense.

\begin{proposition}Let $\caC$ be a simplicial space with duality.
If $\caC$ is a Segal space, the map $\H(\caC) \to \caC$ is a right fibration of simplicial spaces.
\end{proposition}

\begin{proof}
This is direct consequence of \cref{rightfib} and the fact that right fibrations are closed under small limits.
\end{proof}

As announced before, we now proceed to give a different perspective on hermitian objects when $\caC$ is in fact an $\infty$-category with duality. Recall that if $\caC$ is an $\infty$-category with duality, then its twisted arrow simplicial space $\Tw(\caC)$ \textemdash and in consequence $\caH(\caC)$\textemdash is also an $\infty$-category (see \cref{rem:Tw_cat_is_inftycat}).

Observe now that the fiber of the map $\H(\caC) \to \caC $ over some object $X$ of $\caC$
is the space of homotopy $C_2$-fixed points of the fiber of the $C_2$-equivariant map $\Tw(\caC) \to \caC^\op \times \caC$ over the fixed point $(X,X^\dual), $ which is the space $\map_\caC(X, X^\dual) $ with the following $C_2$-action:
A morphism $f\colon X \to X^\dual $ is sent to $f^\dual\colon  X \simeq (X^\dual)^\dual \to X^\dual.$ 

Therefore, a hermitian object in $\caC$ is roughly an object $X$ of $\caC$ together with a morphism $ f\colon  X \to X^\dual$ in $\caC$ that is homotopic to $f^\dual\colon  X \simeq (X^\dual)^\dual \to X^\dual$ and is a non-degenerated object if the morphism  $ f\colon X \to X^\dual$ is an equivalence.

Using the perspective just given, we present now examples of calculations of hermitian objects.

\begin{example}\label{fhjhbffcc}
The $\infty$-category $\caH(\caC \coprod \caC^\op)$ is empty
as there is no duality preserving functor $[0] \ast [0]^\op \to \caC \coprod \caC^\op.$
\end{example}

\begin{example}\label{laxherm}
	Consider the $\infty$-category $[n] \in \Delta^{hC_2} \subset \caP(\Delta)^{hC_2} $, then its hermitian objects are given by
    $$\caH([n]) = \begin{cases} [\frac{n-1}{2}] \hspace{2mm} \text{for odd} \hspace{2mm} n \\
    [\frac{n}{2}] \hspace{2mm} \text{for even} \hspace{2mm} n \end{cases} $$

    Indeed, we observe that there are canonical equivalences $$ \caH([n])_k \simeq \map_{\Delta}([k]\ast[k]^\op, [n])^{hC_2} \simeq  \begin{cases} \map_{\Delta}([k], [\frac{n-1}{2}]) \hspace{2mm} \text{for odd} \hspace{2mm} n \\
    \map_{\Delta}([k], [\frac{n}{2}]) \hspace{2mm} \text{for even} \hspace{2mm} n \end{cases} $$
\end{example}

\begin{example}
Denote $\mathcal{J} \in \caP(\Delta)^{hC_2} $
the nerve of the groupoid with $C_2$-action that has two objects and one isomorphism between them,
equipped with the $C_2$-action
sending the unique isomorphism to its inverse. Then,
$$ \caH(\mathcal{J}) \simeq \mathcal{J}.$$

This follows from the canonical, natural in $[k] \in\Delta$, equivalence $$ \caH(\mathcal{J})_k \simeq \map_{\sSet}([k]\ast[k]^\op, \mathcal{J})^{hC_2} \simeq \mathcal{J}_k.$$ 
\end{example}

We wrap up this subsection with the following remark concerning the zero hermitian structure on a given object in a preadditive $\infty$-category.

\begin{remark}
When we endow a preadditive $\infty$-category $\caC$ with the cartesian symmetric monoidal structure, the right fibration $\caH(\caC) \to \caC $ is a symmetric monoidal functor. This yields a right fibration $$\theta\colon \Calg(\caH(\caC) ) \to \Calg(\caC) \simeq \caC$$ on commutative algebras, whose fiber over an object $X $ of $\caC$ is the contractible space $\Calg(\caH(\caC)_X)$, and thus $\theta$ is an equivalence.
The composition $$\caC \simeq \Calg(\caH(\caC) ) \to \caH(\caC)$$ is a symmetric monoidal section of $\caH(\caC) \to \caC $ that sends an object $X $ of $\caC$ to the unit of the $E_\infty$-space $\map_\caC(X, X^\dual)^{hC_2}$ that lies over the zero morphism $X \to X^\dual.$

\end{remark}

\subsection{$\infty$-categories with genuine duality} \label{subsec:GD}
In this subsection we introduce an enhacement of the main notions of duality in the previous subsection; namely, we define $\infty$-category with genuine duality and simplicial space with genuine duality.

\begin{definition}A simplicial space with genuine duality is a pair $(\caC, \phi)$ consisting of a simplicial space with duality $\caC$ and a map of simplicial spaces $\phi\colon  H \to \caH(\caC)$ for some $H\in\caP(\Delta)$.
\end{definition}

\begin{definition}\label{def:infty_category_wGD}
An $\infty$-category with genuine duality is a simplicial space with genuine duality $(\caC,\phi)$ such that $\caC$ is a complete Segal space and $\phi\colon H \to\caH(\caC)$ is a right fibration.
\end{definition}

Given a simplicial space (or $\infty$-category) with genuine duality, we call the map $\phi\colon  H \to \caH(\caC)$ the genuine refinement of the duality on the simplicial space (resp. $\infty$-category) $\caC$.

As is often the case, we are not only interested in the objects themselves but rather the $\infty$-category they conform. With this in mind, we make the following definitions that show that simplicial spaces (or $\infty$-categories) with genuine duality form naturally an $\infty$-category. 

Consider the evaluation at the target functors $$ \mathrm{ev}_1\colon \Fun([1], \caP(\Delta)) \to \caP(\Delta),  \ \mathrm{ev}_1\colon\Fun([1], \Cat_\infty) \to \Cat_\infty.$$ 

\begin{definition}\label{fgbnmjkk}\label{def:GDpre}
We define the $\infty$-category $\caP(\Delta)^\gd$ of simplicial spaces with genuine duality by the pullback square

\[
\begin{tikzcd}
\caP(\Delta)^\gd \ar[r] \ar[d]      &\caP(\Delta)^{hC_2} \ar[d, "\caH"]\\
 \Fun([1],\caP(\Delta)) \ar[r, "\mathrm{ev}_1"']     & \caP(\Delta).
\end{tikzcd}
\]
\end{definition}

Similarly, for $\infty$-categories with genuine dualilty, we make the following definition.

\begin{definition}\label{def:GD}
We define the $\infty$-category $\Cat_\infty^\gd$ of small $\infty$-categories with genuine duality by the pullback square
\[
\begin{tikzcd}
\Cat_\infty^\gd \ar[r] \ar[d]       & \Cat_\infty^{hC_2} \ar[d, "\caH"]\\
\caR \ar[r, "\mathrm{ev}_1"']     & \Cat_\infty,
\end{tikzcd}
\]
where $\caR \subset \Fun([1], \Cat_\infty) $ denotes the full reflexive subcategory spanned by the right fibrations.
\end{definition}

\begin{remark}\label{rmk:left_adj_GDpre_to_PDelta}
The pullback squares of \cref{def:GDpre} and  \cref{def:GD} live in the $\infty$-category of presentable $\infty$-categories and right adjoint functors and so are presentable. The left adjoint of the projection $ \caP(\Delta)^\gd \to  \caP(\Delta)^{hC_2}$ sends $\caC$ to $ (\caC, \emptyset \to \H(\caC))$ and is thus fully faithful. It restricts to a left adjoint of the projection $\Cat_\infty^\gd \to \Cat_\infty^{hC_2}.$
\end{remark}

\begin{remark}\label{rmk:canonical_gd_Cat}
Both top horizontal functors in \cref{def:GDpre} and \cref{def:GD} admit right adjoints. Indeed, the projection $ \caP(\Delta)^\gd \to  \caP(\Delta)^{hC_2}$ has a fully faithful right adjoint that sends $\caC$ to $ (\caC, \id \colon  \H(\caC) \to \H(\caC)).$ It restricts to a right adjoint of the projection $\Cat_\infty^\gd \to \Cat_\infty^{hC_2}.$

 Therefore, via these right adjoint embeddings $\caP(\Delta)^{hC_2}\to \caP(\Delta)^\gd$ and $\Cat_\infty^{hC_2}\to \Cat_\infty^\gd$, we can view simplicial spaces with duality and $\infty$-categories with duality as their genuine versions and say that objects in the essential image of the embeddings carry a standard genuine refinement.
\end{remark}

\begin{remark}\label{rmk:ff_gd_to_hC2_is_cartesian}
The functor $ \Fun([1],\Cat_\infty) \to  \Cat_\infty$ given by evaluating at the target
is a cartesian fibration. Indeed, since right fibrations are stable under pullback, and the restriction $\caR \to \Cat_\infty$ is a cartesian fibration, then the pullback $\Cat_\infty^{gd} \to \Cat_\infty^{hC_2}$ is a cartesian fibration.
\end{remark}

\subsubsection{Genuine $C_2$-spaces}\label{subsubsec:genuine_C2_spaces}
We now devote ourselves to motivate the definition of $\infty$-category with genuine duality by seeing that they encompass the classical notion of genuine $C_2$-space. For this, we begin by recalling the latter.

\begin{definition}
A genuine $C_2$-space is a pair $Y= (X, \varphi)$ consisting of a space $X$ with a $C_2$-action together with a map $\varphi\colon Z \to X^{hC_2}$.
\end{definition}

This can be seen as a refinement of the original $C_2$-action on the space $X$, and we use the notation below to illustrate this.

\begin{notation}\label{not:fixed_points} Given a genuine $C_2$-space $Y=(X,\varphi)$, we call $X$ the underlying $C_2$-space of $Y$ and write $Y^{hC_2}$ for $X^{hC_2}$ and we write $Y^{C_2}$ for $Z$ and call $Y^{C_2}$ the space of $C_2$-fixed points of $Y.$
\end{notation}

We now see that when we think of a space as an $\infty$-category, the notions of genuine $C_2$-space and genuine duality coincide. Precisely, a genuine duality on a space seen as an $\infty$-category is the same as a genuine $C_2$-space.

Indeed, a $C_2$-space $\caC$ is an $\infty$-category with duality, whose underlying $\infty$-category is a space. In this case the functor $\Tw(\caC) \to \caC$ is an equivalence and so yields an equivalence $\caH(\caC) \simeq \caC^{hC_2}$ on homotopy $C_2$-fixed points. So $\caH(\caC) $ is a space. As right fibrations are conservative functors and every map of spaces is a right fibration, a functor $H \to \caH(\caC)$ is a right fibration if and only if $H$ is a space.

More is true, we can compare not only the objects themselves but their $\infty$-categories. In order to do this, we define the $\infty$-category of genuine $C_2$-spaces.

\begin{definition}\label{def:Spc^C2}
We define the $\infty$-category of genuine $C_2$-spaces, that we denote $\Spc^{C_2}$, by the following pullback
\[
\begin{tikzcd}
\Spc^{C_2} \ar[r] \ar[d] & \Spc[C_2] \ar[d, "(-)^{hC_2}"] \\
\Fun([1], \Spc) \ar[r, "\mathrm{ev}_1"']      & \Spc.
\end{tikzcd}
\]
\end{definition}

From \cref{def:GD} and \cref{def:Spc^C2}, together with the discussion above, we get a canonical embedding $\Spc^{C_2} \subset \Cat_\infty^\gd$ that identifies genuine $C_2$-spaces with $\infty$-categories with genuine duality, whose underlying $\infty$-category is a space.

\begin{remark}
There is a canonical equivalence
$$\Spc^{C_2} \simeq \Fun(B(C_2)^{\triangleleft},\Spc)$$
sending $(X, Y \to X^{hC_2})$ to the functor $B(C_2)^{\triangleleft}\to\Spc$ that sends the object of $B(C_2)$ to $X$, the added terminal object to $Y$ and the unique non trivial map to the composite $Y \to X^{hC_2} \to X$.
\end{remark}

\begin{remark}
The pullback square of \cref{def:Spc^C2} lives in the $\infty$-category of presentable $\infty$-categories and right adjoint functors and so is itself presentable. The left adjoint of the projection $\Spc^{C_2} \to \Spc[C_2]$ sends $X$ to $ (X, \emptyset \to X^{hC_2})$, and is thus fully faithful. 
Moreover, the canonical embedding $\Spc^{C_2} \subset \Cat_\infty^\gd$
admits a left adjoint as the embeddings inducing the embedding admit left adjoints.

The right adjoint embedding $ \Spc^{C_2} \subset \Cat_\infty^\gd$ admits itself a right adjoint
that sends an $\infty$-category with genuine duality $(\caC, \phi)$ to
the genuine $C_2$-space $(\caC^\simeq, \psi), $ where $\psi$ is the pullback of the map $\phi^\simeq \colon   H^\simeq  \to \caH(\caC)^\simeq $ to
$(\caC^\simeq)^{hC_2} \simeq \caH(\caC^\simeq) \subset \caH(\caC)^\simeq. $
\end{remark}

\begin{remark}\label{rmk:canonical_gd_Spc}
The fully faithful right adjoint of the projection $\Cat_\infty^\gd \to \Cat_\infty^{hC_2}$ restrict to a right adjoint of the projection $\Spc^{C_2} \to \Spc[C_2]$. Therefore, we view $C_2$-spaces as genuine $C_2$-spaces via the right adjoint embedding (see similarity with \cref{rmk:canonical_gd_Cat}). In particular, we view $C_2$-sets, i.e. discrete $C_2$-spaces, as genuine $C_2$-spaces via the right adjoint embedding $\Set[C_2] \subset \Spc[C_2] \subset \Spc^{C_2}$.
\end{remark}

Further in the spirit of \cref{rmk:canonical_gd_Spc}, we note that we can see spaces as genuine $C_2$-spaces via the embedding $\Spc \to \Spc^{C_2}$ that sends a space $X$ to the genuine $C_2$-space given by the pair $(X,X \to X^{BC_2})$, where the first coordinate is $X$ with the trivial $C_2$-action and we use that the homotopy $C_2$-fixed points of the trivial action are $X^{BC_2}.$ This embedding is left adjoint to taking $C_2$-fixed points.

In the next remark we show the relation of genuine $C_2$-spaces with genuine simplicial spaces. 

\begin{remark}
There is a canonical pullback square
\[
\begin{tikzcd}
\caP(\Delta)^\gd\ar[r]\ar[d]  &\caP(\Delta)^{hC_2}\ar[d,"\Tw"] \\
\Fun(\Delta^\op,\Spc^{C_2}) \ar[r,"\mathrm{ev}_1"'] &\Fun(\Delta^\op, \Spc[C_2]).
\end{tikzcd}
\]
\end{remark}

\subsubsection{The $\infty$-categories $\Cat_\infty^\gd$ and $\caP(\Delta)^\gd$}\label{subsubsec:GD_is_cartesian_closed}
Now that we have justified why we care about the $\infty$-categories $\caP(\Delta)^\gd$ and $\Cat_\infty^\gd$, we proceed to show further good properties of them. In the next proposition we prove that both of them are cartesian closed and give a description of their internal homs. 

\begin{proposition}\label{htwdzfd}
Let $\mathcal{B}$ and $\mathcal{C}$ be small monoidal closed $\infty$-categories, $H\colon \mathcal{B}\to \caC$ a lax monoidal functor, and $R \subset \Fun([1], \caC)$ a full subcategory satisfying the following conditions.

\begin{enumerate}
\item $R \subset \Fun([1], \caC)$ is closed under the tensor product.

\item The base change of an object in $R$ along any map in $C$ exists and is again in $R$.
\item For all $c \in \caC$ and $\alpha \in R$ we have $\uHom_C(c,\alpha) \in R$.
\end{enumerate}

Then, every $X,Y\in P=\caB\times_{\caC} R$ admit an internal hom in $\caB\times_{\caC}\Fun([1], C)$ that belongs to $P$. Moreover, there is a natural equivalence $\uHom_P(X,Y)_\caB \simeq \uHom_\caB(X_\caB,Y_\caB)$ and a pullback square
\[
\begin{tikzcd}
s(\uHom_P(X,Y)_R) \ar[r] \ar[d]         & H(\uHom_\caB(X_\caB,Y_\caB)) \ar[d] \\
\uHom_\caC(s(X_R),s(Y_R)) \ar[r]        & \uHom_\caC(s(X_R),t(Y_R)),
\end{tikzcd}
\]
where $s$ and $t$ denote source and target.
\end{proposition}

\begin{proof}
Consider the pullbacks in the statement, as depicted below
\[
\begin{tikzcd}
Q \ar[r]\ar[d]      & \caB \ar[d, "H"]      &P\ar[r]\ar[d]        &\caB \ar[d, "H"]\\
\Fun([1], \caC) \ar[r, "\mathrm{ev}_1"']    &\caC           &R\ar[r, "\mathrm{ev}_1"']      &\caC.
\end{tikzcd}
\]
Observe that $Q$ is a monoidal $\infty$-category as $\ev_1 \colon  \Fun([1], \caC) \to \caC$ is a cocartesian fibration of monoidal $\infty$-categories.

We define $\uHom_P(X,Y) \in P$ by $\uHom_P(X,Y)_\caB\simeq \uHom_\caB(X_\caB,Y_\caB)$, and in what follows show that it satisfies the universal property of the internal hom.

We have a commutative square
\[
\begin{tikzcd}
s(\uHom_P(X,Y)_R)\otimes s(X_R)\ar[r]\ar[d]     &\uHom_\mathcal{C}(s(X_R),s(Y_R))\otimes s(X_R)\ar[d] \\
H(\uHom_\mathcal{B}(X_\mathcal{B},Y_\mathcal{B})) \otimes H(X_\mathcal{B}) \ar[d]       &s(Y_R) \ar[d] \\
H(\uHom_\mathcal{B}(X_\mathcal{B},Y_\mathcal{B}) \otimes X_\mathcal{B}) \ar[r]          & H(Y_\mathcal{B})
\end{tikzcd}
\]
that together with the map 
$ \uHom_\mathcal{B}(X_\mathcal{B},Y_\mathcal{B}) \otimes X_\mathcal{B}  \to Y_\mathcal{B}$ define a map
$ \uHom_P(X,Y) \otimes X  \to Y $ in $Q.$

We need to show that for every $Z \in Q$ the canonical map
$$ \map_Q(Z, \uHom_P(X,Y)) \to \map_Q(Z \otimes X , \uHom_P(X,Y) \otimes X) \to \map_Q(Z \otimes X , Y)$$ is an equivalence.

We have a commutative square
\[
\begin{tikzcd}
\map_Q(Z, \uHom_P(X,Y))\ar[r]\ar[d]       &\map_Q(Z \otimes X , Y)\ar[d]\\ 
\map_\mathcal{B}(Z_\mathcal{B}, \uHom_\mathcal{B}(X_\mathcal{B},Y_\mathcal{B})) \ar[r, "\simeq"]        & \map_\mathcal{B}(Z_\mathcal{B} \otimes X_\mathcal{B} , Y_\mathcal{B}).
\end{tikzcd}
\]

Consequently, it is enough to check that this square yields on the fiber over every $\phi\colon  Z_\mathcal{B} \to \uHom_B(X_\mathcal{B},Y_\mathcal{B}) $ adjoint to a morphism $Z_\mathcal{B} \otimes X_\mathcal{B} \to Y_\mathcal{B}$ an equivalence.   

On the fiber over $\phi$ the last square induces the map
$$ \map_{H(\uHom_B(X_\mathcal{B},Y_\mathcal{B})) }(s(Z_R), s(\uHom_Q(X,Y)_R)) \to \map_{H(Y_\mathcal{B})}(s(Z_R) \otimes s(X_R) , s(Y_R)), $$ 
which is equivalent to the equivalence
$$ \map_{\uHom_\mathcal{C}(s(X_R),t(Y_R)) }(s(Z_R), \uHom_\mathcal{C}(s(X_R),s(Y_R)) ) \to $$$$ \map_{t(Y_R)}(s(Z_R) \otimes s(X_R) , s(Y_R)). $$ 
\end{proof}

\begin{corollary}\label{GD(pre)_cartesian_closed}
The $\infty$-categories $\caP(\Delta)^\gd$ and $\Cat_\infty^\gd$ are cartesian closed.
\end{corollary}
\begin{proof}
For $\caP(\Delta)^\gd$ we consider $H\colon \caB\to \caC $ to be $\caH\colon \caP(\Delta)^{hC_2} \to \caP(\Delta)$; see \cref{def:GDpre}. Similarly, for the case of $\Cat_\infty^\gd$ we take $H\colon \caB\to \caC$ to be $\caH\colon \Cat_\infty^{hC_2} \to  \Cat_\infty$ and $R=\caR$; see \cref{def:GD}.
\end{proof}

After having shown that any two $\infty$-categories $\caC,\caD$ with genuine duality have an internal hom of the form $H \to \caH(\Fun(\caC,\caD))$, a natural question is to understand the genuine refinement of this internal hom in terms of the genuine refinements of $\caC,\caD$. The following remark addresses this question.

\begin{remark}\label{htedud}
Let $X,Y \in \Cat_\infty^\gd$ and $X',Y'$ be their images in $\Cat_\infty^{hC_2}$. Let $$H \to \caH(\Fun(X',Y'))$$ be the right fibration corresponding to the internal hom $\uHom_{\Cat_\infty^\gd}(X,Y)$ in $\Cat_\infty^\gd$. Then for $\varphi \in \caH(\Fun(X',Y'))$ the fiber $H_\varphi$ is canonically equivalent to 
the space $$\Fun(X',Y') \times_{\Fun(H_X,\caH^\mathrm{nd}(Y'))} \Fun(H_X,H_Y),$$ which is the limit of the functor
$$H_X^\op \to \caH(X')^\op\overset{\caH(\varphi)^\op}{\longrightarrow} \caH(Y')^\op \to \Spc,$$ where the last functor corresponds to the right fibration $H_Y \to \caH(Y')$.

\end{remark}
\section{Real structures}\label{sec:real_structures}

In this section we present a real notion of complete Segal spaces, which we later identify with $\infty$-categories with genuine duality in \cref{sec:GD_complete_rSegalspaces}. We dedicate \cref{subsec:real_infty_categories} to define real $\infty$-categories and study their main features. In \cref{subsec:real_Segal_objects} we introduce real spine inclusions and relying on them we define (complete) real Segal objects; see \cref{def:complete_real_SegSpc}. Together with real Segal spaces, we introduce an important variant of them, that of \emph{balanced} Segal spaces, that play a pivotal role in the aforementioned identification of \cref{sec:GD_complete_rSegalspaces}.

\subsection{Real $\infty$-categories}\label{subsec:real_infty_categories}

We showed in \cref{unicity} that the $\infty$-category $\Cat_\infty$ of small $\infty$-categories carries precisely two $C_2$-actions, the trivial one and the one sending an $\infty$-category to its opposite $\infty$-category, whose homotopy $C_2$-fixed points are, respectively, the $\infty$-categories with $C_2$-action and the $\infty$-categories with duality.

In \cref{subsec:GD} we generalized $\infty$-categories with duality to $\infty$-categories with genuine duality, whose $\infty$-category $\Cat_\infty^\gd$ we defined as the pullback
\[
\begin{tikzcd}
\Cat_\infty^\gd \ar[r] \ar[d]       & \Cat_\infty^{hC_2} \ar[d, "\caH"] \\
\caR \ar[r, "\mathrm{ev}_1"']   & \Cat_\infty,
\end{tikzcd}
\]
where $\caR \subset \Fun([1], \Cat_\infty) $ is the full subcategory spanned by the right fibrations.

Similarly, $\infty$-categories with $C_2$-action
generalize to $\infty$-categories with genuine $C_2$-action,
whose $\infty$-category, that we denote by $ \Cat_\infty^{C_2}$, can be defined via the pullback below 
\[
\begin{tikzcd}
\Cat_\infty^{C_2} \ar[r] \ar[d]         & \Cat_\infty[C_2] \ar[d, "(-)^{hC_2}"] \\
\caR \ar[r, "\mathrm{ev}_1"']             & \Cat_\infty.
\end{tikzcd}
\]

For our purpose, it will be more useful to consider a full subcategory of $\Cat^{C_2}_\infty$, that of $\Spc^{C_2}$-enriched $\infty$-categories. Indeed, by \cref{Venr:lem:restriction_of_G,Venr:Fun_enriched_is_Fun_to_enriched} we know that there exists a left adjoint full embedding
$$\Cat_\infty^{\Spc^{C_2}}  \subset \Cat_\infty^{C_2}$$
with essential image the $\infty$-categories with $C_2$-action, that we think of as pairs $(\caC, \phi\colon H \to \caC^{hC_2})$, such that the functor $\phi$ is essentially surjective.

\begin{Terminology}\label{terminology:real}
In this section, we will use the adjective \emph{real} to mean $\Spc^{C_2}$-enriched. In particular, a real $\infty$-category is a $\Spc^{C_2}$-enriched $\infty$-category, a real functor between real $\infty$-categories is a $\Spc^{C_2}$-enriched functor, and so on.
\end{Terminology} 

Given two real $\infty$-categories $\caC$ and $\caD$, we
denote by $\Fun_{\Spc^{C_2}}(\caC, \caD) $
the real $\infty$-category of real functors $\caC \to \caD$,
which is by definition the internal hom in the 
$\infty$-category of real $\infty$-categories with respect to the cartesian product.

There are different ways to extract related $\infty$-categories from a real $\infty$-category, we explain some natural ways to do so below. Let $\caC$ be a real $\infty$-category.

\begin{itemize}
\item We can consider the $\infty$-category arising from $\caC$ by taking $C_2$-fixed points (recall that these are the source of the genuine refinement) on all mapping genuine $C_2$-spaces. As the functor given by taking $C_2$-fixed points is corepresented by the terminal object, i.e.\ the tensorunit of $\Spc^{C_2},$ this is the $\infty$-category carrying the enrichment. For this reason we will also denote it by $\caC$.

\item The $\Spc[C_2]$-enriched $\infty$-category arising from $\caC$ by forgetting the genuine refinement on all mapping genuine $C_2$-spaces but keeping the $C_2$-action, which we denote by $\caC_{\mid C_2}$. 

\item The $\infty$-category arising from $\caC$ by forgetting the genuine $C_2$-action on all mapping genuine $C_2$-spaces (underlying plain $C_2$-action included). We call this $\infty$-category the underlying $\infty$-category of $\caC$, and we denote it by $\caC^u.$
\end{itemize}

Note that the natural transformation $(-)^{C_2} \to (-)^u$ of functors $\Spc^{C_2} \to \Spc$ gives rise to a functor $ \caC \to \caC^u$.

In the same way that we obtain $\infty$-categories from real $\infty$-categories, we can obtain underlying functors from real ones. Given two real $\infty$-categories $\caC$ and $\caD$, the functor $(-)^{C_2}\colon  \Spc^{C_2} \to \Spc$
yields a functor $$  \Fun_{ \Spc^{C_2}}(\caC, \caD) \to \Fun(\caC, \caD).$$

Note also that the forgetful functor $\Spc^{C_2} \to \Spc[C_2]$ yields a real functor 
$$\Fun_{ \Spc^{C_2}}(\caC, \caD) \to \Fun_{\Spc[C_2]}(\caC_{\mid C_2}, \caD_{\mid C_2}).$$

Similarly, the forgetful functor $\Spc^{C_2} \to \Spc$, yields a functor
$$ \Fun_{ \Spc^{C_2}}(\caC, \caD) \to \Fun(\caC^u, \caD^u).$$

We present canonical examples of real $\infty$-categories and their associated $\infty$-categories that will be useful to have at hand.

\begin{example} 
\begin{enumerate}
\item The $\infty$-category $\Spc^{C_2}$ is a real $\infty$-category, for which $(\Spc^{C_2})_{\mid C_2} \simeq \Spc[C_2]$ and $(\Spc^{C_2})^u \simeq \Spc.$ 

\item The $\Spc^{C_2}$-enrichment on $\Spc^{C_2}$ restricts to a $\Spc^{C_2}$-enrichment on $\Spc[C_2]$, with $\Spc[C_2]^u \simeq \Spc.$ One can arrive to this same $\Spc^{C_2}$-enrichment on $\Spc[C_2]$ from the closed action induced by the cartesian structure on $\Spc[C_2]$ and the forgetful functor $\Spc^{C_2} \to \Spc[C_2]$.
\end{enumerate}
\end{example}

To achieve our objective, we need to see the category $\Delta$ as a $\Spc^{C_2}$-enriched category. This follows from \cref{rmk:Delta_as_Spc[C2]_enriched}, and we will use the same notation $\underline{\Delta}$.

\begin{definition}\label{def:real_sC}
Let $\caC$ be a real $\infty$-category. We call real functors $ \underline{\Delta}^\op \to \caC$ 
real simplicial objects in $\caC.$ We denote by $\mathrm{rs}\caC$ the $\infty$-category $\Fun_{\Spc^{C_2}}(\underline{\Delta}^\op, \caC).$
\end{definition}

To define a real version of the geometric realization as below, we use weighted limits as presented in \cref{Appx:enriched}. 

\begin{definition}\label{def:real_geometric_realization}
Let $\caC$ be a real $\infty$-category and $F\colon \underline{\Delta}^\op \to \caC$ a real simplicial object. Consider $H\colon \underline{\Delta} \to \Spc^{C_2}$ to be the constant real functor with image the final object. We define the real geometric realization of $F$ as the $H$-weighted colimit of $F$.
\end{definition}

\begin{notation}\label{not:underlying_simpl_object}
Given $X$ a real simplicial object in $\caC$,
we write $X^u\colon  \Delta^\op \to \caC^u$ for its image under the functor
$$ \mathrm{rs}\caC =\Fun_{\Spc^{C_2}}(\underline{\Delta}^\op, \caC) \to  \Fun((\underbrace{\underline{\Delta}^\op)^u}_{\Delta^\op}, \caC^u).$$
Note that with this notation, for $\caC=\Spc^{C_2}$ we have $\caC^u=\Spc$, so that $X^u$ is a simplicial space.
\end{notation}

\begin{notation}For those $\infty$-categories of real simplicial objects that we use regularly, we slightly change the notation of \cref{def:real_sC} as described here.
\begin{itemize}
\item We set $\rsSpc\coloneqq \mathrm{rs} \Spc^{C_2}= \Fun_{\Spc^{C_2}}(\underline{\Delta}^\op, \Spc^{C_2}) $ and call its objects real simplicial spaces.
		
\item We set $\rsSet\coloneqq \Fun_{\Set[C_2]}(\underline{\Delta}^\op, \Set[C_2]) \simeq \Fun_{\Spc^{C_2}}(\underline{\Delta}^\op, \Set[C_2])$ and call its objects real simplicial sets. 
\end{itemize}
\end{notation}

The embedding $\Set[C_2] \subset \Spc[C_2] \subset \Spc^{C_2}$ yields an embedding $\rsSet \subset \rsSpc.$ Thus, examples of real simplicial sets will canonically give examples of real simplicial spaces.

An important class of the first ones are given by observing that, by \cref{gghhjjn}, there exists a canonical equivalence $\rsSet \simeq \sSet^{hC_2}$ together with noticing that the nerve functor $\Cat_1 \to \sSet $ from the category of small categories to simplicial sets gives rise to an embedding
$$\Cat_1^{hC_2} \hookrightarrow \sSet^{hC_2} \simeq \rsSet.$$
This allow us to see the nerve of any category with strict duality as a real simplicial set.
	
\begin{example}\label{fhnnnhffg}\label{ex:real_spine_inc}
\begin{enumerate}
\item\label{Deltan_as_rsset} For every $n \geq 0 $ the category $[n]$ has a unique strict duality that gives its nerve $\Delta^n \in \sSet$ the structure of a real simplicial set.

\item\label{real_spine_inc} For every $n \geq 2$ the functors $[1] \to [n]$, $0 \mapsto i-1$, $1 \mapsto i$ for $1 \leq i \leq n$ yield the spine inclusion $$ \Lambda_n\coloneqq \Delta^1 \sqcup_{\Delta^0 } \Delta^1 \sqcup_{ \Delta^0 } \ldots \sqcup_{\Delta^0 } \Delta^1 \subset \Delta^n.$$

The structure of real simplicial set on $\Delta^n$ given in \cref{Deltan_as_rsset} restricts to $ \Lambda_n$, as we will see in \cref{Appx:real_spine_inc}. We call $ \Lambda_n \subset \Delta^n $ the  n-th real spine inclusion.

\item\label{boundary} The structure of real simplicial set on $\Delta^n$ given in \cref{Deltan_as_rsset} restricts to $\partial\Delta_n$.

\item The contractible category with two objects carries a unique strict duality permuting the two objects. Therefore, its nerve $ \mathcal{J}$ admits a canonical structure of real simplicial set. 
\end{enumerate}
\end{example}		

We finish this subsection with a result that will appear several times in key parts of the paper, that allows us to compute real geometric realizations as a non-enriched one.

\begin{proposition}\label{prop:geom_realization_with_edgewise}Let $\caC$ be a real $\infty$-category and $X\colon\underline{\Delta}^\op \to\caC$ a real simplicial object in $\caC$. Then there is an equivalence
$$\vert X\vert\simeq \vert X\circ e^\op\vert,$$
where the left-most geometric realization is as in \cref{def:real_geometric_realization}, and $e\colon\Delta \to \Delta^{hC_2}$ is the edgewise subdivision.
\end{proposition}

\begin{proof}
First note that, for any $Y \in \caC$, it follows from the corresponding definitions of geometric realizations that the canonical maps 
$$\map_\caC(\vert X\vert,Y) \to \map_{\Fun_{\Spc^{C_2}}(\underline{\Delta}^\op, \caC)}( \delta(\vert X\vert),  \delta(Y)) \to \map_{\Fun_{\Spc^{C_2}}(\underline{\Delta}^\op, \caC)}(X,\delta(Y)) $$
and 
$$\map_\caC(\vert X \circ e^\op \vert,Y) \to \map_{\Fun(\Delta^\op, \caC)}(\delta(X \circ e^\op),\delta(\Y)) \to $$$$ \map_{\Fun(\Delta^\op, \caC)}(X \circ e^\op,\delta(\Y))$$
are equivalences. 

Consequently it suffices to check that the canonical map
$$\theta\colon  \map_{\Fun_{\Spc^{C_2}}(\underline{\Delta}^\op, \caC)}(X,\delta(Y)) \to \map_{\Fun(\Delta^\op, \caC)}(X \circ e^\op,\delta(\Y))$$
induced by the real functor $e\colon \Delta \to \underline{\Delta}$
adjoint to the edgewise subdivision $\Delta \to \Delta^{hC_2}$
is an equivalence. For this, we use \cref{Venr:lem:restriction_of_G} to obtain an adjunction
 $$L\colon \Cat_\infty^{\Spc^{C_2}} \subset \Cat_\infty^{C_2}\colon R,$$
where the right adjoint $R$ takes $C_2$-fixed points.
Moreover by \cref{equivalence from adjunction} for any $\caC \in \Cat_\infty^{\Spc^{C_2}}, \caD \in \Cat_\infty^{C_2}$ there is a canonical equivalence
$$ \Fun(L(\caC),\caD)^{hC_2} \times_{\Fun(L(\caC)^{C_2},\caD^{hC_2})  } \Fun(L(\caC)^{C_2},\caD^{C_2}) \simeq \Fun_{\Spc^{C_2}}(\caC,R(\caD))$$
of $\infty$-categories and so for any $\caB, \caC \in \Cat_\infty^{\Spc^{C_2}}$ a canonical equivalence
$$ \Fun(L(\caB),L(\caC))^{hC_2} \times_{\Fun(L(\caB)^{C_2},L(\caC)^{hC_2})  } \Fun(L(\caB)^{C_2},L(\caC)^{C_2}) \simeq \Fun_{\Spc^{C_2}}(\caB,\caC)$$
of $\infty$-categories.
We look at this equivalence for $\caB= \underline{\Delta}^\op, \Delta^\op.$
Hence by functoriality in $\caB$, the map $\theta$ is equivalent to the following map.
$$ \map_{\Fun(\Delta^\op,\caC^u)^{hC_2}}(X,\delta(Y)) \times_{\map_{\Fun((\Delta^{hC_2})^\op,\caC^{hC_2})}(X,\delta(Y)) } \map_{\Fun((\Delta^{hC_2})^\op,\caC^{C_2})}(X,\delta(Y)) $$$$  \to \map_{\Fun(\Delta^\op,\caC^u)^{hC_2}}(X \circ e^\op,\delta(Y)) \times_{\map_{\Fun(\Delta[C_2]^\op,\caC^{hC_2})}(X\circ e^\op,\delta(Y)) } $$$$ \map_{\Fun(\Delta^\op,\caC^{C_2})}(X\circ e^\op,\delta(Y)).$$

In turn, the latter is equivalent to the following one:
$$ \map_{\caC^{hC_2}}(\colim(X^u),Y) \times_{\map_{\caC^{hC_2}}(\colim(X),Y) } \map_{\caC^{C_2}}(\colim(X),Y) $$$$  \to \map_{\caC^{hC_2}}(\colim(X^u \circ e^\op),Y) \times_{\map_{\caC^{hC_2}}(\colim(X \circ e^\op),Y)} \map_{\caC^{C_2}}(\colim(X \circ e^\op \circ \nu^\op),Y), $$
where $\nu \colon \Delta \to \Delta[C_2]$ is the diagonal functor,
which is an equivalence.

Now, that the latter map is an equivalence follows from the fact that the functor edgewise subdivision $e^\op\colon \Delta^\op \to \Delta^\op$ is cofinal \textemdash hence so is the induced functor $\Delta[C_2]^\op \to (\Delta^{hC_2})^\op$\textemdash and that the functor $\nu^\op$ is cofinal as well.
\end{proof}

\subsection{Real Segal objects}\label{subsec:real_Segal_objects}
In this subsection we introduce a real (i.e.\ a $\Spc^{C_2}$-enriched) version of (complete) Segal spaces; see \cref{def:complete_real_SegSpc}. In analogy to usual Segal spaces, characterized by being local with respect to the spine inclusions, we define real Segal spaces as real simplicial spaces that are local with respect to the real version of spine inclusions $\Lambda_n \to \Delta^n $ as presented in \cref{def:real_spine_inc} (see also \cref{real_spine_inc} of \cref{ex:real_spine_inc}). 

More generally, we define Segal objects in any real $\infty$-category 
$\caC$ as those real simplicial objects $X$ in $\caC$ such that for every $Z \in \caC$ the composition $\underline{\Delta}^\op \xrightarrow{X} \caC \xrightarrow{ \map_\caC(Z,-)} \Spc^{C_2}$ is a real Segal space. This give us an uncomplicated way to define real monoids, which play an important role in \cref{sec:universal_prop}.

If $\caC$ admits finite limits and cotensors with $C_2$, we construct real Segal maps for any real simplicial object $X$ in $\caC$, and we note that we can equivalently define real Segal objects in $\caC$ as real simplicial objects, whose real Segal maps are equivalences (Corollary \ref{reali}). From this, we deduce a recognition principle for real simplicial spaces looking separately at their underlying spaces and $C_2$-fixed points spaces (see \cref{cor:recognition_pple_rSeg}).

We start by defining real spine inclusions, which we show in \cref{Appx:real_spine_inc} that it gives meaning to the description of real spine inclusions stated in \cref{real_spine_inc} of \cref{ex:real_spine_inc}. Having this presentation will allow us to prove a recognition principle for real Segal spaces; see \cref{cor:recognition_pple_rSeg}. For this, we introduce first the following notation.

\begin{notation}\label{not:tilde}Let $\caC$ be a real $\infty$-category that admits cotensors with $C_2$, and $Z \in \caC$, then we denote the cotensor $Z^{C_2}$ by $\widetilde{Z \times Z}$.
\end{notation}

This notation may seem unconventional, so we offer a brief justification for this choice. When we consider $\caC=\Spc^{C_2}$ and $Z\in \Spc^{C_2}$, we know that the cotensor $Z^{C_2}$ is $Z \times Z$ with the $C_2$-action given by switching the factors and the standard genuine refinement, usually denoted by $ \widetilde{Z \times Z} $. Further, this notation has symbolic advantages; for example, for any $Y \in \Spc^{C_2}$ there is a canonical equivalence
$$ \caC(Y, \widetilde{Z \times Z}) \simeq \widetilde{\caC(Y, Z) \times \caC(Y,Z)} $$
of genuine $C_2$-spaces.

\begin{definition}\label{def:real_spine_inc}We define the spine inclusions $\Lambda^n \to \Delta^n$ for $n \geq 1$. For $n$ even, we set $\Lambda^n$ to be the coproduct of $[0]$ with $\frac{n}{2}$ copies of $\widetilde{[1] \coprod [1]}$, constructed by using alternatively the two maps $[0]\to[1]$. We have
$$\Lambda^n \coloneqq [0] \coprod_{\widetilde{[0] \coprod [0]}} \widetilde{[1] \coprod [1]} \coprod_{\widetilde{[0] \coprod [0]}}\dots \coprod_{\widetilde{[0] \coprod [0]}} \widetilde{[1] \coprod [1]}. $$
 
Given $\Lambda^n$, we define the $n-th$ spine inclusion $\Lambda^n\to\Delta^n$ by the universal property of the coproduct, where the compatible maps considered are the ones induced by considering for every $1 \leq i \leq \frac{n}{2}$, the maps $ [1] \to [n]$ in $\Delta$ sending $0 $ to $i-1$ and $1$ to $i$, that yield maps $\widetilde{[1] \coprod [1]} \to [n]$ in $\rsSpc$; and the map $[0] \to [n]$ in $\Delta^{hC_2}$ sending $0$ to $\frac{n}{2}$.

For $n \geq 3 $ odd, we set
$$\Lambda^n \coloneqq [1] \coprod_{\widetilde{[0] \coprod [0]}} \widetilde{[1] \coprod [1]} \coprod_{\widetilde{[0] \coprod [0]}}\dots \coprod_{\widetilde{[0] \coprod [0]}} \widetilde{[1] \coprod [1]}.$$

In this case, the spine inclusion $\Lambda^n\to \Delta^n$ is induced by considering, for every $1 \leq i \leq \frac{n-1}{2}$, the maps $ [1] \to [n]$ sending $0 $ to $i-1$ and $1$ to $i$ in $\Delta$; and the map $[1] \to [n]$ in $\Delta^{hC_2}$, sending $0$ to $ \frac{n-1}{2}$ and $1$ to $\frac{n-1}{2}+1$ in $\Delta^{hC_2}$.
\end{definition}

\begin{remark}
We show in \cref{Appx:real_spine_inc} that this definition refines the usual spine inclusions in a canonical way.
\end{remark}

We now introduce \emph{real} notions of Segal spaces and complete Segal spaces.

\begin{definition}\label{Seg}
A real Segal space is a real simplicial space $X$, which is local, in the $\Spc^{C_2}$-enriched sense,
with respect to all spine inclusions $\Lambda^n \subset \Delta^n$ for $n \ge 2$. This is, for every $n \ge 2$ the restriction map $$\theta_n\colon X_n \simeq \mathrm{map}_{\rsSpc}(\Delta^n, X) \to  \mathrm{map}_{\rsSpc}(\Lambda^n, X)$$ of genuine $C_2$-spaces is an equivalence.
\end{definition}

Given a real simplicial space $X$, we call $\theta_n$ the $n$-th real Segal map of $X$.

\begin{definition}\label{def:complete_real_SegSpc}
A real Segal space is called \emph{complete} if it is local, in the $\Spc^{C_2}$-enriched sense, with respect to the map $ \mathcal{J} \to \ast$. This is, the induced map $$ X_0\simeq\mathrm{map}_{\rsSpc}(\ast, X) \to \mathrm{map}_{\rsSpc}(\mathcal{J}, X)$$ of genuine $C_2$-spaces is an equivalence.
\end{definition}

Since we have defined (complete) real Segal spaces via localizations, we can consider the full real reflexive subcategories spanned by them, that we denote by $\crSegSpc \subset \rSegSpc \subset \rsSpc$.

\begin{remark}
Since $\crSegSpc$ and $\rSegSpc$ are localizations, they are closed under small limits and cotensors in $ \rsSpc, $ which are formed objectwise in $\Spc^{C_2}$. In particular, for every $ X \in \rsSpc$ and $Z \in \Spc^{C_2}$ the cotensor
$X^Z \simeq \map_{\Spc^{C_2}}(Z,-) \circ X$ is a (complete) real Segal space if $X$ is.
\end{remark}

The above remark motivates the following extension of \cref{Seg}.

\begin{definition}\label{fhjhbf}
Let $\caC$ be a real $\infty$-category. We call a real simplicial object $X\colon  \underline{\Delta}^\op \to \caC$ a real Segal object in $\caC$ if for every $Z \in \caC$ the composition
$ \underline{\Delta}^\op \xrightarrow{X} \caC \xrightarrow{ \map_\caC(Z,-)} \Spc^{C_2}$ is a real Segal space.
\end{definition}

The following special case of \cref{fhjhbf} will be fundamental in \cref{sec:universal_prop}.

\begin{definition}
We call a real Segal object $X$ in $\caC$ a real monoid in $\caC$ if $X_0$ is a final object in $\caC.$
\end{definition}

We write $\mathsf{r}\Mon(\caC) \subset \mathsf{r}\Seg(\caC) \subset \mathsf{rs}\caC$ for the full subcategories of real simplicial objects on $\caC$ spanned, respectively, by real monoids and real Segal objects in $\caC.$

In most of what remains of this section, we explore the extent of the analogous behavior of the definitions above with respect to their usual counterparts.

\begin{remark}
Note that given a real simplicial space $X$, for every $n \ge 2$, the maps of genuine $C_2$-spaces as below left induce on underlying spaces the maps as below right.
$$\mathrm{map}_{\rsSpc}(\Delta^n, X) \to  \mathrm{map}_{\rsSpc}(\Lambda_n, X),\hspace{5mm} \mathrm{map}_{\caP(\Delta)}(\Delta^n, X^u) \to  \mathrm{map}_{\caP(\Delta)}(\Lambda_n, X^u).$$ 
Similarly, the map below left induces on underlying spaces the map below right.
$$\mathrm{map}_{\rsSpc}(\ast, X) \to \mathrm{map}_{\rsSpc}(\mathcal{J}, X), \hspace{5mm}\mathrm{map}_{\caP(\Delta)}(\ast, X^u) \to \mathrm{map}_{\caP(\Delta)}(\mathcal{J}, X^u).$$

Therefore, the underlying simplicial space of a (complete) real Segal space is a (complete) Segal space.
\end{remark}

We proceed to give a more familiar description of the Segal maps of \cref{Seg}, for which we use \cref{not:tilde}.

Let $X$ be a real simplicial object on $\caC$, where $\caC$ is a real $\infty$-category that admits finite limits and cotensors with $C_2$. Then for $n = 2k +j \geq 2 $, with $j=0,1$, we will define a map in $\caC$,
$$\rho\colon  X_n \to \widetilde{X_1 \times X_1} \times_{\widetilde{X_0 \times X_0}}\widetilde{X_1 \times X_1} \times_{\widetilde{X_0 \times X_0}} \cdots \times_{\widetilde{X_0 \times X_0}} \widetilde{X_1 \times X_1} \times_{\widetilde{X_0 \times X_0}} X_j.$$

In order to do so, we will separate in cases. Let us begin with $n\geq 2$ even. For every $1 \leq i \leq \frac{n}{2}$, we consider the maps $ [1] \to [n]$ in $\Delta$ sending $0 $ to $i-1$ and $1$ to $i$ that yield maps $[1] \coprod [1]^\op \to [n],$
which in turn induce compatible maps 
$$ X_n \simeq \map_{\rsSpc } ([n],X) \to \map_{  \rsSpc } ([1] \coprod [1]^\op, X) \simeq \widetilde{X_1 \times X_1}. $$ 

The map $[0] \to [n]$ in $\Delta^{hC_2}$ sending $0$ to 
$ \frac{n}{2}$ yields a map $\X_n \to \X_0.$ Note that here we are using $\Delta^{hC_2}$ because it's the underlying category of $\underline{\Delta}$, which is the source of the real simplicial space $X$.

All the above allow us to construct the following map in $\caC$
$$\rho_n\colon  X_n \to \widetilde{X_1 \times X_1} \times_{\widetilde{X_0 \times X_0}}\widetilde{X_1 \times X_1} \times_{\widetilde{X_0 \times X_0}} \cdots \times_{\widetilde{X_0 \times X_0}} \widetilde{X_1 \times X_1} \times_{\widetilde{X_0 \times X_0}} X_0,$$
where the pullbacks are formed
via the canonical maps $\widetilde{X_1 \times X_1}\to\widetilde{X_0 \times X_0}$ induced by alternating
evaluation at the source and target (starting with evaluation at the target, as in the usual case) and the canonical map $X_0\to\widetilde{X_0\times X_0}$ in the last one.

Similarly, given $n \geq 3 $ odd, for every $1 \leq i \leq \frac{n-1}{2}$ the maps $ [1] \to [n]$ sending $0 $ to $i-1$ and $1$ to $i$ in $\Delta$ induce compatible maps
$$ X_n \simeq \map_{ \rsSpc} ([n],X) \to \map_{  \rsSpc } ([1] \coprod [1]^\op, X) \simeq \widetilde{X_1 \times X_1}.$$

In this case, we consider the map $X_n\to X_1$, induced by the map $[1] \to [n]$ in $\Delta^{hC_2}$, sending $0$ to $ \frac{n-1}{2}$ and $1$ to $\frac{n-1}{2}+1$ in $\Delta^{hC_2}$. 

All the above allow us to construct the following map in $\caC$
$$\rho_n\colon  X_n \to \widetilde{X_1 \times X_1} \times_{\widetilde{X_0 \times X_0}}\widetilde{X_1 \times X_1} \times_{\widetilde{X_0 \times X_0}} \cdots \times_{\widetilde{X_0 \times X_0}} \widetilde{X_1 \times X_1} \times_{\widetilde{X_0 \times X_0}} X_1,$$

where the pullbacks are again formed via the canonical maps $\widetilde{X_1 \times X_1}\to\widetilde{X_0 \times X_0}$ induced by alternating evaluation at the source and target (starting with evaluation at the target) and the canonical map $X_1\to X_0\to\widetilde{X_0\times X_0}$ in the last one.

In the following lemma we show that, when $\caC=\Spc^{C_2}$, the maps constructed above effectively describe the $n$-th real Segal maps of real simplicial spaces.

\begin{lemma}\label{reali}
Let $X $ be a real simplicial space and $n \geq 2$. The $n$-th real Segal map of $X$, $\theta_n\colon \mathrm{map}_{\rsSpc}(\Delta^n, X) \to  \mathrm{map}_{\rsSpc}(\Lambda_n, X)$ admits the following description.

\begin{itemize}
\item
If $n$ is even, $\theta_n$ is equivalent to the map
$$\rho_n\colon  X_n \to \widetilde{X_1 \times X_1} \times_{\widetilde{X_0 \times X_0}}\widetilde{X_1 \times X_1} \times_{\widetilde{X_0 \times X_0}} \cdots \times_{\widetilde{X_0 \times X_0}} \widetilde{X_1 \times X_1} \times_{\widetilde{X_0 \times X_0}} X_0.$$

\item If $n$ is odd, $\theta_n$ is equivalent to the map
$$\rho_n\colon  X_n \to \widetilde{X_1 \times X_1} \times_{\widetilde{X_0 \times X_0}}\widetilde{X_1 \times X_1} \times_{\widetilde{X_0 \times X_0}} \cdots \times_{\widetilde{X_0 \times X_0}} \widetilde{X_1 \times X_1} \times_{\widetilde{X_0 \times X_0}} X_1.$$
\end{itemize}
\end{lemma}

\begin{proof}
The result follows directly from \cref{def:real_spine_inc} and \cref{Yoneda_lemma}.
\end{proof}

\cref{reali} shows that the real Segal maps of $X$ induce on underlying spaces the Segal maps of $X^u$. As announced, we give in the following corollary a recognition principle for real Segal spaces. 

\begin{corollary}\label{htt4sf}\label{cor:recognition_pple_rSeg}
A real simplicial space $X$ is a real Segal space if and only if the following conditions are satisfied:
	
\begin{enumerate}
\item\label{cor:recog_pple_item1}the underlying simplicial space of $X$, this is $X^{u}$ as in \cref{not:underlying_simpl_object}, is a Segal space,
		
\item\label{cor:recog_pple_item2} for every even $n \ge 2$ the natural map
		
$$X_n^{C_2} \to X_{\frac{n}{2}}^u \times_{X_0^u} X_0^{C_2}$$ is an equivalence, and
		
\item\label{cor:recog_pple_item3} for every odd $n \ge 3$ the natural map
		
$$X_n^{C_2} \to X_{\frac{n-1}{2}}^u \times_{X_0^u} X_1^{C_2}$$ is an equivalence.
\end{enumerate}
\end{corollary}

\begin{proof}
This follows directly from \cref{reali} using that $\widetilde{Z\times Z}^{C_2}\simeq Z^{u}$ for any genuine $C_2$-space $Z$.
\end{proof}

We introduce now a special class of real simplicial spaces that will be key in the next section.

\begin{definition}\label{fghjjnngl}\label{def:balanced}
We say that a real simplicial space $X$ is balanced if for every even $n \geq 0 $ 
the commutative square
\[
\begin{tikzcd}
X_n^{C_2} \ar[r] \ar[d]         & X_{n+1}^{C_2} \ar[d]\\
X_n^{hC_2} \ar[r] & X_{n+1}^{hC_2},
\end{tikzcd}
\]
induced by the unique surjective map $[n+1] \to [n]$ in $\Delta^{hC_2}$, is a pullback square.

\end{definition}

\begin{remark}
A real Segal space is balanced if and only if the commutative square of \cref{def:balanced} is a pullback square for the case $n=0$. Indeed, to see this we use \cref{htt4sf} to factor the map $$X_n^{C_2} \to X_n^{hC_2} \times_{X_{n+1}^{hC_2}} X_{n+1}^{C_2} $$
as the following composition of equivalences 
$$X_n^{C_2} \simeq X_{\frac{n}{2}}^u \times_{X_0^u} X_0^{C_2} \simeq X_{\frac{n}{2}}^u \times_{X_0^u} ( X_{0}^{hC_2} \times_{X_{1}^{hC_2 }}X_1^{C_2}) \simeq $$$$ X_{0}^{hC_2} \times_{X_{1}^{hC_2 }} ( X_{\frac{n}{2}}^u \times_{X_0^u} X_1^{C_2}) \simeq X_n^{hC_2} \times_{(X_{n}^{hC_2} \times_{X_{0}^{hC_2}}X_{1}^{hC_2 })} ( X_{\frac{n}{2}}^u \times_{X_0^u} X_1^{C_2}) $$$$ \simeq X_n^{hC_2} \times_{X_{n+1}^{hC_2}} X_{n+1}^{C_2}. $$

\end{remark}
\begin{remark}
If $X$ is a real Segal space, whose underlying Segal space $X^u$ is complete,
then $X$ is balanced if and only if the commutative square
\[
\begin{tikzcd}
X_0^{C_2} \ar[r] \ar[d]         & X_1^{C_2} \ar[d] \\
X_0^u \ar[r]                    & X_1^u
\end{tikzcd}
\]
induced by the map $[1] \to [0]$ is a pullback square. This follows from the pasting law for pullbacks as the commutative square
\[
\begin{tikzcd}
 X_0^{hC_2}\ar[r]\ar[d]         &X_1^{hC_2}\ar[d] \\
 X_0^u  \ar[r] &  X_1^u
\end{tikzcd}
\]
induced by the map $[1] \to [0] $ is a pullback square, since $X^u$ is a complete Segal space and therefore the map of spaces $X^u_0 \to X_1^u$.
\end{remark}

\section{\texorpdfstring{Complete real Segal spaces are $\infty$-categories with genuine duality}{Complete real Segal spaces are Infinity-categories with genuine duality}}\label{sec:GD_complete_rSegalspaces}

The main result of this section is \cref{equiv_crSegSpc_GD}, where we prove that $\infty$-categories with genuine duality (see \cref{def:infty_category_wGD}) can be interpreted as complete real Segal spaces (see \cref{def:complete_real_SegSpc}). Contrary to \cref{dhfvghjjk}, where we identified $\infty$-categories with duality as $\Spc[C_2]$-enriched functors $\underline{\Delta}^{op}\to \Spc[C_2]$ whose underlying functors are complete real Segal spaces, it is not enough to require the complete Segal condition to the underlying functor of $\Spc^{C_2}$-enriched functors $\underline{\Delta}^{op}\to \Spc^{C_2}$ to obtain $\infty$-categories with \emph{genuine} duality. This is due to the fact that the forgetful functor  $\Spc[C_2] \to \Spc$ is conservative whereas the forgetful functor $\Spc^{C_2} \to \Spc$ is not.

In order to prove \cref{equiv_crSegSpc_GD}, we part from and adjunction
\[
\begin{tikzcd}[column sep = huge]
	{\rsSpc} & {\caP(\Delta)^\gd}
	\arrow[""{name=0, anchor=east, inner sep=0}, "{\Phi}", bend left=30, from=1-1, to=1-2]
	\arrow[""{name=1, anchor=center, inner sep=0}, "{\Psi}", bend left=30, from=1-2, to=1-1]
	\arrow["\dashv"{anchor=center, rotate=-90}, draw=none, from=0, to=1]
\end{tikzcd}
\]
which we successively restrict until we obtain an equivalence of $\infty$-categories $\crSegSpc\simeq\Cat_\infty^\gd$. We describe the process below. 

\textbf{Step 1: (\cref{step1})} We construct the adjunction $\Psi\vdash\Phi$ (\cref{original_adj}), and show that it is actually a localization whose local objects are the balanced real simplicial spaces (\cref{loc_and_local_objects}); therefore, there exists an equivalence $\mathrm{b} \rsSpc\simeq\caP(\Delta)^\gd$ (see \cref{loc_and_local_objects}).

\textbf{Step 2:} This step consists of imposing a Segal-type condition on both sides of the equivalence obtained in Step 1, and show that those concepts correspond to each other via such equivalence giving us a new equivalence, between balanced real \emph{Segal} spaces and \emph{Segal} spaces with genuine duality, $\mathrm{b}\rSegSpc\simeq\SegSpc^\gd$ (\cref{prop:equiv_brSegSpc_GDseg}).

{\textbf{Step 3:}} Finally, we show that such equivalence also behaves well when imposing the completeness, and therefore obtain an equivalence $\crSegSpc\simeq \Cat_\infty^\gd$, between complete real Segal spaces and $\infty$-categories with genuine duality (\cref{equiv_crSegSpc_GD}).

We walk towards the proof of \cref{equiv_crSegSpc_GD}.

\subsection{Step 1}\label{step1}

As announced before, we first embark on proving the result below.
\begin{theorem}\label{Phi_Psi_loc} There exists a localization

\[
\begin{tikzcd}[column sep = huge]
	{\rsSpc} & {\caP(\Delta)^\gd}
	\arrow[""{name=0, anchor=east, inner sep=0}, "{\Phi}", bend left=30, from=1-1, to=1-2]
	\arrow[""{name=1, anchor=center, inner sep=0}, "{\Psi}", bend left=30, from=1-2, to=1-1]
	\arrow["\dashv"{anchor=center, rotate=-90}, draw=none, from=0, to=1]
\end{tikzcd}
\]

whose local objects are the balanced real simplicial spaces as defined in \cref{fghjjnngl}.
\end{theorem}

The proof of this theorem involves several steps, that we present as separate results. We being by constructing the adjunction $\Psi\vdash\Phi$ (\cref{fghjjl}), and then show that it is actually a localization whose local objects are the balanced real simplicial spaces (\cref{loc_and_local_objects}).

Before starting with the construction of the functors involved in the adjunction, we present the following lemma that will be necessary for it.

\begin{lemma}\label{dfhjkkkl}
	
The functor $\Tw\colon  \caP(\Delta)^{hC_2} \to \caP(\Delta)[C_2]$
factors, up to canonical identifications, as 
\[
\begin{tikzcd}[column sep=tiny, row sep=large]
\caP(\Delta)^{hC_2} \ar[rr, "\Tw"]\ar[rd]       &               &\caP(\Delta)[C_2]\\
\       &\Fun((\Delta^{hC_2})^\op, \Spc[C_2]),\ar[ru, "e^\ast"']      &
\end{tikzcd}
\]
where the map on the left is given by forgetting the enrichment, and $e^\ast$ is induced by the lift $e\colon \Delta\to \Delta^{hC_2}$ of the edgewise dubdivision map.
\end{lemma}

\begin{proof}
	
The functor $e\colon\Delta \to \Delta^{hC_2}$ is transpose to a $\Spc[C_2]$-enriched functor $e'\colon \Delta \to \underline{\Delta}$, where $\Delta$ has the trivial enrichment.
Since the map $e'$ is sent by $L\colon  \Cat_\infty^{\Spc[C_2]} \to \Cat_\infty[C_2] $ to the $C_2$-equivariant functor $e\colon  \Delta \to \Delta$, by naturality we have a commutative square below
\[
\begin{tikzcd}
\caP(\Delta)^{hC_2}\ar[r, "\simeq"] \ar[d, "\Tw"']    & \Fun_{\Spc[C_2]}(\underline{\Delta}^\op, \Spc[C_2])\ar[d, "e'^\ast"]\\
\caP(\Delta)[C_2] \ar[r, "\simeq"]                         & \Fun(\Delta^\op, \Spc[C_2]).
\end{tikzcd}
\]

Note that the right vertical functor admits the following factorization
$$ \Fun_{\Spc[C_2]}(\underline{\Delta}^\op, \Spc[C_2]) \to
\Fun((\Delta^{hC_2})^\op, \Spc[C_2]) \xrightarrow{e^\ast} \Fun(\Delta^\op, \Spc[C_2]).$$
\end{proof}

We are now in condition of constructing the functor $\Phi\colon\rsSpc\to \caP(\Delta)^\gd$. For this, we will use the universal property of $\caP(\Delta)^\gd$ as seen in \cref{def:GDpre} and therefore we are after compatible maps $\varphi_1\colon\rsSpc\to \caP(\Delta)^{hC_2}$ and $\varphi_2\colon\rsSpc\to \Fun([1],\caP(\Delta))$ . For the first one, we observe that the forgetful functor $\Spc^{C_2} \to \Spc[C_2]$ yields a real functor 
$$\rsSpc=\Fun_{\Spc^{C_2}}(\underline{\Delta}^\op, \Spc^{C_2})\to \Fun_{\Spc[C_2]}(\underline{\Delta}^\op, \Spc^{C_2}_{\mid C_2}).$$
Since $\Spc^{C_2}_{\mid C_2}\simeq \Spc[C_2]$, \cref{dhfvghjjk} gives a canonical map 
$$\rsSpc\xrightarrow{\varphi_1}\caP(\Delta)^{hC_2}.$$

The construction of the second map we are searching for is a bit more involved. Firstly, we note that the functor $(-)^{C_2}\colon  \Spc^{C_2} \to \Spc$ yields a functor $$ \rsSpc \to \Fun((\Delta^{hC_2})^\op, \Spc^{C_2}).$$

Secondly, we recall that the edgewise subdivision $e\colon  \Delta \to  \Delta^{hC_2}$, that sends  $[n] \mapsto [n] \ast [n]^\op \simeq [2n+1] $, induces a functor
$$ \Fun((\Delta^{hC_2})^\op, \Spc^{C_2}) \to \Fun(\Delta^\op, \Spc^{C_2}).$$

Thirdly, the functor $\Spc^{C_2} \to \Fun([1], \Spc) $, which sends a genuine $C_2$-space $X$ to its genuine refinement $X^{C_2} \to X^{hC_2}$, gives rise to a functor
$$ \Fun(\Delta^\op, \Spc^{C_2}) \to \Fun(\Delta^\op,  \Fun([1], \Spc)) \simeq  \Fun([1], \caP(\Delta)).$$

Finally, by composing these three functors we obtain the map
$$\varphi_2\colon  \rsSpc \to \Fun([1], \caP(\Delta)).$$

By \cref{dfhjkkkl}, we know that there is a commutative diagram as below

\begin{equation}\label{diag:Phi_GDpre}
\begin{tikzcd}
\rsSpc\ar[rdd, bend right=30, "\varphi_2"']\ar[rrd, bend left=30, "\varphi_1"]\ar[rd, dashed, "\Phi"]     &       &\\
\       &\caP(\Delta)^\gd\ar[r]\ar[d]     &\caP(\Delta)^{hC_2}\ar[d, "\caH"]\\
\       &\Fun([1], \caP(\Delta))\ar[r, "\mathrm{ev}_1"]      &\caP(\Delta).
\end{tikzcd}
\end{equation}

Using that $\caP(\Delta)^\gd$ was defined as a pullback, we obtain a functor $\Phi\colon  \rsSpc \to \caP(\Delta)^\gd$ as depicted above. Explicitly, the functor $\Phi\colon  \rsSpc \to \caP(\Delta)^\gd$ sends a real simplicial space
$X$ to the simplicial space with genuine duality
consisting of the simplicial space with duality
$X^u$ and the genuine refinement 
$$([n] \mapsto X^{C_2}_{[n] \ast [n]^\op }) \to \caH(X^u) \simeq ([n] \mapsto X^{hC_2}_{[n] \ast [n]^\op }). $$

\begin{proposition}\label{fghjjl}\label{original_adj}The functor $\Phi\colon\rsSpc\to \caP(\Delta)^\gd$ is part of an adjunction as below 

\[
\begin{tikzcd}[column sep = huge]
	{\rsSpc} & {\caP(\Delta)^\gd.}
	\arrow[""{name=0, anchor=east, inner sep=0}, "{\Phi}", bend left=30, from=1-1, to=1-2]
	\arrow[""{name=1, anchor=center, inner sep=0}, "{\Psi}", bend left=30, from=1-2, to=1-1]
	\arrow["\dashv"{anchor=center, rotate=-90}, draw=none, from=0, to=1]
\end{tikzcd}
\]
\end{proposition}

\begin{proof}
We will see that the functor $\Phi\colon\rsSpc\to\caP(\Delta)^\gd$ preserves small colimits and so admits a right adjoint $\Psi.$ 

Let $\caB$ be a small $\infty$-category and $\phi\colon \caB \to \caP(\Delta)^\gd$ a functor with components $\alpha\colon \caB \to \caP(\Delta)^{hC_2}$ and $ \caB \to \Fun([1], \caP(\Delta) ) $ corresponding to a natural transformation
$\beta \to \caH \circ \alpha $ of functors $ \caB \to \caP(\Delta) $.
The colimit of $\phi$ is given by $$(\colim(\alpha),  \colim(\beta) \to \colim(\caH \circ \alpha) \to \caH(\colim(\alpha)) ). $$ In particular, the projection
$\caP(\Delta)^\gd \to \caP(\Delta)^{hC_2}$ and the composite $\caP(\Delta)^\gd \to \caP(\Delta)$ of projecting to $\Fun([1], \caP(\Delta))$ followed by evaluation at $0$ preserve small colimits.

The proof concludes by observing that $\varphi_1$ and $\varphi_2$ preserve small colimits.
\end{proof}

With the following remarks we show that the adjunction of \cref{original_adj} actually is more structured. 

\begin{remark}\label{rsSpc_has_SpcC2_action}
The $\infty$-category of real simplicial spaces, $\rsSpc$, admits a canonical closed left $\Spc^{C_2}$-action induced by the finite product preserving real diagonal embedding $\Spc^{C_2} \to  \rsSpc.$ This left $\Spc{C_2}$-action on $\rsSpc$ retrieves the canonical $\Spc^{C_2}$-enrichment on $\rsSpc$ coming from seeing it as the real $\infty$-category of real functors $\Fun(\underline{\Delta}^{op}, \Spc^{C_2})$.
\end{remark}

\begin{remark}\label{GDpre_has_SpcC2_action}
By \cref{GD(pre)_cartesian_closed}, the $\infty$-category $\caP(\Delta)^\gd$ is cartesian closed and therefore it carries a canonical closed left action over itself. Restricting this left action along the canonical left and right adjoint embedding $\Spc^{C_2} \subset \Cat_\infty^\gd\subset \caP(\Delta)^\gd $, we get a closed left action of $\caP(\Delta)^\gd$ over 
$\Spc^{C_2}.$
\end{remark}

\begin{remark}When considering the left $\Spc^{C_2}$-actions of \cref{rsSpc_has_SpcC2_action,GDpre_has_SpcC2_action}, the functor $\Phi\colon\rsSpc\to\caP(\Delta)^\gd$ is $\Spc^{C_2}$-linear. This follows directly by observing that the triangle
\[
\begin{tikzcd}[column sep=tiny]
\rsSpc\ar[rr, "\Phi"]       &       &\caP(\Delta)^\gd\\
\       &\Spc^{C_2}\ar[ul, hook, "\delta"]\ar[ur, hook, "\iota"']
\end{tikzcd}
\]
is a commutative diagram of finite products preserving functors, where $\delta\colon \Spc^{C_2}\to \rsSpc$ is the real diagonal embedding and $\iota\colon\Spc^{C_2}\to \caP(\Delta)^\gd$ is induced by seeing genuine $C_2$-spaces as spaces with genuine duality as discussed after \cref{not:fixed_points}.

Since $\Phi$ is $\Spc^{C_2}$-linear, it is $\Spc^{C_2}$-enriched and it preserves tensors with genuine $C_2$-spaces; see \cref{Venr:equiv_left_modules_enriched_cats}. This, together with the fact that the functor $\Phi$ preserves small colimits (as seen in the proof of \cref{fghjjl}), guarantees, by \cref{Venr:F_admits_right_adjoint} that it has a real right adjoint.

In this framework, we can give an explicit description of the real right adjoint $\Psi$ as the composite
$$ \Psi\colon  \caP(\Delta)^\gd \subset \Fun_{ \Spc^{C_2} }( (\caP(\Delta)^\gd)^\op, \Spc^{C_2}) 
\to \rsSpc = \Fun_{ \Spc^{C_2} }( \underline{\Delta}^\op, \Spc^{C_2}),$$
where the right-most map is the Yoneda embedding of \cref{Venr:def:Yoneda_embedding} followed by the restriction along the composition $ \underline{\Delta} \subset \rsSpc \xrightarrow{\Phi} \caP(\Delta)^\gd.$
\end{remark}

The following proposition, together with \cref{original_adj}, completes the proof of \cref{Phi_Psi_loc}.

\begin{proposition}\label{loc_and_local_objects}\label{ffzjnbfg}
The adjunction

\[
\begin{tikzcd}[column sep = huge]
	{\rsSpc} & {\caP(\Delta)^\gd}
	\arrow[""{name=0, anchor=east, inner sep=0}, "{\Phi}", bend left=30, from=1-1, to=1-2]
	\arrow[""{name=1, anchor=center, inner sep=0}, "{\Psi}", bend left=30, from=1-2, to=1-1]
	\arrow["\dashv"{anchor=center, rotate=-90}, draw=none, from=0, to=1]
\end{tikzcd}
\]
is a localization with local objects the balanced real simplicial spaces. Consequently, it restricts to an equivalence $$\mathrm{b} \rsSpc \simeq \caP(\Delta)^\gd. $$

\end{proposition}

In order to prove this proposition, we will use two lemmas whose precise statements and proofs we provide immediately after this proof, but we present here for readability purposes. Namely, 
\begin{enumerate}
\item\label{re4zwej_item} The right adjoint $\Psi$ is conservative and sends object of $ \caP(\Delta)^\gd $ to balanced real simplicial spaces, as seen in \cref{re4zwej} and \cref{rmk:Psi_conservative}.
\item\label{unit_item} A real simplicial space $X$ is balanced if and only if the unit $X \to \Psi\Phi(X)$ is an equivalence, as shown in \cref{unit}.
\end{enumerate}

\begin{proof}[Proof of \cref{loc_and_local_objects}.]
With these lemmas in mind, we observe that for any $Y \in \caP(\Delta)^\gd$, its image $\Psi(Y)$ is balanced by \cref{re4zwej_item}; so the unit component $\Psi(Y)\to \Psi\Phi\Psi(Y)$ is an equivalence by \cref{unit_item}. Therefore, we know by the triangular identities that $\Psi$ sends the counit component $\Phi\Psi(Y) \to Y$ to an equivalence. Given that, by \cref{re4zwej_item}, the functor $\Psi$ is conservative, we deduce that the counit component $\Phi\Psi(Y) \to Y$ must be an equivalence itself. This, together with \cref{unit_item}, concludes the proof.
\end{proof}

We now proceed, as promised, to prove the results stated without proof so far and used for the proof of \cref{loc_and_local_objects}. In order to see that $\Psi$ is conservative, we describe levelwise the image of the functor $\Psi\colon\caP(\Delta)^\gd\to\rsSpc$ for a given $X\in \caP(Delta)$ in \cref{lem:desc_Psi_underlying,lem:desc_Psi_fixed_pts}. Precisely, we know that $\Psi(X)_n$ is a genuine $C_2$-space
 for every $n$, and therefore it is given by its underlying $C_2$-space (which is described in \cref{lem:desc_Psi_underlying}) and its genuine refinement (which we describe in \cref{lem:desc_Psi_fixed_pts}).

The following lemma will be helpful when showing that the functor $\Psi\colon\caP(\Delta)^\gd\to\rsSpc$ is conservative. First, we recall that the forgetful functors $$\caP(\Delta)^\gd \to \caP(\Delta)^{hC_2} \hspace{0.5em}\text{ and }\hspace{0.5em} \rsSpc \to \Fun_{\Spc^{C_2}}(\underline{\Delta}^\op, \Spc[C_2])$$ admit left and right adjoints. Indeed, by \cref{rmk:left_adj_GDpre_to_PDelta} and \cref{rmk:canonical_gd_Cat} the forgetful functor $\caP(\Delta)^\gd \to \caP(\Delta)^{hC_2}$ admits a left (resp. right) adjoint sending a simplicial space $X$ with duality to $(X, \emptyset \to \caH(X))$ (resp. $(X, \id\colon \caH(X) \to \caH(X))$). Also, the real forgetful functor $\Spc^{C_2} \to \Spc[C_2]$ admits a real left (right) adjoint equipping a $C_2$-space $X$ with the genuine refinement $\emptyset \to X^{hC_2}$ respectively the identity of $X^{hC_2}$. This implies that the induced real functor
$$\rsSpc= \Fun_{\Spc^{C_2}}(\underline{\Delta}^\op, \Spc^{C_2}) \to \Fun_{\Spc^{C_2}}(\underline{\Delta}^\op, \Spc[C_2])\simeq \Fun_{\Spc[C_2]}(\underline{\Delta}^\op, \Spc[C_2])$$
admits a real left (right) adjoint with the natural description.

\begin{lemma}\label{lem:desc_Psi_underlying}\label{unit_counit_underlying}
Consider the adjunction

\[
\begin{tikzcd}[column sep = huge]
	{\rsSpc} & {\caP(\Delta)^\gd.}
	\arrow[""{name=0, anchor=east, inner sep=0}, "{\Phi}", bend left=30, from=1-1, to=1-2]
	\arrow[""{name=1, anchor=center, inner sep=0}, "{\Psi}", bend left=30, from=1-2, to=1-1]
	\arrow["\dashv"{anchor=center, rotate=-90}, draw=none, from=0, to=1]
\end{tikzcd}
\]

Then the following hold:
\begin{enumerate}
\item The unit and counit component of any object are sent to an equivalence in the $\infty$-categories
$\Fun_{\Spc[C_2]}(\underline{\Delta}^\op, \Spc[C_2]) \simeq \caP(\Delta)^{hC_2}$ via the corresponding forgetful functors as depicted below. 

\[
\begin{tikzcd}[row sep=huge, column sep=small]
\rsSpc\rar["\Phi", bend left=20]\rar[phantom, "\bot"]\ar[d] & \caP(\Delta)^\gd\lar["\Psi", bend left=20]\ar[d]\\
\Fun_{\Spc[C_2]}(\underline{\Delta}^\op, \Spc[C_2])\rar["\simeq", bend left]\rar[phantom, "\bot"] & \caP(\Delta)^{hC_2}\lar["\simeq" pos=0.55, bend left]
\end{tikzcd}
\]

\item There exists a commutative diagram as below
\[
\begin{tikzcd}[row sep=huge, column sep=small]
\Fun_{\Spc[C_2]}(\underline{\Delta}^\op, \Spc[C_2])\rar[bend left, "\simeq"]\rar[phantom, "\bot"]\ar[d, hook] & \caP(\Delta)^{hC_2}\lar[bend left, "\simeq"]\ar[d, hook]\\
\rsSpc\rar["\Phi", bend left=20]\rar[phantom, "\bot"] & \caP(\Delta)^\gd,\lar["\Psi", bend left=20]\
\end{tikzcd}
\]
where we embed 
$\Fun_{\Spc[C_2]}(\underline{\Delta}^\op, \Spc[C_2])$ into $\rsSpc$ and
$\caP(\Delta)^{hC_2}$ into $\caP(\Delta)^\gd$ via the left (resp. right) adjoints of the forgetful functors.
\end{enumerate}
\end{lemma}

\begin{proof}We will deal with both statements simultaneously. By construction the left adjoint $\Phi\colon\rsSpc \to \caP(\Delta)^\gd$
forgets to the canonical equivalence
$$\Fun_{\Spc[C_2]}(\underline{\Delta}^\op, \Spc[C_2]) \simeq \caP(\Delta)^{hC_2}.$$
Moreover, again by construction, the functor $\Phi$ also restricts to the canonical equivalence
$$\Fun_{\Spc[C_2]}(\underline{\Delta}^\op, \Spc[C_2]) \simeq \caP(\Delta)^{hC_2},$$
where we embed $\Fun_{\Spc[C_2]}(\underline{\Delta}^\op, \Spc[C_2])$ into $\rsSpc$ and
$\caP(\Delta)^{hC_2}$ into $\caP(\Delta)^\gd$ as described in item (2).
By adjointness, the functor $\Psi\colon\caP(\Delta)^\gd \to \rsSpc$ forgets to such a canonical equivalence as well, and restricts to the equivalence of full subcategories 
$$\Fun_{\Spc[C_2]}(\underline{\Delta}^\op, \Spc[C_2]) \simeq \caP(\Delta)^{hC_2},$$
where we now embed 
$\Fun_{\Spc[C_2]}(\underline{\Delta}^\op, \Spc[C_2])$ into $\rsSpc$ and
$\caP(\Delta)^{hC_2}$ into $\caP(\Delta)^\gd$
via the right adjoint (both times) of the forgetful functors.
Then \cref{htwdzfd,htedud} implies that $\Psi\colon\caP(\Delta)^\gd \to \rsSpc$
restricts to the full subcategory 
$$\Fun_{\Spc[C_2]}(\underline{\Delta}^\op, \Spc[C_2]) \simeq \caP(\Delta)^{hC_2},$$
where we embed 
$\Fun_{\Spc[C_2]}(\underline{\Delta}^\op, \Spc[C_2])$ into $\rsSpc$ and
$\caP(\Delta)^{hC_2}$ into $\caP(\Delta)^\gd$
via the left adjoint of the forgetful functors.
\end{proof}

Now we analyze the space $\Psi(X)^{C_2}_n$. 

\begin{remark}\label{rmk:C2_fixed_pts_of_Psi}Note that, by adjointness, the fixed points $\Psi(X)_n^{C_2}$ can be expressed in terms of mapping out of $\Phi([n])$. Indeed, for every $X \in \caP(\Delta)^\gd $ and $[n] \in \underline{\Delta} $ 
we have a natural equivalence
$$ \Psi(X)_n^{C_2} \simeq \map_{\rsSpc}([n], \Psi(X)) \simeq \map_{\caP(\Delta)^\gd}(\Phi([n]), X),$$
where the first equivalence is given by the enriched Yoneda lemma \cref{Yoneda_lemma}.
\end{remark}

The following lemma makes use of \cref{rmk:C2_fixed_pts_of_Psi} to compute $\Psi(X)_n^{C_2}$.

\begin{lemma}\label{re4zwej}\label{lem:desc_Psi_fixed_pts}
	
Let $X= (\caC, \phi\colon H \to \caH(\caC))$ be an object of $\caP(\Delta)^\gd$, i.e., a simplicial space with genuine duality. Then:

\begin{enumerate}
\item\label{lem:desc_Psi_fixed_pts_item1} For odd $n \geq 0$ there exists a canonical equivalence over $ \caH(\caC)_{ \frac{n-1}{2} }\coloneqq\caC_n^{hC_2}$ as below 
$$ \Psi(X)_n^{C_2} \simeq H_{\frac{n-1}{2}}.$$

\item\label{lem:desc_Psi_fixed_pts_item2} For even $n \geq 0$ there exists a canonical equivalence  over $ \caC_n^{hC_2}$ as below
$$\Psi(X)_n^{C_2} \simeq H_{\frac{n}{2}} \times_{\caC_{ n+1}^{hC_2}} \caC_n^{hC_2},$$
where the pullback is formed over the maps 
$\phi_{\frac{n}{2}}\colon H_{\frac{n}{2}} \to \caH(\caC)_{\frac{n}{2}} \coloneqq\caC_{ n+1}^{hC_2}$ and $ \caC_{ n}^{hC_2} \to \caC_{n+1}^{hC_2} $
induced by the map $\caC_n \to \caC_{ n+1}$ of $C_2$-spaces, induced in turn by the unique surjective duality preserving map $[n+1] \to [n]$. 
\end{enumerate}

In particular, the real simplicial space $\Psi(X)$ is balanced.
\end{lemma}

\begin{proof}
In order to prove \cref{lem:desc_Psi_fixed_pts_item1}, we consider the following pullback square induced by the pullback definition of $\caP(\Delta)^\gd$.

\[
\begin{tikzcd}
\underbrace{\map_{\caP(\Delta)^\gd} (\Phi([n]), X)}_{\simeq\Psi(X)_n^{C_2}  } \ar[r] \ar[d] & \underbrace{\map_{ \Fun_{\Spc[C_2]}(\underline{\Delta}^\op, \Spc[C_2])   }([n], \caC)}_{ \simeq \caC_n^{hC_2} } \ar[d, "\simeq"] \\
\underbrace{\map_{\Fun([1], \caP(\Delta)) }(\id_{[\frac{n-1}{2}]}, \phi)}_{\simeq H_{ \frac{n-1}{2} }  }  \ar[r] & \underbrace{\map_{\caP(\Delta) }([\frac{n-1}{2}], \caH(\caC))}_{ \simeq \caH(\caC)_{ \frac{n-1}{2}}}. 
\end{tikzcd}
\]
Since the right vertical map is an equivalence and pullbacks are stable under equivalence, this concludes the proof of the statement.

For \cref{lem:desc_Psi_fixed_pts_item2}, we consider the pullback square

\[
\begin{tikzcd}
\underbrace{\map_{\caP(\Delta)^\gd} (\Phi([n]), X)}_{ \simeq\Psi(X)_n^{C_2}  } \ar[r] \ar[d] & \underbrace{\map_{\Fun_{\Spc[C_2]}(\underline{\Delta}^\op, \Spc[C_2])}([n], \caC)}_{ \simeq \caC_n^{hC_2} } \ar[d] \\
\underbrace{\map_{ \Fun([1], \caP(\Delta)) }(\id_{[\frac{n}{2}] }, \phi)}_{\simeq H_{ \frac{n}{2} }  }  \ar[r] & \underbrace{\map_{\caP(\Delta) }([\frac{n}{2}], \caH(\caC))}_{ \simeq \caH(\caC)_{ \frac{n}{2} } \coloneqq \caC_{n+1}^{hC_2}} .
\end{tikzcd}
\]
\end{proof}

\begin{remark}\label{rmk:Psi_conservative}\cref{lem:desc_Psi_underlying,lem:desc_Psi_fixed_pts} imply that the functor $\Psi\colon\caP(\Delta)^\gd\to\rsSpc$ is conservative. Indeed, consider the two functors $\alpha,\beta\colon  \caP(\Delta)^\gd \to \caP(\Delta)$ that send a real simplicial space with genuine duality $X= (\caC, \phi\colon H \to \caH(\caC))$ to $\caC$ and $H$ respectively, which result in the conservative functor $(\alpha,\beta)\colon  \caP(\Delta)^\gd \to \caP(\Delta) \times \caP(\Delta)$. With this, the conservativity of $\Psi$ follows from the existence of the following commutative squares for $n \geq 0$
\[
\begin{tikzcd}[column sep=large]
\caP(\Delta)^\gd \ar[r, "\Psi"] \ar[d, "\alpha" ]  &\rsSpc \ar[d, "(-)_n " ] \\
\caP(\Delta)  \ar[r, "(-)_n"']     & 
\Spc\\
\end{tikzcd}
\text{ and }
\begin{tikzcd}[column sep=large]
\caP(\Delta)^\gd \ar[r, "\Psi"] \ar[d, "\beta" ]  &\rsSpc \ar[d, "(-)_{2n+1}^{C_2} " ] \\
\caP(\Delta)  \ar[r, "(-)_n"']     & 
\Spc,\\
\end{tikzcd}
\]
which is guarateed by \cref{lem:desc_Psi_underlying,lem:desc_Psi_fixed_pts} respectively. 
\end{remark}

The second fact that we stated without proof in order to prove \cref{loc_and_local_objects} is that the component unit $\eta_X\colon C\to\Psi\Phi (X)$ of the adjunction involved is a levelwise equivalence precisely when $X$ is balanced. This is the content of the following lemma.

\begin{lemma}\label{unit}Consider the adjunction 
\[
\begin{tikzcd}[column sep = huge]
	{\rsSpc} & {\caP(\Delta)^\gd.}
	\arrow[""{name=0, anchor=east, inner sep=0}, "{\Phi}", bend left=30, from=1-1, to=1-2]
	\arrow[""{name=1, anchor=center, inner sep=0}, "{\Psi}", bend left=30, from=1-2, to=1-1]
	\arrow["\dashv"{anchor=center, rotate=-90}, draw=none, from=0, to=1]
\end{tikzcd}
\]
For any given real simplicial space $X$, the unit component $\eta_X\colon X\to \Psi\Phi(X)$ is an equivalence on all odd levels. Moreover, $\eta_X$ is also an equivalence on even levels if and only if $X$ is balanced. 
\end{lemma}

\begin{proof}
By \cref{unit_counit_underlying} the unit component $\eta_X\colon X\to \Psi\Phi(X)$ induces an equivalence on underlying simplicial spaces. Consequently to show that the unit component $\eta_X\colon X\to \Psi\Phi(X)$ is a levelwise equivalence, it is enough to show that for any $n$ the following composite $\rho$ is an equivalence 
$$\rho\colon \map_\rsSpc([n],\X) \to \map_{\rsSpc}([n],\Psi(\Phi(\X))\simeq\map_{\caP(\Delta)^\gd}(\Phi([n]),\Phi(\X)).$$

For $n$ odd, the diagram (\ref{diag:Phi_GDpre}) yields a diagram on mapping spaces as below, where all vertical maps are equivalences.
\[
\begin{tikzcd}
\underbrace{\map_{\rsSpc} ([n], X)}_{\simeq X_n^{C_2}  } \ar[rdd, bend right=20,"\simeq"']\ar[rrd, bend left=20]\ar[rd, "\rho"]     &       &\\
   &\map_{\caP(\Delta)^\gd}(\Phi([n]),\Phi(X)) \ar[r]\ar[d, "\simeq"]     &\underbrace{\map_{\caP(\Delta)^{hC_2}}([n], X^u)}_{ \simeq X_n^{hC_2} } \ar[d, "\simeq"] \\
\ & \underbrace{\map_{\Fun([1], \caP(\Delta)) }(\id_{[\frac{n-1}{2}]}, \phi) }_{\simeq X_n^{C_2}  }  \ar[r]   &\underbrace{\map_{\caP(\Delta) }([\frac{n-1}{2}], \caH(X^u))}_{ \simeq X_n^{hC_2} }.
\end{tikzcd}
\]
This shows that the unit component $\eta_X$ is an equivalence when evaluated at $[n]$ odd.

For $n$ even, we can consider a similar pullback diagram, where the $[\frac{n-1}{2}$ of the bottom line is replaced by $[\frac{n}{2}]$. From it, we deduce that the $C_2$-fixed points of the unit component $\eta_X\to \Psi(\Phi(X))$ evaluated at $[n]$ is given by the map 
$$X_n^{C_2}  \to X_{n+1}^{C_2}  \times_{X_{n+1}^{hC_2}} X_n^{hC_2}.$$
Therefore, this is an equivalence precisely when $X$ is a balanced real simplicial space. 

\end{proof}

Now that we have given all the components of the proof of \cref{loc_and_local_objects}, we are ready to further restrict the equivalence $\mathrm{b} \rsSpc\simeq \caP(\Delta)^\gd$ pursuing our aim of seeing $\infty$-categories with genuine duality as complete real Segal spaces.

\subsection{Step 2}

We devote this part to show that the adjoint equivalence of $\mathrm{b} \rsSpc\simeq\caP(\Delta)^\gd$ of \cref{loc_and_local_objects} restricts to full subcategories verifying certain corresponding Segal conditions. For real simplicial spaces, we made such ``Segal condition'' precise in \cref{Seg}. For simplicial spaces with genuine duality we mean the following.

\begin{definition}
We call a simplicial space with genuine duality $(\caC, \phi)$ a Segal space with genuine duality if $\caC$ is a Segal space and $\phi\colon  H \to \caH(\caC)$ is a right fibration of simplicial spaces.
\end{definition}

\begin{definition}\label{def:GDseg}
We define the $\infty$-category $\SegSpc^\gd$ of Segal spaces with genuine duality by the pullback square
\[
\begin{tikzcd}
\SegSpc^\gd\ar[r]\ar[d]            &\SegSpc^{hC_2}\ar[d, "\caH"]\\
\mathfrak{R}\ar[r, "\mathrm{ev}_1"']          &\SegSpc,
\end{tikzcd}
\]
where $\mathfrak{R} \subset \Fun([1], \SegSpc)$ is the full reflexive subcategory spanned by the right fibrations and $\ev_1\colon \Fun([1], \SegSpc) \to \SegSpc, \ \Fun([1], \Cat_\infty) \to \Cat_\infty$ the evaluation at the target functor. 
\end{definition}

It is worth noticing the similitude between this definition and \cref{def:GD,def:GDpre}; the $\infty$-category of Segal spaces with genuine duality is in fact a pivot concept between those two, much as Segal spaces are a concept between simplicial spaces and $\infty$-categories. The pullback square of \cref{def:GDseg} also lives in the $\infty$-category of presentable $\infty$-categories and right adjoint functors and so $\SegSpc^\gd$ is presentable.

\begin{proposition}\label{rifib}\label{prop:equiv_brSegSpc_GDseg}\label{ggbbhm}
The equivalence $\mathrm{b} \rsSpc\simeq\caP(\Delta)^\gd$ of \cref{loc_and_local_objects} restricts to an equivalence $$\brSegSpc\simeq\SegSpc^\gd. $$
\end{proposition}

\begin{proof}
We first show that for a real Segal space $X$, the simplicial space with genuine duality $\Phi(X)$ is in fact a Segal space with genuine duality. Let us name $\Phi(X)$ by $(\caC, \phi\colon H\to \caH(\caC))$. We then need to verify that the composition $H \xrightarrow{\phi} \caH(\caC) \to \caC$ is a right fibration. To obtain this, it is enough to observe that for every $n\geq 0$, the natural map 
$$H_n \to \caC_n \times_{\caC_0} H_0$$
induced by $\{n\} \subset [n]$ is $$X([n]*[n]^\op)^{C_2} \to X_n^u \times_{X_0^u} X_1^{C_2},$$ which is an equivalence by item 3 of \cref{cor:recognition_pple_rSeg}.

Take now a Segal space with genuine duality $Y=(\caC, \phi\colon  H \to \caH(\caC))$. We wish to show that the real simplicial space $\Psi(Y)$ is a real Segal space. We use the recognition principle of \cref{cor:recognition_pple_rSeg}. As $\caC$ is a Segal space, we proceed to show that \cref{cor:recog_pple_item2,cor:recog_pple_item3} of \cref{cor:recognition_pple_rSeg} hold.

By \cref{lem:desc_Psi_fixed_pts} the map in \cref{cor:recog_pple_item2}, this is, when $n$ is even, is given by 
$$\theta\colon  H_{\frac{n}{2}} \times_{\caC_{n+1}^{hC_2}} \caC_n^{hC_2} \to \caC_{\frac{n}{2}} \times_{\caC_0} H_0 \times_{\caC_1^{hC_2}} \caC_0^{hC_2}.$$

Note that the embedding $$\caP(\Delta)^{hC_2} \simeq \Fun_{\Spc[C_2]}(\underline{\Delta}^\op, \Spc[C_2]) \subset \rsSpc $$
sends a simplicial space with duality, whose underlying simplicial space is a Segal space, to a real Segal space. So $\caC$ may be viewed as a real Segal space with $\caC_n^{C_2} = \caC_n^{hC_2}$ for any $n \geq 0$. Therefore, there exist canonical equivalences 
$$\caC_{n+1}^{hC_2} \simeq \caC_{\frac{n}{2}} \times_{\caC_0} \caC_1^{hC_2} \text{ and }\ \ \caC_n^{hC_2} \simeq \caC_{\frac{n}{2}} \times_{\caC_0} \caC_0^{hC_2}.$$

Using them, we obtain that the map $\theta$ is equivalent to the map below, which is an equivalence because $H \to \caC$ is a right fibration.
$$H_{\frac{n}{2}} \times_{\caC_1^{hC_2} } \caC_0^{hC_2} \to \caC_{\frac{n}{2}} \times_{\caC_0} H_0 \times_{\caC_1^{hC_2}} \caC_0^{hC_2}$$

Again by \cref{lem:desc_Psi_fixed_pts}, when $n$ is odd the map in \cref{cor:recog_pple_item3} of \cref{cor:recognition_pple_rSeg} is
$$ H_{\frac{n-1}{2}} \to \caC_{\frac{n-1}{2}} \times_{\caC_0} H_0,$$
which is an equivalence since $H \to \caC$ is a right fibration.
\end{proof}

\subsection{Step 3}\label{step3} In this last step, we make a further restriction, obtaining the main theorem of this section. We first state such result, and present the main components of the proof as a separate result that we prove immediately afterwards.

\begin{theorem}\label{ggbbh}\label{equiv_crSegSpc_GD}The equivalence $\brSegSpc \simeq \SegSpc^\gd $ of \cref{prop:equiv_brSegSpc_GDseg} restricts to an equivalence $\crSegSpc \simeq \Cat_\infty^\gd$.
\end{theorem}
\begin{proof}
Follows directly from \cref{X_complete_iff_X_balanced}.
\end{proof}

\begin{proposition}
\label{artru6}
Let $X$ be a real Segal space whose underlying Segal space is complete. Then the commutative square below is a pullback square.
\begin{equation}\label{fghbnv}
\begin{tikzcd}
\map_{\rSegSpc}(\mathcal{caJ}, X)  \ar[r] \ar[d]        &\map_{\rSegSpc}(\Delta^1, X) \simeq X_1^{C_2} \ar[d] \\ \map_{\caP(\Delta)} (\mathcal{caJ}, X^u) \ar[r]         &\map_{\caP(\Delta)}(\Delta^1, X^u) \simeq X_1^u
\end{tikzcd}
\end{equation}

In other words, the mapping space $\map_{\rSegSpc}(\mathcal{caJ}, X) $ is the full subspace of $ X_1^{C_2}$ consisting of those objects whose image in $X_1^u$ is an equivalence.
\end{proposition}

\begin{proof}
By pasting of pullback squares, the commutative square (\ref{fghbnv}) is a pullback if and only if the commutative below square is so.
\begin{equation}\label{fghjh}
\begin{tikzcd}
\map_{\rSegSpc}(\mathcal{caJ}, X)  \ar[r] \ar[d] &  \map_{\rSegSpc}(\Delta^1, X) \simeq X_1^{C_2} \ar[d] \\ \map_{\caP(\Delta)} (\mathcal{caJ}, X^u)^{hC_2} \ar[r]  & \map_{\caP(\Delta)}(\Delta^1, X^u)^{hC_2} \simeq X_1^{hC_2}
\end{tikzcd}
\end{equation}

The commutative square (\ref{fghjh}), in turn, is a pullback square if and only if the square obtained from it by replacing $X$ with $\Psi(\Phi(X))$ is a pullback square. We begin by justifying this. 

In order to see this, since the unit $X \to \Psi(\Phi(X))$ yields an equivalence on underlying spaces (see \cref{unit_counit_underlying}), it is enough to prove that it also induces equivalences 
$$\map_{\rSegSpc}(\Delta^1, X) \simeq \map_{\rSegSpc}(\Delta^1, \Psi\Phi(X)) \text{ and}$$
$$\map_{\rSegSpc}(\mathcal{J}, X) \simeq \map_{\rSegSpc}(\mathcal{J}, \Psi\Phi(X)).$$

We see this by arguing that the full subcategory $\caW \subset \rSegSpc$ spanned by all $Z$ such that the induced map $$\map_{\rSegSpc}(Z, X) \to \map_{\rSegSpc}(Z, \Psi(\Phi(X))) $$ is an equivalence, contains both $\Delta^1$ and $\mathcal{caJ}.$ The case of $\Delta^1$ follows from \cref{unit} that states that $\caW$ contains $\Delta^n$ for $n \geq 1$ odd. Since the unit $X \to \Psi(\Phi(X))$ yields an equivalence on underlying spaces, the subcategory $\caW$ contains the free real simplicial spaces $C_2 \times \Delta^n \simeq  \Delta^n \coprod (\Delta^n)^\op $ for all $n \geq 0$. Furthermore,  $\caW$ is closed under small colimits. Therefore, as $\mathcal{caJ} $ can be built by gluing cells of the form $\Delta^n$ for $n \geq 1$ odd, and $C_2 \times \Delta^n $ for all $n \geq 0$, we have $\mathcal{caJ} \in \caW.$

Having replaced $X$ by $\Psi(\Phi(X))$ in square (\ref{fghjh}), by adjointness of $\Phi$ and $\Psi$, the new square is equivalent to the following commutative square 
\begin{equation}\label{fghj}
\begin{tikzcd}
\map_{\caP(\Delta)^\gd}(\mathcal{caJ}, \Phi(X))  \ar[r] \ar[d] &  \map_{\caP(\Delta)^\gd}(\Delta^1, \Phi(X)) \ar[d] \\ \map_{\caP(\Delta)} (\mathcal{caJ}, X^u)^{hC_2} \ar[r]  & \map_{\caP(\Delta)}(\Delta^1, X^u)^{hC_2}.
\end{tikzcd}
\end{equation}

Square (\ref{fghj}) is, in turn, the base change of the commutative square 
\begin{equation}\label{fgkhlj}
\begin{tikzcd}
\map_{\caP(\Delta)}(\caH(\mathcal{caJ}), H)  \ar[r] \ar[d] &  \map_{\caP(\Delta)}(\caH(\Delta^1), H) \ar[d] \\ \map_{\caP(\Delta)} (\caH(\mathcal{caJ}), \H(X)) \ar[r]  & \map_{\caP(\Delta)}(\caH(\Delta^1), \H(X)).
\end{tikzcd}
\end{equation}

By \cref{laxherm}, the functor $\caH(\Delta^1) \to \caH(\mathcal{caJ})$
is equivalent to both of the functors $ [0] \to \mathcal{J} $ and is thus a local equivalence with respect to complete Segal spaces. Now, by \cref{completeright} and \cref{rifib}, since $X^u$ is a complete Segal space, so are $\caH(X)$ and $H$. In consequence, both horizontal maps of square (\ref{fgkhlj}) are equivalences, which concludes the proof.
\end{proof}

\begin{corollary}\label{X_complete_iff_X_balanced} Let $X$ be a real Segal space whose underlying Segal space is complete. Then $X$ is complete (as real space) if and only if $X$ is balanced.
\end{corollary}

\begin{proof}
By \cref{artru6}, the commutative square below is a pullback square.
\[
\begin{tikzcd}
\map_{\rSegSpc}(\mathcal{caJ}, X)  \ar[r] \ar[d] &  \map_{\rSegSpc}(\Delta^1, X) \simeq X_1^{C_2} \ar[d] \\ \map_{\caP(\Delta)} (\mathcal{caJ}, X^u) \ar[r]  & \map_{\caP(\Delta)}(\Delta^1, X^u) \simeq X_1^u
\end{tikzcd}
\]

By pasting of pullback squares, the commutative square
\[
\begin{tikzcd}
X_0^{C_2} \ar[r] \ar[d] & \map_{\rSegSpc}(\mathcal{caJ}, X)   \ar[d] \\ 
X_0^u  \ar[r] & \map_{\caP(\Delta)}( \mathcal{caJ}, X^u)
\end{tikzcd}
\]
induced by the map $\mathcal{J} \to \Delta^0 $
is a pullback square if and only if the commutative square 
\[
\begin{tikzcd}
X_0^{C_2} \ar[r] \ar[d] & \map_{\rSegSpc}( \Delta^1, X) \simeq X_1^{C_2}  \ar[d] \\ 
X_0^u  \ar[r] & \map_{\caP(\Delta)}( \Delta^1, X^u) \simeq X_1^u.
\end{tikzcd}
\]
\end{proof}

\part{\texorpdfstring{Waldhausen $\infty$-categories with genuine duality}{Waldhausen Infinity-categories with genuine duality}}\label{part:2}

\section{\texorpdfstring{(Pre)additive and stable $\infty$-categories with genuine duality}{(Pre)additive and stable Infinity-categories with genuine duality}}\label{sec:genuine_add_preadd_stable}

In this section we develop a theory of preadditive, additive, and stable $\infty$-categories with genuine duality. In this genuine (i.e.\ enriched) context, as in the classical one, we fully expect the theory of stable $\infty$-categories with genuine duality to be rich and well-behaved. As evidence of this, we dedicate \cref{sec:quadratic_functors} to show an equivalence between stable $\infty$-categories with genuine duality and quadratic functors on a stable $\infty$-category.

\subsection{Preadditive $\infty$-categories with genuine duality}\label{subsec:preadd_with_gd}
Let $\caC$ be an $\infty$-category with duality that is preadditive. The cartesian monoidal structure on $\caC$ yields a symmetric monoidal structure on $\caH(\caC)$ such that $\caH(\caC) \to \caC $ is a symmetric monoidal right fibration, which we define to be a symmetric monoidal functor whose underlying functor is a right fibration.

Note that the tensor unit of $\caH(\caC)$ is the initial object, as it lies over the zero object of $\caC$ and $\caH(\caC) \to \caC$ is a right fibration whose fiber over the zero object of $\caC$ is contractible. As the tensorunit of $\H(\caC)$ is initial, every functor $\H(\caC)^\op \to \Spc $ admits a unique lift to an oplax symmetric monoidal functor (see \cite[Proposition 2.4.3.9]{lurie.higheralgebra}). We are now in conditions to define preadditive $\infty$-category with genuine duality.

\begin{definition}\label{preadd}
Let $(\caC, \phi\colon  H \to \caH(\caC))$ be an $\infty$-category with genuine duality. We call $(\caC, \phi)$ a preadditive $\infty$-category with genuine duality if
\begin{enumerate}
\item $\caC$ is preadditive.
\item The unique oplax symmetric monoidal functor
lifting the functor $$\H(\caC)^\op \to \Spc$$
classified by $\phi$ is symmetric monoidal.
\end{enumerate}
\end{definition}

We denote by $$\gdpreadd \subset \Cat_\infty^\gd$$ the subcategory
whose objects are the preadditive $\infty$-categories with genuine duality and whose morphisms are the maps of $\infty$-categories with genuine duality, whose underlying functor preserves finite products. 

\begin{remark}\label{fibers_HtoC_are_Einfty}
Given an $\infty$-category with duality $\caC$ such that $\caC$ is preadditive, by taking the cartesian monoidal structure on $\caC$, we obtain that $\caH(\caC) \to \caC $ is a symmetric monoidal right fibration so that its fibers admit canonical structures of $\mathrm{E}_\infty$-spaces.

If we actually have a preadditive $\infty$-category with genuine duality $(\caC,\phi)$, by \cref{preadd} we know that the functor $\phi\colon H\to \H(\caC)$ classifies a symmetric monoidal functor $\H(\caC)^\op \to \Spc $, and it is therefore a symmetric monoidal right fibration. From this we conclude that the composite $$H \to \caH(\caC) \to \caC$$
is a symmetric monoidal right fibration. Therefore, it classifies a lax symmetric monoidal functor $\caC^\op \to \Spc$.  Since $\caC$ is preadditive, this corresponds to a functor $\caC^\op \to \mathrm{CMon}(\Spc)$ (see \cite[Theorem 2.4.3.18]{lurie.higheralgebra}), in particular, the fibers of $H \to \caC $ admit canonical structures of $\mathrm{E}_\infty$-spaces.
\end{remark}

Recall that the $\infty$-category $\Pre\Add$ of preadditive $\infty$-categories is preadditive (as opposed to fact that $\Add$ is not additive). Thus the $\infty$-category $\Pre\Add^{hC_2}$ of preadditive $\infty$-categories with duality is preadditive.

We now prepare the ground to show that $\gdpreadd$ is a preadditive $\infty$-category, which means that the forgetful functor $ \Cmon(\gdpreadd) \to \gdpreadd$ is an equivalence so that every preadditive $\infty$-category with genuine duality $(\caC, \phi)$ has the canonical structure of a commutative monoid in $ \gdpreadd $ corresponding to a canonical structure of a symmetric monoidal functor on $ \phi\colon  H \to \caH(\caC)$, where $\caC$ carries the cartesian structure. We begin with the following result.

\begin{lemma}\label{54fkow} Let $\caC$ be a symmetric monoidal $\infty$-category whose tensor unit is initial. Then the forgetful functor $$\{\caC \} \times_{\Cmon(\Cat_\infty)} \Cmon(\caR) \to \{\caC \} \times_{ \Cat_\infty } \caR$$
becomes an equivalence after pulling back to the full subcategory of 
$ \{\caC \} \times_{ \Cat_\infty } \caR$ spanned by the right fibrations $H \to \caC$ such that $H(1)$ is contractible and for every $A,B \in \caC$ the canonical map $H(A \otimes B) \to H(A) \times H(B)$ is an equivalence.
\end{lemma}

\begin{proof} For notational simplicity, we prove the dual statement. Namely, that given $\caC$ a symmetric monoidal $\infty$-category whose tensor unit is final, the forgetful functor $$\Cmon(\Cat_{\infty})_{ / \caC} \to \Cat_{\infty/\caC}$$
becomes an equivalence after pulling back to the full subcategory of 
$\Cat_{\infty/\caC}$ spanned by the left fibrations $H \to \caC$ such that $H(1)$ is contractible and for every $A,B \in \caC$ the canonical map $H(A \otimes B) \to H(A) \times H(B)$ is an equivalence.

By \cite[Proposition 2.4.3.9]{lurie.higheralgebra}, the forgetful functor $$ \Fun^{\otimes, \lax}(\caC^\op, \Spc^\op) \to \Fun(\caC^\op, \Spc^\op) $$ is an equivalence, since the tensor unit of $\caC^\op$ is initital and $\Spc^\op$ is a cocartesian symmetric monoidal $\infty$-category. So we get an equivalence 
$$ \Fun^{\otimes, \lax}(\caC^\op, \Spc^\op)^\op \simeq \Fun(\caC^\op, \Spc^\op)^\op \simeq \Fun(\caC, \Spc), $$ 
under which $ \Fun^{\otimes}(\caC^\op, \Spc^\op)^\op \simeq \Fun^{\otimes}(\caC, \Spc) $ corresponds to the full subcategory of $ \Fun(\caC, \Spc) $ spanned by the functors classified by a left fibration over $\caC$ such that for every $A,B \in \caC$ the canonical map $H(A \otimes B) \to H(A) \times H(B)$ is an equivalence and $H(1) $ is contractible.

Note that there is a canonical equivalence between $\Fun^{\otimes, \lax}(\caC, \Spc) $ and the full subcategory of $\Cmon(\Cat_{\infty})_{ / \caC}$ spanned by the symmetric monoidal functors over $\caC$, whose underlying functor is a left fibration. Under this equivalence $\Fun^{\otimes}(\caC, \Spc) $ corresponds to those symmetric monoidal functors over $\caC$, whose underlying functor $H \to \caC$ is a left fibration such that for every $A,B \in \caC$ the canonical map $H(A \otimes B) \to H(A) \times H(B)$ is an equivalence and $H(1) $ is contractible.
\end{proof}

\begin{proposition}\label{fdghjkkkl}\label{gdpreadd_is_preadditive} The $\infty$-category $\gdpreadd$ is preadditive.
\end{proposition}

\begin{proof}
There exists a commutative square
\[
\begin{tikzcd}
\Cmon(\gdpreadd)\ar[r,"\theta"] \ar[d]      &\gdpreadd \ar[d]\\
\Cmon(\Pre\Add^{hC_2})\ar[r, "\simeq"]     &\Pre\Add^{hC_2}
\end{tikzcd}
\]
which induces on the fiber over any 
$\caC \in \Cmon(\Pre\Add^{hC_2})$ the functor $\theta_\caC$ resulting from pulling back the forgetful functor
$$\{\caH(\caC) \} \times_{\Cmon(\Cat_\infty)} \Cmon(\caR) \to \{\caH(\caC) \} \times_{ \Cat_\infty } \caR$$ to the full subcategory spanned by the right fibrations $H \to \caH(\caC)$ such that for every $A,B \in \caH(\caC)$ the canonical map $H(A \otimes B) \to H(A) \times H(B)$ is an equivalence and $H(1) $ is contractible. By \cref{54fkow}, the functors $\theta_\caC$ are equivalences, and therefore so is the forgetful functor
$$ \Cmon(\gdpreadd) \to \gdpreadd.$$
\end{proof}

With \cref{gdpreadd_is_preadditive} in hand, we will often see $ \gdpreadd \simeq \Cmon(\gdpreadd )$ as a full subcategory of $\Cmon(\Cat_{\infty }^\gd  )$
and identify a preadditive $\infty$-category with genuine duality
with its canonical commutative monoid in $ \Cat_{\infty}^\gd.$

\subsection{Additive $\infty$-categories with genuine duality}\label{subsec:additive_with_gd}
Let $(\caC,\phi)$ be a preadditive $\infty$-category with genuine duality. If $\caC$ is additive, the canonical $E_\infty$-structures on the fibers of $ \caH(\caC) \to \caC$ are grouplike. Indeed, for any $\X \in \caC$
the space $\map_\caC(\X,\X^\dual)$ carries a canonical structure of a grouplike $E_\infty$-space, and therefore so does the fiber $\caH(\caC)_\X = \map_\caC(\X,\X^\dual)^{hC_2}$. In this subsection, we introduce the notion of additive $\infty$-category with genuine duality in such a way that we have the same for the $E_\infty$-structures on the fibers of the functor $H \to \caC$ (see \cref{fibers_HtoC_are_Einfty}). 

We begin by introducing the following concept.
\begin{definition} Let $\caC,\caD$ be symmetric monoidal $\infty$-categories. A symmetric monoidal right fibration $q\colon  \caD \to \caC$ is called an additive right fibration if $\caC$ is additive and carries the cartesian monoidal structure, and $\caD \to \caC$ classifies a functor $\caC^\op \to \Sp_{\geq 0}.$
\end{definition}

Now we can define the main objects of this subsection.

\begin{definition}
We call a preadditive $\infty$-category with genuine duality $(\caC, \phi)$ 
an additive $\infty$-category with genuine duality if $ H \to \caC$
is an additive right fibration.
\end{definition}

We denote by $$\gdadd\subset \gdpreadd$$ 
the full subcategory spanned by the additive $\infty$-categories with genuine duality.

\begin{remark}\label{rmk:Add_gd_preadditive}
Note that $\gdadd$ is closed under small limits in $\gdpreadd$, and it is therefore preadditive, too.
\end{remark}

In the following lemma, that will be of use later, we describe the extra condition that a preadditive $\infty$-category $(\caC,\phi)$ whose underlying $\infty$-category $\caC$ is additive must verify for it to be an \emph{additive $\infty$-category}. 

\begin{lemma}\label{dghjkkkjg}
Let $(\caC, \phi)$ be a preadditive $\infty$-category with genuine duality such that $\caC$ is additive. Let us consider $X \in \caC$ and $Y,Z \in \caH(\caC)_X$, and denote $YZ \to Y \otimes Z$ the lift in $ \caH(\caC) $ of the diagonal map $X \to  X \oplus X$ in $\caC$.

Then, $(\caC, \phi)$ is an additive $\infty$-category with genuine duality if and only if for every $X \in \caC$ and $Y,Z \in \caH(\caC)_X$
the maps $YZ \to Y \otimes Z$ and $Y \simeq Y \otimes 1 \to Y \otimes Z$
yield an equivalence $\rho\colon  H_{ Y \otimes Z} \to H_Y \times H_{YZ}.$

\end{lemma}

\begin{proof}We know that the pair $(\caC, \phi)$ is an additive $\infty$-category with genuine duality if and only if for every $X \in \caC$ the canonical $E_\infty$-space structure on the fiber $H_X$ is grouplike.	By definition, the $E_\infty$-space $H_X$ is grouplike if and only if the shear map $\alpha\colon  H_X \times H_X \to H_X \times H_X $ is an equivalence, where the shear map is the projection to the first factor on the first factor and the multiplication of $H_X$ on the second factor.
	
As the canonical $E_\infty$-space structure on $\caH(\caC)_X$ is grouplike, the shear map $\beta$ for $\caH(\caC)_X$ is an equivalence.

There exists a commutative square
\[
\begin{tikzcd}
 H_X \times H_X\ar[r, "\alpha"] \ar[d] &  H_X \times H_X \ar[d] \\
\caH(\caC)_X \times \caH(\caC)_X \ar[r, "\beta"]  & \caH(\caC)_X \times \caH(\caC)_X.
\end{tikzcd}
\]

Therefore, in order to prove that $\alpha$ is an equivalence, it is enough to check that for every $Y,Z \in \caH(\caC)_X$ the induced map on the fiber $\theta\colon  H_Y \times H_Z \to H_Y \times H_{YZ}$ is an equivalence.
But $\theta$ factors as $H_Y \times H_Z \simeq H_{Y \otimes Z} \xrightarrow{\rho} H_Y \times H_{YZ}$, as the multiplication in the fiber $H_X$ of any $X\in\caC$ factors as $H_X \times H_X \to H_{X \oplus X} \to H_X$, induced by the tensor product of $H$ and the diagonal $X \to  X \oplus X$ in $\caC$.
\end{proof}

\subsection{Stable $\infty$-categories with genuine duality}\label{subsec:stable_with_gd}

Let $(\caC,\phi)$ be an additive $\infty$-category with genuine duality.
If $\caC$ is stable, the canonical $E_\infty$-structures on the fibers of 
$ \caH(\caC) \to \caC$ are $\Omega^\infty (A)$ for some (not necessarily connective) spectrum $A$. In this subsection, we introduce the notion of stable $\infty$-category with genuine duality in such a way that we have the same for the $E_\infty$-structures on the fibers of the functor $H \to \caC$ (see \cref{fibers_HtoC_are_Einfty}).

Let $(\caC, \phi)$ be an $\infty$-category with genuine duality such that $\caC$ is stable. As $\caC$ admits pullbacks, also $\H(\caC)$ does and the right fibration $\H(\caC) \to \caC$ preserves and detects pullback squares. In consequence, given objects $A,B \in \H(\caC)$ there exists a pullback square in $\caH(\caC)$ as below left, lying over any pushout (and therefore also pullback) square in $\caC$ as below right.
\begin{equation*}\label{fhjlkgd}
\begin{tikzcd}
1 \ar[r] \ar[d] & A \otimes 1 \ar[d]        &       &0 \ar[r] \ar[d] & A \ar[d]\\
1 \otimes B \ar[r]  & A \otimes B           &       &B \ar[r]  & A \oplus B
\end{tikzcd}
\end{equation*}

Condition 2 of \cref{preadd} says that $\phi$ classifies a reduced functor $\alpha\colon  \H(\caC)^\op \to \Spc $ that sends all squares as above left to pullback squares. If $\caC$ is stable, it is natural to ask that $\alpha\colon  \H(\caC)^\op \to \Spc $ sends not only those pullback squares but all pullback squares in $\H(\caC)$ to pullback squares. In other words, we would be asking $\alpha\colon  \H(\caC)^\op \to \Spc $ to send pushout squares in $\H(\caC)^\op$ to pullback squares, i.e. $\alpha$ to be an excisive functor.

The discussion above motivates the following definition.

\begin{definition}\label{def:stable_gd}We call an additive $\infty$-category with genuine duality $(\caC, \phi)$ a stable $\infty$-category with genuine duality if
\begin{enumerate}
\item $\caC$ is stable.
\item $\phi$ classifies an excisive functor $\H(\caC)^\op \to \Spc.$
\end{enumerate}
\end{definition}

We will identify in \cref{sec:quadratic_functors} these as quadratic functors.

So an $\infty$-category with genuine duality $(\caC, \phi)$ 
is a stable $\infty$-category with genuine duality if and only if
$\caC$ is stable, $\phi$ classifies a reduced and excisive functor $\H(\caC)^\op \to \Spc $
and the canonical $E_\infty$-structures on the fibers of $H \to \caC $ are grouplike.

\begin{remark}
Note that $\phi$ classifies an excisive functor $\H(\caC)^\op \to \Spc $
if and only if for every pullback square $[1]^2\coloneqq [1] \times [1] \to \caC$
the pullback of $ [1]^2 \times_\caC H \to [1]^2 \times_\caC  \caH(\caC)$ along any section of $[1]^2 \times_\caC  \caH(\caC) \to [1]^2$
classifies a pullback square of spaces.
\end{remark}
 
We denote by $$ \St^\gd\subset \St \times_{\Add} \Add^\gd $$ the full subcategory spanned by the stable $\infty$-categories with genuine duality.

\begin{remark}\label{rmk:St_gd_preadditive}
Note that the $\infty$-category $\St^\gd $ is closed under small limits in $\St \times_{\Add} \Add^\gd$, and thus $\St^\gd $ is preadditive. This give us an inclusion $\St^\gd \simeq \Cmon(\St^\gd ) \subset \Cmon(\Cat_{\infty}^\gd).$
\end{remark}

\subsection{Structural results about (pre)additive $\infty$-categories with genuine duality}\label{subsec:structural_results_of_additive_infty_categores_w_gd}
In this subsection we prove the following three structural results about (pre)additive and stable $\infty$-categories with genuine duality, which we will apply later to our theory of Waldhausen $\infty$-categories with genuine duality:

\begin{itemize}
\item $\Pre\Add^\gd,\ \Add^\gd$ and $\St^\gd$ are presentable (\cref{Prese});
\item $\Pre\Add^\gd,\ \Add^\gd$ and $\St^\gd$ are genuine preadditive (see \cref{fghhjjkmnl}) and so in particular preadditive (\cref{fghhgkkkn});

\item $\Pre\Add^\gd,\ \Add^\gd$ and $\St^\gd$ admit canonical closed symmetric monoidal structures (\cref{Symmon}).
\end{itemize}

\subsection*{Presentability}
In this subsection we prove that the $\infty$-categories $\Pre\Add^\gd, \Add^\gd$ and $\St^\gd$ are presentable. 

We note that the $\infty$-categories $\Pre\Add, \Add$ and $\St$ are presentable. A proof of this for $\Pre\Add$ can be found in \cite[Proposition 3.6]{hls}; the cases of $\Add$ and $\St$ admit a similar proof. From this, we can deduce that the fiber product $\Theta^{hC_2}\times_{\Cat_\infty^{hC_2}}\Cat_\infty^\gd$ is a presentable $\infty$-category for $\Theta \in \{\Pre\Add,\Add,\St\}$ using that the canonical inclusion $\Theta \subset \Cat_\infty$ admits a left adjoint.

The above observation, together with the next proposition, imply the presentability of $\Pre\Add^\gd, \Add^\gd$ and $\St^\gd$.

\begin{proposition}\label{Prese}Let $\Theta\in \{\Pre\Add,\Add,\St \}$. Then $\Theta^\gd$ is an accessible localization of $\Theta^{hC_2}\times_{\Cat_\infty^{hC_2}}\Cat_\infty^\gd$ and any local equivalence is sent to an equivalence by the projection $\Theta^{hC_2}\times_{\Cat_\infty^{hC_2}}\Cat_\infty^\gd \to \Theta^{hC_2}$.
\end{proposition}
\begin{proof}

Since the functor $\caH\colon  \Cat_\infty^{hC_2} \to  \Cat_\infty$ is accessible, we deduce that the embedding $\Theta^\gd \hookrightarrow \Theta^{hC_2} \times_{(\Cat_{\infty})^{hC_2}} \Cat_{\infty}^\gd$ is also accessible.

We know that the forgetful functor $\Cat_{\infty}^\gd \to \Cat_{\infty}^{hC_2} $ is a cartesian fibration by \cref{rmk:ff_gd_to_hC2_is_cartesian}. Therefore so is the base change $\Theta^{hC_2} \times_{(\Cat_{\infty})^{hC_2}} \Cat_{\infty}^\gd \to \Theta^{hC_2}$. Moreover, the full subcategory $\Theta^\gd$ is stable under the fiber transports of the cartesian fibration $\Theta^{hC_2} \times_{(\Cat_{\infty})^{hC_2}} \Cat_{\infty}^\gd \to \Theta^{hC_2}$.

Then, by \cite[Proposition 7.3.2.6]{lurie.higheralgebra}, in order to prove the statement it is enough to show that for every $\caC \in \Theta^{hC_2}$ the fiber $\{\caC\} \times_{ \Theta^{hC_2}} \Theta^\gd$ is a localization of $ \{\caH(\caC)\} \times_{ \Cat_\infty} \caR$, where $\caR$ is the $\infty$-category of right fibrations. We devote the rest of the proof to this.

Denote by $\caK_\Theta$ the set of commutative squares in $\caH(\caC)$ of the following form:
\begin{itemize}
\item For $\Theta \in \{ \Pre\Add, \Add\},$ for every $A, B \in \caH(\caC) $ we consider the canonical commutative square in $\caH(\caC)$ as below.
\[
\begin{tikzcd}
1 \ar[r] \ar[d] & A \otimes 1 \ar[d] \\
1 \otimes B \ar[r]  & A \otimes B 
\end{tikzcd}
\]

\item For $\Theta= \St$, we consider the collection of all pullback squares in $\caH(\caC).$
\end{itemize}

Now, we consider $$\caB^\Theta \subset \{\caH(\caC)\} \times_{ \Cat_\infty} \caR \simeq \Fun(\caH(\caC)^\op, \Spc)$$ the full subcategory spanned by the reduced functors $\caH(\caC)^\op \to \Spc$ that send all squares of $\caK_\Theta$ to pullback squares. Then, the embedding $\caB^\Theta \hookrightarrow \{\caH(\caC)\} \times_{ \Cat_\infty} \caR$ is an accessible localization, as $\caK_\Theta$ is a set. Especially, $\caB^\Theta$ is presentable.

The proof now trifurcates. When $\Theta=\Pre\Add$, for every $\caC$ in $\Pre\Add^{hC_2}$ we have $\{\caC\} \times_{\Pre\Add^{hC_2} } \Pre\Add^\gd = \caB^{\Pre\Add}$ as full subcategories of $\{\caH(\caC)\} \times_{ \Cat_\infty} \caR$, which allows us to conclude the first part.

We proceed to discuss the case $\Theta=\Add$. We will prove that for every $\caC$ in $\Add^{hC_2}$ the full subcategory $\{\caC\} \times_{\Add^{hC_2}} \Add^\gd  \subset  \caB^\Add \subset \{\caH(\caC)\} \times_{ \Cat_\infty} \caR $ is an accessible localization of $\caB^\Add $ (and so also of 
$ \{\caH(\caC)\} \times_{ \Cat_\infty} \caR $). 

Denote $\caP \to \caB^\Add$ the base change map of the forgetful functor $$\{\caH(\caC) \}  \times_{\Cmon(\Cat_\infty)} \Cmon(\caR) \to \{\caH(\caC) \} \times_{ \Cat_\infty} \caR$$ along the embedding $\caB^\Add \subset \{\caH(\caC) \} \times_{ \Cat_\infty} \caR$. Then, by \cref{54fkow} the projection $\caP \to \caB^\Add$ is an equivalence.

For every $X \in \caC$ we have a right adjoint functor
$$ \caB^\add \simeq \caP \subset \{\caH(\caC)\} \times_{ \Cmon(\Cat_\infty)} \Cmon(\caR) \simeq $$$$ \{\caH(\caC)^\op\} \times_{ \Cmon(\Cat_\infty)} \Cmon(\caL) \to  \Cmon(\Spc)_{/ \caH(\caC)_X^\op}$$ induced by pulling back along the functor $\caH(\caC)_X \to \caH(\caC)$ and the equivalence $\caR \simeq \caL$ that takes the opposite $\infty$-category, where $\caL$ is the $\infty$-category of left fibrations. This functor sends $(\caC, H \to \caH(\caC)) \mapsto H_X$. Since $\caB^\Add$ is presentable, these maps induce a right adjoint functor $$\Psi\colon \caB^\Add \to \mathrm{\prod}_{X \in \caC} \Cmon(\Spc)_{/ \caH(\caC)_X^\op}.$$

We conclude this case by observing that the fiber $(\Add^\gd)_\caC \subset  \caB^\add$ is the pullback of $\Psi$ along the accessible localization of presentable $\infty$-categories as below $$\mathrm{\prod}_{X \in \caC} \Grp_{E_\infty}(\Spc)_{/ \caH(\caC)_X^\op} \hookrightarrow \mathrm{\prod}_{X \in \caC} \Cmon(\Spc)_{/ \caH(\caC)_X^\op}.$$

Finally, when $\Theta=\St$, it is enough to observe that the full subcategory $$\{\caC\} \times_{\St^{hC_2}} \St^\gd \subset \{\caH(\caC)\} \times_{ \Cat_\infty} \caR$$ is the intersection of $\caB^\St $ and $\{\caC\} \times_{\Add^{hC_2}} \Add^\gd$.
\end{proof}

\subsection*{Genuine preadditivity} We introduce the notion of \emph{genuine preadditivity}, which is the $\Spc^{C_2}$-enriched version of usual preadditivity. Then we show that the categories $\Pre\Add^\gd,$ $\Add^\gd,$ and $\St^\gd$ are genuine preadditive $\infty$-categories.

One of the conditions for this real notion of preadditivity is a comparison between tensors and cotensors with $C_2$. In order to state this, we present the following construction. 

Let $\caC$ be a real $\infty$-category that admits cotensors with $C_2$ and such that $\caC$ has finite products and a zero object. For any $X, Y \in \caC$, the map $X \to X \times X$ in $\caC^u$ that is the identity on the first factor and the zero map on the second factor corresponds to a map $\xi\colon C_2 \to \caC(X,X^{C_2}) $ of genuine $C_2$-spaces under the equivalence $$\Spc^{C_2}(C_2,\caC(X,X^{C_2})) \simeq  \map_{\caC^u}(X, X \times X). $$ 

The map $\xi$ exhibits $X^{C_2}$ as the tensor of $C_2$ with $X$ if for any $Y \in \caC$ the induced map of genuine $C_2$-spaces below is an equivalence.
$$ \theta\colon \caC(X^{C_2},Y) \to \Spc[C_2](C_2, \map_{\caC}(X,Y))$$

If $\caC$ admits tensors with $C_2$, the map $\xi$ corresponds to a map $C_2 \otimes X \to \X^{C_2}$ in $\caC$.

\begin{definition}\label{fghhjjkmnl}\label{def:genuine_preadd}
We call a real $\infty$-category $\caC$ genuine preadditive if 
\begin{enumerate}
\item the underlying $\infty$-category of $\caC$ is preadditive (we mean underlying in the sense of the first item below \cref{terminology:real}),
\item $\caC$ admits a zero object, finite products and finite coproducts (in the $\Spc^{C_2}$-enriched sense) and tensors and cotensors with $C_2$, and
\item for every $X \in \caC$, the natural map $C_2 \otimes X \to X^{C_2}$ is an equivalence.
\end{enumerate}
\end{definition}

Note that the underlying $\infty$-category of $\caC$ and $\caC^u$ do not mean the same.

\begin{remark}\label{rmk:genuine_preadd_char}
Let $\caC$ be a real $\infty$-category that
admits cotensors with genuine $C_2$-spaces. Then $\caC$ is genuine preadditive if and only if the following conditions hold.
\begin{enumerate}
\item The underlying $\infty$-category of
$\caC$ is preadditive.
\item The canonical map
$$ \map_\caC(X^{C_2},Y) \to \map_{\Spc^{C_2}}(C_2, \map_{\caC}(X,Y)) \simeq \map_{\caC^u}(X,Y) $$
is an equivalence.
\end{enumerate}
\end{remark}

\begin{notation} Given objects $A, B \in \caC$ we have a canonical equivalence $A \coprod B \simeq A \times B$ in $\caC$ and write $A \oplus B $ for $A \coprod B \simeq A \times B$.
\end{notation}

\begin{notation}Given an object $X \in \caC$ we have a canonical equivalence $C_2 \otimes X \simeq X^{C_2}$ in $\caC$ and write $\widetilde{X \oplus X} $ for $ C_2 \otimes X \simeq X^{C_2}$; see \cref{not:tilde}.
\end{notation}

\begin{remark}
Let $\caC$ be a real $\infty$-category that admits tensors with $C_2$. Then for every $X, Y \in \caC$ there exists a canonical equivalence of genuine $C_2$-spaces
\begin{equation*}\label{eqq}
\caC(C_2 \otimes X ,Y) \simeq \Spc^{C_2}(C_2, \caC(X,Y)) 
\end{equation*}
that yields on $C_2$-fixed points an equivalence 
\begin{equation*}\label{eqqq}
\map_\caC(C_2 \otimes X ,Y) \simeq \map_{\Spc^{C_2}}(C_2, \caC(X,Y)) 
\simeq \map_{\caC^u}(X,Y)
\end{equation*}
since $C_2$ is the free genuine $C_2$-space on a point, and on underlying spaces an equivalence 
$$ \map_{\caC^u}(C_2 \otimes X ,Y) \simeq \map_{\Spc}(C_2, \map_{\caC^u}(X,Y)) \simeq \map_{\caC^u}(X,Y) \times \map_{\caC^u}(X,Y)$$
that exhibits $C_2 \otimes X $ as the coproduct of $X$ with $X$ in $\caC^u$.

Dually, there is a canonical equivalence 
$$ \map_{\caC^u}(Y,X^{C_2}) \simeq \map_{\Spc}(C_2, \map_{\caC^u}(Y,X)) \simeq \map_{\caC^u}(Y,X) \times \map_{\caC^u}(Y,X)$$
that exhibits $X^{C_2} $ as the product of $X$ with $X$ in $\caC^u$. When $\caC$ is genuine preadditive, this can be seen as a justification for the notation $\widetilde{X \oplus X}$.

\end{remark}

\begin{proposition}\label{genpre}
Let $\caC$ be a $\Spc[C_2]$-enriched $\infty$-category. Then $\caC$ is genuine preadditive if and only if
$\caC^u$ is preadditive and $\caC$ has cotensors with $C_2$.
\end{proposition}

\begin{proof}
Since $\caC^u$ is preadditive, also the real $\infty$-category $\caC$ has a zero object and finite (co)products in the enriched sense, which are especially the zero object and finite (co)product in the underlying $\infty$-category of $\caC$. 

Moreover, for any $A,B \in \caC$
the canonical map $\caC(A\times B,C) \to \caC(A,C) \times \caC(B,C)$
of $C_2$-spaces induces on underlying spaces the canonical equivalence $$\map_{\caC^u}(A\times B,C) \to \map_{\caC^u}(A,C) \times \map_{\caC^u}(B,C)$$ and is thus itself an equivalence. This implies that the underlying $\infty$-category of $\caC$ is preadditive.

The map $\theta$ from the introduction of this subsection induces on underlying spaces the map $$\rho\colon  \map_{\caC^u}(X^{C_2},Y) \to \map_{\Spc}(C_2, \map_{\caC^u}(X,Y)) \simeq \map_{\caC^u}(X,Y) \times \map_{\caC^u}(X,Y).$$
If $\caC^u$ is preadditive, $X^{C_2}$ 
satisfies the universal property of the sum $X \oplus X $ in $\caC^u$ so that the map $\rho$ is an equivalence. Hence also $\theta$ is an equivalence and in consequence $\xi$ exhibits $X^{C_2}$ as the tensor of $C_2$ with $X$. Therefore, $\caC$ is genuine preadditive. 
\end{proof}

\begin{corollary}\label{cor:hC2_gives_genuine_preadditive}
Let $\caC$ be a preadditive $\infty$-category with a $C_2$-action. Them the real $\infty$-category $\caC^{hC_2}$ is genuine preadditive.
\end{corollary}

\begin{proof}
The $\infty$-category $\caC^{hC_2}$ is canonically $\Spc[C_2]$-enriched.
So, by \cref{genpre}, it is enough to see that $\caC^{hC_2}$ has cotensors with $C_2$ and $(\caC^{hC_2})^u$ is preadditive.
The $\infty$-category $(\caC^{hC_2})^u$ 
is the essential image of the forgetful functor $\caC^{hC_2} \to \caC$, which is closed under finite products in $\caC$, and thus preadditive.
As $\caC$ has finite products, the forgetful functor $\caC^{hC_2} \to \caC$ admits a right adjoint (sending an object $\X$ to
a $C_2$-action on $\X \times \X$).

For any object $\X\in\caC$, the cofree $C_2$-action $\X'$ generated by $\X$ enjoys the universal property of the cotensor with $C_2$: for any $\Y \in \caC$ we have a canonical equivalence
$$ \map_{\caC^{hC_2}}(\Y, \X') \simeq \map_\caC(\Y,\X) \simeq \Spc^{C_2}(\C_2,\map_{\caC^{hC_2}}(\Y, \X)). $$
\end{proof}

\begin{proposition}\label{fghhgkkkn}\label{PreAdd_gd_genuine_preadditive}
The real $\infty$-category $\Pre\Add^\gd $ is genuine preadditive.
\end{proposition}

\begin{proof}
We will use the characterization of genuine preadditive presented in \cref{rmk:genuine_preadd_char}. 
We have already shown in \cref{fdghjkkkl} that the $\infty$-category $\Pre\Add^\gd $ is preadditive. 

The real $\infty$-category $\Pre\Add^\gd$ has cotensors with all genuine $C_2$-spaces. Indeed, the real forgetful functor $\Cat_{\infty}^\gd \to \Cat_\infty^{hC_2}$ preserves cotensors and the real $\infty$-category $\Pre\Add^{hC_2}$ is closed under taking cotensors with genuine $C_2$-spaces in $ \Cat_{\infty }^{hC_2}$, therefore the pullback $\Pre\Add^{hC_2} \times_{\Cat_\infty^{hC_2}} \Cat_{\infty }^\gd$ has cotensors with genuine $C_2$-spaces. To conclude, we note that the real $\infty$-category $\Pre\Add^\gd$ is closed under taking cotensors with genuine $C_2$-spaces in $\Pre\Add^{hC_2} \times_{\Cat_\infty^{hC_2}} \Cat_{\infty }^\gd$.

We dedicate the rest of the proof to show that the canonical map 
$$\theta\colon \map_{\Pre\Add^\gd }(\caC^{C_2}, \caD) \to \map_{\Spc^{C_2}}(C_2, \map_{\Pre\Add^\gd }(\caC,\caD)) $$
is an equivalence for every $\caC,\caD\in\Pre\Add^\gd$ .

The previous discussion implies that the real forgetful functor $\Pre\Add^\gd \to \Pre\Add^{hC_2}$ preserves cotensors with genuine $C_2$-spaces. Therefore, we can consider the following commutative square

\[
\begin{tikzcd}
\map_{\Pre\Add^\gd }(\caC^{C_2}, \caD) \ar[r, "\theta"] \ar[d]   & \map_{\Spc^{C_2}}(C_2, \map_{\Pre\Add^\gd }(\caC,\caD)) \ar[d] \\
\map_{\Pre\Add^{hC_2} }(\caC^{C_2}, \caD) \ar[r, "\gamma"]    & \map_{\Spc[C_2]}(C_2, \map_{\Pre\Add^{hC_2} }(\caC,\caD)). 
\end{tikzcd}
\]
We show that $\theta$ is an equivalence by showing that $\gamma$ is an equivalence and that the square is a pullback square. To see that $\gamma$ is an equivalence, it is enough to see that the real $\infty$-category $\Pre\Add^{hC_2}$ is genuine preadditive.

From above we know that the real $\infty$-category $\Pre\Add^{hC_2}$ has cotensors with all $C_2$-spaces.
Thus in view of Proposition \ref{genpre} the real $\infty$-category
$\Pre\Add^{hC_2}$ is genuine preadditive if and only if its underlying $\infty$-category is preadditive. The underlying $\infty$-category of $\Pre\Add^{hC_2}$ is the essential image of the forgetful functor $\Pre\Add^{hC_2} \to \Pre\Add$ and is thus preadditive since $ \Pre\Add$ is so and the functor $\beta$ preserves finite products.

It remains to prove that the square above is a pullback square. For this it is enough to check that for every map $\f\colon  \caC \to \caD$ in $\Pre\Add$ corresponding to a
map $\f'\colon \caC^{\caC_2} \to \caD$ in $\Pre\Add^{hC_2}$ the map $$\{ \f'\} \times_{\map_{\Pre\Add^{hC_2} }(\caC^{C_2}, \caD) } \map_{\Pre\Add^\gd }(\caC^{C_2}, \caD) \to $$$$ \{ \f'\} \times_{\map_{\Spc[C_2]}(C_2, \map_{\Pre\Add^{hC_2} }(\caC,\caD)) }\map_{\Spc^{C_2}}(C_2, \map_{\Pre\Add^\gd }(\caC,\caD)) \simeq $$$$ \{ \f\} \times_{\map_{\Pre\Add}(\caC,\caD)} \map_{\Pre\Add}(\caC,\caD) \simeq \ast$$
induced by $\theta$ is an equivalence. 
So we need to check that the space 
$$\{ \f'\} \times_{\map_{\Pre\Add^{hC_2} }(\caC^{\caC_2}, \caD) } \map_{\Pre\Add^\gd }(\caC^{C_2}, \caD)$$ 
is contractible.
Note that for every $\infty$-category with genuine duality its cotensor with $C_2$ carries the standard genuine refinement. So the space
$$\{ \f'\} \times_{\map_{\Pre\Add^{hC_2} }(\caC^{C_2}, \caD) } \map_{\Pre\Add^\gd }(\caC^{C_2}, \caD)$$ 
is canonically equivalent to the space $$\{\caH(\f)\} \times_{\Fun(\caH(\caC^{C_2}), \caH(\caD))} \Fun(\caH(\caC^{C_2}), H)$$
$$\simeq \Fun_{\caH(\caC^{C_2})}(\caH(\caC^{C_2}),H'), $$
where $H \to \caH(\caD)$ is the genuine refinement of $\caD$ and $H'\to\caH(\caC^{C_2}) $ is the pullback of $H \to \caH(\caD)$ along the functor $\caH(\f')\colon  \caH(\caC^{C_2}) \to \caH(\caD).$

Given a right fibration $A \to B$, denote by $\Gamma(A)\coloneqq\Fun_B(B,A)$ its space of sections. With this notation we can write $\Gamma(H')=\Fun_{\caH(\caC^{C_2})}(\caH(\caC^{C_2}),H')$. Therefore we can reformulate what we need to prove by saying that $\Gamma(H')$ must be contractible.

Recall that the underlying $\infty$-category of $\caC^{C_2}$ is $\caC \oplus \caC$
and that $\f' \colon \caC^{C_2} \to \caD$ sends $(x,y)$ to $ \f(x) \oplus \f(y^\vee)^\vee$. 

Let $\gamma\colon \caD \to \caH(\caD)$ the section of the right fibration $\caH(\caD) \to \caD$ which equips each object with the zero hermitian structure. Let $\mathrm{pr}_1, \pr_2\colon  \caC \oplus \caC \to \caC $ the projections. 
We define functors 
$$\pi_1 \colon \caH(\widetilde{\caC \oplus \caC}) \to \caC \oplus \caC \overset{\mathrm{pr}_1}{\longrightarrow} \caC \overset{f}{\longrightarrow}
\caD \overset{\gamma}{\longrightarrow} \caH(\caD),$$
$$\pi_2 \colon \caH(\widetilde{\caC \oplus \caC}) \to \caC \oplus \caC \overset{\mathrm{pr}_2}{\longrightarrow} C \xrightarrow{(-)^\dual} C^\op \overset{f^\op}{\longrightarrow} \caD^\op \xrightarrow{(-)^\dual} \caD \overset{\gamma}{\longrightarrow} \caH(\caD).$$

Denote $H_i \to \caH(\widetilde{\caC \oplus \caC})$ the pullback of $H \to \caH(\caD)$ along $\pi_i$ for $i=1,2$.

There are canonical natural transformations $\pi_i \to \caH(f')$
that yield a map of right fibrations $H' \to H_i$ over $ \caH(\widetilde{\caC \oplus \caC}) $ and so a map $\kappa\colon H' \to  H_1 \times_{\caH(\widetilde{\caC \oplus \caC})} H_2$ of right fibrations over $ \caH(\widetilde{\caC \oplus \caC}) $ that induces on the fiber over
an object $(x,y, \tau\colon x \to y^\dual) \in \caH(\widetilde{\caC \oplus \caC})$
the canonical map $$\rho\colon  H_{ f(x) \oplus f(y^\vee)^\vee} \to  H_{f(x)} \times H_{f(y^\dual)^\dual}.$$

The map $\rho$ is an equivalence, as we have a commutative square
\[
\begin{tikzcd}
0 \ar[r] \ar[d] & f(x) \ar[d] \\
f(y^\dual)^\dual \ar[r] & f(x) \oplus f(y^\vee)^\vee
\end{tikzcd}
\]
in $\caH(\caD)$ lying over a pushout square in $\caD$ and $\caD$ is a preadditive $\infty$-category with genuine duality. Thus $\kappa\colon  H' \to  H_1 \times_{\caH(\widetilde{\caC \oplus \caC})} H_2$ is an equivalence.

Consequently, it is enough to check that $\Gamma(H_i)$ is contractible for $i=1,2$. We show that for $H_1$, the other case is similar.

Define the functors $$p\colon \caH(\widetilde{\caC \oplus \caC}) \to \caC \oplus \caC \xrightarrow{\pr_1} \caC, $$
$$ s\colon \caC \xrightarrow{i_1 } \caC \oplus \caC \xrightarrow{\gamma} \caH(\widetilde{\caC \oplus \caC}),$$
where $i_1$ denotes the inclusion into the first summand.

Note that $p \circ s $ is the identity and there is a canonical natural transformation $\psi\colon s \circ p \to \id$. The transformation $\psi$ induces a map $\Psi\colon H_1 \to p^\ast(s^\ast(H_1))$ of right fibrations over $\caH(\widetilde{\caC \oplus \caC})$, which is an equivalence as $\pi_1 \circ \psi \simeq \id.$

Given a right fibration $A \to B$ and a functor $\alpha\colon  Z \to B$
denote $\alpha^\ast$ the functor $\Gamma(A) \to \Gamma(\alpha^\ast(A))$
that sends a section $ \omega$ of $A$ to the section of $\alpha^\ast(A)$
adjoint to the functor $\omega \circ \alpha\colon  Z \to B \to A$ over $B.$
Note that given a functor $\beta\colon W \to Z$ we have $(\alpha \circ \beta)^\ast \simeq \beta^\ast \circ \alpha^\ast.$
The composition 
$$ p^\ast \circ s^\ast \simeq (s \circ p)^\ast \colon  \Gamma(H_1) \to \Gamma(s^\ast(H_1)) \to \Gamma(p^\ast(s^\ast(H_1))) $$
is the map induced by $\Psi\colon H_1 \to p^\ast(s^\ast(H_1))$ and so an equivalence. Thus $\Gamma(H_1)$ is a retract of $\Gamma(s^\ast(H_1))$, which leads us to verify that $ \Gamma(s^\ast(H_1))$ is contractible.

Note that $ \Gamma(s^\ast(H_1))$ is the limit of the functor
$F\colon \caC^\op \to \Spc$ classified by the right fibration $s^\ast(H_1) \to \caC.$
As $\caC^\op$ has an initial object $0$, the limit of $F$ is canonically equivalent to
$$F(0) \simeq s^\ast(H_1)_0 \simeq (H_1)_{s(0)}  \simeq H_{\pi_1(s(0))} \simeq \{0\} \times_\caC H, $$ which is contractible.
\end{proof}

\begin{corollary}
The real $\infty$-categories $\Add^\gd$ and $\St^\gd $ are genuine preadditive.
\end{corollary}

\begin{proof}
The real $\infty$-categories $\Add^\gd$ and $\St^\gd$ are both closed under small limits and cotensors with small genuine $C_2$-spaces in $ \Pre\Add^\gd$.
\end{proof}

\subsection*{Symmetric monoidal structures}

As announced, the main result of this subsection is \cref{Theta_gd_symmetric_monoidal}, where we prove that the $\infty$-categories $\Pre\Add^\gd$, $\Add^\gd$ and $\St^\gd$ are symmetric monoidal. We also show in \cref{proposition:free} that the respective free functors $\Cat_\infty^\gd\to \Pre\Add^\gd,\Add^\gd,\St^\gd$ is a symmetric monoidal functor.

We begin with the following intermediate step. 

\begin{lemma}\label{dfghjlojh}\label{Theta_hC2_symmetric_monoidal}
Let $\Theta$ denote one of the $\infty$-categories in $\{\Pre\Add, \Add, \St \}$. Then  $\Theta$ is a closed symmetric monoidal $\infty$-category with $C_2$-action and the free functor $\Cat_\infty \to  \Theta$ is a $C_2$-equivariant symmetric monoidal functor.
	
\end{lemma}

\begin{proof}
Recall that the $\infty$-category $\widehat{\Cat}_\infty^{\prod}$ of small
$\infty$-categories with finite products and finite products preserving functors is preadditive. Then we know that the forgetful functor $$\Cmon(\widehat{\Cat}_\infty^{\prod})[C_2] \to \widehat{\Cat}_\infty^{\prod}[C_2] $$ is an equivalence, which allows us to see $\Cat_\infty$ as a symmetric monoidal $\infty$-category with $C_2$-action.
	
The canonical subcategory inclusion $\Cmon(\widehat{\Cat}_\infty) \subset \widehat{\mathrm{Op}}$ that considers a symmetric monoidal $\infty$-category as an
$\infty$-operad gives rise to a subcategory inclusion $$\Cmon(\widehat{\Cat}_\infty^{\prod})[C_2] \subset \Cmon(\widehat{\Cat}_\infty)[C_2]
\subset \widehat{\mathrm{Op}}_\infty[C_2]$$ that endows
$\Cat_\infty$ with the structure of an $\infty$-operad with $C_2$-action
(classifying a $B(C_2)$-family of $\infty$-operads).
	
As next we define a suboperad with $C_2$-action of the $\infty$-operad with $C_2$-action $\Cat_\infty$ by specifying colors and multimorphisms that are preserved by the operadic composition and the $C_2$-action. For this, we take the $\infty$-category of colors as $\Theta \subset \Cat_\infty$ and choose the following multi-morphism spaces: for $X_1,\dots,X_n, Y \in \Theta$ the space $$ \mathrm{Mul}_{\Theta}(X_1,\dots,X_n;Y) \subset \mathrm{map}_{\Cat_\infty}(X_1 \times \dots \times X_n,Y) $$ is the full subspace spanned by the functors that are maps in $\Theta$ in each component.
	
We will complete the proof by showing that the defined $\infty$-operad structure on $\Theta$ is a closed symmetric monoidal $\infty$-category
(which, then, will automatically be compatible with the $C_2$-action).
Denote $\caB_\Theta \subset \Cat_\infty$ the subcategory with objects the small $\infty$-categories that admit finite coproducts, and morphisms the functors preserving finite coproducts if $\Theta \in \{\Pre\Add, \Add\}$; and the subcategory with objects the small $\infty$-categories that admit finite colimits and morphisms the functors preserving finite colimits if $\Theta=\St.$ Note that $ \Theta$ is a localization of $\caB_\Theta$.

We define another suboperad of the $\infty$-operad $\Cat_\infty$ by specifying colors and multimorphisms that are preserved by the operadic composition:
We take the $\infty$-category of colors as $\caB_\Theta \subset \Cat_\infty$ and choose the following multi-morphism spaces: for $X_1,...,X_n, Y \in \caB_\Theta$ the space $$ \mathrm{Mul}_{\caB_\Theta}(X_1,...,X_n;Y) \subset \mathrm{map}_{\Cat_\infty}(X_1 \times ... \times X_n,Y) $$ is the full subspace spanned by the functors that are maps in $\caB_\Theta$ in each component.
By definition the $\infty$-operad structure on $ \Theta$ is a full suboperad of the $\infty$-operad structure on $ \caB_\Theta$.
By \cite[Corollary 4.8.1.4]{lurie.higheralgebra}, the $\infty$-operad structure on $ \caB_\Theta$ is a closed symmetric monoidal $\infty$-category, where the internal hom of $C,D \in \caB_\Theta$ is the
full subcategory of $\Fun(C,D)$ spanned by the finite coproducts respectively finite colimits preserving functors. Moreover, the lax symmetric monoidal functor $ \caB_\Theta \subset \Cat_\infty$ admits a symmetric monoidal left adjoint (see \cite[Remark 4.8.1.8]{lurie.higheralgebra}).

For every $C, D \in \caB_\Theta$ the internal hom of $C,D $ in $\caB_\Theta$ belongs to $ \Theta$ if $D$ does. So the localization $\caB_\Theta \subset \Theta$ is symmetric monoidal so that $ \Theta$ is a closed symmetric monoidal $\infty$-category and the free functor $  \Cat_\infty \to \Theta$ is a symmetric monoidal functor.
\end{proof}

\begin{corollary}\label{cor:Theta_hC2_symmetric_monoidal}Let $\Theta$ denote one of the $\infty$-categories in $\{\Pre\Add, \Add, \St \}$. Then  $\Theta^{hC_2}$  admits a closed symmetric monoidal $\infty$-structure and the free functor $\Cat_\infty^{hC_2}\to \Theta^{hC_2}$ is a symmetric monoidal functor.
\end{corollary}

We use the monoidal structures of \cref{cor:Theta_hC2_symmetric_monoidal} to construct symmetric monoidal structures on the categories $\Pre\Add^{\gd},\ \Add^\gd$ and $\St^\gd$. With this propose we make the following disgression.

Denote by $\caR \subset \Fun([1], \Cat_\infty)$ the full subcategory of right fibrations and $\xi\colon  \caR \to \Cat_\infty$ the cocartesian fibration given by evaluation at the target. Since $\xi$-cocartesian morphisms are closed under finite products, the functor $\xi$ yields on cartesian structures a cocartesian fibration of symmetric monoidal $\infty$-categories.

By Corollary \ref{cor:Theta_hC2_symmetric_monoidal}, the
$\infty$-category $\Theta^{hC_2}$ carries a canonical closed symmetric monoidal structure such that the inclusion
$\Theta^{hC_2} \hookrightarrow \Cat_\infty^{hC_2}$ is lax symmetric monoidal.
Hence the pullback $$\Theta^{hC_2} \times_{\Cat_\infty} \caR \to \Theta^{hC_2}$$ of $\xi$ along the composition
$\Theta^{hC_2} \hookrightarrow \Cat_\infty^{hC_2} \xrightarrow{\caH}\Cat_\infty$
of (lax) symmetric monoidal functors promotes to a cocartesian fibration of symmetric monoidal $\infty$-categories. Especially, $ \Cat_\infty^{hC_2} \xrightarrow{\caH}\Cat_\infty $ carries a symmetric monoidal structure, which is closed by \cref{htwdzfd}.

Before proceeding to endow $\Theta^\gd$ with a symmetric monoidal structure, we will show a result that will be paramount for this.

\begin{lemma}\label{fhjkkllk}
Let $\Theta\in \{\Pre\Add, \Add, \St \}$ and let $\caC, \caD \in \Theta^{hC_2} \times_{\Cat_\infty} \caR$. The internal hom of $\caC,\caD$ in $ \Theta^{hC_2} \times_{\Cat_\infty} \caR $ lies in $\Theta^\gd$ if $\caD$ does.
\end{lemma}

\begin{proof}
Let $(\caC, \phi)$ and $(\caD, \psi)$ be objects in  $\Theta^{hC_2} \times_{\Cat_\infty} \caR \simeq \Theta^{hC_2}\times_{\Cat_\infty^{hC_2}} \Cat_\infty^\gd$.

We can describe the internal hom in $\Theta^{hC_2} \times_{\Cat_\infty} \caR$ from the respective internal homs in $\Cat_\infty^\gd$ and $\Theta^{hC_2}$ (see \cref{GD(pre)_cartesian_closed}). Indeed, if $\theta\colon  T \to \caH(\Fun(\caC, \caD))$ denotes the genuine refinement of the internal hom of $(\caC, \phi), (\caD, \psi) $ in $\Cat_\infty^\gd$ and $\Fun^\mathrm{\Theta}(\caC, \caD)$ denotes the
internal hom in $\Theta^{hC_2}$, then the internal hom in $\Theta^{hC_2} \times_{\Cat_\infty} \caR$ consists of $ \Fun^\mathrm{\Theta}(\caC, \caD) $ and the pullback $\zeta\colon  H \to \caH(\Fun^\mathrm{\Theta}(\caC, \caD))$ of $\theta\colon  T \to \caH(\Fun(\caC, \caD))$ along $ \caH(\Fun^\mathrm{\Theta}(\caC, \caD)) \to \caH(\Fun(\caC, \caD)).$

By \cref{htedud}, for every $\varphi \in \caH(\Fun^{\Theta}(\caC,\caD))$ the fiber $H_\varphi$
is canonically equivalent to the limit of the functor
$$H''^\op \xrightarrow{\phi} \caH(\caC)^\op \overset{\caH(\varphi)^\op}{\longrightarrow} \caH(\caD)^\op \xrightarrow{\nu} \Spc,$$
where the last functor $\nu$ is classified by the right fibration $\psi\colon  H' \to \caH(\caD)$.

Using this description, we first show that $H_\varphi$ is contractible when $\varphi$ lies over the zero exact functor $\caC \to \caD$. The functor 
$$H''^\op \xrightarrow{\phi} \caH(\caC)^\op \overset{\caH(\varphi)^\op}{\longrightarrow} \caH(\caD)^\op \to \caD^\op $$ factors as $$H''^\op \xrightarrow{\phi} \caH(\caC)^\op \to \caC^\op \overset{\varphi^\op}{\longrightarrow} \caD^\op $$ and is thus constant with value the zero object. So the functor $$H''^\op \xrightarrow{\phi} \caH(\caC)^\op \overset{\caH(\varphi)^\op}{\longrightarrow} \caH(\caD)^\op \xrightarrow{\nu} \Spc $$
is constant with value the contractible space as, by hypothesis, $(\caD,\psi) \in \Theta^\gd.$

Now we check that the right fibration $H \to \caH(\Fun^\Theta(\caC,\caD))$ classifies a symmetric monoidal functor $\caH(\Fun^\Theta(\caC,\caD))^\op \to \Spc$ for $\Theta \in \{\Pre\Add, \Add \}$, and an excisive functor for $\Theta = \St$. We only show the case $\Theta= \St$, the other case is similar.

We want to prove that for every commutative square
\[
\begin{tikzcd}
A \ar[r, "f"] \ar[d, "\alpha"'] & B \ar[d, "\beta"] \\
C \ar[r, "g"']  & D
\end{tikzcd}
\]
in $\caH(\Fun^\Theta(\caC,\caD))$ lying over a pushout square in $\Fun^\Theta(\caC,\caD)$, the induced square of spaces
\[
\begin{tikzcd}
H_D \ar[r] \ar[d] & H_B \ar[d] \\
	H_C \ar[r] & H_A
\end{tikzcd}
\]
is a pullback square.
By the description of the fibers we gave above, we can see that the last square arises by taking the limit over $H''^\op$ from the following commutative square in $\Fun(H''^\op, \Spc):$ 
\[
\begin{tikzcd}
\nu \circ \caH(D)^\op \circ \phi \ar[r] \ar[d] & \nu \circ \caH(B)^\op \circ \phi  \ar[d] \\
\nu \circ \caH(C)^\op \circ \phi  \ar[r] & \nu \circ \caH(A)^\op \circ \phi.
\end{tikzcd}
\]

Consequently, it is enough to check that the last square induces a pullback square when evaluated at any $Z \in H''$. Given $Z \in H''$, this evaluation is the image under $\nu$ of the following commutative square
\[
\begin{tikzcd}
\caH(A)(\phi(Z)) \ar[r] \ar[d] & \caH(B)(\phi(Z))  \ar[d] \\
\caH(C)(\phi(Z)) \ar[r] & \caH(D)(\phi(Z)) 
\end{tikzcd}
\]
in $\caH(\caD)$.
The claim follows from the fact that the last square lies over the pushout square 
\[
\begin{tikzcd}[column sep=large, row sep=large]
A(\phi(Z)) \ar[r, "f_{\phi(Z)}"] \ar[d, "\alpha_{\phi(Z)}"']         & B(\phi(Z)) \ar[d, "\beta_{\phi(Z)}"] \\
C(\phi(Z)) \ar[r, "g_{\phi(Z)}"'] & D(\phi(Z))
\end{tikzcd}
\]
in $\caD$ and that by hypothesis $(\caD,\psi) \in \St^\gd.$

This concludes the proof for $\Theta=\Pre\Add$. When $\Theta \in \{\Add, \St\} $ we still need to show that the internal hom $(\Fun^\Theta(\caC, \caD), \zeta)$ is an additive $\infty$-category with genuine duality. By \cref{dghjkkkjg} this is equivalent to verify that for every $ X \in \Fun^\Theta(\caC, \caD)$ and $Y,Z \in \caH(\Fun^\Theta(\caC, \caD))_X$
the maps $$YZ \to Y \otimes Z, \ Y \simeq Y \otimes 1 \to Y \otimes Z$$ in $\caH(\Fun^\Theta(\caC, \caD))$ lying over the maps $X \to X \oplus X$ and $X \simeq X \oplus 0 \to X \oplus X$ in $ \Fun^\mathrm{\Theta}(\caC, \caD)$ yield an equivalence $\rho\colon  H_{ Y \otimes Z} \to H_Y \times H_{YZ}.$

The functor $\rho$ is the limit of the map
$$ \nu \circ \caH(Y \otimes Z)^\op \circ \phi \to (\nu \circ \caH(Y)^\op \circ \phi) \times (\nu \circ \caH(YZ)^\op \circ \phi) $$
in $\Fun(H''^\op, \Spc) $. The component of this map at any given $K \in H''$ is the map
$$\xi\colon  H'_{\caH(Y \otimes Z)(\phi(K))}  \to H'_{\caH(Y)(\phi(K))} \times H'_{\caH(YZ)(\phi(K))} $$
induced by the maps $$ \caH(Y)(\phi(K)) \to \caH(Y \otimes Z)(\phi(K)) \simeq \caH(Y)(\phi(K)) \otimes \caH(Z)(\phi(K)), $$$$ \caH(YZ)(\phi(K)) \to \caH(Y \otimes Z)(\phi(K))\simeq \caH(Y)(\phi(K)) \otimes \caH(Z)(\phi(K))$$ in 
$\caH(\caD)$ that lie, respectively, over the maps 
$$X(\phi(K)) \to X(\phi(K)) \oplus X(\phi(K))$$ 
$$X(\phi(K)) \simeq X(\phi(K)) \oplus 0 \to X(\phi(K)) \oplus X(\phi(K))$$ in $\caD.$ So $\xi$, and in consequence $\rho$, is an equivalence.
\end{proof}

\begin{proposition}\label{Symmon}\label{Theta_gd_symmetric_monoidal}
Let $\Theta$ denote one of the $\infty$-categories in $\{\Pre\Add, \Add, \St \}$. Then $ \Theta^\gd \subset \Theta^{hC_2} \times_{\Cat_\infty} \caR $
is a symmetric monoidal localization.
\end{proposition}

\begin{proof}
By \cref{Prese} the $\infty$-category $ \Theta^\gd \subset \Theta^{hC_2} \times_{\Cat_\infty} \caR $ 
is a localization relative to $\Theta^{hC_2}$. As the internal hom of $\caC, \caD \in \Theta^{hC_2} \times_{\Cat_\infty} \caR$ belongs to $ \Theta^\gd $ if $ \caD$ does (\cref{fhjkkllk}), $ \Theta^\gd \subset  \Theta^{hC_2} \times_{\Cat_\infty} \caR $ is a symmetric monoidal localization relative to $ \Theta^{hC_2}.$ 
Hence the restriction $\Theta^\gd \to\Theta^{hC_2}$ is a cocartesian fibration of symmetric monoidal $\infty$-categories and the localization functor
$\Theta^{hC_2} \times_{\Cat_\infty} \caR \to \Theta^\gd $
is a symmetric monoidal functor.

\end{proof}

Note that, since the symmetric monoidal structure on $\Theta^{hC_2} \times_{\Cat_\infty} \caR$ is closed, the above result implies that the symmetric monoidal structure on $\Theta^\gd$ is also closed.

We conclude this section by proving that the free functor $\Cat_\infty^\gd\to \Theta^\gd$ is symmetric monoidal. In order to do this, we will need the following technical lemma.

\begin{lemma}\label{fdghjkllu}
Consider the following pullback square of $\infty$-categories.
\[
\begin{tikzcd}
\caA \ar[r, "g"] \ar[d, "\psi"'] & \caC \ar[d, "\phi"] \\
\caB \ar[r, "G"'] & \caD
\end{tikzcd}
\]
Let $X\in\caC$, and $Y\in\caB$. Let also $(Y, \alpha\colon  \phi(X) \to G(Y)) $ be a pair that corepresents the functor $\map_\caD(\phi(X),-) \circ G \colon  \caB \to \Spc $ and such that $\alpha\colon  \phi(X) \to G(Y)$ admits a $\phi$-cocartesian lift $ \beta\colon  X \to \alpha_\ast(X)$.
Then the pair $(Z, \beta\colon  X \to g(Z))$, where  $Z\coloneqq (Y, \alpha_\ast(X)) \in \caA$, corepresents the functor $\map_\caC(X,-) \circ g \colon  \caA \to \Spc$. 

\end{lemma}

\begin{proof}
For every $T \in \caA$ we have a commutative square of spaces
\[
\begin{tikzcd}
\map_\caA(Z,T) \ar[r] \ar[d] & \map_\caC(\alpha_\ast(X),g(T))\ar[d] \ar[r] \ar[d] & \map_\caC(X,g(T))\ar[d]\\
\map_\caB(Y,\psi(T)) \ar[r] & \map_\caD(G(Y),G(\psi(T))) \ar[r] & \map_\caD(\phi(X),G(\psi(T))),
\end{tikzcd}
\]
whose bottom horizontal map is an equivalence by assumption. It follows from the pullback square in our hypotheses that the leftmost square is a pullback square. On the other hand, the right-most square is a pullback square as $\beta$ is $\phi$-cocartesian. 
\end{proof}

\begin{remark}\label{rmk:preserving_cocart_fib}
Consider the following pullback square of $\infty$-categories.
\[
\begin{tikzcd}
\caA \ar[r, "g"] \ar[d, "\psi"'] & \caC \ar[d, "\phi"] \\
\caB \ar[r, "G"'] & \caD
\end{tikzcd}
\]
If $\phi$ is a cocartesian fibration and $G$ admits a left adjoint $F$, then also $g$ admits a left adjoint $f$. Moreover, $\phi$ preserves the unit of the adjunctions componentwise. Thus $f$ covers $F$,
and $f$ sends $\phi$-cocartesian morphisms to $\psi$-cocartesian morphisms.
\end{remark}

\begin{proposition}\label{proposition:free}
Let $\Theta$ denote one of the $\infty$-categories in $\{\Pre\Add, \Add, \St \}$. Then the free functor $ \Cat_\infty^\gd \to \Theta^\gd $ is a symmetric monoidal functor.
\end{proposition}

\begin{proof}
By \cref{rmk:preserving_cocart_fib}, the projection $ \Theta^{hC_2} \times_{\Cat_\infty^{hC_2}} \Cat_\infty^\gd \to \Cat_\infty^\gd$ admits a left adjoint $\Cat_\infty^\gd \to \Theta^{hC_2} \times_{\Cat_\infty^{hC_2}} \Cat_\infty^\gd$ that preserves $\phi$-cocartesian morphisms, where $\phi$
is the cocartesian fibration $\Cat_\infty^\gd \to \Cat_\infty^{hC_2}.$
Hence such left adjoint is lax symmetric monoidal. Since this left adjoint covers the symmetric monoidal free functor $ \Cat_\infty^{hC_2} \to \Theta^{hC_2},$ it is itself a symmetric monoidal functor.

Now, by \cref{Theta_gd_symmetric_monoidal}, the embedding $\Theta^\gd \subset \Theta^{hC_2} \times_{\Cat_\infty^{hC_2}} \Cat_\infty^\gd $
is a symmetric monoidal localization. Hence the free functor $ \Cat_\infty^\gd \to \Theta^\gd $ is a symmetric monoidal functor as well.
\end{proof}

Note that the closed symmetric monoidal structure on $\Theta^\gd$ makes $\Theta^\gd$ left tensored (and cotensored) over itself and so especially left tensored
and cotensored over $\Cat_\infty^\gd $
using the free functor $ \Cat_\infty^\gd \to \Theta^\gd $.
\section{Stable genuine dualities and quadratic functors}\label{sec:quadratic_functors}

The aim of this section is to prove that
\begin{center}
    \emph{Stable $\infty$-categories with genuine duality are precisely quadratic functors.}
\end{center}
This result supports the definition of Waldhausen $\infty$-categories with genuine duality presented in this paper.

We begin with a speedy introduction to quadratic functors in \cref{subsec:recap_qu}, where we also give a new definition of quadratic functors that better adjust to our purposes, \cref{def:quadratic_functor}. We prove in \cref{lem:defs_qu_are_equivalent} that this definition coincides with Lurie's when both are available. \cref{subsec:dualities_as_symmetric_monoidal_righ_fibs} is devoted to present dualities as symmetric monoidal right fibrations, with the objective of using such description for the proof of the main result of the section. 

If the reader believes \cref{thm:duality_as_sym_monoidal_right_fibrations}, they can hop to \cref{subsec:stable_wgd_are_quadratic_functors} where we conclude the section by proving our sought result, \cref{equiv_qu_gd}.

\subsection{Recap on quadratic functors}\label{subsec:recap_qu}
Recall that a cube in an $\infty$-category $\caC$ is a functor $[1]^3\coloneqq [1] \times [1] \times [1] \to \caC$. Since the object $(0,0,0) \in [1]^3$ is initial, there exists an equivalence $[1]^3 \simeq ([1]^3 \setminus \{(0,0,0)\})^\triangleleft$.

\begin{definition}Let $\caC$ be an $\infty$-category, and $X\colon [1]^3 \simeq ([1]^3 \setminus \{(0,0,0)\})^\triangleleft\to \caC$ a cube in it. We say that
\begin{itemize}
\item the cube $X\colon [1]^3 \simeq ([1]^3 \setminus \{(0,0,0)\})^\triangleleft\to \caC$ is cartesian if $X$ is a limit diagram, and that
\item the cube $X\colon  [1]^3 \to \caC$ strongly cartesian if $X$ sends all six faces of $[1]^3$ to a pullback square.
\end{itemize}
Dually, we define cocartesian and strongly cocartesian cubes in $\caC$.
\end{definition}
Note that every strongly (co)cartesian cube is in particular a (co)cartesian cube. 

The following functors, singled out by its behavior with respect to (co)cartesian cubes, will play a fundamental role in this section. 

\begin{definition}
A functor $\caC \to \caD$ is 2-excisive if it sends strongly cocartesian cubes to cartesian cubes. 
\end{definition}

\begin{definition}\label{def:quadratic_functor}
A functor $\caC \to \caD$ is quadratic if it is 2-excisive and reduced, i.e.\  preseserves the final object if it exists.
\end{definition}

Denote by
$$\Fun^{\mathrm{qu}}(\caC, \caD) \subset \Fun(\caC, \caD)^{2-\mathrm{exc}} \subset \Fun(\caC, \caD)$$
the full subcategories spanned respectively by the
quadratic functors, 2-excisive functors.

Note that a functor $\caC \to \caD$ is excisive, i.e. sends pushout squares to pullback squares, if and only if it sends strongly cocartesian cubes to strongly cartesian cubes. So every excisive functor is 2-excisive.

\begin{example}\label{fhjjjjkk}
Let $\caC, \caC'$ be $\infty$-categories which admit finite colimits and
$\caD$ a $\infty$-category which admits finite limits. Every functor $\beta\colon \caC \times \caC' \to \caD$ that is excisive in both variables is 2-excisive as functor of one variable. Observe that the diagonal functor $\delta\colon \caC \to \caC \times \caC$ is 2-excisive. So if $\caD$ admits limits indexed by $B(C_2)$ and $\beta\colon \caC \times \caC \to \caD$ is a $C_2$-equivariant functor that is excisive in both variables, the functor $(\beta \circ \delta)^{hC_2}\colon \caC \to \caD$ is 2-excisive.
\end{example}

\begin{remark}
Note that if $\caD$ is stable and we have a fiber sequence $F \to G \to H$ in
$\Fun(\caC, \caD)$, where $H$ is 2-excisive,
then the functor $F$ is 2-excisive if and only if the functor $G$ is.
\end{remark}

We state a direct consequence of \cite[Lemma 6.1.1.35]{lurie.higheralgebra}.
\begin{lemma}\label{exc_exc2_left_localizations}
Let $\caC$ be an $\infty$-category that admits finite colimits and a final object
and $\caD$ an $\infty$-category that admits colimits indixed by
the poset of natural numbers and finite limits such that the formation of colimits over $\mathbb{N}$ preserves finite limits. Then by  the full subcategories $\Fun(\caC, \caD)^\exc \subset \Fun(\caC, \caD)^{2-\exc} \subset \Fun(\caC, \caD)$ spanned by the excisive and 2-excisive functors respectively, are left exact localizations.
\end{lemma}

We denote the left adjoint of these embeddings into $\Fun(\caC,\caD)$ by $P_1, P_2$ and call them excisive, respectively 2-excisive, approximation.
Note that the functors $P_1, P_2$ send reduced functors $\caC \to \caD$ to reduced functors.

\begin{definition}
A functor $\caC \to \caD$ is 2-homogenous if it is 2-excisive and
its excisive approximation vanishes.
\end{definition}

\begin{example}\label{exam}
Given a 2-excisive functor $q\colon \caC \to \caD$ the fiber of the canonical map
$q \to P_1(q)$ in $\Fun(\caC, \caD)$ is 2-homogenous.
\end{example}

When we take $n=2$ in\cite[Proposition 6.1.4.14]{lurie.higheralgebra}, we obtain the following.

\begin{proposition}\label{quad_lurie_h}
Let $\caC$ be an $\infty$-category that
admits finite colimits and a final object
and $\caD$ a stable $\infty$-category that admits colimits indixed by
the poset of natural numbers. Then the functor $$ \Fun( \caC \times \caC, \caD)^{hC_2} \to \Fun(\caC, \caD)[C_2] \xrightarrow{(-)_{hC_2}} \Fun(\caC, \caD) $$ induced by composition with the $C_2$-equivariant diagonal functor $\delta\colon  \caC \to \caC \times \caC$
restricts to an equivalence between the symmetric functors $ \caC \times \caC \to \caD$ that preserve finite colimits in both components and the 2-homogenous functors.
\end{proposition}

We proceed to introduce a notion associated to quadratic functors, its polarization, that will allow us to connect \cref{def:quadratic_functor} with Lurie's definition of quadratic functor (see \cite{lurie.L4}). 

In linear algebra there is a close relation between quadratic and symmetric functors: given a quadratic functor $q\colon  \caV \to K$
on a $K$-vector space $\caV$
one can associate a symmetric functor, called polarization, by the formula 
$B(X,Y)\coloneqq q(X+Y)-q(X)-q(Y).$
The same is true for quadratic and symmetric functors. Any quadratic functor gives rise to a symmetric functor called polarization by 
categorifying the formula above:
One replaces the sum by the direct sum and difference by the cofiber to get to the following definition.

\begin{definition}\label{def:polarization}
Let $q\colon  \caC \to \caD$ be a functor between $\infty$-categories, where $\caD$ is preadditive and has cofibers. We call polarization of $q$ to the cofiber in $ \Fun(\caC \times \caC, \caD)^{hC_2} $ of the $C_2$-equivariant natural transformation $ (X,Y) \mapsto q(X)\oplus q(Y) \to q(X \oplus Y)$ of $C_2$-equivariant functors  $\caC \times \caC \to  \caD.$ 
\end{definition}

\begin{lemma}
Let $\caC$ be an $\infty$-category that
admits finite colimits and a final object
and $\caD$ a stable $\infty$-category that admits colimits indixed by
the poset of natural numbers. Given a quadratic functor $q\colon\caC \to \caD$,
the polarization $\alpha$ of $q$ preserves finite colimits in both variables and fits into a fiber sequence in $\Fun(\caC, \caD)$ as below 
$$(\alpha \circ \delta)_{hC_2} \to q \to P_1(q),$$
with $P_1$ the left adjoint of the embedding $\Fun(\caC,\caD)^{\mathrm{exc}}\subset\Fun(\caC,\caD)$ and $\delta$ the diagonal functor.
\end{lemma}

\begin{proof}
Given any 2-excisive functor $q\colon\caC \to \caD$, the fiber of the canonical map $q \to P_1(q)$ in $\Fun(\caC, \caD)$ is 2-homogenous (\cref{exam}) and so by \cref{quad_lurie_h} of the form
$(\beta \circ \delta)_{hC_2}$ for a unique symmetric functor $\beta\colon \caC \times \caC \to \caD$ that preserves finite colimits in both components. If $q\colon\caC \to \caD$ is moreover reduced, also $P_1(q)$ is reduced and so preserves finite colimits, and therefore the polarization of $P_1(q)$ vanishes. The fiber sequence $$(\beta \circ \delta)_{hC_2} \to q \to P_1(q)$$ in $\Fun(\caC, \caD)$ yields a fiber sequence on polarizations that identifies the polarization
of $(\beta \circ \delta)_{hC_2}$ being $\beta$ with the polarization of $q$.
\end{proof}

\begin{corollary}\label{reco}
Let $\caC$ be an $\infty$-category that
admits finite colimits and a final object
and $\caD$ a stable $\infty$-category that admits colimits indexed by
the poset of natural numbers. Then the following conditions on a functor $F\colon\caC \to \caD$ are equivalent.
\begin{itemize}
\item $F\colon  \caC \to \caD$ is reduced and excisive.
\item $F\colon  \caC \to \caD$ is 2-excisive and the polarization $\alpha\colon  \caC \times \caC \to \caD$ of $F$ vanishes.
\end{itemize}
\end{corollary}

So given a fiber sequence $F \to G \to H$ of 2-excisive functors in
$\Fun(\caC, \caD)$ the functor $F$ is reduced and excisive if and only if 
the map $G \to H$ induces an equivalence on polarizations.

As announced, the following result shows the equivalence of \cref{def:quadratic_functor} with Lurie's, when $\caC$ is stable.

\begin{lemma}\label{lem:defs_qu_are_equivalent}
Let $\caC$ be a pointed $\infty$-category that admits finite colimits
and $\caD$ a stable $\infty$-category. A functor $q\colon\caC \to \caD $ is quadratic if and only if its polarization $\alpha$ preserves finite colimits in both variables
and the fiber of the canonical map $q \to (\alpha \circ \delta)^{hC_2}$, with $\delta$ the diagonal functor, preserves finite colimits.
\end{lemma}

\begin{proof}
Since $q$ is a quadratic functor, its polarization $\alpha$ preserves finite colimits in both variables. Thus $(\alpha \circ \delta)^{hC_2}$ is quadratic and $\alpha$ is the polarization of $(\alpha \circ \delta)^{hC_2}$. So the map $q \to (\alpha \circ \delta)^{hC_2}$ induces an equivalence on polarizations and hence its fiber is exact. 

On the other hand, if $\alpha$  preserves finite colimits in both variables and the fiber $F$ of the map $q \to (\alpha \circ \delta)^{hC_2}$ preserves finite colimits, the functor $q$ sits in a fiber sequence $F \to q \to (\alpha \circ \delta)^{hC_2}$ between quadratic functors and is thus itself quadratic.
\end{proof}

The next result shows that taking the connective cover yields an equivalence
$$ \Fun^{\mathrm{qu}}(\caC, \Sp) \simeq \Fun^{\mathrm{qu}}(\caC, \Sp_{\geq 0})$$
between quadratic functors. In consequence, we can identify quadratic functors $\caC \to \Sp $ with their connective cover $\caC \to \Sp_{\geq 0} $. For that, let us consider $\caC$ to be a small $\infty$-category	that admits finite colimits and a final object. Then the adjunctions $ \iota\colon  \Sp_{\geq 0} \rightleftarrows \Sp \colon R $ and $\mathrm{P}_2\colon  \Fun(\caC, \Sp) \rightleftarrows \Fun(\caC, \Sp)^{2-\exc} $ give rise to an adjunction 
$$\Psi\colon  \Fun(\caC, \Sp_{\geq 0})^{2-\exc} \rightleftarrows \Fun(\caC, \Sp)^{2-\exc} $$ 
that restricts to an adjunction
$$ \Fun^{\mathrm{qu}}(\caC, \Sp_{\geq 0}) \rightleftarrows \Fun^{\mathrm{qu}}(\caC, \Sp) \colon  \Phi $$ 
on quadratic functors. 

\begin{proposition}\label{hjjkkhbvc}\label{proposition_cover}

Let $\caC$ be a small $\infty$-category	that admits finite colimits and a final object. Then the adjunctions above verify the following properties.
\begin{enumerate}
\item The functor $\Psi$ is fully faithful.
\item If $\caC$ is stable, the functor $\Phi$ is an equivalence. 
\item A quadratic functor $\caC^\op \to \Sp$ is excisive if and only if
$RF$ preserves finite products. In other words, the equivalence of (2)  restricts to an equivalence
$$ \Fun(\caC, \Sp_{\geq 0})^{\exc,\ast} \simeq \Fun(\caC, \Sp)^{\exc,\ast}.$$ 
\end{enumerate}
\end{proposition}

\begin{proof}
For the proof of (1), we observe that since the right adjoint $R$ preserves small limits and filtered colimits, it commutes with 2-excisive approximations. Given $F \in \Fun(\caC, \Sp_{\geq 0})^{2-\exc}$, the unit component $F \simeq R \iota F \to R \mathrm{P}_2 (\iota F) \simeq \mathrm{P}_2(F)  $ is the canonical equivalence, as $ R \iota F \simeq F$ is $2$-excisive.

For (2), we will show that $\Phi$ is conservative, which by (1) implies that it is an equivalence. Since $\Fun^{\mathrm{qu}}(\caC, \Sp)$ is stable, it is enough to prove that a functor $G \in \Fun^{\mathrm{qu}}(\caC, \Sp)$ is zero if $RG$ vanishes:
Note that an exact functor $H\colon  \caC \to \Sp $ vanishes if $RH\colon  \caC \to \Sp_{\geq 0} $ vanishes and similarly for a functor $\caC \times \caC \to \Sp$ that is exact in both variables.

Since $G$ is quadratic, its polarization $\alpha$ is exact in both variables. We know that $R \alpha$ is the polarization of $RG \simeq 0.$ Therefore, $R \alpha $ vanishes, and so $\alpha$ vanishes as well. As $G$ has vanishing polarization, by \cref{reco} $G$ is exact. So $RG \simeq 0 $ implies $G \simeq 0.$ 

Finally, we prove (3). By \cref{reco}, a quadratic functor $F\colon  \caC \to \Sp$ is excisive if and only if its polarization $\alpha$ vanishes, where $\alpha$ is exact in both variables. If $RF$ preserves finite products, its polarization, which is given by $R\alpha$, vanishes. Consequently, $\alpha$ vanishes as $\alpha$ is exact in both variables. Thus $F$ is excisive.
\end{proof}

\subsection{Dualities as symmetric monoidal right fibrations}\label{subsec:dualities_as_symmetric_monoidal_righ_fibs}

In this subsection we give a description of dualities on a preadditive $\infty$-category in terms of symmetric monoidal right fibrations. The trustful reader can skip this subsection if they are willing to believe \cref{thm:duality_as_sym_monoidal_right_fibrations}. 

We present the proof of \cref{thm:duality_as_sym_monoidal_right_fibrations} forth referencing results that appear later within the subsection.

\begin{theorem}\label{thm:duality_as_sym_monoidal_right_fibrations}
There exists a subcategory inclusion
$$\Pre\Add^{hC_2} \subset \Cmon(\Pre\Add \times_{\Cat_\infty} \caR)$$ over $ \Pre \Add$. Its essential image is the $\infty$-category whose objects are the symmetric monoidal right fibrations $\phi\colon \caD \to \caC$, where $\caC$ carries the cartesian monoidal structure,
that classify functors $\psi\colon  \caC^\op \to \Spc_{\geq0},$ whose polarization $\alpha$ is non-degenerate, and such that the canonical map $\psi \to \alpha^{hC_2}$ is an equivalence; and its morphisms are those maps of symmetric monoidal right fibrations 
that induce strict maps on polarizations.
\end{theorem}

\begin{proof}
Denote by $$ (\Pre\Add\times \Pre\Add) \times^{\prod}_{\Cat_\infty} \caR \subset (\Pre\Add \times \Pre\Add) \times_{  \Cat_\infty } \caR$$ the full subcategory (closed under finite products and the $C_2$-action) spanned by the triples $(\caC,\caD,\psi)$, where $\caC,\caD$ are preadditive $\infty$-categories and $\psi$ is a right fibration $\caB \to \caC \times \caD$ that classifies a functor $\caC^\op \times \caD^\op \to \Spc$ which preserves finite products in both components. 
By \cite[Theorem 5.9]{hls}, there is a (non-full) subcategory inclusion over $\Cat_\infty $
$$\Cat_\infty^{hC_2} \subset \Cat_\infty \times_{\Cat_\infty[C_2]} \caR[C_2]$$ 
sending an $\infty$-category $\caC$ with duality to $(\caC, \Tw(\caC) \to \caC \times \caC)$. Pulling back this inclusion to $\Pre\Add$ we get an inclusion
$$\Pre\Add^{hC_2} \subset \Pre\Add \times_{\Cat_\infty[C_2]} \caR[C_2]$$ over $ \Pre\Add$ that factors through $\Pre\Add \times^{\prod}_{\Cat_\infty[C_2]} \caR[C_2]$,
where we set 
$$ \Pre\Add \times^{\prod}_{\Cat_\infty[C_2]} \caR[C_2] \coloneqq ((\Pre\Add\times \Pre\Add) \times^{\prod}_{\Cat_\infty} \caR)^{hC_2} $$$$ \subset\Pre\Add \times_{\Cat_\infty[C_2]} \caR[C_2].$$

We will complete the proof by constructing
an embedding 
$$\psi\colon \Pre\Add \times^{\prod}_{\Cat_\infty[C_2]} \caR[C_2] \subset \Cmon(\Pre\Add\times_{\Cat_\infty} \caR)$$ over $ \Pre\Add$
that sends $H \to \caC $ to $H^{hC_2} \to \caC^{hC_2}.$

In order to construct $\psi$, we first prove that the $\infty$-category $ \Pre\Add \times^{\prod}_{\Cat_\infty[C_2]} \caR[C_2]$ is preadditive. Note that the $\infty$-category $ (\Pre\Add \times \Pre\Add) \times^{\prod}_{\Cat_\infty} \caR$ is preadditive. Indeed, the forgetful functor
$$\Cmon((\Pre\Add \times \Pre\Add) \times^{\prod}_{\Cat_\infty} \caR)
\to (\Pre\Add\times \Pre\Add) \times^{\prod}_{\Cat_\infty} \caR$$
is a map of cartesian fibrations over $ \Pre\Add \times \Pre\Add$ that yields on the fiber over $\caC^\op,\caD^\op \in  \Pre\Add$ the forgetful functor
$$ \Fun^{\prod}(\caC \times \caD, \Cmon(\Spc)) \to \Fun^{\prod}(\caC \times \caD, \Spc),$$ that is in turn an equivalence. Therefore, also the subcategory $ \Pre\Add \times^{\prod}_{\Cat_\infty[C_2]} \caR[C_2]$ is preadditive.

Now, since the functor $ \caR[C_2] \xrightarrow{(-)^{hC_2}} \caR $
yields a finite products preserving functor 
$$\Pre\Add\times^{\prod}_{\Cat_\infty[C_2]} \caR[C_2] \subset \Pre\Add \times_{\Cat_\infty[C_2]} \caR[C_2] \to \Pre\Add \times_{\Cat_\infty} \caR $$
over $\Pre\Add $ we obtain a functor 
$$ \psi\colon \Pre\Add \times^{\prod}_{\Cat_\infty[C_2]} \caR[C_2] \to \Cmon(\Pre\Add \times_{\Cat_\infty} \caR)$$
over $\Pre\Add$. By Corollary \ref{rhjkklj}, we know that $\psi$ is fully faithful.

Finally, composing both inclusions we obtain the desired functor
$$\Pre\Add^{hC_2} \subset \Cmon(\Pre\Add\times_{\Cat_\infty} \caR)$$ over $ \Pre\Add$
that sends $\caC $ to $ \caH(\caC)= \Tw(\caC)^{hC_2} \to (\caC \times \caC)^{hC_2} \simeq \caC.$

The description of the essential image of this inclusion, follows directly from its construction.
\end{proof}

\begin{corollary}There exists a subcategory inclusion over $\Add$,
$$\Add^{hC_2} \subset \Cmon(\Add \times_{\Cat_\infty} \caR),$$  whose essential image admits an analogous description to that of \cref{thm:duality_as_sym_monoidal_right_fibrations}.
\end{corollary}
\begin{proof}
It follows from \cref{thm:duality_as_sym_monoidal_right_fibrations} by pulling back to $\Add$.
\end{proof}

\begin{corollary}There exists a subcategory inclusion over $\Pre\Add$
$$\Pre\Add^{hC_2} \times_{\Cmon(\Cat_\infty^{hC_2}) } \Cmon(\Cat_\infty^\gd)\subset\Cmon(\Pre\Add \times_{\Cat_\infty} \caR)^{[1]}$$ that sends $(\caC, H \to \caH(\caC))$ to $\caH(\caC) \to \caC, H \to \caH(\caC).$ A similar statement holds for $\Add$ by pulling back to it.
\end{corollary}

\begin{proof}It follows from the following chain of inclusions
$$\Pre\Add^{hC_2} \times_{\Cmon(\Cat_\infty^{hC_2}) } \Cmon(\Cat_\infty^\gd) \simeq \Pre\Add^{hC_2} \times_{\Cmon(\Cat_\infty) } \Cmon(\caR)$$$$  \subset \Cmon(\Pre\Add\times_{\Cat_\infty} \caR \times_{\Cat_\infty} \caR) \simeq \Cmon( \Pre\Add\times_{\Cat_\infty} \caR^{[1]} ) \simeq $$
$$  \Cmon(\Pre\Add \times_{\Cat_\infty} \caR)^{[1]}.$$
\end{proof}

It remains to prove that the functor over $\Pre\Add$
$$\psi\colon \Pre\Add \times^{\prod}_{\Cat_\infty[C_2]} \caR[C_2] \to \Cmon(\Pre\Add \times_{\Cat_\infty} \caR)$$  
that we constructed in the proof of \cref{thm:duality_as_sym_monoidal_right_fibrations} is fully faithful; this is the content of \cref{rhjkklj}. In order to show this, we first need to prove the result that we embark to prove next.

Let $\caC$ be a preadditive $\infty$-category and $\caD$ a preadditive $\infty$-category that admits limits over $B(C_2)$ and cofibers. Given a functor $q\colon \caC \to \caD$ we have a $C_2$-equivariant natural transformation of $C_2$-equivariant functors 
$\caC \times \caC \to  \caD$:
$$(X,Y) \mapsto q(X)\oplus q(Y) \to q(X \oplus Y).$$ We will denote by $ q'\in \Fun(\caC \times \caC, \caD)^{hC_2}$ its cofiber, i.e.\ its polarization. 

Denote by $\Fun^{\prod}(\caC \times \caC, \caD)^{hC_2} \subset \Fun(\caC \times \caC, \caD)^{hC_2}$ the full subcategory spanned by the $C_2$-equivariant functors $\caC \times \caC \to \caD$ that preserve finite products in both components; and by $\Fun'(\caC,\caD) \subset \Fun(\caC,\caD)$ the full subcategory spanned by the functors
$q\colon \caC \to \caD$ such that $q'$ belongs to $\Fun^{\prod}(\caC \times \caC, \caD)^{hC_2}.$

The composition $$\Fun(\caC \times \caC, \caD)^{hC_2} \xrightarrow{\delta^\ast} \Fun(\caC, \caD)^{hC_2} \simeq \Fun(\caC, \caD[C_2]) \xrightarrow{(-)^{hC_2}} \Fun(\caC, \caD)$$
with $\delta\colon \caC \to \caC \times \caC$ the $C_2$-equivariant diagonal functor
restricts to a functor $$ \beta\colon \Fun^{\prod}(\caC \times \caC, \caD)^{hC_2} \to \Fun'(\caC, \caD) .$$

\begin{lemma}\label{dghjkkl} 
Let $\caC$ be a preadditive $\infty$-category and $\caD$ a preadditive $\infty$-category that admits limits over $B(C_2)$ and cofibers. The functor $\beta$ is fully faithful and admits a left adjoint $L$ that sends 
$q \in \Fun'(\caC, \caD)$ to its polarization. 
\end{lemma}

\begin{proof}
The canonical natural transformation $\lambda\colon \id \to \oplus \circ \delta$ of functors
$\caC \to \caC$ exhibits $\oplus\colon \caC \times \caC \to \caC$ as a right adjoint to $\delta$. So $\oplus$ is a $C_2$-equivariant functor and $\lambda\colon\id \to \oplus \circ \delta$ is a $C_2$-equivariant natural transformation.

Let $q \in \Fun(\caC, \caD)$. By definition of $q'$, the polarization of $q$, we have a $C_2$-equivariant natural transformation $ q \circ \oplus \to q'$ of $C_2$-equivariant functors $\caC \times \caC \to \caD.$
So we obtain a $C_2$-equivariant natural transformation
$$q \xrightarrow{q \circ \lambda} q \circ \oplus  \circ \delta \to q' \circ \delta $$ of $C_2$-equivariant functors $\caC \to \caD$, which corresponds to a natural transformation
$q \to q' \circ \delta$ of functors $\caC \to \caD[C_2]$. The latter is adjoint to a natural transformation
$h\colon q \to (q' \circ \delta)^{hC_2}= \beta(q')$ of functors $\caC \to \caD$.

Denote by  $\Fun^\mathrm{pol}(\caC, \caD) \subset  \Fun'(\caC, \caD)$ the full subcategory spanned by the functors $q\colon \caC \to \caD$ such that the natural transformation $h$ is an equivalence.

We will show that the left map of the factorization of $\beta$ below, is an equivalence.
\[
\begin{tikzcd}
\Fun^{\prod}(\caC \times \caC, \caD)^{hC_2}\ar[r]       &\Fun^\mathrm{pol}(\caC, \caD)\subset\Fun'(\caC, \caD)
\end{tikzcd}
\]

For every $\gamma \in \Fun^{\prod}(\caC \times \caC, \caD)^{hC_2}$, the polarization of $\beta(\gamma)$, this is, the cofiber $ \beta(\gamma)' $ of the $C_2$-equivariant natural transformation $$ (X,Y) \mapsto \gamma(X,X)^{hC_2} \oplus \gamma(Y,Y)^{hC_2} \to \gamma(X \oplus Y,X \oplus Y)^{hC_2} \simeq $$$$ \gamma(X,X)^{hC_2} \oplus \gamma(Y,Y)^{hC_2} \oplus (\gamma(X,Y) \oplus \gamma(Y,X))^{hC_2} \simeq $$$$ \gamma(X,X)^{hC_2} \oplus \gamma(Y,Y)^{hC_2} \oplus \gamma(X,Y) $$ of $C_2$-equivariant functors $\caC \times \caC \to  \caD$ is $ \gamma$. Moreover, the map $$h\colon(\gamma \circ \delta)^{hC_2} =\beta(\gamma) \to \beta(\beta(\gamma)') \simeq \beta(\gamma)= (\gamma \circ \delta)^{hC_2} $$ is the identity.
Note that for every $ q \in \Fun'(\caC, \caD)$ the induced map $h'\colon  q' \to \beta(q')' \simeq q' $ is the identity.

As next we prove that $\beta\colon  \Fun^{\prod}(\caC \times \caC, \caD)^{hC_2} \to \Fun^\mathrm{pol}(\caC, \caD)$ is an equivalence. 
Note that $\beta$ is an equivalence if and only if for every $\infty$-category $T$ the induced functor $$\Fun(T, \beta)\colon \Fun(T,  \Fun^{\prod}(\caC \times \caC, \caD)^{hC_2}) \to  \Fun(T,\Fun^\mathrm{pol}(\caC, \caD)) $$ yields a bijection on equivalence classes. We reduce this to the case when $T$ is contractible.

Note that a functor $\omega\colon  T \to  \Fun(\caC, \caD)$ factors through $\Fun^\mathrm{pol}(\caC, \caD)$ if and only if its adjoint functor $q\colon \caC \to \Fun(T,\caD)$ belongs to $\Fun^\mathrm{pol}(\caC, \Fun(T,\caD)) $. This follows from the fact that the map $h\colon q \to (q' \circ \delta)^{hC_2} $ formed with respect to $\Fun(T,\caD)$ yields after evaluation at $t \in T$ the corresponding map $h\colon  \omega(t) \to (\omega(t)' \circ \delta)^{hC_2}$ for $ \omega(t) $ formed with respect to $\caD$.

The functor $\Fun(T, \beta)$ is equivalent to the functor
$$\Fun(T, \beta)\colon \Fun(T, \Fun^{\prod}(\caC \times \caC, \caD)^{hC_2}) \simeq \Fun^{\prod}(\caC \times \caC, \Fun(T,\caD))^{hC_2} \xrightarrow{\beta} $$$$ \Fun^\mathrm{pol}(\caC, \Fun(T,\caD)) \simeq \Fun(T,\Fun^\mathrm{pol}(\caC, \caD)). $$

We conclude that we can indeed reduce our statement to show that $\beta$ induces a bijection on equivalence classes.

Now, we know already that for every $\gamma \in \Fun^{\prod}(\caC \times \caC, \caD)^{hC_2}$ we have $\beta(\gamma)' \simeq \gamma.$ On the other hand, for every functor $q \in \Fun^\mathrm{pol}(\caC, \caD)$ the map $h\colon  q \to \beta(q')$ is an equivalence. So $\beta$ is fully faithful with essential image $ \Fun^\mathrm{pol}(\caC, \caD)$.

Moreover, for $T= \Fun'(\caC,\caD)$ and $q\colon  \caC \to \Fun(T,\caD)$ adjoint to 
the subcategory inclusion $T \to \Fun(\caC,\caD)$ the map $h\colon  q \to \beta(q') $ in 
$\Fun(\caC,\Fun(T,\caD)) \simeq \Fun(T, \Fun(\caC,\caD)) $ provides a natural transformation $\lambda\colon  \id \to \beta \circ L$, where $L$ is the functor $$\Fun(\caC,\caD)' \to \Fun^{\prod}(\caC \times \caC, \caD)^{hC_2}$$ adjoint to $q' \in \Fun^{\prod}(\caC \times \caC, \Fun(T,\caD))^{hC_2}. $  So for every $q \in \Fun(\caC,\caD)' $ the component $\lambda(q) \colon  q \to \beta(L(q)) $ is $h.$ Thus the map $$L(\lambda(q)) =L(h) \colon  L(q) \to L(\beta(L(q))) $$ is an equivalence and $\lambda(q)$ is an equivalence if $q \in \Fun^\mathrm{pol}(\caC, \caD).$ 
\end{proof}

\begin{corollary}\label{rhjkklj}
The functor over $\Pre\Add$ below (as in the proof of \cref{thm:duality_as_sym_monoidal_right_fibrations}) is fully faithful and admits a left adjoint.
$$\psi\colon\Pre\Add \times^{\prod}_{\Cat_\infty[C_2]} \caR[C_2] \to \Cmon(\Pre\Add\times_{\Cat_\infty} \caR).$$
\end{corollary}

\begin{proof}
Denote by 
$$\Cmon'(\Pre\Add \times_{\Cat_\infty} \caR) \subset  \Cmon(\Pre\Add \times_{\Cat_\infty} \caR)$$
the full subcategory spanned by the symmetric monoidal right fibrations $\caD \to \caC$ for some $\caC \in \Cmon(\Pre\Add)$ that classify a functor $q\colon  \caC^\op \to \Cmon(\Spc)$ such that
the polarization of $q$, $q' \colon \caC^\op \times \caC^\op \to \Spc$, preserves finite products in both components.

The functor $\psi $ induces a functor 
$$ \Pre\Add \times^{\prod}_{\Cat_\infty[C_2]} \caR[C_2] \to \Cmon'(\Pre\Add \times_{\Cat_\infty} \caR).$$

As $\psi$ is a map of cartesian fibrations over $\Pre\Add$, it is enough to check that for every $\caC \in \Pre\Add$ the induced functor $$\psi_{\caC^\op}\colon\{\caC^\op \} \times_{\Pre\Add} (\Pre\Add \times^{\prod}_{\Cat_\infty[C_2]} \caR[C_2])  \to $$$$ \{\caC^\op\} \times_{\Pre\Add} \Cmon'(\Pre\Add \times_{\Cat_\infty} \caR) $$ 
on the fiber over $\caC^\op$ is fully faithful and admits a left adjoint.
But $\psi_{\caC^\op} $ factors as $$\{\caC^\op \} \times_{\Pre\Add} (\Pre\Add \times^{\prod}_{\Cat_\infty[C_2]} \caR[C_2]) \simeq \Fun^{\prod}(\caC \times \caC, \Spc)^{hC_2}$$$$\simeq \Fun^{\prod}(\caC \times\caC, \Cmon(\Spc))^{hC_2} \xrightarrow{\beta} \Fun'(\caC, \Cmon(\Spc)) $$$$ \simeq \{\caC^\op\} \times_{\Pre\Add} \Cmon'(\Pre\Add \times_{\Cat_\infty} \caR), $$
where $\beta $ is defined in \cref{dghjkkl}, so the result follows from \cref{dghjkkl}.
\end{proof}

\subsection{Stable $\infty$-categories with genuine duality are quadratic functors}\label{subsec:stable_wgd_are_quadratic_functors}
In this subsection we prove that for any stable $\infty$-category with genuine duality $(\caC, \phi\colon H\to \caH(\caC))$, the symmetric monoidal right fibration $H \to \caH(\caC)$ classifies a quadratic functor $\caC^\op \to \Sp_{\geq 0}$ \textemdash and then by \cref{hjjkkhbvc} a quadratic functor $\caC^\op \to \Sp$. This is the content of \cref{equiv_qu_gd}.

\begin{proposition}\label{fjjjkljh} Let $(\caC, \phi)$ be a stable $\infty$-category with genuine duality. Then the symmetric monoidal right fibration $H \to \caC$ classifies a quadratic functor
$\caC^\op \to \Sp_{\geq 0}.$
\end{proposition}

\begin{proof}
By \cref{fhjjjjkk}, we know that the right fibration $\caH(\caC) \to \caC$ classifies a 2-excisive functor $\caC^\op \to \Spc.$  So it is enough to see that for every strongly cartesian cube $X\colon  [1]^3 \to \caC$, the pullback $\caZ$ of $ H \to \caH(\caC)$ along any lift $A \colon  [1]^3 \to \caH(\caC)$ of $X$ classifies a cartesian cube of spaces.

Denote by $\tau\colon [1]^2\coloneqq [1] \times [1] \to [1]^3$ the inclusion of one of the six faces. The pullback of $\caZ \to [1]^3$ along $\tau\colon [1]^2 \to [1]^3 $, i.e.\ the pullback of $ H \to \caH(\caC)$ along $A \circ \tau \colon  [1]^2 \to \caH(\caC) $, classifies a cartesian square of spaces. Indeed, $X \circ \tau\colon  [1]^2 \to \caC $ is a cartesian square by hypotheis and $\phi$ classifies an excisive functor $\caH(\caC)^\op \to \Spc.$ We conclude then that $\caZ$ classifies a strongly cartesian cube and so especially a cartesian cube.
\end{proof}

\begin{definition}\label{definition_quadratic}
We call an additive right fibration $q\colon\caD \to \caC$
a quadratic right fibration if $\caC$ is stable and $\caD \to \caC$ classifies a quadratic functor $\caC^\op \to \Sp_{\geq 0}.$
\end{definition}

The next two propositions address situations in which we have good local-global properties for a commutative monoid in $\Cat_\infty^\gd$, say $(\caC, \phi)$. This is, we describe properties of $\phi$ in terms of its fibers. To highlight their, perhaps subtle at first sight, difference, we present both statements before writing their respective proof.

\begin{proposition}\label{dfjklkh_1}Let $(\caC, \phi)$ be a commutative monoid in $\Cat_\infty^\gd$. If $H \to \caC$ is an additive right fibration, the following conditions are equivalent.
\begin{enumerate}
\item $\phi$ classifies a symmetric monoidal functor $\theta\colon\H(\caC)^\op \to \Spc. $
\item The fiber of $\phi$ classifies a functor that preserves finite products.
\end{enumerate}
\end{proposition}

\begin{proposition}\label{dfjklkh_2}Let $(\caC, \phi)$ be a commutative monoid in $\Cat_\infty^\gd$. If $H \to \caC$ is a quadratic right fibration, the following conditions are equivalent.
\begin{enumerate}
\item $\phi$ classifies an excisive functor $\theta\colon\H(\caC)^\op \to \Spc. $
\item The fiber of $\phi$ classifies an exact functor.
\end{enumerate}
\end{proposition}

\begin{proof}[Proof of \cref{dfjklkh_1}]
The fiber of $\phi$ classifies the functor $\caC^\op \xrightarrow{\nu^\op} \H(\caC)^\op \xrightarrow{\theta} \Spc, $ where $\nu\colon \caC \to \H(\caC)$ is the canonical symmetric monoidal section of $ \H(\caC) \to \caC.$ Therefore (1) implies (2).

Assume now that (2) holds. As $\nu$ is symmetric monoidal, the fiber of $\phi\colon H\to\H(\caC)$ over the tensor unit $\mathbbm{1}\in\H(\caC)$ is contractible. It remains to see that for every $A,B \in \H(\caC)$ the canonical map $H_{A \otimes B} \to H_A \times H_B$ is an equivalence.
When $H_A \times H_B$ is empty, it is trivially satisfied. We can then assume that
$H_A \times H_B$ is non-empty.

Let $X \in H_A, Y \in H_B.$ Denote by $U,V$ the image of $A,B$, respectively, in $\caC$.
Take $S = [1]$ and $\gamma\colon  X \simeq X \otimes 1 \to X \otimes Y$
the canonical map in $H$ lying over $\beta\colon  A \simeq A \otimes 1 \to A \otimes B$
in $ \caH(\caC) $, which in turn lies over $\alpha\colon  U \simeq U \oplus 0 \to U \oplus V $
in $\caC$.

By \cref{lem:for_local_global_phiprop}, we get an equivalence $(\nu \circ \alpha)^\ast(H) \simeq \beta^\ast(H) $ over $[1]$, where $\nu\colon\caC\to\H(\caC)$ is as before, classifying a commutative square of spaces
\[
\begin{tikzcd}
H_{A \otimes B} \ar[r] \ar[d, "\simeq"'] & H_A \ar[d, "\simeq"] \\
H_{\nu(U \oplus V) } \ar[r] &  H_{\nu(U)}. 
\end{tikzcd}
\]

Similarly, for $\gamma\colon Y \simeq 1 \otimes Y \to X \otimes Y$
the canonical map in $H$ lying over $\beta\colon  B \simeq 1 \otimes B \to A \otimes B$
in $ \caH(\caC) $, which in turn lies over $\alpha\colon  V \simeq 0 \oplus V \to U \oplus V $
in $\caC$ we get an equivalence $(\nu \circ \alpha)^\ast(H) \simeq \beta^\ast(H) $ over $[1]$ classifying a commutative square of spaces
\[
\begin{tikzcd}
H_{A \otimes B} \ar[r] \ar[d, "\simeq"'] & H_B \ar[d, "\simeq"] \\
H_{\nu(U \oplus V) } \ar[r] &  H_{\nu(V)}.
\end{tikzcd}
\]

From those squares, we obtain a commutative square of spaces as below.
\[
\begin{tikzcd}
H_{A \otimes B} \ar[r] \ar[d, "\simeq"'] & H_A \times H_B \ar[d, "\simeq"] \\
H_{\nu(U \oplus V) } \ar[r] &  H_{\nu(U)} \times  H_{\nu(V)}
\end{tikzcd}
\]
Now the claim follows from (2).
\end{proof}

For the next result, that we have used in the proof of \cref{dfjklkh_1}, consider a map of symmetric monoidal right fibrations over $\caC$,  $\rho\colon  H \to H'$. We know that it classifies a natural transformation of functors $\caC^\op \to \Alg_{E_\infty}(\Spc)$, that is itself classified by a map of commutative monoids in the $\infty$-category of right fibrations over $\caC$. This makes $\rho$ a map of commutative monoids in the $\infty$-category of right fibrations over $\caC$. We denote by $\mu\colon  H \times_\caC H \to H, \ \mu'\colon  H' \times_\caC H' \to H'$ the multiplications and $\eta\colon  \caC \to H, \eta' \colon  \caC \to H'$ the units of $H, H'$ respectively. 

\begin{lemma}\label{lem:for_local_global_phiprop}
Let $\caC$ be an additive $\infty$-category, $\rho\colon  H \to H'$ a map of additive right fibrations over $\caC$, and $\alpha\colon S \to \caC, \gamma\colon S \to H$ functors over $\caC$. There is a canonical equivalence $(\eta' \circ \alpha)^\ast(H) \simeq \beta^\ast(H) $ over $S$, where $\beta\coloneqq \rho \circ \gamma\colon S \to H'$ and $\eta'$ is the unit of $H' \to \caC.$
	
\end{lemma}
\begin{proof}
	
Let us assume first that $\alpha$ is the identity. Then the composition $\psi\colon  H \xrightarrow{\gamma \times_\caC H} H \times_\caC H \xrightarrow{\mu} H$ of functors over $\caC$ is an equivalence, as it induces an equivalence on the fiber over every object $X$ of $\caC$, $ H_X \xrightarrow{\gamma_X \times H_X} H_X \times H_X \xrightarrow{\mu_X} H_X$, using that the $E_\infty$-structure on $H_X$ is grouplike. Similarly, the composition $\psi'\colon  H' \xrightarrow{\beta \times_\caC H'} H' \times_\caC H' \xrightarrow{\mu'} H'$ is an equivalence. 

The commutative square of right fibrations over $\caC$ below
\[
\begin{tikzcd}
H \ar[r, "\psi"] \ar[d]     & H \ar[d] \\
H' \ar[r, "\psi'"'] &  H'
\end{tikzcd}
\]

yields an equivalence $\eta'^\ast(H) \simeq \beta^\ast(H) $ over $\caC$ by pulling back along $\eta'\colon  \caC \to H'$, since $\beta$ factors as $\caC \xrightarrow{\eta'} H' \xrightarrow{\psi'} H'.$

Now consider functors $\alpha\colon S \to \caC$ and $\gamma\colon S \to H$ over $\caC$, with $\gamma$ corresponding to a section of $S \times_\caC H \to S$ as in the hypotheses. We can pull-back $\rho$ along $\alpha\colon  S \to \caC$ to a map of commutative monoids of right fibrations over $S$, and obtain this way the desired equivalence $(\eta' \circ \alpha)^\ast(H) \simeq \beta^\ast(H) $ over $S$ with $\beta\coloneqq \rho \circ \gamma.$
\end{proof}	
 
For the stable case, that we prove in what follows, we will also need a lemma, which appears immediately after said proof.

\begin{proof}[Proof of \cref{dfjklkh_2}]
The fiber of $\phi$ classifies the functor $\caC^\op \xrightarrow{\nu^\op} \caH(\caC)^\op \xrightarrow{\theta} \Spc, $ where $\nu\colon \caC \to \caH(\caC)$ is the canonical symmetric monoidal section of $ \caH(\caC) \to \caC.$ A square in $\caH(\caC)$ is a pullback square if it lies over a pushout square in $\caC$, hence $\nu$ is excisive.
Moreover, $\caH(\caC)$ has an initial object lying over the zero object of $\caC$. So $\nu$ sends the zero object to the initial object. The functor $\theta$ sends the initial object of $\caH(\caC)$ (lying over the zero object in $\caC$)
to the contractible space, as both the fiber of $\caH(\caC) \to \caC$ and the fiber of $H \to \caC$ are contractible (the latter needs $H \to \caC$ to be reduced, which it is since it classifies a quadratic functor). Therefore (1) implies (2).

Let us assume now that (2) holds. By \cite[Proposition 1.4.2.13]{lurie.higheralgebra} it is enough to check that the functor $\caH(\caC)^\op \to \Spc$ classified by $\phi$ sends pullback squares in $\caH(\caC)$ as below left lying over some pullback square in $\caC$ as below right
\begin{equation}\label{ffhjknbfd}
\begin{tikzcd}
S \ar[r] \ar[d] & 1 \ar[d] \\
1 \ar[r] &  T 
\end{tikzcd}
\hspace{3em}
\begin{tikzcd}
A \ar[r] \ar[d] & 0 \ar[d] \\
0 \ar[r] &  B
\end{tikzcd}
\end{equation}
to pullback squares of spaces as below.
\begin{equation}\label{fghjkll}
\begin{tikzcd}
H_T \ar[r] \ar[d] & \ast \ar[d] \\
\ast \ar[r] &  H_S  
\end{tikzcd}
\end{equation}

By \cref{dfghjjjb}, the fiber $H_T \simeq \{T\} \times_{ \caH(\caC)_B } H_B$ is empty if and only if the pullback $ \ast \times_{H_S} \ast \simeq \{ T' \} \times_{\Omega(\caH_A) } \Omega(H_A) $ is so, where $T'$ denotes the image of $T$ in $\Omega(\caH_A).$ When this is the case, there is nothing left to show.

We can then assume that $H_T$ is not empty. Let $Z$ be in $H_T$, and let us consider the unique square in $H$ lying  over the left square in (\ref{ffhjknbfd}); we depict it as below.
\[
\begin{tikzcd}
W \ar[r] \ar[d] & 1 \ar[d] \\
1 \ar[r] &  Z
\end{tikzcd}
\]

Taking $\gamma\colon  [1] \times [1] \to H$ as the functor specified by this square
we get a canonical equivalence $(\nu \circ \alpha)^\ast(H) \simeq \beta^\ast(H) $ over $[1] \times [1]$
classifying an equivalence of commutative squares between square
(\ref{fghjkll}) and the square 
\[
\begin{tikzcd}
H_{\nu(B)} \ar[r] \ar[d] & H_{\nu(0)} \ar[d] \\
H_{\nu(0)} \ar[r] &  H_{\nu(A)}.
\end{tikzcd}
\]
So the claim follows from item (2). 
\end{proof}

\begin{lemma}\label{dfghjjjb}
Consider a commutative square of spectra with exact vertical columns as below.
\[
\begin{tikzcd}
X \ar[r, "\simeq"] \ar[d] & \Omega(F) \ar[d]\\
Y \ar[r] \ar[d] & \Omega(G) \ar[d] \\
Z \ar[r, "\varphi"] & \Omega(H)
\end{tikzcd}
\]

Then, given $T \in Z$, the fiber $\{T\} \times_Z Y $ is nonempty if and only if the fiber $\{\varphi(T)\} \times_{\Omega(H)} \Omega(G)$ is non-empty \textemdash for notational simplicity, we denote the underlying grouplike $E_\infty$-space of a spectrum with the same symbol.
\end{lemma}

\begin{proof}
	
The commutative square in the statement yields a new commutative square of spectra whose vertical columns are also exact, as depicted below.	
\[
\begin{tikzcd}
Y \ar[r] \ar[d] & \Omega(G) \ar[d]\\
Z \ar[r] \ar[d] & \Omega(H) \ar[d] \\
\Sigma(X) \ar[r, "\simeq"] & F,
\end{tikzcd}
\]

Therefore, $T$ is in the image of the map $Y \to Z$ if and only if $T$ vanishes in $\Sigma(X)$. Similarly for $\varphi(T)$.
\end{proof}

We now prove our main result of this section, in order to do so, we introduce one more notion. 

\begin{definition}
We call an additive (resp. quadratic) right fibration $q\colon\colon\caD\to\caC$ non-degenerate if its polarization is non-degenerate as a symmetric functor.
\end{definition}

Denote by 
$$\mathrm{QuR} \subset \mathrm{AddR} \subset \Cmon(\caR)$$
the subcategories with objects the non-degenerate additive
(respectively quadratic) right fibrations and with morphisms those maps that yield on polarizations a map of dualities (and yield exact functors after evaluation at the target $\caR \to \Cat_\infty$).

\begin{theorem}\label{fghkklll}\label{equiv_qu_gd}
There is a localization $\mathrm{AddR} \subset \Add^\gd$, whose essential image consists of those additive $\infty$-categories with genuine duality $(\caC, \phi)$ such that $\phi$ induces an equivalence on polarizations. Moreover, this localization induces an equivalence $\mathrm{QuR} \simeq \St^\gd.$
\end{theorem}

\begin{proof}
By \cref{rhjkklj}, the functor $ \Cmon(\Cat_\infty^\gd) \to \Cmon(\caR) $ that sends $(\caC, \phi)$ to $ H \to \H(\caC) \to \caC$ restricts to an equivalence between 
\begin{itemize}
\item 
$\mathrm{AddR}$ and the subcategory of $\Cmon(\Cat_\infty^\gd) $ with objects the commutative monoids $(\caC, \phi)$ in $\Cat_\infty^\gd $ such that $H \to \caC$ is an additive right fibration and $ \phi\colon  H \to \H(\caC)$ induces an equivalence on polarizations, and with morphisms the maps
in $\Cmon(\Cat_\infty^\gd) $, whose underlying functor is additive.

\item $ \mathrm{QuR}$ and the subcategory of $\Cmon(\Cat_\infty^\gd) $ with objects the commutative monoids $(\caC, \phi)$  in $\Cat_\infty^\gd$ such that $H \to \caC$ is a quadratic right fibration
and $\phi$ induces an equivalence on polarizations, and with morphisms those maps, whose underlying functor is exact.
\end{itemize}

Via these equivalences we view $ \mathrm{QuR}, \mathrm{AddR} $ as subcategories of $\Cmon(\Cat_\infty^\gd).$

Since $ \Add^\gd, \St^\gd $ are preadditive (see \cref{rmk:Add_gd_preadditive} and \cref{rmk:St_gd_preadditive}), we can see these $\infty$-categories as subcategories of $\Cmon(\Cat_\infty^\gd)$. More precisely, $\Add^\gd$ is the subcategory of $\Cmon(\Cat_\infty^\gd) $ with objects the commutative monoids $(\caC, \phi)$ in $\Cat_\infty^\gd $ such that $H \to \caC$ is an additive right fibration and 
$\phi\colon  H \to \caH(\caC) $ classifies a symmetric monoidal functor $\caH(\caC)^\op \to \Spc$, and with morphisms the maps in $\Cmon(\Cat_\infty^\gd) $, whose underlying functor is additive.
And $ \St^\gd$ is the subcategory of $\Cmon(\Cat_\infty^\gd) $ with objects the commutative monoids $(\caC, \phi)$ in $\Cat_\infty^\gd $ such that $H \to \caC$ is an additive right fibration and 
$\phi\colon H \to \caH(\caC) $ classifies an excisive functor $\caH(\caC)^\op \to \Spc$, and with morphisms the maps in $\Cmon(\Cat_\infty^\gd)$, whose underlying functor is exact.

Consequently, for the first part of the statement, it is enough to see that for any commutative monoid $(\caC, \phi)$ in $\Cat_\infty^\gd $ such that $H \to \caC$ is an additive right fibration and such that $\phi$ induces an equivalence on polarizations, the functor $\phi\colon  H \to \caH(\caC) $ classifies a symmetric monoidal functor $\caH(\caC)^\op \to \Spc.$

We show that what is written above is true in the following way:
Since $\phi$ induces an equivalence on polarizations, the cofiber of $\phi$, i.e.
the cofiber of the map of functors $\caC^\op \to \Sp_{\geq 0 }$ classified by $\phi$,
has vanishing polarization (as taking polarizations preserves cofibers).
Hence the cofiber of $\phi$ preserves finite products. Thus also the fiber of $\phi$
preserves finite products. Invoking \cref{dfjklkh_1} we find that $\phi\colon  H \to \caH(\caC)$ indeed classifies a symmetric monoidal functor.

Now we deal with the restriction to the $\infty$-category $\mathrm{QuR}$. By \cref{fjjjkljh}, for every stable $\infty$-category $(\caC, \phi)$ with genuine duality the functor $H \to \caC$ is a quadratic right fibration.
Consequently, we need to see that for any commutative monoid $(\caC, \phi)$ in $\Cat_\infty^\gd $ such that $H \to \caC$ is a quadratic right fibration, the functor $\phi\colon  H \to \caH(\caC)$ classifies an excisive functor if and only if $\phi$ induces an equivalence on polarizations.
By \cref{dfjklkh_1}, the functor $\phi\colon  H \to \caH(\caC)$ classifies an excisive functor if and only if the fiber of $\phi$
is exact. As $H \to \caC$ is a quadratic right fibration and so classifies a functor $\caC^\op \to \Sp$ (\cref{hjjkkhbvc}), this is equivalent to say that the cofiber of $\phi$ is exact.
By \cref{reco} the cofiber of $\phi$ is exact if and only if its polarization vanishes, which is equivalent to say that $\phi$ induces an equivalence on polarizations.
\end{proof}

\section{\texorpdfstring{Waldhausen $\infty$-categories with genuine duality}{Waldhausen Infinity-categories with genuine duality}}\label{sec:Wald_gd}

In this section we will introduce Waldhausen $\infty$-categories with genuine duality, which will be the input of our real K-theory functor. We pursued in \cref{sec:genuine_add_preadd_stable} a theory of (pre)additive and stable $\infty$-categories with genuine duality, which set the stage for the work on this section: defining Waldhausen $\infty$-categories with genuine duality. These categories interpolate between additive and stable $\infty$-categories with genuine duality, in a similar fashion as the theory of exact $\infty$-categories embeds into the more general world of additive $\infty$-categories but contains the very rich theory of stable $\infty$-categories as a special case.

\subsection{Waldhausen $\infty$-categories with duality}\label{subsec:wald_duality}

In this subsection we make a brief recap of Waldhausen $\infty$-categories as presented by Barwick (see \cite{barwick.wald} for details), and define exact $\infty$-categories (\cref{def:exact_infty_cats}) tha twill allow us to introduce the notion of duality. In later sections we will build on this to obtain our input categories.

\begin{definition}
A Waldhausen $\infty$-category is a pair $(\caD, \caC) $ with $\caD$ an  $\infty$-category with zero object and $\caC$ a subcategory of $\caD$, whose morphisms are called cofibrations, subject to the following conditions.

\begin{itemize}
\item Every object $X$ of $\caD$ is cofibrant, i.e. the unique map $0 \to X$ is a cofibration.
	
\item The cobase change in $\caD$ of any cofibration exists and is again a cofibration.
\end{itemize}

A map of Waldhausen $\infty$-categories $(\caD, \caC) \to (\caD', \caC') $
is a functor $\caD \to \caD'$ that preserves the initial object, cofibrations and pushouts along cofibrations.
\end{definition}

We will denote by $\Wald_\infty \subset \Fun([1], \Cat_\infty)$ the subcategory spanned by the Waldhausen $\infty$-categories and maps between them as defined above. Note that the evaluation at the target defines a faithful functor $\Wald_\infty \to \Cat_\infty$, i.e.\ a functor that yields embeddings of spaces on mapping spaces.

It is interesting to observe that the notion of Waldhausen $\infty$-category is not self dual. We call a pair $(\caC,\caF)$ a co-Waldhausen $\infty$-category if $(\caC^\op,\caF^\op)$ is a Waldhausen $\infty$-category. 

Next we define exact $\infty$-categories, which is a self-dual notion. An exact $\infty$-category will consist of an $\infty$-category $\caD$ compatibly carrying a Waldhausen and co-Waldhausen structure, where the compatibility is encoded in the following notion introduced in \cite{barwick.exact}.

\begin{definition}Let $(\caD,\caC,\caF)$ be a triple where $\caD$ is an $\infty$-category and $\caC,\caF$ are subcategories, whose morphisms we call cofibrations and fibrations respectively. A commutative square $\sigma$ in $\caD$
\[
\begin{tikzcd}
A \ar[r, "f"] \ar[d, "\alpha"'] & B \ar[d, "\beta"] \\
C \ar[r, "g"']  & D
\end{tikzcd}
\]

\begin{itemize}
\item is an ambigressive pushout square if $\sigma$ is a pushout square, $f$ is a cofibration and $\alpha$ is a fibration, and 
\item an ambigressive pullback square if $\sigma$ is a pullback square, $g$ is a cofibration and $\beta$ is a fibration.
\end{itemize}

We call a square as above exact when is both an ambigressive pushout and an ambigressive pullback.
\end{definition}

\begin{definition}\label{def:exact_infty_cats}
An exact $\infty$-category is a triple $(\caD, \caC, \caF)$ with $\caD$ an additive $\infty$-category and $\caC, \caF$ subcategories of $\caD$, whose morphisms are called respectively cofibrations and fibrations, subject to the following conditions.

\begin{enumerate}
\item The pair $(\caD, \caC) $ is a Waldhausen $\infty$-category.
\item The pair $(\caD, \caF) $ is a co-Waldhausen $\infty$-category.
\item A commutative square in $\caD$ is an ambigressive pushout square if and only if it is an ambigressive pullback square.
\end{enumerate}
\end{definition}

The morphisms between them are the natural ones, as given in the following definition.

\begin{definition}
Given exact $\infty$-categories $(\caD, \caC, \caF), (\caD', \caC', \caF') $
a functor $F\colon \caD \to \caD'$ is called exact (with respect to the exact structures on $\caC$ and $\caD$) if $F$ preserves the zero object, cofibrations, fibrations, pushouts along cofibrations and pullbacks along fibrations.

\end{definition}

In light of condition 3, in an exact $\infty$-category we call ambigressive pushout squares respectively ambigressive pullback squares simply ambigressive squares.

\begin{remark}
A functor $F\colon \caD \to \caD'$ is exact if and only if $F$ preserves the zero object, cofibrations and cofibers of cofibrations or dually the zero object, fibrations and fibers of fibrations.
\end{remark}

Denote by $\Exact_\infty \subset \Fun(\Lambda^2_2, \Cat_\infty)$ the subcategory spanned by the exact $\infty$-categories and exact functors.
Note that valuation at $2 \in \Lambda^2_2$ defines a faithful functor $\Exact_\infty \to \Cat_\infty$.\ Restriction along the embedding $[1] \subset \Lambda^2_2$ that sends 0 to 0 and 1 to 2 defines a functor $$\Fun(\Lambda^2_2, \Cat_\infty) \to \Fun([1], \Cat_\infty)$$ that restricts to a fully faithful functor $$\Exact_\infty \to \Wald_\infty.$$

In the next result, we characterize the additivity condition on an $\infty$-category satisfying all conditions of being exact, except being additive.

\begin{proposition}
Let $\caD$ be an $\infty$-category that admits a zero object and $\caC,\caF$ subcategories of $\caD$, whose morphisms we call cofibrations and  fibrations respectively. If conditions (1)-(3) of \cref{def:exact_infty_cats} hold, then $\caD$ is additive if and only if the following condition holds:
a commutative square 
\begin{equation*}\label{jjj1}
\begin{tikzcd}
X \ar[r] \ar[d]     & Y \ar[d]  \\
Z \ar[r]            & W
\end{tikzcd}
\end{equation*}
in $\caD$, where both vertical maps are fibrations, is a pullback square if it induces an equivalence on fibers of the vertical morphisms.
\end{proposition}

\begin{proof}
The `only if' direction follows from \cite[Proposition A.1]{Hull}.

For the `if' direction we introduce the following terminlogy: a sequence $A\to B\to \C$ in $\caD$ is exact when the commutative square below is exact.

\begin{equation*}\label{jjj2}
\begin{tikzcd}
A \ar[r] \ar[d]         & B \ar[d]  \\
0 \ar[r]                & C
\end{tikzcd}
\end{equation*}

We start by proving that $\caD$ is preadditive. For any $A,B \in \caD$
the canonical morphism $\theta\colon A \coprod B \to A \times B$ is an equivalence, since the outer sequences of the commutative diagram below are exact.

\[
\begin{tikzcd}[column sep=small]
& A\ar[dr]\ar[dl] &\\ 
A\coprod B\ar[rr, "\theta"]\ar[dr]       &       &A\times B\ar[dl]\\
 &A     &\,
\end{tikzcd}
\]

Now, we can define the shearing map $\rho\colon A \times A \to A \times A$ for any  $A \in \caD$,
which is the projection to the first factor on the first factor and the codiagonal $A \times A \simeq A \coprod A \to A$ on the second factor. To show additivity, we want this map to be an equivalence. As before, this follows by noticing that the outer sequences of the diagram below are exact.
\[
\begin{tikzcd}[column sep=small]
& A\ar[dr]\ar[dl] &\\ 
A\coprod B\ar[rr, "\rho"]\ar[dr, ]&&A\times B\ar[dl, "\pr_1"]\\
 &A     &\,
\end{tikzcd}
\]
\end{proof}

\begin{remark}The embedding $\Exact_\infty \subset \Wald_\infty$ is a localization. Indeed, we know that $\Exact_\infty$ is presentable, and that the $\infty$-categories $\Wald_\infty \subset \Fun([1], \Cat_\infty), \ \Exact_\infty \subset \Fun(\Lambda^2_2, \Cat_\infty)$ are closed under small limits and small filtered colimits, and therefore $\Exact_\infty$ is closed  under small limits and small filtered colimits in $\Wald_\infty$. 
\end{remark}

In order to define a duality on Waldhausen $\infty$-categories, we note that the non-trivial $C_2$-actions on $\Lambda^2_2$ and $\Cat_\infty$ induce a 
$C_2$-action on $\Fun(\Lambda^2_2, \Cat_\infty)$ that restricts to $\Exact_\infty$. This allows us to introduce the following definition.

\begin{definition}
We define the $\infty$-category of small Waldhausen $\infty$-categories with duality as the homotopy fixed points of the action on $\Exact_\infty$ described above. We write $$\Wald_\infty^\d\coloneqq \Exact_\infty^{hC_2}.$$
\end{definition}

In other words, a Waldhausen $\infty$-category with duality is a Waldhausen $\infty$-category, whose underlying $\infty$-category carries a duality such that the cofibrations together with the opposites of the cofibrations
form an exact $\infty$-category. 

Relevant examples of these kind of categories are given by stable $\infty$-categories with duality.

\begin{example}
Every stable $\infty$-category $\caC$ gives rise to an exact structure on $\caC$ with cofibrations and fibrations all morphisms. More importantly for us, if $E$ is an $\infty$-category with duality whose underlying $\infty$-category is stable, then $E$ is a Waldhausen $\infty$-category with duality in which every morphism is a cofibration.
\end{example}

We conclude this subsection by presenting the monoidal structure on Waldhausen $\infty$-categories with duality. The idea of the proof is the same as the one for \cref{Theta_hC2_symmetric_monoidal}, we write the changed details below.

\begin{proposition}\label{dfghjkgfl}
The $\infty$-category $\Exact_\infty$ is a closed symmetric monoidal $\infty$-category with $C_2$-action. Moreover, the free functor $\Cat_\infty \to \Exact_\infty$ is a $C_2$-equivariant symmetric monoidal functor.
\end{proposition}

\begin{proof}
Recall that the $\infty$-category $\widehat{\Cat}_\infty^{\prod}$ of small
$\infty$-categories with finite products and finite products preserving functors is preadditive. Then we know that the forgetful functor $$\Cmon(\widehat{\Cat}_\infty^{\prod})[C_2] \to \widehat{\Cat}_\infty^{\prod}[C_2] $$ is an equivalence, which allows us to see $\Cat_\infty$ as a symmetric monoidal $\infty$-category with $C_2$-action.
	
The canonical subcategory inclusion $\Cmon(\widehat{\Cat}_\infty) \subset \widehat{\mathrm{Op}}$ that considers a symmetric monoidal $\infty$-category as an
$\infty$-operad gives rise to a subcategory inclusion $$\Cmon(\widehat{\Cat}_\infty^{\prod})[C_2] \subset \Cmon(\widehat{\Cat}_\infty)[C_2]
\subset \widehat{\mathrm{Op}}_\infty[C_2]$$ that endows
the functor $\Fun(\Lambda^2_2, \Cat_\infty) \to \Cat_\infty$
with the structure of a $C_2$-equivariant map of $\infty$-operads classified by a 
map of $B(C_2)$-families of $\infty$-operads \cite[Definition 2.3.2.10.]{lurie.higheralgebra}.
	
Next we define a suboperad with $C_2$-action of the $\infty$-operad with $C_2$-action $\Fun(\Lambda^2_2, \Cat_\infty)$ by specifying colors and multimorphisms that are preserved by the operadic composition and the $C_2$-action. For this, we take the $\infty$-category of colors as $\Exact_\infty \subset \Fun(\Lambda^2_2, \Cat_\infty)$ and choose the following multi-morphism spaces: for exact $\infty$-categories $X_1,\dots,X_n, Y $ the space 
$$ \mathrm{Mul}_{\Exact_\infty}(X_1,\dots,X_n;Y) \subset \mathrm{map}_{\Fun(\Lambda^2_2, \Cat_\infty)}(X_1 \times \dots \times X_n,Y) $$ is the full subspace spanned by the functors that are exact functors in each component.

In the following we will prove that the defined $\infty$-operad structure on $\Exact_\infty$ is a closed symmetric monoidal $\infty$-category (and thus automatically compatible with $C_2$-actions) and the free functor $\Cat_\infty \to \Exact_\infty$ is symmetric monoidal 
(and thus automatically compatible with $C_2$-actions).

For that we define a similar $\infty$-operad structure on $\Wald_\infty$
(that is not compatible with the $C_2$-action).
The functor $\Fun([1], \Cat_\infty) \to \Cat_\infty$ evaluating at the target 
preserves finite products and so extends to a symmetric monoidal functor on cartesian structures. We define a suboperad of the $\infty$-operad $\Fun([1], \Cat_\infty)$ whose colors are $\Wald_\infty \subset \Fun([1], \Cat_\infty)$ and whose space of multi-morphisms
$X_1 \dots X_n \to Y$ for Waldhausen $\infty$-categories $X_1,\dots,X_n, Y $ is the space 
$$ \mathrm{Mul}_{\Wald_\infty}(X_1,\dots,X_n;Y) \subset \mathrm{map}_{\Fun([1], \Cat_\infty)}(X_1 \times \dots \times X_n,Y) $$ spanned by the functors that preserve cofibrations in each component.

The embedding $[1] \subset \Lambda^2_2$ that sends 0 to 0 and 1 to 2 defines a functor $$\Fun(\Lambda^2_2, \Cat_\infty) \to \Fun([1], \Cat_\infty)$$ that preserves finite products and so lifts to a symmetric monoidal functor on cartesian structures that restricts to a full embedding of $\infty$-operads $ \Exact_\infty \subset \Wald_\infty$.

The latter embedding and the forgetful functor $\Wald_\infty \to \Cat_\infty$ admit left adjoints by the adjoint functor theorem since $ \Wald_\infty, \Exact_\infty$ are presentable $\infty$-categories. 

By \cite[1.9. Proposition]{barwick-mult} the $\infty$-operad structure on $\Wald_\infty$ is a closed symmetric monoidal $\infty$-category: for any $\caC, \caD \in \Wald_\infty$ the internal hom $\caE$ is the full subcategory of $\Fun(\caC, \caD)$ spanned by the maps of Waldhausen $\infty$-categories with appropriate Waldhausen structure.
If $\caD$ is an exact $\infty$-category, also $\caE$ is an exact $\infty$-category.
This guarantees that the localization $ \Exact_\infty \subset \Wald_\infty$ is a symmetric monoidal localization. Hence the $\infty$-operad structure on $ \Exact_\infty $ is a closed symmetric monoidal $\infty$-category, too.

The resulting lax symmetric monoidal forgetful functor $\Wald_\infty \to \Cat_\infty$ 
admits a left adjoint $\varphi\colon \Cat_\infty \to \Wald_\infty$, which is automatically an oplax symmetric monoidal functor. By \cite[1.7. Lemma]{barwick-mult} the functor $\varphi$ is symmetric monoidal. Consequently, the forgetful functor $\Exact_\infty \subset \Wald_\infty \to \Cat_\infty$ is a lax symmetric monoidal functor and its left adjoint $ \Cat_\infty \to \Exact_\infty$ is a symmetric monoidal functor.
\end{proof}

\begin{corollary}
The $\infty$-category $\Exact_\infty^{hC_2} $ admits a canonically closed symmetric monoidal structure and the free functor $ \Cat_\infty^{hC_2} \to \Exact_\infty^{hC_2} $ is a symmetric monoidal functor.
\end{corollary}

\subsection{Waldhausen $\infty$-categories with genuine duality}\label{subsec:wald_gd}
As expected, in this subsection we present the concept of Waldhausen $\infty$-categories with genuine duality, a refinement of that of \cref{subsec:wald_duality}. Moreover, this notion generalizes that of \emph{stable} $\infty$-categories with genuine duality introduced in \cref{subsec:stable_with_gd}.

\begin{definition}\label{Waldgd}
A small Waldhausen $\infty$-category with genuine duality is 
a pair $(\caE, \phi)$, with $E$ a small Waldhausen $\infty$-category with duality
and $\phi\colon  H \to \caH(\caE)$ a right fibration enjoying the following conditions.

\begin{enumerate}
\item $(\caE, \phi)$ is an additive $\infty$-category with genuine duality.
\item for every commutative square
\begin{equation*}
\begin{tikzcd}
A \ar[r, "f"] \ar[d, "\alpha"']         &B \ar[d, "\beta"] \\
	C \ar[r, "g"']                      &D
\end{tikzcd}
\end{equation*}
in $\caH(\caE)$ lying over a pushout square in $\caE$ such that $f$ lies over a cofibration, the induced square of spaces
\[
\begin{tikzcd}
H(D) \ar[r] \ar[d]          &H(B) \ar[d] \\
H(C) \ar[r]                 &H(A)
\end{tikzcd}
\]
is a pullback square.
\end{enumerate}
\end{definition}

\begin{definition}
We define the $\infty$-category of small Waldhausen $\infty$-categories with genuine duality  $\Wald_\infty^\gd$ as the full subcategory of $\Wald_\infty^\d \times_{\Add^{hC_2}} \Add^\gd$ spanned by the small Waldhausen $\infty$-categories with genuine duality. 
\end{definition}

\begin{example}
Every stable $\infty$-category with genuine duality is a Waldhausen $\infty$-category with genuine duality, where every map is a cofibration. This example is explored in more detail in \cref{sec:quadratic_functors}.
\end{example}

From the embedding $\Wald_\infty^\gd \subset \Wald_\infty^\d \times_{\Cat_\infty^{hC_2}} \Cat_\infty^\gd$ we infer that $\Wald_\infty^\gd$ is a real $\infty$-category. Moreover, we can see that it is genuine preadditive. Indeed, we can consider the embedding
$$\Wald_\infty^\gd \subset \Wald_\infty^\d \times_{\Add^{hC_2}} \Add^\gd$$ 
and use that the real $\infty$-categories defining the pullback are genuine preadditive. Indeed, it follows from \cref{cor:hC2_gives_genuine_preadditive} that $\Wald_\infty^\d$ and $\Add^{hC_2}$ are genuine preadditive, and it is the content of \cref{PreAdd_gd_genuine_preadditive} that $\Add^\gd$ is so. Therefore, $\Wald_\infty^\gd$ is genuine preadditive.

\begin{remark}
We can pullback any $\infty$-category with genuine duality along a map of Waldhausen 
$\infty$-categories with duality. Indeed, since the full subcategory $\Wald_\infty^\gd$ is stable under the fiber transports of the cartesian fibration $ \Wald_\infty^\d \times_{\Cat_\infty} \caR \to \Wald_\infty^\d$, the forgetful functor $\Wald_\infty^\gd \to \Wald_\infty^\d$ is a cartesian fibration.
\end{remark}

\begin{proposition}\label{dfghjlok}
The full subcategory $$ \Wald_\infty^\gd \subset \Wald_\infty^\d \times_{\Add^{hC_2}} \Add^\gd $$ is an accessible localization relative to $ \Wald_\infty^\d.$
\end{proposition}

\begin{proof}
The $\infty$-category $\Exact_\infty$ is presentable
and the forgetful functor $\Exact_\infty \to \Cat_\infty$ admits a left adjoint.
So the pullback $ \Wald_\infty^\d \times_{\Cat_\infty^{hC_2}} \Cat_\infty^\gd$ is presentable.

Let $\kappa$ be a regular cardinal such that $\caH\colon \Cat_\infty^{hC_2} \to  \Cat_\infty$ preserves small $\kappa$-filtered colimits. Then $ \Wald_\infty^\gd $ is closed in $ \Wald_\infty^\d \times_{\Cat_\infty} \caR $ under small $\kappa$-filtered colimits.

The full subcategory $\Wald_\infty^\gd$ is stable under the fiber transports of the cartesian fibration $\Wald_\infty^\d \times_{\Cat_\infty^{hC_2}} \Cat_\infty^\gd \to \Wald_\infty^\d$.
So to prove that $ \Wald_\infty^\gd \subset \Wald_\infty^\d \times_{\Cat_\infty^{hC_2}} \Cat_\infty^\gd$ is a localization relative to $\Wald_\infty^\d $, it is enough to check this fiberwise. This is, that for every $\caC \in \Wald_\infty^\d$ the fiber $(\Wald_\infty^\gd)_\caC $ is a localization of $  \{\caH(\caC)\} \times_{ \Cat_\infty} \caR.$

For every $\caC \in \Wald_\infty^\d$ denote $\caK$ the set of commutative squares in $\caH(\caC)$ that lie over pushout squares in $\caC$ with two parallel morphisms cofibrations as in \cref{Waldgd}. To verify this, denote by $$\caB \subset \{\caH(\caC)\} \times_{ \Cat_\infty} \caR \simeq \Fun(\caH(\caC)^\op, \Spc)$$ the full subcategory spanned by the reduced functors $\caH(\caC)^\op \to \Spc$ that send diagrams of $\caK$ to pullback squares. Then $\caB \subset \{\caH(\caC)\} \times_{ \Cat_\infty} \caR $ is an accessible localization as $\caK$ is a set.
The fiber $(\Wald_\infty^\gd)_\caC$ is the intersection within $\{\caH(\caC)\} \times_{ \Cat_\infty} \caR $ of $\caB$ and $(\Add^\gd)_\caC.$ Then, by \cref{Prese}, the fiber $(\Add^\gd)_\caC$  is an accessible localization of $ \{\caH(\caC)\} \times_{ \Cat_\infty} \caR. $
\end{proof}

\begin{corollary}
$\Wald_\infty^\gd$ is presentable and the cartesian fibration
$ \Wald_\infty^\gd \to \Wald_\infty^\d$ is a bicartesian fibration.
\end{corollary}

We conclude this subsection by presenting a closed symmetric monoidal structure on exact $\infty$-categories with genuine duality, which leads to an enrichment over $\Cat_\infty^\gd$.

\begin{proposition}\label{fhbhnj}
The $\infty$-category $\Wald^\gd_\infty $ carries a closed symmetric monoidal structure such that the forgetful functor $\Wald^\gd_\infty \to \Exact_\infty^{hC_2}$ is a cocartesian fibration of symmetric monoidal $\infty$-categories. Futhermore, the free functor $ \Cat_\infty^\gd \to \Wald^\gd_\infty $ is symmetric monoidal.
\end{proposition}

\begin{proof}
By \cref{dfghjkgfl} the $\infty$-category $\Exact_\infty$ is a closed symmetric monoidal $\infty$-category with $C_2$-action and the free functor $\Cat_\infty \to \Exact_\infty$ is a $C_2$-equivariant symmetric monoidal functor.
So the result follows as in the proof of \cref{Symmon}, this is, analogously to \cref{fhjkkllk}, one checks that the internal hom of $\caC, \caD \in  \Exact_\infty^{hC_2} \times_{\Cat_\infty} \caR$ belongs to $ \Wald^\gd_\infty $ if $ \caD$ does.
\end{proof} 

\begin{corollary}\label{cor:Waldgd_is_cotensored_with_Catgd}
The $\infty$-category $\Wald^\gd_\infty $ is left tensored, cotensored, and enriched over $\Cat_\infty^\gd $. Moreover, the free functor $ \Cat_\infty^\gd \to \Wald^\gd_\infty $ is a $\Cat_\infty^\gd $-linear functor.
\end{corollary}

\part{\texorpdfstring{Real $K$-theory}{Real K-theory}}\label{part:3}

\section{\texorpdfstring{The additivity theorem for Waldhausen $\infty$-categories with genuine duality}{The additivity theorem for Waldhausen Infinity-categories with genuine duality}}\label{sec:additivity_thm}

Certainly one of the most important features of algebraic K-theory is its 
property to split short exact sequences in the following way:
the induced functor 
$$ \gamma_2\colon \{  A \to B \to C \} \to \caC \times \caC, A \to B \to C \mapsto (A,C) $$
from the $\infty$-category of exact sequences in $\caC$ to the $\infty$-category of pairs in $\caC$
induces an isomorphism on $K$-groups.
From that point of view K-theory does not distinguish between short
exact sequences and split short exact sequences. This seems already reasonable from the definition of $K_0(\caE)$ of an exact category $\caE$
as the free abelian group on the set of equivalence classes of $\caE$ modulo 
the relationship $B=A +C$ for any short exact sequence $A \to B \to C.$
In this section we prove an analogous additivity property for real algebraic K-theory.
Noe that the $\infty$-category of short exact sequences can be identified with $S(\caC)_2$.
Different from the classical additivity theorem in this real addivity theorem 
the role of the $\infty$-category $S(\caC)_2$ is taken by $S(\caC)_3$
that consists of diagrams
\begin{equation}\label{diagram}
\begin{tikzcd}
A \ar[r] & B \ar[r] \ar[d] & C \ar[d] \\
&  B/A \ar[r] & C/A \ar[d] \\
&& C/B
\end{tikzcd}
\end{equation}
In this case the functor $\gamma_2$ is replaced by the functor $\gamma_3$ 
that sends a diagram \ref{diagram} to the triple $(A, C/B, B/A).$
Assuming $K$-theory splits short exact sequences, we can build
$B$ from $A$ and $B/A$ and we can build $C/A$ from $B/A$ and $C/B$
to build $C$ from $A$ and $C/A$.

To prove this additivity theorem, we need to get some control over the class of morphisms in $\Wald_\infty^\gd$ that yield an equivalence on geometric realization in $\D\Wald_\infty^\gd$.
We separate a subclass of these morphisms, the (relative) real simplicial homotopy equivalences, that enjoys good stability properties (see \cref{subsec:real_s_htpies}).

\subsection{The statement and directions to read the proof}\label{subsec:add_statement_and_directions}

Although we do not yet count with all the necessary ingredients for it, we present the statement of the main theorem of this section. 

\begin{theorem}[Additivity]\label{thm:add}\label{gghjgfdd}
For every $\caD \in \Wald_\infty^\gd$
the map
$$S(\gamma_3)\colon S(S(\caD)_3) \to S(\widetilde{\caD \times \caD}) \times S(\caD)$$
in $\rs\Wald_\infty^\gd$ becomes an equivalence after geometric realization in $\D\Wald_\infty^\gd$, where $\widetilde{\caD\times \caD}$ is the $\Spc^{C_2}$-enriched cotensor $\caD^{C_2}$.
\end{theorem}

In what follows we describe the elements of this statement and/or where to find them. 

Let us define first the map $\gamma_3$. For every $\caD \in \Wald^\gd_\infty$ the morphism $[1] \simeq \{0 \to 1\} \subset [3]$ in $\Delta$ yields a map $$\ev_{01}\colon  S(\caD)_3 \to S(\caD)_1 \simeq \caD, \ (A \to B \to C) \mapsto A$$ in $\Exact_\infty,$ whilst the morphism $[1] \simeq \{1 \to 2 \} \subset [3]$ in $\Delta^{hC_2}$ 
yields a map $$\ev_{12}\colon  S(\caD)_3 \to S(\caD)_1 \simeq \caD, \ (A \to B \to C) \mapsto B/A$$ in $ \Wald^\gd_\infty.$ The map $\ev_{12}\colon S(\caD)_3 \to \caD $ 
in $\Wald^\gd_\infty $ and $\ev_{01}'\colon S(\caD)_3 \to \widetilde{\caD \times \caD}$ in $\Wald^\gd_\infty$, transpose to $\ev_{01}\colon S(\caD)_3 \to \caD$
in $\Wald_\infty $
give rise to the map 
$$\gamma_3\colon  S(\caD)_3 \to \widetilde{\caD \times \caD} \times \caD$$
in $\Wald^\gd_\infty $, which sends $\ (A \to B \to C) \mapsto (A, C/B, B/A)$ .

The $S_\bullet$ construction we refer to here is the real version of the usual $S_\bullet$-construction; we introduce this in \cref{subsec:S_construction}. The $\infty$-category $\D\Wald_\infty^\gd$ is the completion of $\Wald_\infty^\gd$ under sifted colimits. We define this completion here, but refer the reader to \cref{subsec:DWald} for its development.

\begin{defn}
We define the non-abelian derived $\infty$-category of $\Wald_\infty^\gd$, that we denote by  
$$\D\Wald_\infty^\gd \subset \Fun((\Wald_\infty^\gd)^\op, \Spc)$$ 
as the smallest full subcategory containing the representables and closed under small sifted colimits. 
\end{defn}

To prove this additivity theorem, we need to get some control over the class of morphisms in $\Wald_\infty^\gd$ that yield an equivalence on geometric realization in $\D\Wald_\infty^\gd$. We separate a subclass of these morphisms, the (relative) real simplicial homotopy equivalences, that enjoys good stability properties (see \cref{subsec:real_s_htpies}). The skeleton of the proof can be found in \cref{subsec:add_thm}, while many of the technical results and constructions are defered to \cref{subsec:lemmas_add}.

\subsection{The real Waldhausen $S_\bullet$-construction}\label{subsec:S_construction}

In this subsection we define a real version of Waldhausen's $S_\bullet$-construction defined by Barwick in \cite[Section 3]{barwick.wald} and prove some basic properties that will be paramount later on.  

\subsection*{$S_\bullet$-construction for exact $\infty$-categories}
\begin{definition} Let $\caC$ be an exact $\infty$-category and $n \geq 0$.
We define  $$S^e(\caC)_n \subset \Fun(\Ar([n]),\caC)$$ to be the full subcategory
of functors $A\colon \Ar([n]) \to \caC$ such that for every $ 0 \leq i \leq n$
the image $A_{i,i}$ is zero and for every 
$0 \leq i \leq j \leq \ell \leq k \leq n$ the induced square

\begin{equation}\label{pullbu}
\begin{tikzcd}
A_{i,\ell} \ar[r]\ar[d] &  A_{i,k} \ar[d]  \\
A_{j,\ell} \ar[r] &  A_{j,k}
\end{tikzcd}
\end{equation}
in $\caC$ is an exact square.	
\end{definition}

\begin{lemma} Let $\caC$ be an exact $\infty$-category, $n \geq 0$ and $A\colon \Ar([n]) \to \caC$ a functor such that for every  $0 \leq i \leq n$ the image $A_{i,i}$ is zero and for every $0 \leq i \leq j \leq \ell \leq k \leq n$ the square 

\begin{equation*}
\begin{tikzcd}
A_{i,\ell} \ar[r]\ar[d] &  A_{i,k} \ar[d]  \\
A_{j,\ell} \ar[r] &  A_{j,k} 
\end{tikzcd}
\end{equation*}

in $\caC$ is a pushout square and the top horizontal map is a cofibration.
Then for every $0 \leq i \leq j \leq \ell \leq k \leq n$ the left vertical map in the square  is a fibration.

\end{lemma}

\begin{proof}
Consider the commutative diagram 
\begin{equation*}
\begin{tikzcd}
A_{i,j} \ar[r] \ar[d] &A_{i,\ell} \ar[r]\ar[d] &  A_{i,k} \ar[d]  \\
0=A_{j,j} \ar[r]& A_{j,\ell} \ar[r] &  A_{j,k}.
\end{tikzcd}
\end{equation*}
By assumption the left-most square is a pushout square whose top horizontal map is a cofibration. Since the left vertical map of the left-most square is a fibration,
the middle vertical map of the diagram is a fibration, too.
\end{proof}

\begin{remark}
Let $\caC$ be an exact $\infty$-category and $n \geq 0.$
The canonical equivalence $$ \Fun(\Ar([n]),\caC)^\op \simeq \Fun(\Ar([n])^\op,\caC^\op) \simeq \Fun(\Ar([n]),\caC^\op)$$
restricts to an equivalence
$$ S^e(\caC)^\op_n \simeq S^e(\caC^\op)_n.$$
\end{remark}

The next lemma shows how we can see an object of $S^e(\caC)_n$ as a sequence of $n-1$ cofibrations, equivalently fibrations\textemdash meaning the extra quotient data, equivalently fiber data, is redundant. 

\begin{lemma}\label{rema}
Let $\caC$ an exact $\infty$-category and $n \ge 1$. 

\begin{enumerate}
\item The functor $$[n-1] \to \Ar([n]) \text{\ given by} \ i \mapsto(0,i+1) $$ induces an equivalence
$$S^e(\caC)_n \simeq \Fun^\mathrm{cof}([n-1],\caC),$$
where $\Fun^\mathrm{cof}([n-1],\caC) \subset \Fun([n-1],\caC)$ is the full subcategory of sequences of cofibrations.

\item The functor $$[n-1] \to \Ar([n]) \text{\ given by} \ i \mapsto (n-i-1,n)$$ induces an equivalence
$$S^e(\caC)_n \simeq \Fun^\mathrm{fib}([n-1],\caC),$$ where $\Fun^\mathrm{fib}([n-1],\caC) \subset \Fun([n-1],\caC)$ is the full subcategory of sequences of fibrations.
\end{enumerate}	
\end{lemma}

\begin{proof}
For the first statement, we observe that the induced functor $$S^e(\caC)_n \to\Fun^\mathrm{cof}([n-1],\caC)$$ is inverse to the functor that sends a sequence of cofibrations $A_1 \to \cdots \to A_n$ to the functor $\Ar([n]) \to \caC$ sending $(i,j) \mapsto A_j/A_i,$ where $A_0\coloneqq 0$. 

The second statement follows from the first one when considering the canonical equivalence $S(\caC^\op)_n^\op \simeq S(\caC)_n$.
\end{proof}

The $\infty$-category $\Fun(\Ar([n]),\caC)$ is an exact $\infty$-category
with objectwise structure, however this exact structure is not inherited by $S^u(\caC)_n$.
In what follows we endow $S^u(\caC)_n$ with the structure of an exact $\infty$-category with fewer cofibrations and fibrations than what we would obtain from the objectwise one.

\begin{definition}\label{cofibrations_in_Se} Let $\caC$ be an exact $\infty$-category.
\begin{itemize}
\item We set a map $A \to B$ in $S^e(\caC)_n$ to be a cofibration if 
for every $1 \leq i < j \leq n $ the canonical maps
$ A_{0,1} \to B_{0,1} $ and
$ A_{0,j} \coprod_{A_{0,i}} B_{0,i} \to B_{0,j}$ are cofibrations.

\item We set a map $A \to B$ in $S^u(\caC)_n$ to be a fibration if for every $1 \leq i < j \leq n $ the canonical maps $A_{n-1,n} \to B_{n-1,n}$ and
$ A_{i,n} \to B_{i,n} \prod_{B_{j,n}} A_{j,n} $ are fibrations.
\end{itemize}
\end{definition}

\begin{proposition}\label{exact_structure_on_S}
Let $\caC$ be an exact $\infty$-category and $n \geq 0$. The $\infty$-category $S^e(\caC)_n$ with the choice of cofibrations and fibrations as above forms an exact $\infty$-category and the embedding $S^e(\caC)_n \subset \Fun(\Ar([n]),\caC)$ is exact. 
\end{proposition}

\begin{proof}
In view of the canonical equivalence $$S^e(\caC)^\op_n \simeq S^e(\caC^\op)_n$$
it is enough to show that the pushout in $\Fun(\Ar([n]),\caC)$ of a cofibration in $S^e(\caC)_n$ along any map in $S^e(\caC)_n$ belongs to and is a cofibration in $S^e(\caC)_n$, that the pushout in $\Fun(\Ar([n]),\caC) $ of a fibration in $S^e(\caC)_n$ along a cofibration in $S^e(\caC)_n$ is again a fibration in $S^e(\caC)_n$, and that the embedding $S^e(\caC)_n \subset \Fun(\Ar([n]),\caC)$ preserves cofibrations.

We will begin by showing the second part of the statement. First note that we can rephrase the definition of cofibration in $S^e(\caC)_n$ by saying that a map $A \to B$ in $S^e(\caC)_n$ is a cofibration if and only if for every $0\leq i < j \leq n $ the induced map $A_{0,j}\coprod_{A_{0,i}} B_{0,i} \xrightarrow{\alpha} B_{0,j}$ is a cofibration.

Let $A \to B$ be a map in $\Fun(\Ar([n]),\caC)$, instead of $S^e(\caC)_n$, such that the above condition holds. Then for every $0\leq k \leq i$ the pushout $\psi\coloneqq \alpha \coprod_{B_{0,i}} B_{k,i}$ given by 
$$A_{k,j} \coprod_{A_{k,i}} B_{k,i} \simeq A_{0,j} \coprod_{A_{0,i}} B_{k,i} \simeq  A_{0,j} \coprod_{A_{0,i}} B_{0,i} \coprod_{B_{0,i}} B_{k,i} \to $$
$$ B_{0,j} \coprod_{B_{0,i}} B_{k,i} \simeq B_{k,j} $$ is a cofibration. 

Now, for $i=k$ the map $\alpha \coprod_{B_{0,i}} B_{k,i}$ identifies with the map $A_{k,j} \to B_{k,j}$ and the latter is thus a cofibration. This shows that the  embedding $S^e(\caC)_n \subset \Fun(\Ar([n]),\caC)$ preserves cofibrations, since the cofibrations in $\Fun(\Ar([n]),\caC)$ are given by levelwise cofibrations.

We now move to prove that pushouts of cofibrations exist and are cofibrations.
Let
\begin{equation}\label{pullbur}
\begin{tikzcd}
A \ar[r]\ar[d] & A' \ar[d]  \\
B \ar[r] &  B'
\end{tikzcd}
\end{equation}
be a pushout square in $\Fun(\Ar([n]),\caC)$, where $A,A',B \in S^u(\caC)_n$ and the left vertical map is a cofibration in $S^e(\caC)_n$.

Let us consider, for every $0 \leq k \leq i < j \leq n$, the canonical commutative diagram below
\begin{equation}\label{pullbub}
\begin{tikzcd}
A_{k,j} \ar[r]\ar[d] & A'_{k,j} \ar[d]  \\ 
B_{k,i}\coprod_{A_{k,i}}A_{k,j}\ar[r]\ar[d, "\psi"'] &  B_{k,i} \coprod_{A_{k,i}} A'_{k,j}\simeq   B'_{k,i} \coprod_{A'_{k,i}} A'_{k,j}\ar[d, "\psi'"]  \\
B_{k,j}  \ar[r] & B'_{k,j}.
\end{tikzcd}
\end{equation}
Since both the upper and outer square are pushout squares, the lower square is so too. Moreover, since $\psi$ is a cofibration, $\psi'$ is a cofibration. 

Because $A' \in S^u(\caC)_n$, the map $B'_{k,i} \to B'_{k,j}$ is a cofibration: indeed, if $A' \in S^e(\caC)_n$, then $A'_{k,i} \to A'_{k,j}$ is a cofibration in $\caC$ and therefore its pushout $ B'_{k,i} \to A'_{k,j} \coprod_{A'_{k,i}} B'_{k,i}$ is a cofibration. Thus the canonical map $ B'_{k,i} \to B'_{k,j}$ is a cofibration as it factors as $$ B'_{k,i} \to A'_{k,j} \coprod_{A'_{k,i}} B'_{k,i} \xrightarrow{\psi} B'_{k,j} .$$

Moreover square (\ref{pullbu}) for $B'$ is a pushout square, since pushouts commute with pushouts. Thus $B' \in S^e(\caC)_n$. Hence square (\ref{pullbur}) is actually pushout square in
$S^e(\caC)_n$, where the right vertical map is a cofibration in $ S^e(\caC)_n$.

Finally, we work on the third condition of \cref{def:exact_infty_cats}.

We consider the square (\ref{pullbur}), with the additional assumption that the top horizontal map is a fibration in $S^e(\caC)_n$. Dually to what we have shown above, every fibration is objectwise a fibration. Thus the pushout square (\ref{pullbur}) is a pullback square in $\Fun(\Ar([n]),\caC)$ and so a pullback square in $S^e(\caC)_n$.

It remains to show that the map $B \to B'$ is a fibration. By looking at the bottom square of (\ref{pullbub}) we know that the commutative square
\begin{equation*}
\begin{tikzcd}
B_{0,i}\coprod_{A_{0,i}}    
A_{0,j}\ar[r]\ar[d] & B'_{0,i} \coprod_{A'_{0,i}} A'_{0,j}\ar[d]  \\
B_{0,j} \ar[r] & B'_{0,j}
\end{tikzcd}
\end{equation*}
is a pushout square, from which we can deduce that the canonical map
\begin{equation}\label{rho}
B_{0,j}/(B_{0,i}\coprod_{A_{0,i}}    
A_{0,j}) \to B'_{0,j}/ (B'_{0,i} \coprod_{A'_{0,i}} A'_{0,j}) \end{equation}
is an equivalence.
The pushout square
\begin{equation*}
\begin{tikzcd}
B_{0,n} / (A_{0,j} \coprod_{A_{0,i}} B_{0,i}) \ar[r]\ar[d] &  B'_{0,n} / (A'_{0,j} \coprod_{A'_{0,i}} B'_{0,i}) \ar[d]  \\
B_{0,n} / B_{0,j} \ar[r] &  B'_{0,n} / B'_{0,j},
\end{tikzcd}
\end{equation*}
where the vertical maps are fibrations, induces on fibers the equivalence (\ref{rho}). 
We can write the latter square as
\begin{equation*}
\begin{tikzcd}
(A_{0,n} \coprod_{A_{0,n}} B_{0,n}) / (A_{0,j} \coprod_{A_{0,i}} B_{0,i}) \ar[r]\ar[d] &  (A'_{0,n} \coprod_{A'_{0,n}} B'_{0,n}) / (A'_{0,j} \coprod_{A'_{0,i}} B'_{0,i}) \ar[d]  \\
B_{0,n} / B_{0,j} \ar[r] &  B'_{0,n} / B'_{0,j}.
\end{tikzcd}
\end{equation*}
Therefore this square identifies with the pushout square
\begin{equation}\label{uth}
\begin{tikzcd}
A_{j,n} \coprod_{A_{i,n}} B_{i,n} \ar[r]\ar[d] & A'_{j,n} \coprod_{A'_{i,n}} B'_{i,n} \simeq A'_{j,n} \coprod_{A_{j,n}} A_{j,n} \coprod_{A_{i,n}} B_{i,n} \ar[d]  \\
B_{j,n} \ar[r] &  B'_{j,n} \simeq A'_{j,n} \coprod_{A_{j,n}} B_{j,n}.
\end{tikzcd}
\end{equation}
Note that the bottom horizontal map of square (\ref{uth}) is a fibration 
because $A_{j,n} \to A'_{j,n}$ is a fibration and $A_{j,n} \to A_{j,n} \coprod_{A_{i,n}} B_{i,n}$ is a cofibration (as $A_{i,n} \to B_{i,n}$
is a cofibration).
Similarly, the top horizontal map of square (\ref{uth}) is a fibration 
since $A_{j,n} \to A'_{j,n}$ is a fibration and $A_{j,n} \to B_{j,n}$ is a cofibration.
Thus square (\ref{uth}) is a pullback square, therefore the canonical map below is an equivalence.
$$\alpha\colon B'_{i,n} \prod_{(A'_{j,n}\coprod_{A'_{i,n}} B'_{i,n})} (A_{j,n} \coprod_{A_{i,n}} B_{i,n}) \to B'_{i,n} \prod_{ B'_{j,n}} B_{j,n}$$

The canonical map
$$ \beta\colon (A'_{i,n} \prod_{A'_{j,n}} A_{j,n}) \coprod_{A_{i,n}} B_{i,n} \to $$$$ (A'_{i,n}\coprod_{A_{i,n}} B_{i,n}) \prod_{(A'_{j,n}\coprod_{A_{i,n}} B_{i,n})} (A_{j,n} \coprod_{A_{i,n}} B_{i,n}) $$
is an equivalence as well, since the map $A_{i,n} \to B_{i,n}$ is a cofibration, both maps in the pullback $A'_{i,n} \prod_{A'_{j,n}} A_{j,n}$ are fibrations,
and all maps from $A_{i,n}$ to the factors of the pullback are fibrations.

Because the map $A_{i,n} \to B_{i,n}$ is a cofibration, the pushout
$$\gamma\colon B_{i,n} \simeq A_{i,n} \coprod_{A_{i,n}} B_{i,n} \to (A'_{i,n} \prod_{A'_{j,n}} A_{j,n}) \coprod_{A_{i,n}} B_{i,n} $$
is a fibration.

All this proofs that the canonical map
$B_{i,n} \to B'_{i,n} \prod_{ B'_{j,n}} B_{j,n}$
is a fibration, as it factors as
$$B_{i,n} \simeq A_{i,n} \coprod_{A_{i,n}} B_{i,n} \xrightarrow{\gamma} (A'_{i,n} \prod_{A'_{j,n}} A_{j,n}) \coprod_{A_{i,n}} B_{i,n} \xrightarrow{\beta} $$$$ (A'_{i,n}\coprod_{A_{i,n}} B_{i,n}) \prod_{(A'_{j,n}\coprod_{A_{i,n}} B_{i,n})} (A_{j,n} \coprod_{A_{i,n}} B_{i,n}) \simeq $$
$$B'_{i,n} \prod_{(A'_{j,n}\coprod_{A'_{i,n}} B'_{i,n})} (A_{j,n} \coprod_{A_{i,n}} B_{i,n}) \xrightarrow{\alpha} B'_{i,n} \prod_{ B'_{j,n}} B_{j,n}.$$
\end{proof}

\subsection*{The real $S_\bullet$-construction}

We recall that by \cref{cor:Waldgd_is_cotensored_with_Catgd} the real $\infty$-category $\Wald_\infty^\gd$ is cotensored over $\Cat_\infty^\gd$. In fact, for a small $\infty$-category $K$ with genuine duality and a small Waldhausen $\infty$-category $\caC$ with genuine duality, the cotensor $\caC^K$ has underlying $\infty$-category $\Fun(K, \caC)$ and objectwise fibrations and cofibrations.

We begin by endowing the exact $\infty$-category $S^e(\caC)_n$ with a genuine duality.

\begin{proposition}\label{Se_has_gd}
Let $\caC$ be a Waldhausen  $\infty$-category with genuine duality and let $n \geq 0$.
The exact subcategory $S^e(\caC)_n \subset \Fun(\Ar([n]),\caC)$ is closed under 
the genuine duality of the functor $\infty$-category.
\end{proposition}

\begin{proof}
We need to show three things. Firstly, that the functor that gives the duality on $\Fun(\Ar([n]),\caC)$ restricts to  $S^e(\caC)_n$. Secondly, we must check that the genuine duality inherited is compatible with the exact structure on $S^e(\caC)_n$ defined in \cref{cofibrations_in_Se} \textemdash this amounts to prove that the restricted duality functor on  $S^e(\caC)_n$ is exact and that axiom (2) of \cref{Waldgd} holds.

Let us consider $A\in S^e(\caC)_n$ and let $B\in \Fun(\Ar([n]),\caC)$ its dual. Then for any $0\leq i\leq j \leq n$ we have $B_{i, j} \simeq A^\dual_{n-j, n-i}$. Thus $B_{i, i}\simeq A_{n-i, n-i} \simeq 0 $ and for any $0 \leq i \leq j\leq k\leq  \ell \leq n$ the square

\begin{equation}
\begin{tikzcd}
B_{i, \ell} \ar[r]\ar[d] &  B_{i,k} \ar[d]  \\
B_{j, \ell} \ar[r] &  B_{j,k}
\end{tikzcd}
\end{equation}
is an exact square as it identifies with the exact square
\begin{equation}
\begin{tikzcd}
A^\dual_{n-\ell,n-i} \ar[r]\ar[d] &  A^\dual_{n-k,n-i} \ar[d]  \\
A^\dual_{n-\ell, n-j} \ar[r] &  A^\dual_{n-k, n-j}.
\end{tikzcd}
\end{equation}
We conclude from this that $S^e(\caC)$ is closed under the relevant duality.

Now, to prove that the duality functor is exact, first note that by definition given a cofibration $A \to A'$ in $S^e(\caC)_n$, for every $0 \leq i \leq j \leq n $ the canonical map $ A_{0,j} \coprod_{A_{0,i}} A'_{0,i} \to A'_{0,j}$ is a cofibration, too. Let $B' \to B$ be the dual map. Then for every $0 \leq i \leq j \leq n $ the canonical map
$ B'_{i, n} \to B_{i, n} \prod_{B'_{j, n}} B_{j, n} $ is a fibration
as it identifies with the induced map
$$ A'^\dual_{0, n-i} \to A^\dual_{0, n-i} \prod_{A'^\dual_{0, n-j}} A^\dual_{0, n-j} \simeq (A_{0, n-i} \coprod_{A'_{0, n-j}} A_{0, n-j})^\dual.$$

The last part follows from the fact that the inclusion $S^e(\caC)_n \subset \Fun(\Ar([n]),\caC) $ is exact and so preserves exact squares.
\end{proof}

To define the real $S_\bullet$-construction we now want to assemble the categories $S^e(\caC)_n$ into a real simplicial object.

\begin{notation}
Let $\Ar\colon \underline{\Delta} \to \Cat_\infty^\gd$ the real functor 
$$\underline{\Delta} \subset \Cat_\infty^\gd \xrightarrow{\Hom_{ \Cat_\infty^\gd }([1], -)}  \Cat_\infty^\gd.$$

\end{notation}

\begin{corollary}

The real functor $\underline{\Delta}^\op \times \Wald_\infty^\gd  \to \Wald_\infty^\gd$ that sends $ (K,\caC) \mapsto \caC^{\Ar(K)}$ induces a real functor 
\[
\begin{tikzcd}[row sep=tiny]
S\colon \underline{\Delta}^\op \times \Wald_\infty^\gd \ar[r, "\sigma"]     &\Wald_\infty^\gd\\
\end{tikzcd}
\]
that sends $(K,\caC) \in \underline{\Delta}^\op \times \Wald_\infty^\gd$
to the genuine duality on $S^e(\caC)_n$ inherited from $\caC^{\Ar([n])}$
via \cref{Se_has_gd}.
\end{corollary}

Such $\sigma$ is the transpose of a real functor $S_\bullet\colon \Wald_\infty^\gd \to \mathrm{rs}\Wald_\infty^\gd $, which we call the real Waldhausen $S_\bullet$-construction. From now on we will abuse notation by calling by $S$ both the real and usual $S_\bullet$-construction.

We will later use the following lemma.
\begin{lemma}\label{SSSSSS} Let $\caC, \caD \in \Wald_\infty^\gd$ and $[n] \in \underline{\Delta}.$ Then the canonical equivalence $$ \Hom_{\Wald_\infty^\gd }(\caC, \caD^{\Ar(n)}) \simeq \Hom_{\Wald_\infty^\gd }(\caC, \caD)^{\Ar(n)} $$ 
restricts to an equivalence
$$ \Hom_{\Wald_\infty^\gd }(\caC, S(\caD)_n) \simeq S(\Hom_{\Wald_\infty^\gd }(\caC, \caD))_n. $$
\end{lemma}

\begin{proof}
	
It is enough to check the statement after forgetting genuine dualities. We begin by showing the statement after forgetting farther, this is, forgetting even Waldhausen structures, and then 
care for the Waldhausen structures in a second step.

Then we want to start by showing that an exact functor $F\colon \caC \to \caD^{\Ar(n)} $ factors through $S(\caD)_n$ if and only if its corresponding object $A \in \Fun(\Ar(n), \Fun^\mathrm{exact}(\caC, \caD)) $ belongs to $ S(\Fun^\mathrm{exact}(\caC, \caD))_n.$

This follows from the fact that a morphism $U \to V$ in $\Fun^\mathrm{exact}(\caC, \caD)$
is a cofibration if and only if it is objectwise a cofibration and for every cofibration $X \to Y$ in $\caC$ the canonical map $$ U(Y) \coprod_{U(X)} V(X) \to V(Y)$$ is a cofibration in $\caD.$ So we get a canonical equivalence of $\infty$-categories 
$$\Fun^\mathrm{exact}(\caC,  S(\caD)_n) \simeq S(\Fun^\mathrm{exact}(\caC, \caD))_n$$ and now we want to identify both Waldhausen structures. 

A morphism $F \to G$ in $\Fun^\mathrm{exact}(\caC, S(\caD)_n)$
is a cofibration if and only if it is objectwise a cofibration and for every cofibration $X \to Y$ in $\caC$ the canonical map $$ F(Y) \coprod_{F(X)} G(X) \to G(Y)$$ is a cofibration in $S(\caD)_n$, which means that it is objectwise a cofibration and for every $ i \leq j \leq k$ the canonical map $$ G(Y)_{i,j} \coprod_{(F(Y)_{i,j} \coprod_{F(X)_{i,j}} G(X)_{i,j}) } (F(Y)_{i,k} \coprod_{F(X)_{i,k}} G(X)_{i,k}) \to G(Y)_{i,k}$$ is a cofibration in $\caD.$ 

A morphism $A \to B$ in $ S(\Fun^\mathrm{exact}(\caC, \caD))_n$ is a cofibration
if and only if it is objectwise a cofibration and for every $ i \leq j \leq k$
the canonical map $$ B_{i,j} \coprod_{A_{i,j}} A_{i,k} \to B_{i,k}$$ is a cofibration
in $\Fun^\mathrm{exact}(\caC, \caD), $ which means that it is objectwise a cofibration and for every cofibration $X \to Y$ in $C$ the canonical map $$ (B_{i,j}(Y) \coprod_{A_{i,j}(Y)} A_{i,k}(Y)) \coprod_{(B_{i,j}(X) \coprod_{A_{i,j}(X)} A_{i,k}(X))} B_{i,k}(X) \to B_{i,k}(Y)$$ is a cofibration in $\caD.$ So the statement follows.
\end{proof}

An important feature of the real $S_\bullet$-construction is that, after taking its edgewise subdivision,
we obtain a Segal object. We will prove this fact in what follows, along with other results that are fundamental for the proof of the Additivity Theorem.

\begin{notation} We denote by $\mathrm{Amb}(\caC) \subset \caC^{\widetilde{[1] \times [1]}}$ the full exact subcategory with genuine duality spanned by the exact squares.
\end{notation}

The embedding of $\infty$-categories with duality $\widetilde{[1] \times [1]} \subset \Ar([3])$ sending $(i,j) \mapsto (i,j+2)$ induces a map of exact $\infty$-categories with genuine duality as below $$\kappa\colon S(\caC)_3 \to \mathrm{Amb}(\caC) \subset {\caC^{\widetilde{[1] \times [1]}}}.$$

\begin{proposition}\label{ehtru64e}
Let $\caC$ be an exact $\infty$-category with genuine duality and $n \geq 1$. The canonical map $\theta$ given by the following composition
$$S(\caC)_{[n]*[n]^\op} \to S(\caC)_{[1]*[1]^\op} \times_{ S(\caC)_{[0]*[0]^\op} }\cdots  \times_{ S(\caC)_{[0]*[0]^\op} }S(\caC)_{[1]*[1]^\op}\simeq $$
$$ S(\caC)_{[3]} \times_{ \caC} \cdots  \times_{\caC} S(\caC)_{[3]} \xrightarrow{\kappa \times_\caC \kappa \cdots \times_\caC \kappa} \mathrm{Amb}(\caC) \times_{ \caC} \cdots  \times_{ \caC} \mathrm{Amb}(\caC) $$ where the first map is the Segal map, induced by the maps
$[1] \simeq \{i-1,i\} \subset [n]$ for $1 \leq i \leq n$, is an equivalence in $ \Wald^\gd_\infty$.
\end{proposition}

\begin{proof}
We first prove that $\theta$ induces an equivalence on underlying exact $\infty$-categories. One can show that the map $\theta$ sends $A$ to the sequence of exact squares $\sigma_A^i$ in $\caC$, for $0\leq i\leq n$:
\begin{equation*}
\begin{tikzcd}
A_{i-1,2n+1-i} \ar[r] \ar[d] &  A_{i-1,2n+2-i} \ar[d]\\
A_{i,2n+1-i} \simeq A_{i-1,2n+1-i}/ A_{i-1,i} \ar[r] & A_{i,2n+2-i} \simeq A_{i-1,2n+2-i}/ A_{i-1,i}.
\end{tikzcd}
\end{equation*} 

It is easy to check that an inverse of the exact functor underlying $\theta$ is given by the exact functor sending a sequence of exact squares in $\caC$, for $1 \leq i \leq n$, as depicted below
\begin{equation*}
\begin{tikzcd}
A^i \ar[r] \ar[d] & B^i \ar[d]  \ar[d]\\
C^i \ar[r] & D^i
\end{tikzcd}
\end{equation*}  to the commutative diagram in $\caC$ as below.
\begin{equation*}
\begin{tikzcd}
\cdots \ar[r]&\ast \ar[d]\ar[r]& A^1 \ar[r] \ar[d] & B^1 \ar[d]  \ar[d]\\
&A^2 \ar[r] \ar[d] & C^1=B^2 \ar[r]\ar[d] & D^1 \ar[d] \\
& C^2=B^3 \ar[r]&D^2 \ar[r]& \ast \ar[d]\\
&&&\vdots
\end{tikzcd}
\end{equation*}

It remains to show that $\theta$ induces an equivalence on genuine refinements.

Let us now consider the genuine refinements $H' \to \caH^\mathrm{nd}(S(\caC)_{[n]*[n]^\op})$ and  $H'' \to \caH^\mathrm{nd}(\caC^{\widetilde{[1] \times [1]}})$ together with an object $A \in  \caH^\mathrm{nd}(S(\caC)_{[n]*[n]^\op})$. We want to see that the induced map $$\varrho\colon H'_A \to H''_{\sigma^1_A} \times_{H_{A_{1,2n}}} \cdots \times_{H_{A_{n-1,n+2}}} H''_{\sigma^{n}_A}$$ is an equivalence.
For that, using \cref{htwdzfd}, we observe that 
$$H''_{\sigma^{i}_A} \simeq H_{A_{i,2n+1-i} } \times_{H_{ A_{i-1,2n+1-i}} } H_{ A_{i-1,2n+2-i}},$$

$$H'_A \simeq  H_{ A_{0,2n+1}} \times_{ H_{ A_{0,2n}} } H_{ A_{1,2n}}\times_{ H_{ A_{1,2n-1}} } H_{A_{2,2n-1}} \cdots \times_{ H_{A_{n-2,n+2}} } H_{A_{n-1,n+2}} \times_{ H_{A_{n-1,n+1}} } H_{A_{n,n+1}} .$$

The map $\varrho$ identifies with the canonical equivalence
$$  H_{ A_{0,2n+1}} \times_{ H_{ A_{0,2n}} } H_{ A_{1,2n}}\times_{ H_{ A_{1,2n-1}} } H_{A_{2,2n-1}} \times_{H_{A_{2,2n-2}}} \cdots \times_{ H_{A_{n-1,n+1}} } H_{A_{n,n+1}}  \simeq $$
$$ ( H_{A_{1,2n}} \times_{H_{ A_{0,2n}} } H_{ A_{0,2n+1}}) \times_{H_{A_{1,2n}}} \cdots \times_{H_{A_{n-1,n+2}}} (H_{A_{n,n+1} } \times_{H_{ A_{n-1,n+1}} } H_{ A_{n-1,n+2}}).$$
\end{proof}

The previous result, when applied to $n=1$, gives us the following.

\begin{corollary}\label{corqa}
Let $\caC$ be in $\Wald^\gd_\infty$. The map
$$\kappa\colon S(\caC)_3 \to \mathrm{Amb}(\caC)$$ of exact $\infty$-categories with genuine duality is an equivalence.
\end{corollary}

We also deduce the following.

\begin{proposition}
Let $\caC \in \Wald^\gd_\infty$.
Then $S (\caC) \circ e\colon \Delta^\op \to \Wald^\gd_\infty$ is a Segal object.	
\end{proposition}

We now move to introduce a fundamental ingredient of the Additivity Theorem, that will serve as motivation for the coming sections.

For every $\caD \in \Wald^\gd_\infty$ the morphism $[1] \simeq \{0 \to 1\} \subset [3]$ in $\Delta$ yields the following map in $\Exact_\infty$
\begin{equation}
\begin{tikzcd}[row sep=tiny]
S(\caD)_3\ar[r]        & S(\caD)_1 \simeq \caD\\
\ (A \to B \to C)\ar[r, mapsto]     &A.
\end{tikzcd}
\end{equation}

On the other hand, the morphism $[1] \simeq \{1 \to 2 \} \subset [3]$ in $\Delta^{hC_2}$ yields a map in $ \Wald^\gd_\infty$ as below
\begin{equation}
\begin{tikzcd}[row sep=tiny]
S(\caD)_3\ar[r]        & S(\caD)_1 \simeq \caD\\
\ (A \to B \to C)\ar[r, mapsto]     & B/A.
\end{tikzcd}
\end{equation}

The adjoint map $ S(\caD)_3 \to \widetilde{\caD \times \caD}$ in $\Wald^\gd_\infty$ of the first map together with the second map give rise to a map in $\Wald^\gd_\infty $

\begin{equation}\label{const:gamma_3}
\gamma_3\colon S(\caD)_3 \to \widetilde{\caD \times \caD} \times \caD, \ (A \to B \to C) \mapsto (A, C/B, B/A).
\end{equation}

One could naively attempt to show that the map of real simplicial objects in $\Wald^\gd_\infty$ as below
$$S(\gamma_3)\colon S(S(\caD)_3) \to S(\widetilde{\caD \times \caD}) \times S(\caD)$$
becomes an equivalence after geometric realization in $\Wald^\gd_\infty$. This is not what we actually want, since the forgetful functor $\Wald^\gd_\infty\to\Spc^{C_2}$ does not preserve geometric realizations and therefore this would not lead to the usual description of $K$-theory. Such failure is corrected by considering the completion of $\Wald_\infty^\gd$ under sifted colimits.

\subsection{From $\Wald_\infty^\gd$ to $ \D\Wald_\infty^\gd $.}\label{subsec:DWald}

To state the addivity theorem we embed $\Wald_\infty^\gd$ into its derived non-abelian $\infty$-category $\D\Wald_\infty^\gd$ \textemdash this is, its completion under sifted colimits \textemdash and consider the geometric realization therein. The reason for taking this deviation is that the geometric realization in $\D\Wald_\infty^\gd$, which is different from that in $\Wald_\infty^\gd$, is sent to the geometric realization in $\Spc$ under the forgetful functor $\D\Wald_\infty^\gd \to \Spc$, a crucial property when set to deduce additivity for real $K$-theory from the additivity theorem as stated in \cref{thm:add}.

We show in this section that the derived non-abelian $\infty$-category $\D\Wald_\infty^\gd$ is genuine preadditive (\cref{fghfghjj}).

The following nomenclature is motivated by the fact that for any ring R the derived $\infty$-category of a ring $R$, $\Mod_{\mathrm{H}(\mathrm{R})}$, is the smallest full subcategory containing the projective $R$-modules and closed under small sifted colimits.

\begin{definition}Given a small $\infty$-category $\caC$ with finite coproducts, we define its non-abelian derived $\infty$-category $$\D\caC \subset \Fun(\caC^\op, \Spc)$$ 
as the smallest full subcategory containing the representables and closed under small sifted colimits. 
\end{definition}

The next lemma is a direct consequence of \cite[Proposition 2.48]{heine-enriched2}, when considering completion under small sifted colimits.

\begin{lemma}\label{lem:functors_extend_to_DWald}
Let $\caE$ be a real $\infty$-category with sifted colimits and $\caC$ a small real $\infty$-category as in the definition above. Then restriction along the real embedding $\caC \subset D\caC$ induces an equivalence
$$ \Fun^\Sigma_{\Spc^{C_2}}(\D\caC,\caE) \to \Fun_{\Spc^{C_2}}(\caC,\caE),$$
where $\Sigma$ indicates that we are taking sifted colimits preserving functors.
\end{lemma}

By \cref{lem:functors_extend_to_DWald}, for any real $\infty$-category $\caE$ we can uniquely extend real functors $\Wald_\infty^\gd\to \caE$ to sifted colimits preserving real functors $\D\Wald_\infty^\gd\to \caE$.

\begin{remark}\label{fhjkllj}
Note that $\D\caC$ coincides with the full subcategory $\caB$ of $\Fun(\caC^\op, \Spc) $ spanned by the presheaves that preserve finite products. Indeed, $\caB$ is closed under small sifted colimits in $\Fun(\caC^\op, \Spc)$, as such commute with finite products. Therefore we know that $\D\caC \subset \caB$. 

On the other hand, $\caB \subset \Fun(\caC^\op, \Spc)$ is a localization, and then, by \cite[ Lemma 5.5.8.13]{lurie.HTT}, every object of $\caB$ is a geometric realization of a simplicial object in $\caB$ whose terms are small filtered colimits of representables \textemdash since $\caC$ is closed under finite coproducts in $\caP^\Sigma(\caC)$, as for every $F \in \caP^\Sigma(\caC) $ the presheaf $\map_{ \caP(\caC)}(y(-),F) \simeq F$ preserves finite products.
\end{remark}

\begin{remark}
The embedding $\caC \subset \D\caC$ preserves limits and cotensors with $C_2$, see \cref{Venr:Yoneda_embedding_preserves_limits_and_cotensors}.
\end{remark}

Note that $\D\caC $ is generally not small. However, by \cref{fhjkllj}, it is an accessible localization of $ \Fun(\caC^\op, \Spc)$ (and so presentable) if $\caC$ admits finite coproducts. In practice, we will usually form $\D\caC$ after enlarging the universe for a large $\infty$-category $\caC$ and then replace $\Spc$ by $\widehat{\Spc}$ in the definition of $\D(\caC).$ In this case,  $\D\caC $ is very large and a localization of $ \Fun(\caC^\op, \widehat{\Spc}).$

The aim of this subsection is to show the following result. 

\begin{proposition}
Let $\caC$ be a small real $\infty$-category that admits finite coproducts and tensors with $C_2$, the following $\infty$-categories constructed from $\caC$ are equivalent.
\begin{enumerate}
\item $\D(\caC)$
\item The full subcategory of $\Fun(\caC^\op, \Spc) $ spanned by the presheaves that preserve finite products.
\item The smallest full subcategory of $\Fun_{\Spc^{C_2}}(\caC^\op, \Spc^{C_2})$ that contains the representables and is closed under small sifted colimits.
\item The full subcategory of $\Fun_{\Spc^{C_2}}(\caC^{op},\Spc^{C_2})$ spanned by the preasheaves that preserve finite products and cotensors with $C_2$, that we denote by $\caP^\Sigma(\caC)$.
\end{enumerate}
\end{proposition}

\begin{proof}
The equivalence between (1) and (2) is the content of \cref{fhjkllj}. That the $\infty$-categories in (1) and (3) are equivalent is shown in \cref{dghjlkj}. Finally, the equivalence between (3) and (4) is proven in \cref{dfghkkll}.
\end{proof}

We embark to prove such results.

Let $\caV$ be a monoidal $\infty$-category, $\caC$ a $\caV$-enriched $\infty$-category and $\caK\subset\Cat_\infty$ a full subcategory.
We denote by $ \caP_\caV^\caK(\caC)$ the smallest full subcategory of $  \caP_\caV(\caC)\coloneqq \Fun_{\caV}(\caC^\op,\caV)$ containing the representables and closed under colimits indexed by $\infty$-categories that belong to $\caK$. When $\caV=\Spc$, we drop $\caV$ from the notation.

\begin{lemma}\label{dghjlkj}
Let $\caV$ be a monoidal $\infty$-category, $\caC$ a $\caV$-enriched $\infty$-category and $\caK\subset\Cat_\infty$ a full subcategory. If the functor $\map_\caV(1,- ) \colon \caV \to \Spc$ preserves $\caK$-indexed colimits, the forgetful functor $\psi\colon \caP_\caV(\caC) \to \Fun(\caC^\op, \caV) \xrightarrow{ \map_\caV(1,- )_\ast } \caP(\caC)$ restricts to an equivalence
$$ \caP_\caV^\caK(\caC) \simeq \caP^\caK(\caC) .$$
\end{lemma}

\begin{proof}
We know that the forgetful functor $\psi$ commutes with Yoneda embeddings. The functor $ \caP_\caV(\caC) \to \Fun(\caC^\op, \caV)$ preserves small colimits. So by assumption the functor $\psi$ preserves $\caK$-indexed colimits and so restricts to a functor $ \phi\colon \caP_\caV^\caK(\caC) \to \caP^\caK(\caC)$.

The functor $\phi$ is fully faithful as by definition $ \caP_\caV^\caK(\caC)$ is the smallest full subcategory of $ \caP_\caV(\caC)$ containing the representables and closed under $\caK$-indexed colimits. Hence the essential image of $\phi$ is a full subcategory of $ \caP(\caC)$ containing the  representables and closed under $\caK$-indexed colimits. Form this we conclude that $\phi$ is an equivalence.
\end{proof}

The following lemma guarantees that the $\infty$-category 
$\D\caC $ is not only an accessible localization of $ \Fun(\caC^\op, \Spc)$ but also a real accessible localization of $ \Fun_{\Spc^{C_2}}(\caC^\op, \Spc^{C_2}).$ In particular, this implies that $\D\caC$ is presentable.

\begin{lemma}\label{dfghkkll}Let $\caC$ be a small real $\infty$-category that admits finite coproducts and tensors with $C_2$. Then the full subcategory 
$$\caP^\Sigma(\caC) \subset \Fun_{\Spc^{C_2}}(\caC^\op, \Spc^{C_2})$$
 spanned by the presheaves that preserve finite products and cotensors with $C_2$, is the smallest full subcategory of $\Fun_{\Spc^{C_2}}(\caC^\op, \Spc^{C_2})$ that contains the representables and is closed under small sifted colimits.

\end{lemma}
\begin{proof}
Note that the $\infty$-category $\caP^\Sigma(\caC)$ contains the representable presheaves and is closed under small sifted colimits in $ \Fun_{\Spc^{C_2}}(\caC^\op, \Spc^{C_2})$; the latter fact follows from the fact that small sifted colimits commute with finite products and by \cref{dfghjgl} also with cotensors with $C_2.$ We can say then that $\caP^\Sigma(\caC)$ contains the smallest full subcategory of $\Fun_{\Spc^{C_2}}(\caC^\op, \Spc^{C_2}) $ that contains the representables and is closed under small sifted colimits.

The forgetful functor $\Fun_{\Spc^{C_2}}(\caC^\op, \Spc^{C_2}) \to \Fun(\caC^\op, \Spc^{C_2})$ is monadic, thus also its restriction $\caP^\Sigma(\caC) \to \Fun(\caC^\op, \Spc^{C_2})$ is monadic since the full subcategory inclusion $\caP^\Sigma(\caC) \subset \Fun_{\Spc^{C_2}}(\caC^\op, \Spc^{C_2})$ admits a real left adjoint. So every object of $\caP^\Sigma(\caC)$ is the geometric realization of free objects. Consequently, it is enough to check that every free object of $\caP^\Sigma(\caC)$ is a small sifted colimit of representables.
	
In turn, in view of the facts listed below, for this it is enough to prove that every object of $ \Fun(\caC^\op, \Spc^{C_2}) $ is a geometric realization of a simplicial object in $ \Fun(\caC^\op, \Spc^{C_2}) $, whose terms are small coproducts of objects of the form $$ C_2 \otimes (\map_\caC(-, Y) \otimes \ast) \hspace{1em}\text{ or}\hspace{1em}  \ast \otimes (\map_\caC(-, Y) \otimes \ast) \simeq \map_\caC(-, Y) \otimes \ast $$ for some $Y \in \caC$.
	
\begin{itemize}
\item The free functor preserves tensors (being a real left adjoint \cref{Venr:left_adj_preserve_tensors}) and sends objects of the form $ \map_\caC(-, Y) \otimes \ast $ for some $Y \in \caC$ to the real presheaf represented by $Y.$

\item Every small coproduct of objects of the form 
$$ C_2 \otimes (\map_\caC(-, Y) \otimes \ast) \hspace{1em}\text{or} \hspace{1em} \map_\caC(-, Y) \otimes \ast $$ for some $Y \in \caC$ is a small filtered colimit of finite coproducts of objects of the same form.

\item The Yoneda embedding $\caC \subset \caP^\Sigma(\caC)$ preserves finite coproducts and tensors with $C_2$. Indeed, for every $F \in \caP^\Sigma(\caC) $ the presheaf $$\map_{\caP^\Sigma(\caC) }(y(-),F) \simeq \map_{ \caP(\caC)}(y(-),F) \simeq F$$ preserves finite products and cotensors with $C_2.$ 
\end{itemize}
	
By \cite[Lemma 5.5.8.13]{lurie.HTT}, for every small $\infty$-category $\caB$, every presheaf on $\caB$ is the geometric realization of some simplicial object in $ \Fun(\caB^\op, \Spc)$ whose terms are small coproducts of representables.
	
Then, every object of $\Fun(\caC^\op, \Spc^{C_2}) \simeq \Fun( (\caC \times \caO_{C_2})^\op, \Spc) $ is the geometric realization of some simplicial object in $\Fun(\caC^\op, \Spc^{C_2})$, whose terms are small coproducts of objects of the form $$ \map_\caC(-, Y) \otimes \map_{\caO_{C_2}}(-, Z) \simeq \map_{\caO_{C_2}}(-, Z) \otimes (\map_\caC(-, Y) \otimes \ast) $$ for some $Y \in \caC$ and $Z \in \caO_{C_2}.$ This concludes the proof, since $\map_{\caO_{C_2}}(-, Z)= C_2$ if $Z$ is the free orbit, and $\map_{\caO_{C_2}}(-, Z)= \ast $ if $Z$ is the point orbit.
\end{proof}

We conclude this subsection by proving that $\D\caC$ is genuine preadditive whenever $\caC$ is so. We begin with the following short lemma.

\begin{lemma}\label{dfghjgl}
The functor giving by cotensoring with $C_2$, that is, $$ (-)^{C_2} \simeq \uHom_{\Spc^{C_2}}( C_2,-)\colon  \Spc^{C_2}  \to  \Spc^{C_2},$$ preserves small sifted colimits.
\end{lemma}

\begin{proof}
For any genuine $C_2$-space $X$ we have the following equivalences $$\uHom_{\Spc^{C_2}}( C_2, X)^u \simeq \map_{\Spc}( C_2, X^u) \simeq X^u \times X^u,$$ 
$$\uHom_{\Spc^{C_2}}( C_2, X)^{C_2} \simeq \map_{\Spc^{C_2}}( C_2, X) \simeq X^u.$$ 
\end{proof}

\begin{proposition}\label{fghfghjj} If  $\caC$ is a genuine  preadditive $\infty$-category, so is $\D\caC$. In particular, $\D\Wald_\infty^\gd$ is genuine preadditive.
\end{proposition}

\begin{proof}
We know that the Yoneda embedding $\caC \subset \caP(\caC)$ preserves small limits and cotensors with all genuine $C_2$-spaces, by \cref{Venr:Yoneda_embedding_preserves_limits_and_cotensors}. Therefore, the embedding $ \caC \subset \D(\caC)$ also preserves them. Especially, the final object of $\caC$ is mapped to the final object of $ \D\caC$, that we denote by $\ast$. Now, for every $F \in  \D\caC$ the genuine $C_2$-space $\map_{  \D\caC}(\ast, F) \simeq \map_{\caP(\caC)}(\ast, F) \simeq F(\ast)$ is contractible, as evaluation at $\ast$ preserves small colimits, taking sifted colimits preserves final objects and $\map_{\caP(\caC)}(\ast, F)$ is contractible for any $F$ representable. This completes the proof that $\D\caC$ has a zero object.

The real $\infty$-category $ \D\caC$ is closed under finite products and cotensors with $C_2$ in $ \caP(\caC) $, as cotensoring with $C_2$ preserves sifted colimits (see \cref{dfghjgl}) and finite products preserve small colimits in each component. In what follows, we will show that tensors with $C_2$ are cotensors and coproducts are products.

Given $F,G, H \in  \D\caC$, we want to see that the canonical maps below are equivalences.
$$ \map_{ \caP(\caC)}( F \times G, H) \to  \map_{ \caP(\caC)}(F, H) \times \map_{ \caP(\caC)}(G, H)$$ 
$$ \map_{ \caP(\caC)}(F^{C_2}, H) \to \map_{ \caP(\caC)}(F, H)^{C_2} $$ 

We can restrict ourselves to the case that $F,G, H \in \caC$, since finite products preserve sifted colimits, and preserve sifted limits in each component. Now, when $F,G, H \in \caC$ these equivalences hold, as $\caC$ is genuine preadditive by hypothesis.
\end{proof}

\subsection{(Real) simplicial homotopies}\label{subsec:real_s_htpies}

In this subsection we introduce a class of morphisms in $\Wald^\gd_\infty$ that yield an equivalence after geometric realization in $\D\Wald^\gd_\infty$. We also separate a subclass of these morphisms, the (relative) real simplicial homotopy equivalences, that enjoys good stability properties.

We begin by recalling some essential facts about simplicial homotopies. Given a simplicial object $X$ in an $\infty$-category $\caD$, we write $X^{*}$ for the composition $(\Delta_{/[1]})^\op \to \Delta^\op \xrightarrow{X} \caD.$ Similarly, we obtain $f^\ast$ for maps $f\colon X\to Y$ in $\caD$.

\begin{definition}\label{def:simplicial_hpry}
Let $\caD$ be an $\infty$-category, $X,Y,$ and $Z$ objects in $\caD^{\Delta^\op}$, and $\alpha,\beta\colon X\to Y$ morphisms over $Z$ as below.

\[
\begin{tikzcd}[column sep=small]
X\ar[rr, "\alpha"]\ar[dr, "\gamma"']        &       &Y\ar[dl, "\delta"]      &           &X\ar[rr, "\beta"]\ar[dr, "\gamma"']      &       &Y\ar[dl, "\delta"]\\
&       Z       &       &       &       &Z    &
\end{tikzcd}
\]
	
A simplicial homotopy $h$ from $\alpha$ to $\beta$ relative to $Z$ is a commutative triangle in the $\infty$-category $\Fun((\Delta_{/[1]})^\op,\caD)$ as below
\[
\begin{tikzcd}[column sep=small]
X^\ast\ar[dr, "\gamma^\ast"']\ar[rr, "h"]      &       &Y^\ast\ar[dl, "\delta^\ast"]\\
        &Z^\ast,     &
\end{tikzcd}
\]
whose images under the functors $$i_0^*,i_1^*\colon \Fun((\Delta_{/[1]})^\op,\caD) \to \Fun(\Delta^\op,\caD),$$ where $i_0,i_1\colon  \Delta \to \Delta_{/[1]}$ are the functors induced by the maps $\{0\} \subset [1], \{1\} \subset [1]$ in $\Delta$, are the upper-left and upper-right triangle respectively.

When $Z$ is a terminal object, we call a simplicial homotopy over $Z$ simply `simplicial homotopy'.
\end{definition}

In the evident way we define simplicial homotopy inverses and simplicial homotopy equivalences.

In the following proposition we show that the geometric realization canonically identifies simplicially homotopic map; this is, simplicial homotopy equivalences yield equivalences on geometric realizations.

\begin{proposition}\label{grergetz}\label{prop:geom_real_s_htpies}
Let $X,Y \colon \Delta^\op \to \caD$ be two simplicial objects in an $\infty$-category $\caD$ and $h\colon  X^{*} \to Y^{*}$ a simplicial homotopy in $\caD$ from $\alpha$ to $\beta$. Then $\alpha, \beta\colon X \to Y$ yield equivalent maps on geometric realizations.
\end{proposition}

\begin{proof}
By Quillen's Theorem A, the functor $p\colon (\Delta_{/[1]})^\op \to \Delta^\op$ is cofinal since for each $[n] \in \Delta^\op$ the pullback
$$ T\coloneqq (\Delta_{/[1]})^\op \times_{\Delta^\op} (\Delta^\op)_{[n] /} $$ is weakly contractible. Indeed, we have that $T^\op \simeq \Delta_{/[1]} \times_{\Delta} \Delta_{/[n]} $ is the category of simplices of $\Delta^1 \times \Delta^n$, which is known to be weakly contractible.

Therefore, for any simplicial object $Z$ in $\caD$ the canonical map $\colim (Z^*) \to \colim (Z)$ induced by $p$ is an equivalence. 
And so also the canonical map $\colim (Z) \to \colim (Z^\ast)$ induced by $i_j$ for $j=0,1$ is an equivalence, as $\colim (Z) \to \colim (Z^\ast) \to \colim (Z) $ is the identity. We then have a commutative square 
\begin{equation*}
\begin{tikzcd}
\colim(X) \ar[d, "\colim (\alpha)"'] \ar[r, "\simeq"]   & \colim(X^\ast) \ar[d, "\colim (h)"] \\
\colim (Y) \ar[r, "\simeq"] & \colim (Y^\ast),
\end{tikzcd}
\end{equation*}
and similar for $\beta$.

\end{proof}

As we said before, our aim in this subsection is to define a real version of simplicial homotopy and link these notions.

We begin by presenting such definition first, and introducing its ingredients immediately afterwards. We hope that, albeit it will lack formal meaning at the moment, this will give the reader useful intuition when defining its components.

\begin{notation}
For any real simplicial object $X$ in a real $\infty$-category $\caC$
let $X^{**}$ the composition $\underline{\Delta}_{/\widetilde{[1] \times[1]}}^\op \to \underline{\Delta}^\op \xrightarrow{X} \caC.$ 
\end{notation}

\begin{definition}\label{def:real_s_htpy}
Let $\caD$ be a real $\infty$-category, $X,Y,$ and $Z$ real simplicial object in $\caD$, and $\alpha,\beta\colon X\to Y$ morphisms over $Z$ as below.

\[
\begin{tikzcd}[column sep=small]
X\ar[rr, "\alpha"]\ar[dr, "\gamma"']        &       &Y\ar[dl, "\delta"]      &           &X\ar[rr, "\beta"]\ar[dr, "\gamma"']      &       &Y\ar[dl, "\delta"]\\
&       Z       &       &       &       &Z    &
\end{tikzcd}
\]
	
A real simplicial homotopy $h$ from $\alpha$ to $\beta$ relative to $Z$ is a commutative triangle in the real $\infty$-category $\Fun_{\Spc^{C_2}}((\underline{\Delta}_{/\widetilde{[1] \times[1]}})^\op,\caD)$ as below
\[
\begin{tikzcd}[column sep=small]
X^{\ast\ast}\ar[dr, "\gamma^{\ast\ast}"']\ar[rr, "h"]      &       &Y^{\ast\ast}\ar[dl, "\delta^{\ast\ast}"]\\
        &Z^{\ast\ast},     &
\end{tikzcd}
\]
whose images under the real functors $i_0^{\ast},i_1^{\ast}\colon \Fun_{\Spc^{C_2}}(\underline{\Delta}_{/\widetilde{[1] \times[1]}}^\op,\caD) \to \Fun_{\Spc^{C_2}}(\underline{\Delta}^\op,\caD)$, where $i_0,i_1\colon \underline{\Delta} \to \underline{\Delta}_{/\widetilde{[1] \times[1]}}$ are the real functors induced by the maps $\{0\} \subset [1], \{1\} \subset [1]$ in $\Delta$, are the upper-left and upper-right triangle respectively.

When $Z$ is a terminal object, we call a real simplicial homotopy over $Z$ simply `real simplicial homotopy'.
\end{definition}

Let us now present the components of \cref{def:real_s_htpy}.

Given a real $\infty$-category $\caC$ and an object $\X \in \caC$ by \cref{Venr:u_of_slice category}
we have a real $\infty$-category $\caC_{/X}$
with the evident underlying $\infty$-category and $(\caC_{/\X})^u \simeq  \caC^u_{/\X}.$ So for any $\infty$-category $\X$ with duality we get a real structure on the $\infty$-category $(\Cat_\infty^{hC_2})_{/X}$. We can then consider the pullback of real $\infty$-categories 
$$ \underline{\Delta}_{/ X }\coloneqq\underline{\Delta} \times_{\Cat_\infty^{hC_2}} (\Cat_\infty^{hC_2})_{/ X }.$$  

For $X = [0]$ we have $\underline{\Delta}_{/ X } \simeq \underline{\Delta}.$ Note that $(\underline{\Delta}_{/ X })^u$ is the essential image of the forgetful functor $$ \Delta^{hC_2} \times_{\Cat_\infty^{hC_2}} (\Cat_\infty^{hC_2})_{/ X } \to \Delta_{/X}\coloneqq \Delta \times_{\Cat_\infty} \Cat_{\infty / X}.$$

In the following we consider $\underline{\Delta}_{/\widetilde{[1] \times[1]}}$.
The canonical maps $\{0\} \subset [1], \{1\} \subset [1]$ in $\Delta$ yield two maps
$$ [0] \simeq \widetilde{[0] \times [0]} \to  \widetilde{[1] \times [1]}$$
in $\Cat_\infty^{hC_2}$ that give rise to real functors
$i_0, i_1\colon \underline{\Delta} \to \underline{\Delta}_{/\widetilde{[1] \times[1]}}$ over $ \underline{\Delta}.$

Finally, given a real simplicial object $X$ in a real $\infty$-category $\caC$ we write $X^{**}$ for the composition $$\underline{\Delta}_{/\widetilde{[1] \times[1]}}^\op \to \underline{\Delta}^\op \xrightarrow{X} \caC.$$

We will now give descriptions of the underlying categories of $\underline{\Delta}_{/\widetilde{[1] \times[1]}}$, for which we first introduce the following notation.

\begin{notation}
Let $\Delta'_{/[1]} \subset \Delta_{/[1]}$ be the wide subcategory whose morphisms are the functors $\gamma\colon [m] \to [n]$ over $[1]$ such that $\gamma$ does not only commute with the structure morphisms $[m] \to [1], [n]\to [1]$ but also with the opposites of the structure morphisms $$[m] \simeq [m]^\op \to [1]^\op \simeq [1], [n] \simeq [n]^\op \to [1]^\op \simeq [1].$$
\end{notation}

\begin{lemma}
\begin{enumerate}
\item There is a canonical equivalence $$(\underline{\Delta}_{/\widetilde{[1] \times[1]}})^{hC_2} \simeq \Delta^{hC_2} \times_{\Delta} \Delta_{/[1]}, $$
\item There is a canonical equivalence $$ (\underline{\Delta}_{/ \widetilde{[1] \times[1]}})^u \simeq \Delta'_{/[1]}.$$
\end{enumerate}
\end{lemma}

\begin{proof}
The equivalence in (1) is obtained from the chain of canonical equivalences
$$(\underline{\Delta}_{/\widetilde{[1] \times[1]}})^{hC_2} \simeq \Delta^{hC_2} \times_{\Cat_\infty^{hC_2}} (\Cat_\infty^{hC_2})_{/ \widetilde{[1] \times[1]} } \simeq \Delta^{hC_2} \times_{\Delta} \Delta_{/[1]}.$$

For (2) recall that $(\Cat_\infty^{hC_2})^u$ is the essential image of the forgetful functor $\Cat_\infty^{hC_2} \to \Cat_\infty$. This implies that $ (\underline{\Delta}_{/ \widetilde{[1] \times[1]}})^u $ is the essential image $\Xi$ of the forgetful functor $$ \Delta^{hC_2} \times_{\Cat_\infty^{hC_2}} (\Cat_\infty^{hC_2})_{/ \widetilde{[1] \times[1]}} \to  \Delta \times_{\Cat_\infty} \Cat_{\infty / [1] \times[1]}.$$
The functor $\Delta \times_{\Cat_\infty} \Cat_{\infty / [1] \times[1]} \to \Delta_{/[1]} $ induced by projection to the first factor
$[1]\times [1] \to [1]$ induces an equivalence
$\Xi \simeq \Delta'_{/[1]}.$	
\end{proof}

Let us now try to gain some familiarity with the notion of real simplicial homotopy by showing some basic properties of it.

\begin{remark}
Given $X,Y,Z \colon\underline{\Delta}^\op \to \caC$ real simplicial objects in a real $\infty$-category $\caC$ and $$h\colon X^{**} \to Y^{**}, \ h'\colon Y^{**} \to Z^{**}$$ real simplicial homotopies from $\alpha $ to $\beta$ respectively from $\alpha' $ to $\beta',$ the composite $$h' \circ h\colon  X^{**} \to Y^{**} \to Z^{**}$$ is a real simplicial homotopy in $\caC$ from $\alpha' \circ \alpha $ to $\beta' \circ \beta$.
\end{remark}

\begin{remark}
A real functor $F\colon  \caC \to \caD$ sends a real simplicial homotopy in $\caC$ (relative to some object $Z \in \caC$) from $\alpha$ to $\beta$ to a real simplicial homotopy in $\caD$ from $F(\alpha)$ to $F(\beta)$ (relative to $F(Z)$).
\end{remark}

\begin{remark}
If $\caC$ admits the neccessary pullbacks, the pullback of a real simplicial homotopy $h\colon  X^{**} \to Y^{**}$ in $\caC$ relative to $Z$ from $\alpha$ to $\beta$ along a real simplicial map $Z' \to Z$ in $\caC$ is a real simplicial homotopy $$h'\colon  (Z' \times_Z X)^{**} \to (Z' \times_Z Y)^{**}$$ in $\caC$ relative to $Z'$ from $Z' \times_Z \alpha$ to $Z' \times_Z \beta$.
Especially the pullback of a real simplicial homotopy equivalence in $\caC$ relative to $Z$ along a real simplicial map $Z' \to Z$ in $\caC$ is a real simplicial homotopy equivalence relative to $Z'$.
\end{remark}

\begin{lemma}
Let $\caC$ be a real $\infty$-category. A real simplicial homotopy $h\colon  X^{**} \to Y^{**}$ in $\caC$ from $\alpha$ to $\beta$ gives rise to a simplicial homotopy in $\caC$ from $e(\alpha)\colon e(X) \to e(Y)$ to $e(\beta)\colon e(X) \to e(Y)$; where $e$ denotes the edgewise subdivision. 
\end{lemma}

\begin{proof}
Let $I\colon \Delta_{/[1]} \to \underline{\Delta}_{/\widetilde{[1] \times[1]} } \simeq \Delta^{hC_2} \times_{\Delta} \Delta_{/[1]}$ be the functor that sends a map $\varphi\colon [m] \to [1]$ to the map $I(\varphi)\colon [m]^\op * [m] \to [1]$, which sends everything in $[m]^\op$ to $0$ and is $\varphi$ on $[m]$. It follows that $h \circ I^\op$ is a simplicial homotopy from $e(\alpha)$ to $h \circ I^\op \circ i_1^\op$ and $h \circ (-)^\op \circ I^\op \circ (-)^\op $ is a simplicial homotopy from $h \circ I^\op \circ i_1^\op$ to $e(\beta)$.
\end{proof}

\begin{corollary}
Real simplicially homotopic maps from $X$ to $Y$ become equivalent after realization. 
\end{corollary}

\begin{proof}
It follows from \cref{grergetz}, as the realization of a real simplicial object is the realization of its edgewise subdivision (see \cref{prop:geom_realization_with_edgewise}).
\end{proof}

There is another important relationship between simplicial homotopies and real simplicial homotopies, that is the content of \cref{cresimpho}. Before diving in its proof, we follow a digression that will clarify some steps of it.

We know that the real $\infty$-category of functors between two real $\infty$-categories, denoted by $\Fun_{\Spc^{C_2}}(\caC, \caD)$, admits cotensors with $C_2$ whenever $\caD$ does. As we know, these cotensors can be computed objectwise. The cotensor of a real functor $F\colon\caC\to\caD$ with $C_2$ also admits a different description; in the following proposition we present it.

Consider the forgetful-cofree adjunction $\Spc^{C_2} \rightleftarrows\Spc$, whose right adjoint sends a space X to the space $\widetilde{X \times X}$ with a $C_2$-action given by switching the factors. Passing to enriched $\infty$-categories this adjunction gives rise to an adjunction as below
\begin{equation}\label{adj_11}
(-)^u\colon \Cat^{\Spc^{C_2}}_\infty \rightleftarrows \Cat_\infty\colon \widetilde{(-)}.
\end{equation}

\begin{proposition}\label{cotensor_as_adjoint}
Let $\caD$ be a real $\infty$-category that admits cotensors with $C_2$. Then the unit component $\lambda_\caD\colon\caD \to \widetilde{ \caD^u}$ of the adjunction (\ref{adj_11}) admits a real right adjoint $\phi$ such that the composition $\caD \to \widetilde{\caD^u} \xrightarrow{\phi} \caD$ of left and right adjoint is cotensoring with $C_2$.
\end{proposition}

\begin{proof}
We start by recalling that there is a canonical equivalence
$ \widetilde{\caD^u}^{hC_2} \simeq \caD^u$ such that $\lambda_\caD^{hC_2}\colon \caD^{hC_2} \to \widetilde{\caD^u}^{hC_2} \simeq \caD^u$ identifies with the canonical functor forgetting the fixed point data on mapping spaces. 

We need to check that, for every $Y\in \caD$, the real presheaf $$\map_{\widetilde{\caD^u}}(-,\lambda_\caD(Y)) \circ \lambda^\op_\caD\colon \caD^\op \to \Spc^{C_2}$$ is represented by $Y^{C_2}$. This follows from the existence of a chain of canonical equivalences of genuine $C_2$-spaces as below for $\X \in \caD$
$$\map_\caD(X,Y^{C_2}) \simeq \map_\caD(X,Y)^{C_2} \simeq$$$$ \widetilde{\map_\caD(\lambda_\caD(X),\lambda_\caD(Y))^u \times \map_\caD(\lambda_\caD(X),\lambda_\caD(Y))^u} \simeq \map_{\widetilde{\caD^u}}(\lambda_\caD(X),\lambda_\caD(Y)).$$
\end{proof}

\begin{corollary}
Let $\caD$ be a real $\infty$-category that admits cotensors with $C_2$. The
forgetful functor $\caD^{C_2} \to \caD^u$ admits a right adjoint $\phi$ such that the composition $\caD^{C_2} \to \caD^u \xrightarrow{\phi} \caD^{C_2}$ is cotensoring with $C_2$.
\end{corollary}

\begin{proof}

The unit component $\lambda_\caD\colon\caD \to \widetilde{ \caD^u}$ of the adjunction (\ref{adj_11}) induces the forgetful functor $\caD^{C_2} \to \widetilde{\caD^u}^{C_2} \simeq \caD^u.$

\end{proof}

Throughout this paper, for an object $X$ in a real $\infty$-category $\caC$, we have been using the notation $\widetilde{X\times X}$ to denote the cotensor of $X$ with $C_2$. This, together with \cref{cor:tilde_for_functors_is_cotensor_with_C2}, justifies the following notation.

\begin{notation}\label{Not_tilde}
Let $\caC, \caD$ be real $\infty$-categories such that $\caD$ has cotensors with $C_2$, and let $F\colon \caC^u \to \caD^u$ a functor between them.
We write $\widetilde{F\times F} $ for the composition $$\caC \xrightarrow{\lambda} \widetilde{\caC^u} \xrightarrow{\widetilde{F}} \widetilde{\caD^u} \xrightarrow{\phi} \caD.$$	

Similarly, for a real natural transformation $\alpha\colon F \to G$ of real functors $\caC \to \caD$ we write $\widetilde{\alpha \times \alpha}\colon \widetilde{F \times F}
\to\widetilde{G \times G}$ for $\phi \circ \alpha \circ \lambda\colon  \phi \circ F \circ \lambda \to \phi \circ G \circ \lambda.$
\end{notation}

\begin{corollary}\label{cor:tilde_for_functors_is_cotensor_with_C2} Let $\caC, \caD$ be real $\infty$-categories such that $\caD$ has cotensors with $C_2$. The forgetful functor $$ \Fun_{\Spc^{C_2}}(\caC,\caD)^{C_2} \to \Fun(\caC^u,\caD^u)$$
admits a right adjoint that sends $F$ to $\widetilde{F \times F}$.
The resulting adjunction extends the adjunction
$$ \Fun^{\Spc_{C_2}}(\caC,\caD)^{C_2} \rightleftarrows \Fun^{\Spc_{C_2}}(\caC,\caD)^u$$
of the corollary above.
\end{corollary}

\begin{proof}
The real adjunction $\caD \rightleftarrows \widetilde{\caD^u}$ of \cref{cotensor_as_adjoint} gives rise to a real adjunction 
$$\Fun_{\Spc^{C_2}}(\caC,\caD) \rightleftarrows \Fun_{\Spc^{C_2}}(\caC, \widetilde{\caD^u}) \simeq \widetilde{\Fun(\caC^u,\caD^u)},$$ where the left adjoint factors as
$$\Fun_{\Spc^{C_2}}(\caC,\caD) \to \widetilde{\Fun_{\Spc^{C_2}}(\caC,\caD)^u} \hookrightarrow \widetilde{\Fun(\caC^u,\caD^u)}.$$
It can be easily checked that the right adjoint sends $F$ to $\widetilde{F \times F}$, and thus takes objectwise the cotensor with $C_2$.

We conclude then that the induced real adjunction 
$$\Fun_{\Spc^{C_2}}(\caC,\caD)^{C_2} \rightleftarrows  \widetilde{\Fun(\caC^u,\caD^u)}^{C_2} \simeq \Fun(\caC^u,\caD^u)$$ extends the adjunction of the corollary above applied to the real
$\infty$-category $\Fun_{\Spc^{C_2}}(\caC,\caD).$
\end{proof}

We now prove the result we wanted about real simplicial homotopies.

\begin{lemma}\label{cresimpho}
Let $\caD$ be a real $\infty$-category. A simplicial homotopy $h\colon X^{*} \to Y^{*}$ in the $\infty$-category $\caD^u$ (relative to some object $Z \in \caD$) from $\alpha$ to $\beta$ naturally induces a real simplicial homotopy $$\widehat{h}\colon \widetilde{X \times X}^{\ast \ast} \to \widetilde{Y \times Y}^{\ast \ast}$$ in $\caD$ from $\widetilde{\alpha \times \alpha}\colon \widetilde{X \times X} \to \widetilde{Y \times Y}$ to $\widetilde{\beta \times \beta}\colon \widetilde{X \times X} \to \widetilde{Y \times Y}$ (relative to $ \widetilde{Z \times Z}).$
\end{lemma}

\begin{proof}
Consider the forgetful-cofree adjunction $\Spc^{C_2} \rightleftarrows\Spc$, whose right adjoint sends a space $X$ to $\widetilde{X \times X}$. As we saw before, passing to enriched $\infty$-categories this adjunction gives rise to an adjunction $(-)^u\colon \Cat^{\Spc^{C_2}}_\infty \rightleftarrows \Cat_\infty\colon \widetilde{(-)}.$
Next we use the real functor $$\rho\colon (\underline{\Delta}_{/\widetilde{[1] \times [1]}})^\op \longrightarrow \widetilde{ (\Delta_{/[1]})^\op},$$ transpose to the functor 
$$\rho'\colon ((\underline{\Delta}_{/\widetilde{[1] \times [1]}})^\op)^u \subset (\Delta_{/[1] \times [1]})^\op \to (\Delta_{/[1]})^\op$$ induced by the projection $[1] \times [1] \to [1]$ to the first factor. 
The simplicial homotopy $h\colon X^{*} \to Y^{*}$ in $\caD^u$
is given by a functor
$h\colon (\Delta_{/[1]})^\op \to \caD^u$
as in \cref{def:simplicial_hpry}.
We take the composition
$$h'\colon ((\underline{\Delta}_{/\widetilde{[1] \times [1]}})^\op)^u \xrightarrow{\rho'^\op} (\Delta_{/[1]})^\op \xrightarrow{h} \caD^u$$ and have $h' \circ (i_0^\op)^u \simeq \alpha $ and $h' \circ (i_1^\op)^u \simeq \beta$.
Then, calling $\Omega\coloneqq \underline{\Delta}_{/\widetilde{[1] \times [1]}}$, the composition of real functors $\widehat{h}$ given by  
$$\Omega \to \widetilde{\Omega^u} \xrightarrow{\widetilde{h'}} \widehat{\caD^u} \to \caD$$
is the desired real simplicial homotopy 
from $\widetilde{\alpha \times \alpha}\colon \widetilde{X \times X} \to \widetilde{Y \times Y}$ to $\widetilde{\beta \times \beta}\colon \widetilde{X \times X} \to \widetilde{Y \times Y}$ (relative to $ \widetilde{Z \times Z}).$
\end{proof}

\subsection*{Construction of a special real simplicial homotopy}
We construct now a simplicial homotopy that will prove useful later in the proof of the additivity theorem. First, we fix some notation that will also accompany us later on.

\begin{notation}\label{not:maps_jk}
The canonical maps $[m] \to [m]*[n]$ and $[n] \to [m]*[n]$
define, for $k=1,2$, natural transformations $j^k\colon (-)\ast(-)\to \pr_k $ of functors $ \Delta^\op \times \Delta^\op \to \Delta^\op $. For simplicity, we will denote the domain of this map by $T\coloneqq (-)\ast(-)$.
\end{notation}

\begin{construction}\label{fgjhbgch}
Let $n \geq 0$ and $\nu\colon T_{-,n} \times [1] \to T_{-,n}$ the map in $\Fun((\Delta_{/[1]})^\op, \Cat^\op) $ induced by $\{1\}\subset [1].$ In what follows we define $\kappa\colon (T_{-,n})^\ast\to (T_{-,n}\times[1])^\ast$ in a way that it fits in a commutative triangle in $\Fun((\Delta_{/[1]})^\op, \Cat^\op)$ as below.
\begin{equation*}
\begin{tikzcd}[column sep=small]
(T_{-,n})^* \ar[rd, "(j^2_{-,n})^*"'] \ar[rr, "\kappa"]  &   & (T_{-,n} \times [1])^*\ar[ld, "(j^2_{-,n} \circ \nu)^*"] \\   
\       &{[}n{]}           &
\end{tikzcd}
\end{equation*}	

We set the image of $\varphi\colon [m] \to [1]$ under $\kappa$ to be such that it fits in a commutative triangle as below	
\begin{equation*}
\begin{tikzcd}[column sep=small]
{[}m{]}*{[}n{]} \ar[rd, "j^2_{m,n}"'] \ar[rr, "\kappa_{\varphi}"]   &   & {[}m{]}*{[}n{]} \times [1] \ar[ld, "j^2_{m,n} \circ \nu_m"] \\   
\       & {[}n{]}       &
\end{tikzcd}
\end{equation*}	
and verifies the following conditions. The composition below is the identity
$$[m]*[n] \xrightarrow{\kappa_\varphi} [m]*[n] \times [1] \to [m]*[n] \times \{0\}$$ and the composition 
$$[m]*[n] \xrightarrow{\kappa_\varphi} [m]*[n] \times [1] \to [m]*[n] \times \{1\}$$ sends $i$ to 
$$\begin{cases}
i \ \text{if} \ \ i \in [m] \ \text{and} \ \varphi(i)=0, \\
0 \in [n] \ \text{if} \ i \in [m] \ \text{and} \ \varphi(i)=1, \\
i \ \text{if} \ \ i \in [n].
\end{cases}$$

The map $\kappa\colon (T_{-,n})^{*} \to (T_{-,n} \times [1])^{*}$ is a simplicial homotopy in $\Cat^\op$ relative to $[n]$ from the map  $ T_{-,n} \to  T_{-,n} \times [1]$ induced by $[1] \to [0]$ to 
$$\zeta\coloneqq i_1^*(\kappa)\colon T_{-,n} \to  T_{-,n} \times [1].$$
\end{construction}

\begin{notation}\label{not:p_and_q}
The map $[m]*[n] \to [n]$ that is the identity on $[n]$ and sends all of $[m]$ to $0\in [n]$ is natural in $[m] \in \Delta$ (but not in $[n] \in \Delta$!) and so defines a natural transformation $q^n\colon[n] \to T_{-,n} $ of functors $ \Delta^\op \to \Delta^\op.$ Similarly we define another natural transformation $p^n\colon [n]\to T_{n,-}$ of functors $\Delta^\op\to\Delta^\op$, by $p^n\coloneqq (-)^\op \circ q^n \circ (-)^\op$ of functors $\Delta^\op\to\Delta^\op$.
\end{notation}

\begin{remark}\label{reeeem}
The compositions below are the identities on $[n]$, 
$$[n]  \xrightarrow{p^n}T_{n,-} \xrightarrow{j^1_{n,-}} [n], \text{ and } [n] \xrightarrow{q^n} T_{-,n} \xrightarrow{j^2_{-,n} } [n]$$ 
Moreover, the map $[n] \xrightarrow{q^n} T_{-,n} \xrightarrow{j^1_{-,n}  } \id $ factors as $ [n] \xrightarrow{0} [0] \to \id. $
\end{remark}

\begin{remark}\label{rem_fact}
The map $ T_{-,n} \xrightarrow{\zeta}  T_{-,n} \times [1] \to T_{-,n} \times \{1\} $ factors as 
$T_{-,n} \xrightarrow{j^2_{-,n}} [n] \xrightarrow{q^n}  T_{-,n}. $

\end{remark}

\begin{lemma}\label{lem_simhoinv}
With the constructions above, we have the following simplicial homotopies relative to $[n]$.
\begin{enumerate}
 \item The map $j^1_{n,-}\colon T_{n,-} \to [n]$ is a simplicial homotopy inverse of $p^n\colon [n] \to T_{n,-}$.
 \item The map $j^2_{-,n}\colon T_{-,n} \to [n]$ is a simplicial homotopy inverse of $q^n\colon [n] \to T_{-,n}$.
\end{enumerate}
\end{lemma}

\begin{proof}
We will prove the second case, the first one is dual.

In \cref{fgjhbgch} we got a simplicial homotopy $(T_{-,n})^{*} \xrightarrow{\kappa} (T_{-,n} \times [1])^{*}$ in $\Cat^\op$ relative to $[n]$ from the map  $ T_{-,n} \to  T_{-,n} \times [1]$ induced by $[1] \to [0]$ to a map
$$\zeta\coloneqq i_1^*(\kappa)\colon T_{-,n} \to  T_{-,n} \times [1].$$
The induced map $ (T_{-,n})^{*} \xrightarrow{\kappa} (T_{-,n} \times [1])^{*}\to (T_{-,n} \times \{1\})^* $
is a simplicial homotopy in $\Cat^\op$ relative to $[n]$ from 
the identity to the map $ T_{-,n} \xrightarrow{\zeta}  T_{-,n} \times [1] \to T_{-,n} \times \{1\} $, which factors as 
$T_{-,n} \xrightarrow{j^2_{-,n}} [n] \xrightarrow{q^n}  T_{-,n}. $  The other composition follows from \cref{reeeem}.
\end{proof}

\subsection{Some lemmas for the real additivity theorem}\label{subsec:lemmas_add}

In order to make the proof of the main theorem of this section (\cref{thm:add}) easier to read, we have gathered in this subsection several technical lemmas that we will later use. The reader may decide to skip ahead to \cref{subsec:add_thm} and return to this subsection as needed.

Although several things are done on the way, this section is devoted to the construction of the map
$$\omega_1\colon \lambda(\caC)\to S(\caC)_3 $$ in $\Wald^\gd_\infty$ that needed in the proof of \cref{thm:add}. The skeleton of such definition is found in \cref{const_omega1}.

We begin by introducing the domain of such a map. 

\begin{notation}\label{not:lambda(C)}
We define $\lambda(\caC)$ as the pullback in $\Wald^\gd_\infty$ of the square
\[
\begin{tikzcd}
\lambda(\caC) \ar[r, "\beta"]\ar[d, "\omega_0"']    &\widetilde{\caC^{[1]} \times \caC^{[1]}} \times \caC \ar[d, "\widetilde{\ev_0 \times \ev_0}\times \caC"] \\
S(\caC)_3 \ar[r, "\gamma_3"']   &\widetilde{\caC \times \caC} \times \caC.
\end{tikzcd}
\]
So an object of $\lambda(\caC)$ is given by cofibrations $X \to Y \to Z$ in $\caC$ and maps $X \to A, \ B \to Z/Y $ in $\caC.$ Similarly, by taking $S(\caC)'_{3}$ given by the pullback

\[
\begin{tikzcd}
S(\caC)_3' \ar[r, "\varrho"]\ar[d]      &\ast\times\caC\ar[d, "0\times \caC"]  \\
S(\caC)_3 \ar[r, "\gamma_3"']    & \widetilde{\caC \times \caC} \times \caC,
\end{tikzcd}
\]
we define $\lambda'(\caC)$ for any $\caC \in \Wald_\infty$ as the pullback below. 

\[
\begin{tikzcd}
\lambda'(\caC) \ar[r, "\beta"]\ar[d, "\omega_0"']    &\widetilde{\caC^{[1]} \times \caC^{[1]}} \times \caC \ar[d, "\widetilde{\ev_0 \times \ev_0}\times \caC"] \\
S(\caC)'_3 \ar[r, "\gamma_3"']   &\widetilde{\caC \times \caC} \times \caC.
\end{tikzcd}
\]
\end{notation}

In what follows we present some results on cocartesian fibrations that will facilitate the definition of $\omega_1$.

\begin{lemma}\label{fghjkklkl}

Let $\caB, \caD$ be $\infty$-categories such that $\caB$ has an initial object, and let us consider the functor $\gamma\colon \Fun(\caB, \caD) \to \caD$ given by evaluating at the initial object of $\caB$. Let also $H\colon \caB \to \caD$ be a functor and $\theta\colon H(\emptyset) \to X$ a morphism in $\caD$ such that for every morphism $\beta\colon \emptyset \to Z $ in $\caB$ the pushout of the morphisms $\theta, H(\beta)$ exists. Then the following holds:

\begin{enumerate}
    \item There is a $\gamma$-cocartesian lift of $\theta$ starting at $H$.
    \item A lift $\varphi\colon H \to G$ of $\theta$ is $\gamma$-cocartesian if and only if for every morphism $ \emptyset \to Z $ in $\caB$ the induced square
    \[
    \begin{tikzcd}
    H(\emptyset)\ar[r]\ar[d]        &G(\emptyset)\ar[d] \\	
    H(Z)\ar[r]                     &G(Z)
    \end{tikzcd}
    \]
is a pushout square. 
\end{enumerate}
\end{lemma}

\begin{proof}
We start by showing item (2). Let $\delta\colon \caD \to \Fun(\caB, \caD)$ be the diagonal embedding, so we have that $\delta$ is left adjoint to $\gamma$. Observe that a lift $\varphi\colon H \to G$ of $\theta$ satisfies the pushout square condition of the statement if and only if the square
\begin{equation}\label{sqay}
\begin{tikzcd}
\delta(H(\emptyset)) \ar[r, "\delta(\theta)"] \ar[d] & \delta(G(\emptyset)) \ar[d] \\
H \ar[r] & G
\end{tikzcd}
\end{equation} induces a pushout square in $\caD$ at any $Z \in \caB$ .

By assumption for every morphism $\beta\colon \emptyset \to Z $ in $\caB$
the pushout of the morphisms $\theta, H(\beta)$ exists, therefore the pushout of
the morphisms $\delta(\theta)\colon \delta(H(\emptyset)) \to \delta(G(\emptyset))$,
$\delta(H(\emptyset)) \to H$ also exists and is preserved by evaluation at any object of $\caB$. Hence $\varphi$ satisfies the pushout square condition in the statement if and only if square (\ref{sqay}) is a pushout square. This is equivalent to
the condition that for any functor $F\colon \caB \to \caD$ the induced commutative square
\[
\begin{tikzcd}
\map_{\Fun(\caB, \caD)}(G, F) \ar[r] \ar[d] & \map_\caD(G(\emptyset), F(\emptyset)) \simeq \map_{\Fun(\caB, \caD)}(\delta(G(\emptyset)), F) \ar[d] \\
\map_{\Fun(\caB, \caD)}(H, F)  \ar[r] & \map_\caD(H(\emptyset), F(\emptyset)) \simeq \map_{\Fun(\caB, \caD)}(\delta(H(\emptyset)), F)
\end{tikzcd}
\]
is a pullback square, which precisely says that $\varphi\colon  H \to G$ is $\gamma$-cocartesian.

To prove item (1), we observe that the canonical map $H \to \delta(X) \coprod_{\delta(H(\emptyset))} H$ lying over $\theta$ satisfies the pushout square condition of item (2).

\end{proof}

\begin{notation}\label{not:alpha_beta}Let $\caD$ be an exact $\infty$-category and $n \geq 1$. We denote by $$\alpha,\beta\colon S(\caD)_n \to  S(\caD)_1 \simeq \caD$$ the functors induced by the maps $[1] \simeq \{0\to 1\} \subset [n], [1] \simeq \{n-1\to n\} \subset [n]$ respectively.
\end{notation}

In the next lemma, we study the $\alpha,\beta$-cocartesian maps.

\begin{lemma}\label{fjbnjjkk}\label{fjkjjkk}
Let $\caD$ be an exact $\infty$-category and $n \geq 1$. 

\begin{enumerate}
\item\label{alpha_cocart} The functor $\alpha\colon S^e(\caD)_n \to \caD$ is a cocartesian fibration whose cocartesian morphisms are precisely the maps $A \to B$ in $S(\caD)_n$ such that for every $1 \le \ell \le n$ the commutative square
\[
\begin{tikzcd}
A_{0,1} \ar[r] \ar[d] & B_{0,1}\ar[d] \\	
A_{0,\ell} \ar[r] & B_{0,\ell}
\end{tikzcd}
\]
is a pushout square.

\item\label{beta_cocart} The functor $\beta\colon S^e(\caD)_n\to \caD$ is a cartesian fibration whose cartesian morphisms are precisely the maps $A \to B$ in $S(\caD)_n$ such that for every $1 \le \ell \le n$ the commutative square
\[
\begin{tikzcd}
A_{\ell, n} \ar[r] \ar[d] & B_{\ell, n} \ar[d] \\
A_{n-1, n} \ar[r] & B_{n-1, n}
\end{tikzcd}
\]
is a pullback square.
\end{enumerate}

\end{lemma}

\begin{proof}
\cref{alpha_cocart} follows from  \cref{fghjkklkl} and \cref{rema}, and \cref{beta_cocart} follows from \cref{alpha_cocart} via the canonical equivalence $S(\caD^\op)^\op_n \simeq S(\caD)_n$.
\end{proof}

\begin{remark}\label{remref}
Let $\caC$ be an exact $\infty$-category with genuine duality and $n \geq 1.$ The duality on $S(\caC)_n$ 
sends $\alpha$-cocartesian morphisms to $\beta$-cartesian morphisms and vice versa.
\end{remark}

\begin{notation}
For any exact $\infty$-category $\caC$ let $$Q(\caC)\subset S(\caC)_3^{\widetilde{[1] \times [1]}}
$$ be the full subcategory spanned by the commutative squares 
\begin{equation}\label{sqar}
\begin{tikzcd}
A \ar[r] \ar[d] & B \ar[d]  \ar[d]\\
\caC \ar[r] & D
\end{tikzcd}
\end{equation} in $S(\caC)_3$
such that the horizontal maps are $\beta$-cartesian and the vertical maps are $\alpha$-cocartesian, where $\alpha$ and $\beta$ are as in \cref{not:alpha_beta} for $n=3$.
\end{notation}

\begin{remark}
By \cref{remref}, for any exact $\infty$-category $\caC$ with genuine duality the duality on $S(\caC)_3^{\widetilde{[1] \times [1]}}$
restricts to $Q(\caC)$ and makes $Q(\caC)$ into an $\infty$-category with genuine duality.
\end{remark}

\begin{notation}\label{Not:rho}
Let us consider the map $$\rho\colon S(\caC)_3^{\widetilde{[1] \times [1]}} \to \lambda(\caC)= \widetilde{\caC^{[1]} \times \caC^{[1]}} \times_{\widetilde{\caC \times \caC}} S(\caC)_3$$ in $\Wald^\gd_\infty$ that sends a square (\ref{sqar}) to $(B, B_{0,1} \to D_{0,1}, A_{2,3} \to B_{2,3}).$

Precisely, the projection of $\rho$ to $S(\caC)_3$ is evaluation at the non-degenerated object $0,1 \in \widetilde{[1] \times [1]}$, and the projection of $\rho$ to the first factor $\caC^{[1]} $ is the composition $ S(\caC)_3^{{[1] \times [1]}} \to \caC^{{[1] \times [1]}} \to \caC^{[1]}$
of evaluation at $0,1 \in \Ar[3]$ followed by restriction along the functor
$[1] \to [1] \times [1]$ taking $0,1 \to (1,1).$

We consider $$\rho^\prime\colon Q(\caC)\to \lambda(\caC)$$
to be the restriction of $\rho$ to $Q(\caC)$.
\end{notation}

The proof of the following result will occupy the remainder of this subsection. 

\begin{proposition}\label{gr3ehtrer}
The map $\rho^\prime \colon Q(\caC) \to \lambda(\caC)$ 
in $\Wald^\gd_\infty$
is an equivalence.
\end{proposition}

The previous result allows us to finally construct the map $\omega_1$ as follows.

\begin{construction}\label{const_omega1}
We define the desired map $\omega_1\colon \lambda(\caC)\to S(\caC)_3$ as the composite
$$\omega_1\colon \lambda(\caC) \xrightarrow{(\rho^\prime)^{-1}} Q(\caC)\subset S(\caC)_3^{\widetilde{[1] \times [1]}} \xrightarrow{r} S(\caC)_3,$$
where $r$ is evaluation at the non-degenerated object 
$(1,0) \in \widetilde{[1] \times [1]}.$
\end{construction}

\begin{remark}\label{rmk:comp_omega1_sigma}
We observe that  $\omega_1 \circ \sigma = \id.$, fact that will be used later. Indeedn, by construction of $\rho'$ an object of $Q(\caC)$ is a constant square (\ref{sqar}) in $S(\caC)_3$ if and only if its image $(B, B_{0,1} \to D_{0,1}, A_{2,3} \to B_{2,3})$ in $\Lambda(\caC)$ under $\tau'$
has the property that the maps $B_{0,1} \to D_{0,1}, A_{2,3} \to B_{2,3}$
are equivalences.

Let $\delta_\caC \colon \caC \to \caC^{[1]}$ be the diagonal and $\sigma\colon S(\caC)_3 \to \lambda(\caC)$
the pullback along $S(\caC)_3 \xrightarrow{\gamma_3} \widetilde{\caC \times \caC} \times \caC \to \widetilde{\caC \times \caC}$ of the map $\widetilde{\delta_\caC \times \delta_\caC}$.
So $\sigma$ sends $B$ to $(B, id \colon B_{0,1} \to X, id \colon Y \to B_{2,3})$.
Consequently, $(\tau')^{-1} \circ \sigma :S(\caC)_3 \to \caQ(\caC) $ lands in constant squares.
Hence $\omega_1 \circ \sigma = \id.$
\end{remark}

To prove Proposition \ref{gr3ehtrer} we use the following lemma:

\begin{lemma}\label{dfghkhh}\label{hterhedh}
Let $\caD \in \Cat_\infty^{hC_2}$ and $H \to \caH(\caD)$ a genuine refinement.

\begin{enumerate}
\item Let $n \geq 0$ be even, $X \in \caH(\Fun([n],\caD))$ and $H' \to \caH(\Fun([n],\caD))$ the induced genuine refinement.
Then the canonical map $$H'_X \to H_{X(\frac{n}{2})}$$ is an equivalence.

\item Let $H' \to \caH(\Fun(\widetilde{[1] \times [1]},\caD))$ be the induced genuine refinement and an object $X \in \caH(\Fun(\widetilde{[1] \times [1]},\caD))$. Then the square
\[
\begin{tikzcd}
H'_X \ar[r] \ar[d]      & H_{X(10)} \ar[d] \\
H_{X(01)} \ar[r]        & H_{X(00)}
\end{tikzcd}
\]
is a pullback square.
\end{enumerate}

\end{lemma}

\begin{proof}
For (1), we observe that the functor $[0]\simeq \{\frac{n}{2}\} \subset [\frac{n}{2}]$ induces an equivalence
$$[0] \coprod_{\widetilde{[0] \coprod [0]}} \widetilde{[\frac{n}{2}] \coprod [\frac{n}{2}]} \simeq [n]$$ of $\infty$-categories with duality, from which we obtain an equivalence $$ \Fun([n],\caD) \simeq \caD \times_{\widetilde{\caD \times \caD}} \widetilde{\Fun([\frac{n}{2}],\caD) \times \Fun([\frac{n}{2}],\caD)}$$ of $\infty$-categories with genuine duality.
Since the cotensors with $C_2$ carry the standard genuine refinement,
the canonical commutative square 
\[
\begin{tikzcd}
H' \ar[r] \ar[d] & H \ar[d] \\
\caH( \Fun([n],\caD)) \ar[r] & \caH(\caD)
\end{tikzcd}
\]
is a pullback square.

For (2), since $\widetilde{[1] \times [1]} \simeq [2] \coprod_{[1]}[2]$
is the pushout in $\Cat_\infty^{hC_2}$, where $[1] \to [2]$ picks the morphism
$0 \to 2,$ we find that $H'_X \simeq H'_{X_{\mid [2]}} \times_{H'_{X_{\mid [1]}}} H'_{X_{\mid [2]}}.$ Thus the result follows from (1).
\end{proof}

\begin{proof}[Proof of Proposition \ref{gr3ehtrer}]

We begin by proving that the functor $$\rho'\colon Q(\caC) \subset S(\caC)_3^{\widetilde{[1] \times [1]}} \xrightarrow{\rho} \lambda(\caC)= \widetilde{\caC^{[1]} \times \caC^{[1]}} \times_{\widetilde{\caC \times \caC}} S(\caC)_3 $$ 
induces an equivalence on underlying $\infty$-categories, and will later look into its behavior with the genuine refinements involved.

By \cref{beta_cocart} of \cref{fjkjjkk}, evaluation $\beta\colon  S(\caC)_3 \to \caC$ at $(2,3) \in \Ar([3])$ is a cartesian fibration. This implies that $\beta$ and evaluation at the target yield an equivalence 
$$\xi\colon \Fun'([1],S(\caC)_3) \simeq \Fun([1],\caC) \times_{\caC} S(\caC)_3,$$
where the left hand side is the full subcategory of $\beta$-cartesian maps. By \cref{alpha_cocart} of \cref{fjkjjkk}, the functor $\alpha\colon S(\caC)_3 \to \caC$ evaluating at $0,1 \in \Ar([3])$ is a cocartesian fibration.

Let the following be a commutative square in $S(\caC)_3$
\begin{equation*}
\begin{tikzcd}
A \ar[r] \ar[d] & B \ar[d]  \ar[d]\\
C \ar[r] & D,
\end{tikzcd}
\end{equation*} 
where both horizontal maps are $\beta$-cartesian and $B \to D$ is $\alpha$-cocartesian.

Observe that the map $A \to C$ is $\alpha$-cocartesian if and only if it induces an equivalence at $(2,3) \in \Ar([3])$. The `only if' direction is clear. Assume now that the map $A \to C$ induces an equivalence at $(2,3) \in \Ar([3])$. Since every $\beta$-cartesian morphism induces an equivalence on $0,1,(0,2), (1,2)$
and every $\alpha$-cocartesian morphism induces an equivalence on $(1,2),(1,3), (2,3)$, the map $A \to C$ induces an equivalence on $(1,2)$. This, since we assumed that it also induces an equivalence at $(2,3) \in \Ar([3])$, we can conclude that it is $\alpha$-cocartesian.

Consequently, the equivalence $\Fun([1],\xi)$ restricts to an equivalence
$$Q(\caC) \simeq \Fun([1],\mathrm{Eq}(\caC)) \times_{\mathrm{Eq}(\caC)}\Fun''([1],S(\caC)_3),$$ 
where $\Fun''([1],S(\caC)_3) \subset \Fun([1],S(\caC)_3)$ is the full subcategory of $\alpha$-cocartesian morphisms and $\mathrm{Eq}(\caC) \subset \Fun([1],\caC)$ is the full subcategory of equivalences.  

Because the diagonal embedding induces an equivalence $\mathrm{Eq}(\caC) \simeq \caC$ inverse to evaluation at the source,
we obtain an equivalence
$$Q(\caC) \simeq \Fun([1],\caC) \times_{\caC}\Fun''([1],S(\caC)_3)$$ 
induced by evaluation at the target. Moreover, as $\alpha$ is a cocartesian fibration, $\alpha$ and evaluation at the source yield an equivalence 
$$\kappa\colon  \Fun''([1],S(\caC)_3) \simeq S(\caC)_3 \times_{\caC} \Fun([1],\caC).$$

Combining these equivalences, we obtain an equivalence $$Q(\caC) \simeq \Fun([1],\caC) \times_{\caC} S(\caC)_3 \times_{\caC} \Fun([1],\caC) $$ of exact $\infty$-categories that coincides with $\rho'.$

Next we prove that the functor $$\tau\colon \lambda(\caC) \xrightarrow{\rho'^{-1}} Q(\caC) \subset S(\caC)_3^{\widetilde{[1] \times [1]}}$$ is exact.
Because $\tau$ is a section of the exact functor $\rho\colon S(\caC)_3^{\widetilde{[1] \times [1]}} \to  \lambda(\caC) $, this implies that
the fully faithful functor $\tau$ preserves and reflects cofibrations and fibrations. 
So if we transfer the exact structure from $\lambda(\caC) $ to $Q(\caC)$ 
via $\rho'^{-1}$, the functor $\rho'$ is an equivalence of
exact $\infty$-categories and $Q(\caC) \subset S(\caC)_3^{\widetilde{[1] \times [1]}}$ is a full exact subcategory.

The functor $\tau$ factors as 
$$\Fun([1],\caC) \times_{\caC} S(\caC)_3 \times_{\caC} \Fun([1],\caC) \xrightarrow{\Fun([1],\caC) \times_{\caC}\kappa^{-1}} \Fun([1],\caC) \times_{\caC}\Fun([1],S(\caC)_3) $$
$$ \xrightarrow{\delta} \Fun([1],\Fun([1],\caC)) \times_{\Fun([1],\caC)}\Fun([1],S(\caC)_3)$$$$ \xrightarrow{\Fun([1],\xi^{-1})} \Fun([1] \times [1],S(\caC)_3),$$
where $\delta$ is induced by the diagonal embedding
$\caC \to \Fun([1],\caC).$
Because $\delta$ is exact, it is enough to see that
$\kappa^{-1},\xi^{-1}$ are exact functors. Since we have that  $$(\xi^{-1})^\op\colon \Fun([1],\caC^\op) \times_{\caC^\op} S(\caC^\op)_3 \simeq \Fun([1],\caC)^\op \times_{\caC^\op} S(\caC)_3^\op \to $$$$ \Fun( [1],S(\caC)_3)^\op \simeq \Fun( [1],S(\caC^\op)_3) $$ is equivalent to $\kappa^{-1}$ when $\caC$ is replaced by $\caC^\op$, it is enough to show that $$\kappa^{-1}\colon S(\caC)_3 \times_{\caC} \Fun([1],\caC) \to \Fun([1],S(\caC)_3)$$ is exact.

Recall that the exact structure on $\Fun([1],S(\caC)_3)$ is given pointwise.
Therefore since $\kappa^{-1}$ is a section of evaluation at the source
$\Fun([1],S(\caC)_3) \to  S(\caC)_3$, it is enough to show that $\kappa^{-1}$ followed by evaluation at the target, $\ev_1\circ \kappa^{-1}\colon\Fun([1],S(\caC)_3) \to  S(\caC)_3$ is exact.

Let $A \to B$ be a cofibration in $ S(\caC)_3$
and $A' \to B'$ the image under $\ev_1\circ \kappa^{-1}$.
There is a pushout square
\begin{equation*}
\begin{tikzcd}
A_{0,i} \coprod_{A_{0,1}} B_{0,1} \ar[r]\ar[d] &  \ar[d]^{}A'_{0,i} \coprod_{A'_{0,1}} B'_{0,1} \simeq A_{0,i} \coprod_{A_{0,1}} B_{0,1} \coprod_{ B_{0,1 }}  B'_{0,1}  \\ B_{0,i} \ar[r] & B'_{0,i} \simeq B_{0,i} \coprod_{ B_{0,1 }}  B'_{0,1}.
\end{tikzcd}
\end{equation*}
We know that the right map is a cofibration if the left map is. Therefore $\ev_1\circ \kappa^{-1}$ preserves cofibrations. It is clear that $\ev_1\circ \kappa^{-1}$ preserves the zero object and cofibers of cofibrations, and so is an exact functor.

Finally, we turn to genuine refinements. Since the genuine refinements of $\widetilde{\caC^{[1]} \times \caC^{[1]}}$ and  ${\widetilde{\caC \times \caC}}$ are standard, it is enough to see that 
the map $\varrho\colon S(\caC)_3^{\widetilde{[1] \times [1]}} \to S(\caC)_3 $ in $\Wald^\gd_\infty$ evaluating at the non-degenerated object $0,1 \in \widetilde{[1] \times [1]}$
induces an equivalence on the fibers of the genuine refinements over
any hermitian object of $ S(\caC)_3^{\widetilde{[1] \times [1]}}$
that lies over an object of $\caQ(\caC).$

Set $Q\coloneqq \widetilde{[1] \times [1]}.$ By \cref{ehtru64e} the canonical embedding $Q \subset \Ar([3]) $ preserving the duality induces an embedding
$S(\caC)_3 \hookrightarrow \caC^Q $ in $\Wald^\gd_\infty.$ 
Hence the map $\varrho$ is the restriction of the map $\kappa\colon \caC^{Q \times Q} \to \caC^Q $ in $\Wald^\gd_\infty$ evaluating at the non-degenerated object $0,1 \in Q$
and we need to see that $\kappa$ induces an equivalence on the fibers of the genuine refinements over any non-degenerated object $X$ of $ \caC^{Q \times Q}$ that lies over an object of $\caQ(\caC)\subset S(\caC)_3^Q \hookrightarrow \caC^{Q \times Q}. $

Let $H' \to \caH^\mathrm{nd}(\caC^Q)$ and $H'' \to \caH^\mathrm{nd}(\caC^{Q \times Q})$ be the genuine refinements. We want to show that the induced map $H''_X \to H'_{\kappa(X)}$
is an equivalence.

Consider the following commutative square of spaces:
\[
\begin{tikzcd}
H''_X \ar[r] \ar[d] &  H'_{X \mid_{Q \times \{0,1\} }} \times_{H'_{X \mid_{Q \times \{(0,0)\} } }} H'_{X \mid_{Q \times \{(1,0)\} } } \ar[d] \\
H'_{\kappa(X)}\ar[r] & H_{X \mid_{\{0,1\} \times \{0,1\} }} \times_{H_{X \mid_{\{0,1\} \times \{(0,0)\} } }} H_{X \mid_{\{0,1\}\times \{(1,0)\} } }.
\end{tikzcd}
\]
By \cref{hterhedh} the vertical maps are equivalences.
Consequently, it is enough to see that the right vertical map
is an equivalence. This map factors as
$$H'_{X \mid_{Q \times \{0,1\} }} \times_{H'_{X \mid_{Q \times \{(0,0)\} } }} H'_{X \mid_{Q \times \{(1,0)\} } } \to $$$$ H_{X \mid_{\{0,1\} \times \{0,1\} }} \times_{H_{X \mid_{\{(0,0)\} \times \{0,1\} } }} H_{X \mid_{\{(1,0)\}\times \{0,1\} } } \times_{(H_{X \mid_{\{0,1\} \times \{(0,0)\} }} \times_{H_{X \mid_{\{(0,0)\} \times \{(0,0)\} } }} H_{X \mid_{\{(1,0)\}\times \{(0,0)\} } })}$$$$ H_{X \mid_{\{0,1\} \times \{(1,0)\} }} \times_{H_{X \mid_{\{(0,0)\} \times \{(1,0)\} } }} H_{X \mid_{\{(1,0)\}\times \{(1,0)\}}} \to $$
$$H_{X \mid_{\{0,1\} \times \{0,1\} }} \times_{H_{X \mid_{\{0,1\} \times \{(0,0)\}}}} H_{X \mid_{\{0,1\}\times \{(1,0)\} } }.$$
By \cref{hterhedh} the first map in the composition is an equivalence.

So it remains to show that the second map in the composition is an equivalence.
To see that, it is enough to prove that the morphism
$$\nu\colon  H_{X \mid_{\{(1,0)\}\times \{(1,0)\}}} \to H_{X \mid_{\{(0,0)\} \times \{(1,0)\}}}$$ (appearing in the last pullback of the middle term) is an equivalence and the square
\begin{equation}\label{pushhh}
\begin{tikzcd}
H_{X \mid_{\{(1,0)\}\times \{0,1\}}}\ar[r] \ar[d] &  H_{X \mid_{\{(1,0)\}\times \{(0,0)\}}} \ar[d] \\
H_{X \mid_{\{(0,0)\}\times \{0,1\}}} \ar[r] & H_{X \mid_{\{(0,0)\}\times \{(0,0)\}}}
\end{tikzcd}
\end{equation}
is a pullback square.
The morphism $\nu$ takes pullback along a lift in $\caH^\mathrm{nd}(\caC)$ of the morphism
$X((0,0), (1,0)) \to X((1,0), (1,0))$ in $\caC.$
So it is enough to see that this morphism is an equivalence.
The morphism $X((0,0)) \to X((1,0))$ in $\caC^{Q}$ belongs to $ Q(\caC) \subset S(\caC)_3 \hookrightarrow \caC^{Q}$ and thus is $\alpha$-cocartesian.
Hence it induces an equivalence on $(1,0),$ which is equivalent to say that as a morphism of $S(\caC)_3$ it induces an equivalence on $(1,2) \in \Ar([3]).$

The non-degenerated object $X$ of $\caC^{Q \times Q}$ induces a functor
$$[1] \times [1] \subset \Lambda_0^2 \times \Lambda_0^2 \simeq \caH^\mathrm{nd}(Q \times Q) \to \caH^\mathrm{nd}(\caC)$$
corresponding to a commutative square 
\begin{equation}\label{pshhht}
\begin{tikzcd}
X \mid_{\{(0,0)\}\times \{(0,0)\}} \ar[r] \ar[d] & X \mid_{\{(1,0)\}\times \{(0,0)\}} \ar[d] \\
X \mid_{\{(0,0)\}\times \{(0,1)\}} \ar[r] &X \mid_{\{(1,0)\}\times \{(0,1)\}}
\end{tikzcd}
\end{equation} in $\caH^\mathrm{nd}(\caC)$ that lies over the canonical commutative square in $\caC$ below
\begin{equation}\label{pshhh}
\begin{tikzcd}
X((0,0),(0,0))  \ar[r] \ar[d] & X((1,0),(0,0)) \ar[d] \\
X((0,0),(0,1)) \ar[r] & X((1,0),(0,1)).
\end{tikzcd}
\end{equation} 

Since the morphism $X((0,0)) \to X((1,0))$ in $S(\caC)_3 \hookrightarrow \caC^{Q}$is $\alpha$-cocartesian, square (\ref{pshhh}) is exact. We know that square (\ref{pushhh}) is induced by square (\ref{pshhht}) on the fibers. Hence the axioms of an exact $\infty$-category with genuine duality guarantee that square (\ref{pushhh}) is a pullback square.
\end{proof}

\subsection{The real additivity theorem}\label{subsec:add_thm}
In this subsection, we prove \cref{thm:add}. The proof is rather intricate, in order to make it clearer, we tried to group in \cref{subsec:lemmas_add} the proofs of many technical results that are needed. However, we still need to introduce some tedious constructions, we ask to read to bear with us while we do so.

Our aim is to show that the map 
$$S(\gamma_3)\circ \pr_1\colon S(S(\caC)_3)\circ \pr_1 \to (S(\widetilde{\caC \times \caC}) \circ \pr_1) \times (S(\caC) \circ \pr_1)$$ 
induces an equivalence after taking geometric realizations in $D\\Wald^\gd_\infty$. For this, we will consider the pullback of said map along a certain map known to induce an equivalence after geometric realizations in $D\\Wald^\gd_\infty$, and show that the pullback maps also induce an equivalence after geometric realization. For this, we will need to explicitly construct certain real simplicial homotopy inverse.

The reader should notice that we begin by working with (non real) simplicial homotopies. Let us fix notation. 

\begin{notation}\label{not:maps_iota_k}
Using the maps $j^k$ defined in \cref{not:maps_jk}, for every exact $\infty$-category $\caC$ and $k=1,2$, we obtain natural transformations $\iota^k\coloneqq S^e(\caC) \circ (j^k)^\op$ between functors $ \Delta^\op \times \Delta^\op \to \Exact_\infty,$
$$\iota^k\colon \mathfrak{B}(\caC)\to \mathfrak{B}_k(\caC),$$
where $\mathfrak{B}(\caC)\coloneqq S^e(\caC) \circ T^\op$ and $\mathfrak{B}_k(\caC)\coloneqq S^e(\caC) \circ \pr_k^\op$.
\end{notation}

Recall the canonical functor $\Ar\coloneqq\Fun([1],-)\colon \Cat_\infty^\op \xrightarrow{ } \Cat_\infty^\op$.

\begin{remark}\label{rmk:simplicial_rho}
Let $\caC$ be an exact $\infty$-category and $\delta_\caC\colon\caC \to \caC^{[1]}$ the diagonal. When we look at the image under the functor $ \caC^{(-)} \circ \Ar^\op\colon \Cat^\op \to \Cat^\op \xrightarrow{} \Exact_\infty $ of the simplicial homotopy $\kappa$ of construction \cref{{fgjhbgch}} we obtain 
$$ (\caC^{(-)} \circ \Ar \circ T_{-,n})^* \xrightarrow{\caC^{(-)} \circ \Ar \circ \kappa} (\caC^{(-)} \circ \Ar \circ (T_{-,n} \times [1]))^*$$

Combining this with the canonical map $\Ar \circ (T_{-,n} \times [1]) \to (\Ar \circ T_{-,n}) \times [1] $ and the equivalence $(\caC^{(-)} \circ ((\Ar \circ T_{-,n}) \times [1]))^* \simeq ((\caC^{[1]})^{(-)} \circ \Ar \circ T_{-,n})^*$, we obtain a map 

$$(\caC^{(-)} \circ \Ar \circ T_{-,n})^*\to ((\caC^{[1]})^{(-)} \circ \Ar \circ T_{-,n})^*$$
that restricts to a map
$$\mathfrak{B}(\caC)_{-,n}^* \xrightarrow{\rho} \mathfrak{B}(\caC^{[1]})_{-,n}^*.$$

The map $\rho$ is a simplicial homotopy in $\Exact_\infty$ from $$\mathfrak{B}(\delta_\caC)_{-,n}\colon \mathfrak{B}(\caC)_{-,n}\to \mathfrak{B}(\caC^{[1]})_{-,n}$$ to $i_1^*(\rho)\colon  \mathfrak{B}(\caC)_{-,n} \to \mathfrak{B}(\caC^{[1]})_{-,n} $ relative to $\mathfrak{B}(\caC^{\{0\} })_{-,n}.$
\end{remark}

\begin{proposition}\label{prop_square_homotopies}
Let $\caC$ be an exact $\infty$-category and $\lambda\colon \caC \to \caC^{[1]}$ the functor sending $X $ to $X \to 0.$ 
For any $n \geq 1$ there is a canonical commutative square of exact $\infty$-categories as below
\[
\begin{tikzcd}[column sep=large]
\mathfrak{B}(\caC)_{-,n} \ar[r, "i_1^*(\rho)"] \ar[d, "\iota^1_{-,n}"']     & \mathfrak{B}(\caC^{[1]})_{-,n} \ar[d, "\iota^1_{-,n}"]  \\
S^e(\caC) \ar[r, "S^e(\lambda)"']    & S^e(\caC^{[1]}),
\end{tikzcd}
\]
where the vertical maps are as in \cref{not:maps_iota_k} and $i^\ast(\rho)$ as in \cref{rmk:simplicial_rho}.
	
\end{proposition}
\begin{proof}

We want to show that $$ S^e(\lambda) \circ \iota^1_{-,n} \simeq \iota^1_{-,n} \circ i_1^*(\rho).$$

In what follows we show that it is enough to prove such equivalence after post-composing on each side with $S^e(\ev_0)$.

Indeed, the functor $\lambda\colon \caC \to \caC^{[1]}$ is fully faithful and the essential image of $\lambda$ is the fiber of evaluation at the target $\caC^{[1]} \to\caC$. Therefore $S^e(\lambda)$ is objectwise fully faithful, and since $S^e$ preserves fibers, the objectwise essential image of $S^e(\lambda)$, that we are going to call $\Theta$, is the fiber of $S^e(\ev_{1})\colon S^e(\caC^{[1]}) \to S^e(\caC)$. Now, the restriction $\Theta \subset S^e(\caC^{[1]}) \xrightarrow{S^e(\ev_{0})} S^e(\caC)$ is inverse to the induced equivalence $S^e(\caC) \to \Theta$ and thus itself an equivalence.

Moreover, the functor $\iota^1_{-,n} \circ i_1^*(\rho)$ lands in $\Theta$, since
$$S^e(\ev_{1}) \circ \iota^1_{-,n} \circ i_1^*(\rho) \simeq \iota^1_{-,n} \circ \mathfrak{B}(\ev_{1})_{-,n} \circ i_1^*(\rho) \simeq \underbrace{\iota^1_{-,n} \circ q^n}_{0} \circ \iota^2_{-,n} \simeq 0$$ in $\Fun(\mathfrak{B}(\caC)_{-,n},S^e(\caC)).$

We now move on to prove that it is in fact the case that 

$$S^e(\ev_{0})\circ S^e(\lambda) \circ \iota^1_{-,n} \simeq S^e(\ev_{0})\circ \iota^1_{-,n} \circ i_1^*(\rho).$$

This follows from the canonical equivalence 
$$ S^e(\ev_{0}) \circ S^e(\lambda) \circ \iota^1_{-,n} \simeq \iota^1_{-,n} \simeq \iota^1_{-,n} \circ \mathfrak{B}(\ev_{0})_{-,n} \circ i_1^*(\rho) \simeq S^e(\ev_{0}) \circ \iota^1_{-,n} \circ i_1^*(\rho), $$

where the next to last equivalence is deduced from the very last part of \cref{rmk:simplicial_rho}.
\end{proof}

We now promote the previous constructions and results preparring for the case of real simplicial homotopies. 

\begin{notation} Given an exact $\infty$-category with genuine duality $\caC$, for $k=1,2$, we set $i^k\coloneqq \widetilde{\iota^k \times \iota^k}$, where $\iota^k$ is as in \cref{not:maps_iota_k}. Note that $i^k$ are maps of real functors $\underline{\Delta}^\op \times \underline{\Delta}^\op \to \Wald^\gd_\infty$. We will write 
$$i^k\colon B(\caC)\to S(\widetilde{\caC\times\caC})\circ\pr_k,$$
after setting $B(\caC)\coloneqq\widetilde{\mathfrak{B}(\caC) \times \mathfrak{B}(\caC)}$ and noticing that $$\widetilde{S(\caC)^u \circ \pr_k \times S(\caC)^u \circ \pr_k} \simeq \widetilde{S(\caC)^u \times S(\caC)^u} \circ \pr_k \simeq S(\widetilde{\caC \times \caC}) \circ \pr_k.$$
\end{notation}

From the simplicial homotopy $\rho$ of \cref{rmk:simplicial_rho}, we construct the following real simplicial homotopy.

\begin{notation}\label{nota} Given an exact $\infty$-category, and $n\geq 0$, we set $$h\coloneqq\widetilde{\rho \times \rho},$$ which is a real simplicial homotopy by the proof of \cref{cresimpho}. We have that the mao $h\colon B(\caC)_{-,n}^{**} \to B(\caC^{[1]})_{-,n}^{**}$ is a real simplicial homotopy in $\\Wald^\gd_\infty$ from $B(\delta_\caC)_{-,n}$ to $ \widetilde{i_1^*(\rho) \times i_1^*(\rho)}$ relative to $ B(\caC^{\{0\}})_{-,n}$.
\end{notation}

With these new constructions, \cref{prop_square_homotopies} gives us the following result.

\begin{corollary}\label{corrrr}
Let $\caC \in \Wald^\gd_\infty$ and $\lambda\colon  \caC \to \caC^{[1]}$ the functor
sending $X $ to $X \to 0.$  For any $n \geq 1$ there is a canonical commutative square in $\Wald^\gd_\infty$ as below.
\[
\begin{tikzcd}[column sep=huge]
B(\caC)_{-,n} \ar[r, "\widetilde{i_1^*(\rho) \times i_1^*(\rho)}"] \ar[d, "i^1_{-,n}"']       & B(\caC^{[1]})_{-,n} \ar[d, "i^1_{-,n}"]  \\
S(\widetilde{\caC \times \caC}) \ar[r, "S(\widetilde{\lambda \times \lambda})"']    & S(\widetilde{\caC^{[1]} \times \caC^{[1]}}) 
\end{tikzcd}
\]
\end{corollary}

We introduce the real version of \cref{not:p_and_q}.

\begin{notation}\label{not:real_p_and_q}
The natural transformations $q^n\colon[n]\to T_{-,n}$ and $p^n\colon [n]\to T_{n,-}$ of functors $ \Delta^\op \to \Delta^\op$ give rise to natural transformations of functors $ \Delta^\op \to (\Wald^\gd_\infty)^u$ as below
$$S(\caC)^u \circ q^n\colon S(\caC)_n \longrightarrow S(\caC)^u \circ T_{-,n}=\mathfrak{B}(\caC)_{-,n}, $$
$$ S(\caC)^u \circ p^n \colon S(\caC)_n \longrightarrow  S(\caC)^u \circ T_{n,-}=\mathfrak{B}(\caC)_{n,-} $$
which, in turn, yield real natural transformations of real functors $ \underline{\Delta}^\op \to \Wald^\gd_\infty$ as below.
$$\mathfrak{q}^n\coloneqq \widetilde{(S(\caC)^u \circ q^n) \times (S(\caC)^u \circ q^n)} \hspace{1em}\text{ and }\hspace{1em} \mathfrak{p}^n\coloneqq \widetilde{(S(\caC)^u \circ p^n) \times (S(\caC)^u \circ p^n)}$$

Spelling out the domain and codomain of $\mathfrak{q}^n$ and $\mathfrak{p}^n$ we get

$$\mathfrak{q}^n\colon \widetilde{S(\caC)_n \times S(\caC)_n} \simeq S(\widetilde{\caC \times \caC})_n \longrightarrow B(\caC)_{-, n},$$
$$\mathfrak{p}^n\colon \widetilde{S(\caC)_n \times S(\caC)_n}  \simeq S(\widetilde{\caC \times \caC})_n \longrightarrow  B(\caC)_{n,-}.$$
\end{notation}

The following two corollaries are real versions of, and follow from, \cref{reeeem,lem_simhoinv} respectively.

\begin{corollary}\label{corrr}
The compositions $$S(\widetilde{\caC \times \caC})_m  \xrightarrow{\mathfrak{p}^n  }  B(\caC)_{n,-} \xrightarrow{i^1_{n,-}} S(\widetilde{\caC \times \caC})_n, $$$$  S(\widetilde{\caC \times \caC})_n \xrightarrow{\mathfrak{q}^n} B(\caC)_{-, n} \xrightarrow{i^2_{-,n} } S(\widetilde{\caC \times \caC})_n$$ are the identities and the composite $i^1_{-,n}  \circ \mathfrak{q}^n$ factors as 

\[
\begin{tikzcd}
S(\widetilde{C \times C})_n\ar[r, "\mathfrak{q}^n"]\ar[dr]      &B(C)_{-, n}\ar[r, "i^1_{-,n}"] &S(\widetilde{C \times C})\\
\   &S(\widetilde{C \times C})_0 \simeq 0.\ar[ur]
\end{tikzcd}
\]
\end{corollary}

\begin{corollary}\label{coro}
The maps $$ i^1_{n,-} = \widetilde{\iota^1_{n,-} \times \iota^1_{n,-}}\colon B(\caC)_{n,-} 
\xrightarrow{  } S(\widetilde{\caC \times \caC})_n$$
$$i^2_{-,n} = \widetilde{\iota^2_{-,n} \times \iota^2_{-,n}}\colon B(\caC)_{-, n}  
\xrightarrow{  } S(\widetilde{\caC \times \caC})_n, 
$$ are, respectively, real simplicially homotopy inverse relative to $S(\widetilde{\caC \times \caC})_n$
to $$\mathfrak{p}^n\colon S(\widetilde{\caC \times \caC})_n \to B(\caC)_{n,-} \hspace{1em} \text{ and }\hspace{1em} \mathfrak{q}^n\colon S(\widetilde{\caC \times \caC})_n \to B(\caC)_{-,n}.$$
\end{corollary}

We conclude with the following remark, that will be useful in a technical part of the proof of the additivity theorem.

\begin{remark}\label{rem_fact_real} By \cref{rem_fact}, the composition 
 $$B(\caC)_{-,n} \xrightarrow{ \widetilde{i_1^*(\rho) \times i_1^*(\rho)}}  B(\caC^{[1]})_{-,n} \xrightarrow{ B(\ev_1)_{-,n}} B(\caC)_{-,n} $$
 factors as 
$$B(\caC)_{-,n}\xrightarrow{i^2_{-,n}} S(\widetilde{\caC \times \caC})_n \xrightarrow{\mathfrak{q}^n}  B(\caC)_{-,n}. $$
\end{remark}

Before embarking on the proof of the Additivity Theorem, we prove the following lemma.

\begin{lemma}\label{coros}
Let $\caC$ be an exact $\infty$-category with genuine duality, and consider the pullback square in $\Wald^\gd_\infty$ as depicted below.

\begin{equation*}\label{pullb2}
\begin{tikzcd}
S(\caC)_3' \ar[r, "\varrho"]\ar[d]      &\ast\times\caC\ar[d, "0\times \caC"]  \\
S(\caC)_3 \ar[r, "\gamma_3"]    & \widetilde{\caC \times \caC} \times \caC
\end{tikzcd}
\end{equation*}
The map $\varrho$ is an equivalence.	
\end{lemma}
\begin{proof}

The pullback $S(\caC)_3'$ is the full subcategory of $S(\caC)_3$ spanned by
the pairs of cofibrations $A \to B \to C$ such that $A =0$ and $B \to C$ is an equivalence, or equivalently that the maps $B \to C, B \to B/A, C \to C/A$, and  $B/A \to C/A$ are equivalences.
Thus the canonical equivalence $S(\caC)_3 \simeq \mathrm{Amb}(\caC) $
of \cref{corqa} restricts to an equivalence
$S(\caC)_3' \simeq {\caC^{\widetilde{[1] \times [1]}}}'$
of exact $\infty$-categories with genuine duality, where
${\caC^{\widetilde{[1] \times [1]}}}' \subset {\caC^{\widetilde{[1] \times [1]}}}$ is the full exact subcategory with genuine duality spanned by the commutative squares in which all edges are equivalences. 

Using this equivalence, we see that the map $\varrho$ identifies with evaluation $\alpha\colon {\caC^{\widetilde{[1] \times [1]}}}' \to \caC $ at the non-degenerated object $(1,0) \in \widetilde{[1] \times [1]}$. We prove now that $\alpha$ is an equivalence. Firstly, we note that the underlying functor of $\alpha$ is inverse to the diagonal embedding
$\caC \to {\caC^{[1] \times [1]}}'$. In addition, by \cref{htwdzfd}, $\alpha$ induces on the fiber over any non-degenerated object $X $ of ${\caC^{\widetilde{[1] \times [1]}}}'$, the canonical functor
$$\Fun_{\caH^\mathrm{nd}(\caC)}(\caH^\mathrm{nd}(\widetilde{[1] \times [1]}), H) \to H_{X(1,0)},$$ which is in turn an equivalence because $\caH^\mathrm{nd}(\widetilde{[1] \times [1]} \simeq \Lambda_0^2$ is weakly contractible.
\end{proof}

\begin{proof}[\textbf{Proof of \cref{thm:add}}]

Note that the geometric realization of $S(\gamma_3)$ in $\D\Wald^\gd_\infty$ is the two-fold geometric realization of the real bisimplicial map below. 
$$S(\gamma_3)\circ \pr_1\colon S(S(\caC)_3)\circ \pr_1 \to (S(\widetilde{\caC \times \caC}) \circ \pr_1) \times (S(\caC) \circ \pr_1).$$ 

Therefore, we will work with the pullback square in $\Fun_{\Spc^{C_2}}(\underline{\Delta}^\op\times\underline{\Delta}^\op, \Wald^\gd_\infty)$ as below.
\begin{equation}\label{sq:pullback_O}
\begin{tikzcd}
\mathcal{O} \ar[rr, "f"]        \ar[d, "\theta"]       &       &B(C) \times (S(C)\circ \pr_1) \ar[d, "i^1 \times (S(C) \circ \pr_1)"] \\
S(S(C)_3) \circ \pr_1 \ar[rr, "S(\gamma_3) \circ \pr_1"]        &       &(S(\widetilde{C \times C}) \circ \pr_1) \times (S(C)\circ \pr_1)
\end{tikzcd}
\end{equation}

We start by showing that the vertical maps $i^1 \times S(C) \circ \pr_1$ and $\theta$ induce an equivalence on two-fold  geometric realizations in $\D\Wald^\gd_\infty$. To see this, it is enough to check that for any
$n \geq 0$
the vertical maps in the pullback square evaluated at $n$
\[
\begin{tikzcd}
\mathcal{O}_{n,-} \ar[rr, "f_{n,-}"]   \ar[d, "\theta_{n,-}"]  &  &B(\caC)_{n,-} \times S(\caC)_{n} \ar[d, "i^1_{n,-} \times S(\caC)_n"] \\
S(S(\caC)_3)_n \ar[rr, "S(\gamma_3)_n"] &  &S(\widetilde{\caC \times \caC})_n \times S(\caC)_n
\end{tikzcd}
\]
are equivalences in $rs\Wald^\gd_\infty$. This follows from the fact that 
$i^1_{n,-}\colon B(\caC)_{n,-} \to S(\widetilde{\caC \times \caC})_n$
is a real simplicial homotopy equivalence relative to $ S(\widetilde{\caC \times \caC})_n$ by \cref{coro}. Indeed, this implies that 
$$i^1_{n,-} \times S(\caC)_n \colon B(\caC)_{n,-} \times S(\caC)_n\xrightarrow{  } S(\widetilde{\caC \times \caC})_n \times S(\caC)_n $$ 
is a real simplicial homotopy equivalence relative to $S(\widetilde{\caC \times \caC})_n \times S(\caC)_n.$ Which, in turn, implies that $\theta_{n,-}$, its pullback along $S(\gamma_3)_n$ is a real simplicial homotopy equivalence relative to $S(S(\caC)_3)_n.$

Having proven that the vertical maps of the pullback \ref{sq:pullback_O} induce equivalences, it is enough for our purposes to show that the top morphism $f$ also induces an equivalence on two-fold geometric realizations in $\D\Wald^\gd_\infty$. To see this, it is enough to check that for any $n \geq 0$
the real simplicial map $ f_{-,n}$ induces an equivalence on geometric realizations in $\D\Wald_\infty^\gd$.

We show this by constructing a real simplicial homotopy inverse, $\alpha$, of the composition of $f_{-,n}$ with the real simplicial homotopy equivalence $i^2_{-,n} \times S(\caC)$ as depicted below, where we call $P\coloneqq \mathcal{O}_{-,n}$: 
$$g\colon P \xrightarrow{f_{-,n}}  B(\caC)_{-,n} \times S(\caC)\xrightarrow{i^2_{-,n} \times S(\caC)} S(\widetilde{\caC \times \caC})_n \times S(\caC)$$

In order to do this, we consider the pullback square 
\begin{equation}\label{addthm_pb_K}
\begin{tikzcd}
S(\caC)'_3\ar[r]\ar[d]       &\ast\times \caC\ar[d]\\
S(\caC)_3\ar[r, "\gamma_3"] &\widetilde{\caC\times \caC}\times \caC.
\end{tikzcd}
\end{equation}
Note that its top horizontal map is an equivalence by \cref{coros}. Then we consider the following pullback squares in $\rs\Wald^\gd_\infty$, where the outer square is the image under $S$ of (\ref{addthm_pb_K}).

\[
\begin{tikzcd}\label{fdghkkll}
S(S(\caC)'_3)\ar[r, "\simeq"]\ar[d]      &\ast\times S(\caC)\ar[d]\\
P\ar[r, "f_{-,n}"]\ar[d]  &B(C)_{-,n} \times S(C)\ar[d, "i^1_{-,n} \times \id"]\\
S(S(\caC)_3) \ar[r, "S(\gamma_3)"] &S(\widetilde{\caC \times \caC}) \times S(\caC)
\end{tikzcd}
\]

Now, we can take the horizontal map below, given by the composite depicted,

\[
\begin{tikzcd}[column sep=tiny]
S(\widetilde{\caC \times \caC})_n \times S(\caC)\ar[rr, "\varphi"]\ar[dr]      &     &S(S(\caC)_3)\\
\  &* \times S(\caC) \simeq S(S(\caC)'_3)\ar[ru]
\end{tikzcd}
\]
which together with the map 
$$ \mathfrak{q}_n \times S(\caC)\colon S(\widetilde{\caC \times \caC})_n \times S(\caC) \to B(\caC)_{-,n} \times S(\caC)$$
define the map 
\begin{equation}\label{addthm_alpha}
\alpha\colon S(\widetilde{\caC \times \caC})_n \times S(\caC) \to P(\caC),
\end{equation}
as we have, by \cref{corrr}, a commutative square in $\rs\Wald^\gd_\infty$ as below,
\[
\begin{tikzcd}
S(\widetilde{\caC \times \caC} )_n \times S(\caC) \ar[r, "\mathfrak{q}_n \times S(\caC)"] \ar[d]    & B(\caC)_{-,n} \times S(\caC) \ar[d, "i^1_{-,n} \times S(\caC)"] \\
\ast \times S(\caC) \ar[r]     & S(\widetilde{\caC \times \caC} ) \times S(\caC) 
\end{tikzcd}
\]
that induces in turn a commutative square

\[
\begin{tikzcd}
S(\widetilde{\caC \times \caC} )_n \times S(\caC) \ar[r, "\mathfrak{q}_n \times S(\caC)"] \ar[d, "\varphi"']    & B(\caC)_{-,n} \times S(\caC) \ar[d, "i^1_{-,n} \times S(\caC)"] \\
S(S(\caC)_3) \ar[r, "S(\gamma_3)"']     & S(\widetilde{\caC \times \caC} ) \times S(\caC). 
\end{tikzcd}
\]

Note that the map $\alpha$ is a section of $g$, as we have 
$$g \circ \alpha = (i^2_{-,n} \times S(\caC)) \circ f_{-,n} \circ \alpha  =  (i^2_{-,n} \times S(\caC)) \circ (\mathfrak{q}_n \times S(\caC))= \id. $$ 

It then remains to construct a real simplicial homotopy from the identity
to $$\alpha \circ g\colon P(\caC)_{-,n} \to S(\widetilde{\caC \times \caC} )_n \times S(\caC) \to P(\caC)_{-,n}.$$ For this, let us consider now the commutative squares, for $i=0,1$

\begin{equation}\label{addthm_pb_Q}
\begin{tikzcd}
Q \ar[r]\ar[d]           &B(\caC^{[1]})_{-,n} \times S(\caC)\ar[d]\\
S(\lambda(\caC))\ar[r, "S(\beta)"]\ar[d, "S(\omega_i)"']        &S(\widetilde{\caC^{[1]} \times \caC^{[1]}}) \times S(\caC)\ar[d, "S(\widetilde{\ev_i \times \ev_i}) \times \id"]\\
S(S(\caC)_3)\ar[r, "S(\gamma_3)"]      &S(\widetilde{\caC \times \caC}) \times S(\caC),
\end{tikzcd}
\end{equation}
where $\lambda(\caC)$ and $\beta$ are as in \cref{not:lambda(C)}, $\omega_0$ is projection, $\omega_1$ is the map in \cref{const_omega1}, and the top square is by definition a pullback square, whilst the bottom square is a pullback square for $i=0$. Similarly, let us also consider the pullback squares below for $i=0,1$,

\begin{equation}\label{addthm_Pi}
\begin{tikzcd}
P_i \ar[r]\ar[d, "s_i"]            &B(\caC^{[1]})_{-,n} \times S(\caC)\ar[d, "B(\ev_i)_{-,n} \times \id"]\\
P\ar[r, "f_{-,n}"]\ar[d]   &B(\caC)_{-,n} \times S(\caC)\ar[d, "i^1_{-,n} \times \id"]\\
S(S(\caC)_3)\ar[r, "S(\gamma_3)"]          &S(\widetilde{\caC \times \caC}) \times S(\caC),
\end{tikzcd}
\end{equation}
from which we obtain canonical maps $Q \to P_i$ since the right vertical maps in the outer squares of (\ref{addthm_Pi}) are canonically equivalent to the right vertical map in the outer square of (\ref{addthm_pb_Q}) for $i=0,1$.

Since the bottom square of (\ref{addthm_pb_Q}) is a pullback square for $i=0$, so is the outer square, and thus the map $Q \to P_0$ is an equivalence. Let $l$ be the composition 
\begin{equation}\label{addthm_l}
P_0\simeq Q \to P_1 \xrightarrow{s_1} P.
\end{equation}

Following \cref{nota}, there is a real simplicial homotopy 
$$h \colon (B(\caC)_{-,n} \times S(\caC))^{**} \to (B(\caC^{[1]})_{-,n}\times S(\caC))^{**} $$ in $\\Wald^\gd_\infty$ that goes from $B(\delta_\caC)_{-,n}\times S(\caC)$ to  $\widetilde{i_1^*(\rho)\times i_1^*(\rho)} \times S(\caC) $, relative to $B(\caC)_{-,n}\times S(\caC).$

Pulling back  $h$ along $f\colon P \to B(\caC)_{-,n} \times S(\caC)$, we obtain a real simplicial homotopy 
$$H\colon P^{**} \to P_0^{**}$$
that goes from $f^*(B(\delta_\caC)_{-,n}\times S(\caC))\colon P \to P_0 $ to $f^*(\widetilde{i_1^*(\rho) \times i_1^*(\rho)} \times S(\caC))\colon P \to P_0$ relative to $P$, where we view $P_0$ over $P$ via $s_0.$

Finally, considering the composition below 
$$l \circ H\colon P^{**} \to P_0^{**} \to P^{**}$$
we obtain a real simplicial homotopy from $$l \circ f^*(B(\delta_\caC)_{-,n}\times S(\caC))\colon P \to P_0\to P$$ to $$l \circ f^*(\widetilde{i_1^*(\rho) \times i_1^*(\rho)} \times S(\caC))\colon P \to P_0 \to P.$$ 

\cref{old:propos} guarantees that these source and target are respectively the identity and the composition $\alpha \circ g\colon P \to  S(\widetilde{\caC \times \caC})_n \times S(\caC) \to P$.

\end{proof}

To conclude the previous proof, we used that the source and target of the real simplicial homotopy are exactly what we want. We state this precisely below.

\begin{proposition}\label{old:propos}
Let $\caC \in \Wald^\gd_\infty.$ The following holds.
\begin{enumerate}
    \item\label{real_simplicial_htpy_source} The composition $$\sigma_1\coloneqq l \circ f^*(B(\delta_\caC)_{-,n}\times S(\caC))\colon P \to P' \to P $$ is the identity on $P$.
    \item\label{real_simplicial_htpy_target} The composition  $$\sigma_2 \coloneqq l \circ f^*(\widetilde{i_1^*(\rho) \times i_1^*(\rho)} \times S(\caC))\colon P \to P' \to P $$ is
$\alpha \circ g \colon P \to  S(\widetilde{\caC \times \caC})_n \times S(\caC) \to P.$
\end{enumerate}
\end{proposition}

Before presenting the proof of this proposition, we will need to fix some notation and prove two results. 

\begin{notation}\label{not:sections_omega_0}
Let $\caC \in \Wald^\gd_\infty$ and $\lambda\colon  \caC \to \caC^{[1]}$ the functor
sending $X $ to $X\to 0.$ We denote by $\sigma,\chi\colon S(\caC)_3 \to \lambda(\caC)$ the sections of 
$\omega_0\colon \lambda(\caC) \to S(\caC)_3$ that are the pullbacks along $S(\caC)_3 \xrightarrow{\gamma_3} \widetilde{\caC \times \caC} \times \caC \to \widetilde{\caC \times \caC}$, respectively, of the sections $\widetilde{\delta_\caC \times \delta_\caC}, \widetilde{\lambda \times \lambda}\colon \widetilde{\caC \times \caC} \to \widetilde{\caC^{[1]} \times \caC^{[1]}}$ of the map
$\widetilde{ev_0 \times \ev_0}\colon \widetilde{\caC^{[1]} \times \caC^{[1]}} \to \widetilde{\caC \times \caC}$.
\end{notation}

In the following lemma we are using notation introduced in the proof of \cref{thm:add}, specifically, $P$ is as defined before diagram (\ref{addthm_pb_K}) and $P_0$ as defined as in (\ref{addthm_Pi}).

\begin{lemma}\label{bwas}
Let $\caC \in \Wald^\gd_\infty$.
There are canonical commutative squares in $\Wald^\gd_\infty$ as depicted below,
\[
\begin{tikzcd}[column sep=huge]
P\ar[r, "f^*(B(\delta_\caC)_{-,n}\times S(\caC))"]\ar[d]            &P_0\ar[d]\\
S(S(\caC)_3)\ar[r, "S(\sigma)"]                                     &S(\lambda(\caC)),
\end{tikzcd}
\hspace{1em}
\begin{tikzcd}[column sep=huge]
P\ar[r, "f^*(\widetilde{i_1^*(\rho)\times i_1^*(\rho)}\times S(\caC))"]\ar[d]           &P_0\ar[d]\\
S(S(\caC)_3)\ar[r, "S(\chi)"]                                     &S(\lambda(\caC)),
\end{tikzcd}
\]
where $\delta_\caC\colon\caC\to\caC^{[1]}$ is the diagonal map, $\lambda(\caC)$ is as in \cref{not:lambda(C)}, and the map $\rho\colon S(\caC)_3^{\widetilde{[1] \times [1]}} \to \lambda{\caC}$ as in \cref{Not:rho}.
\end{lemma}

\begin{proof}
We will prove both commutativities simultaneously. We will use as an auxiliary tool the pullback square of \cref{not:lambda(C)}.

It is easy to check that all three compositions $$P \xrightarrow{f^*(B(\delta_\caC)_{-,n}\times S(\caC)) }P_0 \to S(\lambda(\caC)) \xrightarrow{S(\omega_0)} S(S(\caC)_3),$$
$$P \xrightarrow{f^*(\widetilde{i_1^*(\rho)\times i_1^*(\rho)} \times S(\caC))}P_0 \to S(\lambda(\caC)) \xrightarrow{S(\omega_0)} S(S(\caC)_3),$$
$$ P \to S(S(\caC)_3) \xrightarrow{S(\sigma), S(\chi)}S(\lambda(\caC))
\xrightarrow{S(\omega_0)}S(S(\caC)_3)$$ are the projection $P \to S(S(\caC)_3)$.

On the other hand, the compositions $$P \xrightarrow{f^*(B(\delta_\caC)_{-,n}\times S(\caC)) }P_0 \to S(\lambda(\caC)) \xrightarrow{S(\beta)} S(\widetilde{\caC^{[1]} \times \caC^{[1]}}) \times S(\caC),$$
$$P \xrightarrow{f^*(\widetilde{i_1^*(\rho)\times i_1^*(\rho)} \times S(\caC))}P_0 \to S(\lambda(\caC)) \xrightarrow{S(\beta)} S(\widetilde{\caC^{[1]} \times \caC^{[1]}}) \times S(\caC)$$
factor, respectively, as 
$$P \xrightarrow{f }B(\caC)_{-,n} \times S(\caC) \xrightarrow{B(\delta_\caC)_{-,n}\times S(\caC) } B(\caC^{[1]})_{-,n} \times S(\caC) \xrightarrow{i^1_{-,n} \times S(\caC) } S(\widetilde{\caC^{[1]} \times \caC^{[1]}}) \times S(\caC),$$
$$P \xrightarrow{f }B(\caC)_{-,n} \times S(\caC) \xrightarrow{\widetilde{i_1^*(\rho)\times i_1^*(\rho)} \times S(\caC) } B(\caC^{[1]})_{-,n} \times S(\caC) \xrightarrow{i^1_{-,n} \times S(\caC) } S(\widetilde{\caC^{[1]} \times \caC^{[1]}}) \times S(\caC).$$

By \cref{corrrr} these, in turn, identify with the following two maps, respectively 
$$P \xrightarrow{}S(S(\caC)_3) \xrightarrow{S(\gamma_3)} S(\widetilde{\caC \times \caC}) \times S(\caC) \xrightarrow{S(\widetilde{\delta_\caC \times \delta_\caC})} S(\widetilde{\caC^{[1]} \times \caC^{[1]}}) \times S(\caC),$$
$$P \xrightarrow{}S(S(\caC)_3) \xrightarrow{S(\gamma_3)} S(\widetilde{\caC \times \caC}) \times S(\caC) \xrightarrow{S(\widetilde{\lambda \times \lambda})} S(\widetilde{\caC^{[1]} \times \caC^{[1]}}) \times S(\caC).$$

Finally, we note that they are, respectively, the following:
$$ P \to S(S(\caC)_3) \xrightarrow{S(\sigma),S(\chi)}S(\lambda(\caC)) \xrightarrow{S(\beta)}S(\widetilde{\caC^{[1]} \times \caC^{[1]}}) \times S(\caC).$$
\end{proof}

\begin{notation}\label{not:tau}
We denote by $\tau\colon S(\caC)_3 \to S(\caC)_1 \simeq \caC$ the map that sends 
$$(A \to B \to C) \mapsto B/A,$$ induced by the map $[1] \simeq \{1\to 2\} \subset [3]$ in $\underline{\Delta}^\op.$
\end{notation}

\begin{lemma}\label{awas}
Let $\caC \in \Wald^\gd_\infty$.
The map $$ S(\caC)_3 \xrightarrow{\chi} \lambda(\caC) \xrightarrow{\omega_1} S(\caC)_3,$$ 
where $\chi$ is as in \cref{not:sections_omega_0} and $\omega_1$ is as in \cref{const_omega1}, is canonically equivalent to the map
$$  S(\caC)_3 \xrightarrow{\tau} \caC \simeq S(\caC)'_3 \to S(\caC)_3,$$ 
where $S(C)_3^\prime$ is as in diagram (\ref{addthm_pb_K}), the right-most morphism is the canonical embedding, and the middle equivalence is provided by \cref{coros}.
\end{lemma}

\begin{proof}
By definition, the map $ S(\caC)_3 \xrightarrow{\chi} \lambda(\caC) \xrightarrow{\omega_1} S(\caC)_3$ sends $(A \to B \to C)$ to $$\omega_1(A \to B \to C, A \to 0, 0 \to C/B),$$ which is $(0 \to B/A \xrightarrow{\id} B/A).$ Hence such composition actually lands in $S(\caC)'_3.$
When we restrict $\tau$ to $S(\caC)'_3$, we obtain the equivalence $\varrho\colon  S(\caC)'_3 \simeq \caC$ of \cref{coros}. Then, it is enough to see that
$$ S(\caC)_3 \xrightarrow{\omega_1\circ \chi} S(\caC)'_3 \xrightarrow{\varrho} \caC$$ 
is canonically equivalent to
$$  S(\caC)_3 \xrightarrow{\tau} \caC \simeq S(\caC)'_3 \xrightarrow{\varrho} \caC,$$ which is $\tau.$

Now, we know by \cref{fjbnjjkk} that both maps $\tau \circ \omega_0, \tau \circ \omega_1$ are canonically equivalent.
Thus the composition $S(\caC)_3 \xrightarrow{\chi} \lambda(\caC) \xrightarrow{\omega_1} S(\caC)_3 \xrightarrow{\tau} \caC$
identifies with the 
composition $S(\caC)_3 \xrightarrow{\chi} \lambda(\caC) \xrightarrow{\omega_0} S(\caC)_3 \xrightarrow{\tau} \caC$
that is $\tau.$
\end{proof}

\begin{proof}[Proof of \cref{old:propos}]
The strategy to prove both equivalences is the same, we will describe it now for the first case. We begin by considering the pullback diagram below.

\[
\begin{tikzcd}
\textcolor{gray}{P}\ar[gray, shift left=0.5ex, "\id_P"']{dr} \arrow[gray, shift right=0.5ex, "\sigma_1"]{dr}           &   & \\
\           &P\ar[r, "f_{-,n}"]\ar[d, "\pi"']          &B(\caC)_{-,n}\times S(\caC)\ar[d, "i^1_{-,n}\times S(\caC)"]\\
\       &S(S(\caC)_3)\ar[r, "S(\gamma_3)"']          &S(\widetilde{\caC\times\caC})\times S(\caC).
\end{tikzcd}
\]
To prove the equivalence in \cref{real_simplicial_htpy_source}, we will show that composing the maps in gray with each leg of the pullback square we obtain compatible equivalent maps. This is, that $$\f_{-,n}\circ \sigma_1\simeq \f_{-,n}\circ \id_P \hspace{1em}\text{and}\hspace{1em} \pi\circ \sigma_1\simeq \pi\circ \id_P,$$ and moreover they fit in a commutative square as below.

\begin{equation}\label{square_Sigma1}
\begin{tikzcd}
(i^1_{-,n} \times S(\caC)) \circ f \circ \sigma_1  \ar[d, "\simeq"'] \ar[r]     & (i^1_{-,n} \times S(\caC)) \circ f \circ \id_P\ar[d, "\simeq"] \\
S(\gamma_3) \circ \pi \circ \sigma_1\ar[r] & S(\gamma_3) \circ \pi \circ \id_P
\end{tikzcd}
\end{equation}

For the equivalence $\f_{-,n}\circ \sigma_1\simeq \f_{-,n}\circ \id_P$, we consider the composition 
$$ P \xrightarrow{f^*(B(\delta_\caC)_{-,n}\times S(\caC))} P' \xrightarrow{l}P \xrightarrow{f} B(\caC)_{-,n} \times S(\caC),$$
where $\delta_C$ is as in \cref{rmk:simplicial_rho}. By construction of $l$, this canonically factors as 
$$ P \xrightarrow{f^*(B(\delta_\caC)_{-,n}\times S(\caC))}P' \to B(\caC^{[1]})_{-,n} \times S(\caC) \xrightarrow{B(\ev_1)_{-,n} \times S(\caC)} B(\caC)_{-,n} \times S(\caC).$$ Which, in turn, factors as $$ P \xrightarrow{f} B(\caC)_{-,n} \times S(\caC) \xrightarrow{B(\delta_\caC)_{-,n}\times S(\caC)} B(\caC^{[1]})_{-,n} \times S(\caC) \xrightarrow{B(\ev_1)_{-,n} \times S(\caC)} B(\caC)_{-,n} \times S(\caC), $$
which is $f$ by functoriality of $B$ and the fact that $\ev_1\circ\delta_\caC=\id$.

For the equivalence $\pi\circ \sigma_1\simeq \pi\circ \id_P$, we first observe that the composition $$ P \xrightarrow{f^*(B(\delta_\caC)_{-,n}\times S(\caC))} P' \xrightarrow{l}P \xrightarrow{\pi} S(S(\caC)_3)$$ 
factors canonically as 
$$ P \xrightarrow{f^*(B(\delta_\caC)_{-,n}\times S(\caC))} P' \simeq Q \xrightarrow{}S(\lambda(\caC)) \xrightarrow{S(\omega_1)} S(S(\caC)_3);$$
to see this factorisation the reader may find useful to look at diagram (\ref{addthm_Pi}). This composition is, by \cref{bwas}, the following $$ P \xrightarrow{\pi} S(S(\caC)_3) \xrightarrow{S(\sigma)} S(\lambda(\caC)) \xrightarrow{S(\omega_1)} S(S(\caC)_3)$$ 
To conclude that the above composition identifies with the map $\pi\colon P \xrightarrow{} S(S(\caC)_3)$, it is enough to show that $\omega_1\circ \sigma=\id$, which is true as we showed in \cref{rmk:comp_omega1_sigma}.

We move now to prove the equivalence in \cref{real_simplicial_htpy_target}. For this, we want equivalences $$f \circ \sigma_2 \simeq f \circ \alpha \circ g \hspace{1em}\text{and}\hspace{1em}\pi \circ \sigma_2 \simeq \pi \circ \alpha \circ g$$ that induce a commutative square as depicted below.
\begin{equation}\label{square_Sigma2}
\begin{tikzcd}
(i^1_{-,n} \times S(\caC)) \circ f \circ \sigma_2 \ar[d, "\simeq"'] \ar[r] & (i^1_{-,n} \times S(\caC)) \circ f \circ \alpha \circ g\ar[d, "\simeq"] \\
S(\gamma_3) \circ\pi \circ \sigma_2\ar[r] &  S(\gamma_3) \circ \pi \circ \alpha \circ g.
\end{tikzcd}
\end{equation}

For the equivalence $f \circ \sigma_2 \simeq f \circ \alpha \circ g$, similarly to what we did before, we consider the composition
 $$ P \xrightarrow{f^*(\widetilde{i_1^*(\rho) \times i_1^*(\rho)} \times S(\caC))} P' \xrightarrow{l}P \xrightarrow{f} B(\caC)_{-,n} \times S(\caC)$$ 
and observe that, by construction of $l$, it canonically factors as
$$ P \xrightarrow{f} B(\caC)_{-,n} \times S(\caC) \xrightarrow{\widetilde{i_1^*(\rho)\times i_1^*(\rho)} \times S(\caC)} B(\caC^{[1]})_{-,n} \times S(\caC) \xrightarrow{B(\ev_1)_{-,n} \times S(\caC)} B(\caC)_{-,n} \times S(\caC).$$ 
By \cref{rem_fact_real}, this composition factors, in turn, as 
$$ P \xrightarrow{f} B(\caC)_{-,n} \times S(\caC) \xrightarrow{i^2_{-,n} \times S(\caC)} S(\widetilde{\caC\times \caC})_{n} \times S(\caC)\xrightarrow{q^n \times S(\caC)} B(\caC)_{-,n} \times S(\caC)$$ 
that can be seen, using the definitions of $\alpha$ and $g$, it identifies with 
$$ P \xrightarrow{g} S(\widetilde{\caC\times \caC})_{n} \times S(\caC) \xrightarrow{\alpha} P \xrightarrow{f} B(\caC)_{-,n} \times S(\caC).$$ 

For the equivalence $\pi\circ\sigma_2\simeq \pi\circ\alpha\circ\ g$, we observe that the composition
$$ P \xrightarrow{f^*(\widetilde{i_1^*(\rho)\times i_1^*(\rho)} \times S(\caC))} P' \xrightarrow{l}P \xrightarrow{\pi} S(S(\caC)_3),$$ 
by definition of $l$, canonically factors as 
$$P \xrightarrow{f^*(\widetilde{i_1^*(\rho)\times i_1^*(\rho)} \times S(\caC))} P' \simeq Q \xrightarrow{}S(\lambda(\caC)) \xrightarrow{S(\omega_1)} S(S(\caC)_3).$$ Which is, by \cref{bwas}, this is the same as 
$$ P \xrightarrow{} S(S(\caC)_3) \xrightarrow{S(\chi)} S(\lambda(\caC)) \xrightarrow{S(\omega_1)} S(S(\caC)_3).$$ 

To conclude, we will show that the composition $\pi\circ\alpha\circ\ g$ can also be identified with the last factorisation above. Looking at the definitions of both $g$ and $\alpha$, we can see that the composition
 $$ P \xrightarrow{g} S(\widetilde{\caC\times \caC})_{n} \times S(\caC) \xrightarrow{\alpha} P \xrightarrow{\pi} S(S(\caC)_3)$$ 
factors as 
$$ P \xrightarrow{f} B(\caC)_{-,n} \times S(\caC) \to S(\caC) \simeq S(S(\caC)'_3) \to S(S(\caC)_3),$$

which simplifies to  
$$ P \xrightarrow{} S(S(\caC)_3) \xrightarrow{S(\tau)} S(\caC) \simeq S(S(\caC)'_3) \to S(S(\caC)_3)$$
by definition of $f$ and $\tau$, where $\tau$ is as in \cref{awas}. Finally, using \cref{awas}, we obtain that such composition identifies with the composition below, as we wanted.

$$ P \xrightarrow{} S(S(\caC)_3) \xrightarrow{S(\chi)} S(\lambda(\caC)) \xrightarrow{S(\omega_1)} S(S(\caC)_3).$$ 

\end{proof}

\section{\texorpdfstring{A universal property of real $K$-theory}{A universal property of real K-theory}}\label{sec:universal_prop}

In this section we characterize the real $K$-theory functor $\KR$ (\cref{def:KR_space}) as initial among certain real functors
$$\D\Wald_\infty^\gd \to \Spc^{C_2}_\ast$$ satisfying the addivity theorem \ref{gghjgfdd}, that we call additive theories. The way we will prove this is by showing a general universal property of the additivization (\cref{subsec:additivization,subsec:universal_add_th}) and defining the real $K$-theory functor as the additivization of the functor 
$$\iota\colon \D\Wald_\infty^\gd \to \Spc^{C_2},$$ 
that sends a Waldhausen $\infty$-category with genuine duality to its maximal subspace.

In \cref{subsec:pre_add_theoroes} we introduce real versions of preadditive, semiadditive and additive theories. In \cref{subsec:additivization} we give a way to obtain, from a preadditive theory $\D\Wald_\infty^\gd$, an additive theory \textemdash we call this process additivization. In order to prove that such associated theory is indeed additive (\cref{additivization_is_additive}), we need two bits of extra theory that we have separated in different subsections. Namely, \cref{subsec:real_bar_construction} that deals with a real version of the bar construction; and \cref{subsec:left_actions} which deals with left actions of Segal spaces \textemdash we defer the latter to the end of the section due to its technical nature.

Finally, in \cref{subsec:universal_add_th} we prove that the additivization functor is a real left adjoint of the embedding 
$$\Add\hookrightarrow \Pre\Add,$$
of additive theories into preadditive theories; this is \cref{univ}. 

It is worth noticing that for this we make use of the real bar construction again. Indeed, the main ingredient in the proof of \cref{univ} is \cref{dfghjhgk}, which shows that every semiadditive theory  identifies the real S-construction with the real bar construction.

Throughout this section we use that a pointed real $\infty$-category is canonically a $\Spc^{C_2}_\ast$-enriched $\infty$-category, as we describe now. The forgetful functor 
$ \Cat_\infty^{\Spc^{C_2}_\ast} \to \Cat_\infty^{\Spc^{C_2}}$ restricts to an equivalence between pointed $\Spc^{C_2}_\ast$-enriched $\infty$-categories and pointed real $\infty$-categories and reduced real functors. Moreover, given pointed real $\infty$-categories $\caC, \caD$ the reduced real functor $ \Fun_{\Spc^{C_2}_\ast  }(\caC, \caD) \to \Fun_{\Spc^{C_2}  }(\caC, \caD)$ is fully faithful with essential image the reduced real functors. Morally, if $\caC$ is a pointed real $\infty$-category with zero object 0, then for any $X,Y \in \caC$ the mapping genuine $C_2$-space $\map_\caC(X,Y)$ is pointed by the zero morphism $X \to 0 \to Y$, which gives the $\Spc^{C_2}_\ast$-enrichment. For a formal proof of this fact we refer the reader to \cref{Venr:pointed_vs_not_pointed}.

\subsection{Preadditive and additive theories}\label{subsec:pre_add_theoroes}

This subsection deals with the definition of preadditive, semiadditive and additive theories; \cref{def:pre_semi_additive}. Preadditive theories are real functors
$\D\Wald_\infty^\gd \to \caD$ \textemdash we care specially for the case when the target category is that of genuine $C_2$-spaces \textemdash preserving finite products and cotensors with $C_2$. The promotion to semiadditive is earned by satisfying additivity in the sense of the real addivity theorem, while it is also required a real version of excision to be plainly additive. 

We begin by recalling \cref{ex:real_spine_inc}, where we show that the category $[1]$ has a unique strict duality that restricts to its full subcategory $\{0,1\}.$ Thus by applying the nerve functor $(\Cat_\infty)^{hC_2}\to \sSet^{C_2} \simeq \rs\Set$, the simplicial set $\Delta^1$ can be seen as a real simplicial set, and $\partial\Delta^1$ as a real simplicial subset.

\begin{definition} We define the pointed genuine $C_2$-space $S^{1,1}$ as the geometric realization of $\mathfrak{S}$ as a real simplicial space (see \cref{def:real_geometric_realization}).
\end{definition}

\begin{remark}
The pointed genuine $C_2$-space $S^{1,1}$
is precisely the canonical genuine $C_2$-space associated to the $C_2$-action on the one point compactification of the real numbers $\mathbb{R}$ endowed with the $C_2$-action given by $X \mapsto -X$.
\end{remark}

\begin{notation}Let $\caC$ be a pointed real $\infty$-category.
\begin{itemize}
\item If $\caC$ admits all cotensors with $S^{1,1}$ as an $\Spc^{C_2}_\ast$-enriched $\infty$-category, we write $ \Omega^{1,1}\colon \caC \to \caC $ for the cotensor with $S^{1,1}$.

\item Dually if $\caC$ admits all tensors with $S^{1,1}$ as an $\Spc^{C_2}_\ast$-enriched $\infty$-category, we write $ \Sigma^{1,1}\colon \caC \to \caC $ for the tensor with $S^{1,1}$.
\end{itemize}
\end{notation}

Before continuing with the topics proper of this subsection, we provide a description of the $C_2$ fixed points of $\Omega^{1,1}(X)$ for $X$ an object of a pointed real $\infty$-category. This will be useful in future sections.

\begin{lemma}\label{lem:fiber_sequence_with_omega11}
Let $\caD$ be a pointed real $\infty$-category that admits cotensors with $C_2$, and let $X$ be an object of $\caD$. Then there exists a canonical fiber sequence in $\caD$ as below:
$$\Omega^{1,1} X \to X \to \widetilde{X \times X}. $$ 
\end{lemma}

\begin{proof}
The canonical fiber sequence in $\caD$ in the statement
is induced by the canonical cofiber sequence
$$C_2 \coprod [0] \to [0] \coprod [0] \to S^{1,1} $$
of pointed genuine $C_2$-spaces.
\end{proof}

\begin{corollary}\label{lem:description_fixed_points_Omega11}
Let $X$ be a genuine $C_2$-space. Then there exists a canonical fiber sequence of spaces 
$$\Omega^{1,1}(X)^{C_2} \to X^{C_2} \to X^u.$$
\end{corollary}

\begin{proof}
By \cref{lem:fiber_sequence_with_omega11} for $\caD= \Spc^{C_2} $ there is a canonical fiber sequence of spaces $$\Omega^{1,1}(X) \to X \to \widetilde{X \times X}$$ and after taking $C_2$-fixed points
$\Omega^{1,1}(X)^{C_2} \to X^{C_2} \to \widetilde{X \times X}^{C_2} \simeq X^u.$
\end{proof}

The following definition provides a $\Spc^{C_2}_\ast$-enriched notion of excisive functors. 

\begin{definition}Let $F\colon\caC\to\caD$ be a reduced real functor between pointed real $\infty$-categories, where $\caC$ admits tensors with $S^{1,1}$ as $\Spc^{C_2}_\ast$-enriched $\infty$-category and $\caD $ admits cotensors with $S^{1,1}$ as $\Spc^{C_2}_\ast$-enriched $\infty$-category. We say that $F\colon \caC \to \caD$ is genuine excisive if for every $X \in \caC$ the natural map below is an equivalence $$ F(X) \to \Omega^{1,1}(F(\Sigma^{1,1}(X))).$$
\end{definition}

\begin{definition}Let $\caD $ be a pointed real $\infty$-category that admits cotensors with $S^{1,1}$ as $\Spc^{C_2}_\ast$-enriched $\infty$-category. We define a theory with values in $\caD$ as a reduced real functor $\D\Wald^\gd_\infty \to \caD$.
\end{definition}

Recall that for every $\caD \in \Wald^\gd_\infty$ we have constructed a map in $\Wald^\gd_\infty$ 
$$\gamma_3\colon S(\caD)_3 \to \widetilde{\caD \times \caD} \times \caD,$$
given by $(A \to B \to C) \mapsto (A, C/B, B/A)$  (see \ref{const:gamma_3}), that extends to a map in $\\Wald^\gd_\infty$ for $\caD \in \\Wald^\gd_\infty$, by \cref{lem:functors_extend_to_DWald}, for which we use the same name.

\begin{definition}\label{def:pre_semi_additive}Let $\caD$ be a be a pointed real $\infty$-category that admits cotensors with $S^{1,1}$ as $\Spc^{C_2}_\ast$-enriched $\infty$-category, and $\phi\colon\D\Wald_\infty^\gd\to \caD$ a theory. We say that $\phi$
\begin{itemize}
\item is preadditive if it preserves finite products and cotensors with $C_2$,
\item is semiadditive if it is preadditive and inverts the map $\gamma_3\colon S(\caD)_3 \to \widetilde{\caD \times \caD} \times \caD$ in $\D\Wald^\gd_\infty $ for every $\caD \in \D\Wald^\gd_\infty$, and \item is additive if it is semiadditive and genuine excisive.
\end{itemize}
\end{definition}

If $\caD$ admits large sifted colimits, we call a real functor $\Wald^\gd_\infty \to \caD$ a preadditive, semiadditive respectively additive theory if its unique sifted colimits preserving extension $\D\Wald^\gd_\infty \to \caD$ (\cref{lem:functors_extend_to_DWald}) is a preadditive, semiadditive, respectively additive theory.

We now construct a functor that will be fundamental for the additivization of a theory, and provide us with an example of semiadditive theory. 

\begin{construction}\label{const:curly_S}
Via the real embedding $\Wald^\gd_\infty \subset \D\Wald_\infty^\gd$ we extend the real functor 
$$\Wald_\infty^\gd \xrightarrow{S} \mathrm{rs}\Wald_\infty^\gd \subset \mathrm{rs}\D\Wald_\infty^\gd \xrightarrow{ |-|} \D\Wald_\infty^\gd $$
to a sifted colimits preserving real endofunctor $\caS^{1,1}$ of $\D\Wald_\infty^\gd$.
\end{construction}

\begin{example}\label{fgghj}
\begin{enumerate}
\item The real functor $\iota\colon \Wald^\gd_\infty \xrightarrow{ } \widehat{\Spc}^{C_2}_\ast $ that sends a Waldhausen $\infty$-category with genuine duality to the maximal subspace
in its underlying $\infty$-category, is a preadditive theory.

\item The functor $\caS^{1,1}\colon \D\Wald_\infty^\gd \to \D\Wald_\infty^\gd $ preserves finite products and cotensors with $C_2$.
So by the addivity theorem (\cref{thm:add}) for every preadditive theory $\phi\colon \D\Wald_\infty^\gd \to \caD$ the composition $$\phi \circ \caS^{1,1}\colon  \D\Wald_\infty^\gd \to \caD$$ is a semiadditive theory.

\item The real functor $\Omega^{1,1}\colon \D\Wald_\infty^\gd \to \D\Wald_\infty^\gd$ preserves finite products and cotensors with $C_2$. So for every semiadditive (preadditive) theory $\varphi\colon\D\Wald_\infty^\gd \to \caD$ the composition $$\Omega^{1,1} \circ \varphi \colon  \D\Wald_\infty^\gd \to \caD$$ is a semiadditive (preadditive) theory.
\end{enumerate}
\end{example}

\subsection{Real bar construction}\label{subsec:real_bar_construction}

We know that given a monoidal $\infty$-category $\caC$, there is a canonical functor $B\colon \Alg(\caC) \to \Fun(\Delta^\op,\caC)$ with source the $\infty$-category of associative algebras in $\caC$, the bar construction. This functor sends an associative algebra $A$ in $\caC$ to a simplicial object $B(A)$ in $\caC$ with $B(A)_n \simeq A^{\otimes n}$. It is a well known fact that when $\caC$ carries the cartesian structure, the functor $B$ is fully faithful with essential image the monoid objects in $\caC$. When $\caC$ carries the cocartesian structure, the forgetful functor $\Alg(\caC) \to \caC$ is an equivalence, and the Bar-construction becomes a functor $\caC \to \Fun(\Delta^\op,\caC)$ sending an object $X$ in $\caC$ to the simplicial object $B(X)$ on $\caC$ with $B(X)_n \simeq X^{\coprod n}.$

With the latter in mind, we define a real version of the bar construction
$B\colon\caC \to \Fun(\Delta^\op,\caC),$ where $\caC$ is a real $\infty$-category and
simplicial objects are replaced by real simplicial objects; see \cref{nth_level_real_bar_construction}.

In this subsection $\caC$ will be a pointed real $\infty$-category that admits tensors with finite pointed $C_2$-sets as $\Spc_\ast^{C_2}$-enriched $\infty$-category. Therefore, if we denote by $\Fin_\ast[C_2]\subset\Spc_\ast^{C_2}$ the full subcategory spanned by finite pointed $C_2$-sets, for every $X\in\caC$ we obtain a real functor $\Fin_\ast[C_2] \xrightarrow{-\wedge X} \caC$ by taking tensors. 

There is still one missing component to define the real bar construction, which we introduce next.
\begin{notation}
Denote by $\mathfrak{S}$ the quotient $\Delta^1 /{\partial\Delta^1}$ in real simplicial spaces, with its canonical basepoint. Note that $\mathfrak{S}$ is in each level a pointed finite $C_2$-set.
\end{notation}

\begin{definition}Let $\caC$ be as in the introduction. Then for every $X \in \caC$ we define the real Bar-construction $B(X)$ as the real simplicial object on $\caC$ given by the composition of real functors
$$ \underline{\Delta}^\op \xrightarrow{ \mathfrak{S} } \Fin_\ast[C_2] \xrightarrow{-\wedge X} \caC.$$
\end{definition}

\begin{remark}
For every $X, Z \in \caC $ we have a canonical equivalence of genuine $C_2$-spaces
$$ \map_{\caC}(|B(X)|, Z) \simeq \map_{\mathrm{rs}\caC}(B(X), \delta(Z)) \simeq \map_{\rsSpc_\ast}(\mathfrak{S}, \delta(\map_\caC(X,Z))) \simeq $$$$ \map_{\Spc^{C_2}_\ast}(S^{1,1}, \map_\caC(X,Z)) $$
that exhibits $|B(X)| \in \caC$ as the tensor $\Sigma^{1,1}(X)= S^{1,1} \wedge X $ of $S^{1,1}$ and $X.$ Then we have seen that $\caC$ admits all tensors with $S^{1,1}$ as $\Spc^{C_2}_\ast$-enriched $\infty$-category and the real functor $\Sigma^{1,1}\colon  \caC \to \caC$ factors as $$ \caC \xrightarrow{B} \mathrm{rs}\caC \xrightarrow{ \mid-\mid} \caC.$$
\end{remark}

As we insinuated before, we are interested in describing the nth-level of the real bar construction.

\begin{remark}\label{nth_level_real_bar_construction}
Let us start by saying that, given an object $X$ of a real $\infty$-category $\caC$, we write
$\widetilde{X \coprod X}$ for the tensor of $X$ with $C_2$; similarly to what we did for the cotensor \cref{not:tilde}.

Note that for every $n \geq 0 $ the set $\mathfrak{S}_n$ has $n$-elements
and the $C_2$-set $\mathfrak{S}_n$ has the following structure:
If $n$ is even, it is $\mathfrak{S}_n \simeq \frac{n}{2} \times C_2$, this the coproduct of $\frac{n}{2}$ copies of $C_2$; if $n$ is odd, it is $\mathfrak{S}_n \simeq (\frac{n-1}{2} \times C_2) \coprod \ast$. Then, we conclude that if $n$ is even

$$B(X)_n \simeq \underbrace{\widetilde{X \coprod X}\coprod\cdots \coprod \widetilde{X \coprod X},}_{\frac{n}{2}-\text{copies}}$$
and if $n$ is odd, 

$$B(X)_n \simeq \underbrace{\widetilde{X \coprod X}\coprod\cdots \coprod \widetilde{X \coprod X},}_{\frac{n-1}{2}-\text{copies}}\coprod X.$$
\end{remark}

We now introduce a real version of the object of monoid, which will be the structure that the real bar construction carries.

\begin{definition}
A real Segal object $X$ is said to be a real monoid if $X_0$ is the final genuine $C_2$-space.
\end{definition}

\begin{lemma}Let $\caC$ be a genuine preadditive real $\infty$-category, and $X$ and object in it. Then the bar construction $B(X)$ is a real monoid. 
\end{lemma}

\begin{proof}
It follows directly from the description of $B(X)_n$ in \cref{nth_level_real_bar_construction}.
\end{proof}

We denote by $\mathrm{rs}_\dag\caC \subset \mathrm{rs}\caC$ the full subcategory spanned by the real simplicial objects $Y$ of $\caC$ with $Y_0$ the zero object.

\begin{lemma}\label{lem:adj_B_ev1}
Let $\caC$ be a real $\infty$-category. Then the following is a real adjunction

\[
\begin{tikzpicture}
\node[](A) {$\caC$};
\node[right of=A,xshift=1.5cm](B) {$\mathrm{rs}_\dag\caC$};
\draw[->] ($(A.east)+(0,.25cm)$) to [bend left=30] node[above]{$B$} ($(B.west)+(0,.25cm)$);
\draw[->] ($(B.west)+(0,-.25cm)$) to [bend left=30] node[below]{$\ev_1$} ($(A.east)-(0,.25cm)$);
\node[] at ($(A.east)!0.5!(B.west)$) {\rotatebox{90}{$\vdash$}};
\end{tikzpicture}
\]

whose unit $\tau\colon \id_\caC \to \ev_1 \circ B $ is an equivalence.
\end{lemma}

\begin{proof}
Since $\mathfrak{S} $ is the quotient in $ \rsSpc$ of $ \Delta^1$ with $\partial\Delta^1 \simeq C_2 \otimes \Delta^0, $ for every real simplicial space
$Y$ we have a canonical equivalence of genuine $C_2$-spaces
$$ \map_{\rsSpc}(\mathfrak{S}, Y) \simeq \map_{\rsSpc}(\Delta^0, Y) \times_{\map_{\rsSpc}( C_2 \otimes \Delta^0 , Y)  } \map_{\rsSpc}(\Delta^1, Y) $$$$ \simeq Y_0 \times_{ \widetilde{Y_0 \times Y_0} }  Y_1. $$

Especially, if $Y_0$ is the final genuine $C_2$-space, we obtain the equivalence of genuine $C_2$-spaces $\map_{\rsSpc}(\mathfrak{S}, Y) \simeq Y_1$. Then for every $X \in \caC, Z \in \mathrm{rs}_\dag\caC$, there is a canonical equivalence of genuine $C_2$-spaces
$$ \map_{\mathrm{rs}\caC}(B(X), Z) \simeq \map_{\rsSpc}(\mathfrak{S}, \map_\caC(X,-) \circ Z) \simeq \map_\caC(X,Z_1).$$

For $Z=B(X)$ the identity corresponds to the unit component $\tau_X\colon X \to B(X)_1,$ which is the identity. Indeed, the map $\mathfrak{S} \to \map_\caC(X,-) \circ \mathfrak{S} \wedge X $
induces on the first level the map
$$\ast \simeq \mathfrak{S}_1 \to \map_\caC(X, \mathfrak{S}_1 \wedge X) \simeq \caC(X,X) $$
selecting the identity.
\end{proof}

\begin{corollary}Let $\caC$ be a genuine preadditive real $\infty$-category. Then the real adjunction of \cref{lem:adj_B_ev1} restricts to a real adjoint equivalence 
\[
\begin{tikzpicture}
\node[](A) {$\caC$};
\node[right of=A,xshift=1.8cm](B) {$\mathsf{r}\Mon(\caC).$};
\draw[->] ($(A.east)+(0,.25cm)$) to [bend left=30] node[above]{$B$} ($(B.west)+(0,.25cm)$);
\draw[->] ($(B.west)+(0,-.25cm)$) to [bend left=30] node[below]{$\ev_1$} ($(A.east)-(0,.25cm)$);
\node[] at ($(A.east)!0.5!(B.west)$) {\rotatebox{90}{$\vdash$}};
\end{tikzpicture}
\]
\end{corollary}

\begin{remark}In particular, from the result above we have that the real monoid structure of $B(X)$ is completely determined by its first term in the same way as a monoid in a preadditive $\infty$-category is determined by its first term.
\end{remark}

\subsection{Additivization}\label{subsec:additivization}
We undertake here the task of constructing, from a preadditive theory $\phi\colon \D\Wald_\infty^\gd \to \caD $, an additive theory 
$$\add(\phi)\colon\D\Wald_\infty^\gd \to \caD,$$ 
that we call its additivization; see \cref{def:additivization}. In \cref{fhjkplk} we provide a way to identify certain genuine excisive funcotrs, that allows us to deduce that the functor $\add(\phi)\colon\D\Wald_\infty^\gd \to \caD$ is indeed an additive theory; see \cref{additivization_is_additive}. Later, in \cref{subsec:universal_add_th}, we show that such construction is universal.

The subsection culminates with \cref{fghkkkkl}, that states that the theory $\add(\phi)$ preserves geometric realizations if $\phi$ does; this property will result crucial in the next section.

\begin{definition}\label{def:additivization}
For every preadditive theory $\phi\colon \D\Wald_\infty^\gd \to \caD $ we define its additivization as the functor $\add(\phi)\colon\D\Wald_\infty^\gd\to\caD$ given by the composite
$$ \add(\phi)\coloneqq \Omega^{1,1} \circ \phi \circ \caS^{1,1}\colon \D\Wald_\infty^\gd \to \caD.$$
\end{definition}

Note that this gives rise to an endofunctor of $\Fun(\D\Wald_\infty^\gd,\caD)$, that can be expressed as $$\add\coloneqq(\caS^{1,1})^\ast \circ (\Omega^{1,1})_\ast \simeq (\Omega^{1,1})_\ast \circ (\caS^{1,1})^\ast.$$

With this, we finally define the real $K$-theory space of a Waldhausen $\infty$-category with genuine duality as follows. 

\begin{definition}\label{def:KR_space}
The real $K$-theory functor
$$\KR\colon\D\Wald_\infty^\gd\to \Spc^{C_2}_\ast$$
is defined as the additivization $\KR\coloneqq\add(\iota)$ of the functor $\iota\colon\D\Wald_\infty^\gd\to\Spc^{C_2}_\ast$, which sends an $\infty$-category with genuine duality $\caC$ to its maximal subspace in $\caC$ equipped with the restricted genuine duality.
\end{definition}

\begin{definition}Let $\caD$ be a real $\infty$-category, and $Z$ an object in it. We say that $Z$ is projective if the real functor $\map_\caD(Z,-)\colon \caD \to \Spc^{C_2}$ preserves geometric realizations.
\end{definition}

\begin{definition}
We call a real $\infty$-category $\caD$ projectively generated if
there is a set $T$ of projective objects of $\caD$ that detect equivalences, i.e.\ a morphism $f$ in $\caD$ is an equivalence if and only if for every $Z \in T$
the map of genuine $C_2$-spaces $\map_\caD(Z,f)$ is an equivalence.
\end{definition}

\begin{lemma}\label{hfh6rere}
Let $\caD$ be a projectively generated real $\infty$-category that admits cotensors with $C_2$ and geometric realizations, and let $X$ be a real Segal object in $\caD$ whose underlying Segal object in $\caD^u$ is a groupoid object. Let us consider the commutative square 
\begin{equation}\label{fdghjkl}
\begin{tikzcd}
X_1 \ar[r]\ar[d]       &\vert X\vert\ar[d]\\
\widetilde{X_0 \times X_0} \ar[r]       &\widetilde{\vert X\vert \times \vert X\vert},
\end{tikzcd}
\end{equation}

where the top map is induced by the map $X \to |X|$ of real Segal objects in $\caD$, where $|X|$ represents the constant real Segal object on the realization $|X|$ and the right vertical map is the diagonal. Then this square is a pullback square in $\caD$.
\end{lemma}

Specializing to the case that 
the zero-th term of $X$ is contractible, we obtain the following result.

\begin{corollary}
Let $\caD$ be a projectively generated real $\infty$-category that admits cotensors with $C_2$ and geometric realizations, and let $X$ be a real monoid in $\caD$, whose underlying monoid in $\caD^u$ is grouplike. Then there exists a canonical equivalence in $\caD$ as below.
$$X_1 \simeq \Omega^{1,1} |X| $$ 
\end{corollary}

\begin{proof}
By \cref{lem:fiber_sequence_with_omega11}, there is a canonical fiber sequence
$$ \Omega^{1,1} |X| \to |X| \to \widetilde{|X| \times |X|}$$
of pointed genuine $C_2$-spaces
induced by the canonical cofiber sequence
$$C_2 \coprod [0] \to [0] \coprod [0] \to S^{1,1} $$
of pointed genuine $C_2$-spaces.
\end{proof}

\begin{proof}[Proof of \cref{hfh6rere}]
Let us show first that we can restrict ourselves to the case $\caD=\Spc^{C_2}$. Given $T$ a set of projective generators of $\caD$, the square (\ref{fdghjkl}) is a pullback square in $\caD$ if and only if it is sent to a pullback square in $\Spc^{C_2}$ by the real functor $\map_\caD(Z,-)\colon \caD \to \Spc^{C_2}$ for every $Z \in T$, for  $\map_\caD(Z,-)$ preserves small limits. Moreover, since for every $Z\in T$, the functor $\map_\caD(Z,-)$ preserves cotensors with $C_2$ \textemdash and therefore real Segal objects in particular \textemdash and geometric realizations, we can replace $X$ by the real Segal space $\map_\caD(Z,-) \circ X\colon \underline{\Delta}^\op \to \caD \to \Spc^{C_2}$, whose underlying Segal space is, too, a groupoid.

We can then exclusively treat the case that $\caD = \Spc^{C_2}$. Forgetting the $C_2$-actions, we get from the commutative square (\ref{fdghjkl}), a commutative square of spaces as below.
\[
\begin{tikzcd}
X_1^u \ar[r] \ar[d]         &\vert X^u\vert \ar[d] \\
X_0^u \times X_0^u \ar[r]   & \vert X^u\vert \times \vert X^u\vert 
\end{tikzcd}
\]
That this is a pullback square follows from the fact that the square below is so by \cite[Theorem 6.1.0.6]{lurie.HTT}.

\[
\begin{tikzcd}
X_1^u \ar[r] \ar[d]         &X_0^u \ar[d] \\
X_0^u \ar[r]   & \vert X^u\vert 
\end{tikzcd}
\]

Taking now $C_2$-fixed points, from the square (\ref{fdghjkl}) we obtain a commutative square of spaces
\[
\begin{tikzcd}
X_1^{C_2} \ar[r] \ar[d]     &\vert X^{C_2}\vert \ar[d] \\
X_0^u \ar[r]                & \vert X^u\vert,
\end{tikzcd}
\]
where $X^{C_2}$ denotes the simplicial space $$\Delta^\op \xrightarrow{e} (\Delta^{hC_2})^\op \to \Spc^{C_2} \xrightarrow{(-)^{C_2}} \Spc$$ induced by $X$. We will now proceed to prove that this is also a pullback square.

Denote by $(\caC, H \to \H(\caC))$ the Segal space with genuine duality associated to $X$. Then the last square is equivalent to the following square
\begin{equation}\label{dfghjkl}
\begin{tikzcd}
H_0 \ar[r] \ar[d]       &\vert H\vert \ar[d] \\
\caC_0 \ar[r]           &\vert\caC\vert 
\end{tikzcd}
\end{equation}
induced by the map of simplicial spaces $ \psi\colon H \to \caH(\caC) \to \caC $. Then it will be enough to prove that square (\ref{dfghjkl}) is a pullback square.

The map $\psi\colon H \to \caH(\caC) \to \caC$ of simplicial spaces is a right fibration of simplicial spaces and so belongs to $\LMod_{\caC}(\Spc)$, see the discussion at the beginning of \cref{subsec:left_actions}.

We have a commutative square 
\begin{equation}\label{fghjkjppp}
\begin{tikzcd}
H \ar[r] \ar[d, "\psi"']     & \delta(|H|) \times \caC \ar[d] \\
\caC \ar[r]                 & \delta(|\caC|) \times \caC 
\end{tikzcd}
\end{equation}
in $ \LMod_{\caC}(\Spc) \subset \caP(\Delta)_{/\caC}$
that yields after evaluation at $0 \in \Delta$ the commutative square
\begin{equation}\label{fjkhkjjp}
\begin{tikzcd}
H_0 \ar[r] \ar[d, "\psi_0"']      &\vert H\vert \times \caC_0 \ar[d] \\
\caC_0 \ar[r]                   &\vert\caC\vert \times \caC_0.
\end{tikzcd}
\end{equation}

By pasting of squares, (\ref{dfghjkl}) is a pullback square if and only if
square (\ref{fjkhkjjp}) is so. Since $ \LMod_{\caC}(\Spc)$ is closed under pullbacks in $\caP(\Delta)_{/\caC}$, it is enough to check that square (\ref{fghjkjppp}) is a pullback square in $\LMod_{\caC}(\Spc)$. Now, under the equivalence $\theta$ of \cref{prop:left_action_equiv_Segal_spaces} square (\ref{fghjkjppp}) corresponds to the pullback square 
\[
\begin{tikzcd}
\theta(H) \ar[r] \ar[d]             & \theta(H) \times \vert\caC\vert \ar[d] \\
\vert\caC\vert \ar[r] & \vert\caC\vert \times \vert\caC\vert.
\end{tikzcd}
\]
\end{proof}

We present now a theorem that allows us to identify certain genuine excisive functors.
\begin{proposition}\label{fhjkplk}
Let $\caC $ be a genuine preadditive real $\infty$-category that admits geometric realizations, and let $\caD$ be a pointed projectively generated real $\infty$-category that admits geometric realizations, finite products and cotensors with $C_2$, and also admits cotensors with $S^{1,1}$ as $\Spc^{C_2}_\ast$-enriched $\infty$-category. Let $F\colon\caC \to \caD$ be a real functor that preserves geometric realizations, finite products and cotensors with $C_2$.

Then, if the canonical lift of $F$ to commutative monoids (in the sense of \cref{Venr:def:monoids}) in $\caD$ takes values in those commutative monoids whose underlying commutative monoid in $\caD^u$ is grouplike, $F$ is genuine excisive.
\end{proposition}

\begin{proof}
We want to prove that for every $Y \in \caC$ the canonical morphism
$$F(Y) \to \Omega^{1,1} F \Sigma^{1,1}(Y) \simeq \Omega^{1,1} F (|B(Y)|) \simeq \Omega^{1,1} |F \circ B(Y)| $$ 
is an equivalence. Set $X= F \circ B(Y) \in \mathrm{rMon}(\caD)$ so that $X_1= F(Y)$ and $X_0 \simeq 0. $

So the claim follows from \cref{hfh6rere} applied to $X, $
whose underlying monoid in $\caD^u$ is grouplike, since the following square commutes:

\[
\begin{tikzcd}
\mathsf{rMon}(\Cmon(\caD)) \ar[r] \ar[d]        & \mathsf{Mon}(\Cmon(\caD^u))  \ar[r, "\simeq"] \ar[d]            & \Cmon(\caD^u) \ar[d] \\
\mathsf{rMon}(\caD) \ar[r]      & \mathsf{Mon}(\caD^u) \ar[r, "="]         & \mathsf{Mon}(\caD^u).
\end{tikzcd}
\]

The right hand square commutes as both diagonals of the square are equivalent when composed with the forgetful functor $\mathsf{Mon}(\caD^u)\to \caD^u$
and $\Cmon(\caD^u)$ is preadditive.
\end{proof}

\begin{corollary}\label{dfghjlk}
Let $\caD$ be a pointed projectively generated real $\infty$-category that admits geometric realizations, finite products and cotensors with $C_2$, and also admits cotensors with $S^{1,1}$ as $\Spc^{C_2}_\ast$-enriched $\infty$-category. Then for every semiadditive theory $\varphi\colon\D\Wald_\infty^\gd \to \caD$, the theory $ \Omega^{1,1} \circ \varphi\colon \D\Wald_\infty^\gd \to \caD \to \caD $ is additive.
\end{corollary}

\begin{proof}
By \cref{fhjkplk}, it is enough to check that the functor $ \D\Wald_\infty^\gd \xrightarrow{\Omega^{1,1} \circ \varphi} \caD \to \caD^u $ lifts to grouplike commutative monoids of $\caD^u$. This follows from the fact that such functor factors as $$\D\Wald_\infty^\gd \xrightarrow{\varphi} \caD \to \caD^u \xrightarrow{\Omega} \caD^u $$ and $\Omega$ lifts to grouplike commutative monoids of $\caD^u$.
\end{proof}

\begin{remark}\label{additivization_is_additive}
Given a theory $\phi\colon\D\Wald_\infty^\gd \to \caD$, applying \cref{dfghjlk} to $\varphi=\phi\circ \caS^{1,1}$, we conclude that $\add(\phi)$ is an additive theory.
\end{remark}

We conclude the section proving that the additivization of a preadditive theory preserves both geometric realizations and filtered colimits, fact that we need in \cref{fhjkkkl}.

\begin{lemma}\label{fghkkkkl}
Let $\caD$ be a pointed projectively generated real $\infty$-category
that admits cotensors with $S^{1,1}$ and $\phi\colon\D\Wald_\infty^\gd \to \caD $ a preadditive theory that preserves geometric realizations (resp.\ filtered colimits). Then the additization of $\phi$, 
$$\add(\phi)\colon \D\Wald_\infty^\gd \to \caD $$ preserves geometric realizations (resp.\ filtered colimits) as well.
\end{lemma}

\begin{proof}
As $\caD$ is projectively generated, we can reduce to the case that
$\caD= \Spc^{C_2}_\ast$.

Set $\psi\coloneqq \phi \circ \caS^{1,1}\colon \D\Wald_\infty^\gd \to \Spc^{C_2}_\ast$, for a reminder of who is $\caS^{1,1}$ see \cref{const:curly_S}. We want to check that both compositions 
\begin{equation*}
\alpha\colon  \D\Wald_\infty^\gd \xrightarrow{\Omega^{1,1} \circ \psi } \Spc^{C_2}_\ast \xrightarrow{(-)^u } \Spc_\ast \hspace{0.8em}     \text{and}\hspace{0.8em} \beta\colon  \D\Wald_\infty^\gd \xrightarrow{\Omega^{1,1} \circ \psi } \Spc^{C_2}_\ast \xrightarrow{(-)^{C_2} } \Spc_\ast
\end{equation*}
preserve geometric realizations (resp. sifted colimits).

The functor $\alpha$ factors as 
$$\D\Wald_\infty^\gd \xrightarrow{\psi } \Spc^{C_2}_\ast \xrightarrow{(-)^u}\Spc_\ast \xrightarrow{\Omega } \Spc_\ast,$$ while the functor $\beta$ factors as 
$$  \D\Wald_\infty^\gd \xrightarrow{\psi } \Spc^{C_2}_\ast \xrightarrow{\rho} \Fun([1], \Spc_\ast) \xrightarrow{\fib} \Spc_\ast,$$
where $\rho$ sends a pointed genuine $C_2$-space $X$ to the map of pointed spaces $X^{C_2} \to X^u$ and $\fib$ takes the fiber (see \cref{lem:description_fixed_points_Omega11}).
This shows that $\add(\phi)$ preserves filtered colimits if $\phi$ does.

By using the following implication of the $\pi_\ast$-kan conditions, it is enough to show that the functor $(-)^u \circ \psi\colon  \D\Wald_\infty^\gd \to \Spc_\ast$ takes values in connected spaces to obtain that both $\alpha$ and $\beta$ preserve geometric realizations. The implication we refer to says: given a fiber sequence of pointed simplicial spaces 
\begin{equation*}\label{fjkhkjp}
\begin{tikzcd}
A \ar[r] \ar[d]         & B \ar[d] \\
\ast \ar[r]             & C
\end{tikzcd}
\end{equation*}
such that $C$ is levelwise connected, the induced commutative square on geometric realizations is a fiber sequence of pointed spaces.

We conclude the proof by showing that $(-)^u \circ \psi\colon  \D\Wald_\infty^\gd \to \Spc_\ast$ takes values in connected spaces. For this, we show that for every $C \in \D\Wald_\infty^\gd$ the pointed genuine $C_2$-space $\phi(\caS^{1,1}(C))$ has an underlying connected space.

We begin by observing that there is a canonical equivalence $\phi(\caS^{1,1}(C)) \simeq |\phi \circ S(C)|$
as $\phi\colon \D\Wald_\infty^\gd \to \caD$ preserves geometric realizations. If $X$ denotes the underlying simplicial space of the real simplicial space $\phi \circ S(C)$, then $|X|$ is the underlying space of the genuine $C_2$-space $|\phi \circ S(C)|.$  Now, as $X_0 =\phi(S(C)_0)^u \simeq \ast$, the simplicial set $\pi_0 \circ X$ is the singleton at level 0. We conclude then that $\pi_0(|X|) = |\pi_0 \circ X| = \ast$.
\end{proof}

\subsection{Universality of additivization}\label{subsec:universal_add_th}

Our aim now is to show that the additivization of a preadditive theory $\phi$, which we call $\add(\phi)$, is the initial additive theory among those additive theories $F$ equipped with a map of theories $\lambda\colon \phi \to F$, this is \cref{univ}.

We have decided to present the proof of \cref{univ} first, and later the results the proof depends on. The main of these results is \cref{dfghjhgk}, where, among other things, we prove that a theory is semiadditive if and only if it identifies the real bar construction and the real $S$-construction. 

Let us proceed to present the main ingredients for the statement of the universal property.

\begin{construction}\label{theta_B_to_S}
For $\caC= \D\Wald_\infty^\gd$ composing the counit $B \circ \ev_1 \to \id$ 
with the $S$-construction we get a map
$$\theta\colon B \simeq B \circ \ev_1 \circ S \to S $$ in $\Fun_{\Spc_\ast^{C_2}}(\D\Wald_\infty^\gd, \mathrm{rs}_\ast\D\Wald_\infty^\gd).$
\end{construction}

\begin{remark}\label{def:xi}
Taking realizations, from the map $\theta\colon B\to S$ we obtain a map 
$$ \xi\colon \Sigma^{1,1} \simeq |-| \circ B \xrightarrow{|-| \circ \theta} |-| \circ S \simeq \caS^{1,1}$$ in $\Fun_{\Spc_\ast^{C_2}}(\D\Wald_\infty^\gd, \D\Wald_\infty^\gd)$ adjoint to a map $\zeta\colon \id \to \Omega^{1,1} \circ \caS^{1,1}$.
\end{remark}

Given a pointed real $\infty$-category $\caD$ that admits cotensors with $S^{1,1}$, there is a canonical real natural transformation 
$ (\Omega^{1,1})^\ast \to (\Omega^{1,1})_\ast $ of real endofunctors of $\Fun_{\Spc_\ast^{C_2}}(\D\Wald_\infty^\gd, \caD).$ This gives rise to a real natural transformation 
$$ \lambda\colon \id \xrightarrow{\zeta^\ast} (\caS^{1,1})^\ast \circ (\Omega^{1,1})^\ast \to \add\coloneqq (\caS^{1,1})^\ast \circ (\Omega^{1,1})_\ast \simeq (\Omega^{1,1})_\ast \circ (\caS^{1,1})^\ast$$ of real endofunctors of $\Fun_{\Spc_\ast^{C_2}}(\D\Wald_\infty^\gd, \caD).$

Denote by $$ \Add_\caD \subset \Pre\Add_\caD \subset \Fun_{ \Spc^{C_2}_\ast }(\D\Wald_\infty^\gd, \caD)$$ the full subcategories spanned by the additive respectively preadditive theories that preserve geometric realizations.

\begin{theorem}[Universal property]\label{univ}
Let $\caD$ be a pointed projectively generated real $\infty$-category
that admits geometric realizations and cotensors with $S^{1,1}$ and $C_2$.
The real natural transformation
$$ \lambda\colon \id \to \add$$ of real endofunctors of $\Fun_{\Spc_\ast^{C_2}}(\D\Wald_\infty^\gd, \caD)$ exhibits
$\add\colon \Pre\Add_\caD \to \Add_\caD $ as a real left adjoint of the embedding $ \Add_\caD \subset \Pre\Add_\caD.$
\end{theorem}
	
\begin{proof}
For readability sake, set $L\coloneqq\add$. We will check that for every additive theory $\phi$ the map $\lambda\colon \phi \to L(\phi)$ is an equivalence, which implies that $L$ is an essentially surjective functor. Then by \cref{Venr:left_adjs_of_embeddings} it is enough to show that $\lambda \circ L\colon L \to L \circ L $ and $L \circ \lambda \colon L \to L \circ L $ are equivalences.

That $\lambda \circ L\colon L \to L \circ L$ is an equivalence follows from the fact that for every additive theory $\phi$ the map $\lambda\colon \phi \to L(\phi)$ is an equivalence, since $L(\phi)$ is additive by \cref{dfghjlk}. 
The map $\lambda \colon\phi \to L(\phi) $ factors as $$ \phi \to \Omega^{1,1} \circ \phi \circ \Sigma^{1,1}  \to \Omega^{1,1} \circ \phi \circ \caS^{1,1} $$
since there is a commutative square

\[
\begin{tikzcd}
\phi \circ \Omega^{1,1} \circ \Sigma^{1,1} \ar[r] \ar[d]    &  \Omega^{1,1} \circ \phi \circ \Sigma^{1,1} \ar[d] \\
\phi \circ \Omega^{1,1} \circ \caS^{1,1} \ar[r]      &  \Omega^{1,1} \circ \phi \circ \caS^{1,1}.
\end{tikzcd}
\]

Since $\phi$ is additive, the first map $ \phi \to \Omega^{1,1} \circ \phi \circ \Sigma^{1,1}$ is an equivalence. By Item (3) of \cref{dfghjhgk}, the map
$$\Omega^{1,1} \circ \phi \circ \Sigma^{1,1}  \to \Omega^{1,1} \circ \phi \circ \caS^{1,1}$$ is an equivalence. 

Then it remains to check that $ (\lambda \circ L)(\phi), (L \circ \lambda)(\phi) $ are both equivalent, in other words that the bottom and top horizontal morphims of the following composed commutative square are equivalent:

\[
\begin{tikzcd}[column sep=large]
\Omega^{1,1} \circ \phi \circ \caS^{1,1} \ar[r, "\Omega^{1,1} \circ \phi \circ \caS^{1,1} \circ \zeta"]\ar[d, "="]      &  \Omega^{1,1} \circ \phi \circ \caS^{1,1} \circ \Omega^{1,1} \circ \caS^{1,1} \ar[r] \ar[d, "\gamma"]       & \Omega^{1,1} \circ  \Omega^{1,1} \circ \phi \circ \caS^{1,1} \circ \caS^{1,1} \ar[d, "="] \\
\Omega^{1,1} \circ \phi \circ \caS^{1,1} \ar[r, "\Omega^{1,1} \circ \phi \circ \zeta  \circ \caS^{1,1}"]            & \Omega^{1,1} \circ \phi \circ \Omega^{1,1} \circ \caS^{1,1}  \circ \caS^{1,1} \ar[r]       & \Omega^{1,1} \circ  \Omega^{1,1} \circ \phi \circ \caS^{1,1} \circ \caS^{1,1}.
\end{tikzcd}
\]

The right hand square commutes since, for reduced real functors $F\colon \caC \to \caD, G\colon \caD \to \caE$ between pointed real $\infty$-categories $\caC,\caD,\caE$ that admit cotensors with $S^{1,1}$ in the pointed sense, the canonical map $GF \Omega^{1,1} \to \Omega^{1,1} GF$ factors as 
$$GF \Omega^{1,1} \to G \Omega^{1,1} F \to \Omega^{1,1} GF$$ and for $F= \Omega^{1,1} \colon \caC \to \caC$ the canonical map $F \Omega^{1,1} \to \Omega^{1,1} F$ is the identity.

To see that the left hand square commutes, it is enough to verify that
the square

\[
\begin{tikzcd}[column sep=large]
\caS^{1,1} \ar[r, "\caS^{1,1} \circ \zeta"] \ar[d, "="]      & \caS^{1,1} \circ \Omega^{1,1} \circ \caS^{1,1} \ar[d] \\
\caS^{1,1} \ar[r, "\zeta \circ \caS^{1,1}"]                 & \Omega^{1,1} \circ \caS^{1,1}  \circ \caS^{1,1}
\end{tikzcd}
\]
commutes. But this square is adjoint to the commutative square 

\[
\begin{tikzcd}
\Sigma^{1,1} \circ \caS^{1,1} \ar[r] \ar[d, "="] & \caS^{1,1} \circ \Sigma^{1,1} \ar[d, "\caS^{1,1} \circ \xi"]\\
\Sigma^{1,1} \circ \caS^{1,1} \ar[r, "\xi \circ \caS^{1,1}"']        & \caS^{1,1}  \circ \caS^{1,1},
\end{tikzcd}
\]
where $\xi$ is as in \cref{def:xi}.
\end{proof}

\begin{remark}
From the universal property of the additivization real functor (see \cref{univ}), we have that for every preadditive theory $\phi$ and additive theory $\psi$ the canonical map $$\lambda\colon \phi \xrightarrow{\phi \circ \zeta } \phi \circ \Omega^{1,1} \circ \caS^{1,1} \to \add(\phi)= \Omega^{1,1} \circ \phi \circ \caS^{1,1} $$
induces an equivalence $$\map_{\Add_\caD}(\add(\phi), \psi) \to \map_{\Pre\Add_\caD}(\phi, \psi)$$ on $\Spc^{C_2}$-enriched mapping spaces.
\end{remark}

\begin{remark}\label{univ_KR}
Since we have defined the real $K$-theory functor as additivization of the preadditive theory $\iota\colon\D\Wald_\infty^\gd\to\Spc_\ast^{C_2}$, \cref{univ} implies a universal property of real $K$-theory.
\end{remark}

\begin{lemma}\label{dfghjkkkk}
Let $\caD$ be in $\D\Wald^\gd_\infty$ and $\phi\colon \D\Wald^\gd_\infty \to \caD$ a real functor that inverts $\gamma_3\colon  S(\caD)_3 \to \widetilde{\caD \times \caD} \times \caD.$ Let $\alpha\colon \widetilde{\caD \times \caD} \to \caD$ be the map in $\Wald^\gd_\infty$ 
corresponding to the identity of $\caD$
under the equivalence $\Wald^\gd_\infty(\widetilde{\caD \times \caD} ,\caD) \simeq \Exact_\infty(\caD,\caD)$. Then the image under $\phi$ of the extension of the composition 
$$S(\caD)_3 \xrightarrow{ \gamma_3 } \widetilde{\caD \times \caD} \times \caD \xrightarrow{\alpha \times \caD} \caD \times \caD \xrightarrow{\oplus} \caD $$
to a map in $ \D\Wald^\gd_\infty$, which we call $\xi$, and the map $ S(\caD)_3 \to \caD$ induced by the map $[1] \simeq \{0 \to 3\} \subset [3]$ in $\Delta^{hC_2}$, is the same.
\end{lemma}

\begin{proof}
We define the following maps in $ \Wald^\gd_\infty$ that extend to similarly defined maps in $\D\Wald_\infty^\gd$ by \cref{lem:functors_extend_to_DWald}:
\begin{itemize}
\item $\beta\colon S(\caD)_3 \to \caD \simeq S(\caD)_1, \ (X \to Y \to Z) \mapsto Z$ the map induced by the map $[1] \to [3]$ in $\Delta$ that sends $0\mapsto 0$ and $1 \mapsto 3$.
\item $\rho\colon \widetilde{\caD \times \caD} \to S(\caD)_3, \ (X,Y) \mapsto (X \to X \to X \oplus Y)$ the map adjoint to the map $ \rho'\colon \caD \simeq S(\caD)_1 \to S(\caD)_3$ in $ \Exact_\infty$
induced by the map $[3] \to [1]$ in $\Delta$ that sends only 0 to 0.
\item $\theta\colon \caD \simeq S(\caD)_1 \to S(\caD)_3, \ X \mapsto (0 \to X \to X) $ the map induced by the map $[3] \to [1]$ in $\Delta$ that sends $0,1\mapsto 0$ and $2,3\mapsto 1$.
\end{itemize} 

The map $\rho + \theta\colon \widetilde{\caD \times \caD} \times \caD \to S(\caD)_3$ in $ \Wald^\gd_\infty$ that sends $\ (X,Y,Z) \mapsto (X \to X \oplus Z \to X \oplus Y \oplus Z) $ is a section of $\gamma_3\colon S(\caD)_3 \to \widetilde{\caD \times \caD} \times \caD $. As $\phi$ inverts $\gamma_3,$ the map $\phi( \rho + \theta)$ is inverse to $\phi(\gamma_3)$ so that $\phi$ identifies the maps 
$$\xi\colon  S(\caD)_3 \xrightarrow{} \caD, $$
$$ S(\caD)_3 \xrightarrow{ \gamma_3 } \widetilde{\caD \times \caD} \times \caD \xrightarrow{ \rho+ \theta}  S(\caD)_3 \xrightarrow{\xi} \caD. $$

The proof concludes by observing that the map
$ \widetilde{\caD \times \caD} \times \caD \xrightarrow{ \rho+ \theta}  S(\caD)_3 \xrightarrow{\xi} \caD $
in $ \Wald^\gd_\infty$ is equivalent to $$ \widetilde{\caD \times \caD} \times \caD \xrightarrow{\alpha \times \caD} \caD \times \caD \xrightarrow{\oplus} \caD. $$
\end{proof}

We are now set to prove characterizations of semiadditive, as stated in \cref{dfghjhgk}. We begin by fixing notation.

Let $\caC \in \Wald^\gd_\infty$ and $n \geq 0.$ If $n$ is even, denote by $\gamma_n$ the canonical map 

$$ S(\caC)_n \longrightarrow (\widetilde{\caC \times \caC})^{\frac{n}{2}}$$ 
$$(X_1 \to \dots \to X_n) \mapsto (X_1, X_n/X_{n-1}, X_2/X_1, X_{n-1}/X_{n-2},\dots, X_{\frac{n}{2} }/X_{\frac{n}{2}-1}, X_{\frac{n}{2}+1 }/X_{\frac{n}{2}} ) $$
in $ \Wald^\gd_\infty$, whose projection to the $i$-th factor for $1 \leq i \leq \frac{n}{2}$ is adjoint to the map $ S(\caC)_n \to S(\caC)_1 \simeq \caC $
in $ \Exact_\infty$ induced by the map $[1] \to [n]$ in $\Delta$ that sends
$0\mapsto i-1$ and $1\mapsto i$.

If $n$ is odd, denote by $\gamma_n$ the canonical map 
$$ S(\caC)_n  \longrightarrow  (\widetilde{\caC \times \caC})^{\frac{n-1}{2}} \times \caC$$
$$(X_1 \to \dots \to X_n) \mapsto $$$$(X_1, X_n/X_{n-1}, X_2/X_1, X_{n-1}/X_{n-2},\dots, X_{\frac{n-1}{2} }/X_{\frac{n-1}{2}-1}, X_{\frac{n+3}{2} }/X_{\frac{n+1}{2}}; X_{\frac{n+1}{2} }/X_{\frac{n-1}{2}} ), $$
whose projection to the $i$-th factor for $1 \leq i \leq \frac{n-1}{2}$ is adjoint to the map $ S(\caC)_n \to S(\caC)_1 \simeq \caC $
in $ \Exact_\infty$ induced by the map $[1] \to [n]$ in $\Delta$ that sends
$0\mapsto i-1$ and $1\mapsto i$, and whose projection to the last factor is the map 
$S(\caC)_n \to S(\caC)_1 \simeq \caC $ in $ \Wald^\gd_\infty$ induced by the map $[1] \to [n]$ in $\Delta^{hC_2}$ that sends $0\mapsto\frac{n-1}{2}$ and $1\mapsto \frac{n-1}{2}+1.$

The maps $\gamma_n$ extend to maps in $\D\Wald^\gd_\infty$ for $\caC \in \D\Wald^\gd_\infty.$

The following result characterizes semiadditive theories. Remarkably, it states that a semiadditive theory identifies the real $S$-construction with the real bar construction, and that it splits all levels of the real $S$-construction coherently.

\begin{proposition}\label{dfghjhgk}
Let $ \caD$ be a pointed real $\infty$-category and $\phi\colon \D\Wald^\gd_\infty \to \caD$ a preadditive theory. The following are equivaelnt:

\begin{enumerate}
\item The theory $\phi$ is semiadditive.

\item For every $\caC \in \D\Wald^\gd_\infty $ the real functor $\phi$ inverts $$\gamma_n\colon S(\caC)_n \to (\widetilde{\caC \times \caC})^{\frac{n}{2}}$$
for $n $ even, and $$ \gamma_n\colon S(\caC)_n  \to  (\widetilde{\caC \times \caC})^{\frac{n-1}{2}} \times \caC $$ for $n $ odd.		

\item The map $\phi \circ \theta\colon \phi \circ B \to \phi \circ S $ is an equivalence, where $\theta$ is as in \cref{theta_B_to_S}.

\item The real simplicial object $ \phi \circ S(\caC) $ is a real monoid in $\caD$.
\end{enumerate}
\end{proposition}

\begin{proof}	
The equivalence between (2) and (4) is easy to see. 

We work now on the equivalence between (2) and (3). For every $\caC \in \D\Wald_\infty^\gd$ we have a map $\theta\colon B(\caC) \to S(\caC)$
in $\mathrm{rs}_\dag\D\Wald_\infty^\gd$ (\cref{theta_B_to_S}), unique with the property that it yields the identity $\caC \simeq B(\caC)_1 \to S(\caC)_1 \simeq \caC $ in degree 1. Consider now the real functor given by the composite  $\phi \circ \theta\colon \phi \circ B(\caC) \to \phi \circ S(\caC) $ in $\mathrm{rs}_\ast\caD$. 

As $B(\caC)$ is a real monoid in $\D\Wald_\infty^\gd$, its image $\phi \circ  B(\caC)$ is a real monoid in $\caD$, since $\phi$ preserves finite products and cotensors with $C_2.$ As a map of real monoids in $\caD$ is an equivalence if it gets an equivalence after evaluation at $1 \in \underline{\Delta},$ the map $\phi \circ \theta\colon \phi \circ B(\caC) \to \phi \circ S(\caC) $ in
$\mathrm{rs}_\ast\caD$ is an equivalence if and only if $ \phi \circ S(\caC) $
is a real monoid in $\caD$, which is equivalent to (2) by \cref{cor:recognition_pple_rSeg}.

Item (1) follows trivially from (2). We show now that (1) implies (2). For $n \geq 0$ we define the following maps in $ \Wald^\gd_\infty$ that extend to similarly defined maps in $\D\Wald_\infty^\gd:$
\begin{itemize}
\item $\rho\colon S(\caC)_{n-2} \to S(\caC)_n, \ (X_1 \to \dots\to X_{n-2}) \mapsto (0 \to X_1 \to\dots\to X_{n-2} \xrightarrow{\id}  X_{n-2}) $ the map in $\Wald_\infty^\gd$ induced by the map $[n] \to [n-2]  $ in $\Delta^{hC_2}$ that identifies $0$ and $1$, identifies $n-1$ and $n$, and is injective when restricted to $ 1 < i < n-1$.

\item $\tau\colon \widetilde{\caC \times \caC} \to S(\caC)_n, \ (A,B) \mapsto (A \to \dots \to A \to A \oplus B) $ the map in $\Wald_\infty^\gd$ adjoint to the map $\delta\colon \caC \simeq S(\caC)_1 \to S(\caC)_n, \ A \mapsto (A \to \dots \to A)$ in $\Exact_\infty$ induced by the map $[n] \to [1] $ in $\Delta$ that sends only $0 $ to $0$.

\item $\sigma\colon S(\caC)_n \to S(\caC)_{n-2} \times \widetilde{\caC \times \caC}$
the map in $\Wald_\infty^\gd$, whose projection to the first factor is
the map $$S(\caC)_n \to S(\caC)_{n-2}, \ (X_1 \to\dots\to X_n) \mapsto (X_2/X_1 \to\dots\to X_{n-1}/X_1) $$ in $\Wald_\infty^\gd$
induced by the map $[n-2] \to [n]  $ in $\Delta^{hC_2}$ that sends $i $ to $i+1$ 
and whose projection to the second factor is the map $$ S(\caC)_n \to \widetilde{\caC \times \caC}, \ (X_1 \to \dots \to X_n) \mapsto (X_1, X_n/X_{n-1}) $$ in $\Wald_\infty^\gd$ adjoint to the map 
$$ \psi\colon S(\caC)_n \to \caC \simeq S(\caC)_1, \ (X_1 \to\dots \to X_n) \mapsto X_1 $$ in $\Exact_\infty$ induced by the map $[1] \to [n]  $ in $\Delta$ that sends $0 $ to $0$ and $1$ to $1$.

\item $ \xi\colon S(\caD)_3 \xrightarrow{ \gamma_3 } \widetilde{\caD \times \caD} \times \caD \xrightarrow{\alpha \times \caD} \caD \times \caD \xrightarrow{\oplus} \caD, $$$ (A \to B \to \caC) \mapsto A \oplus (C/B) \oplus (B/A) $$
the composition of maps in $\Wald^\gd_\infty$,
where $\caD\coloneqq S(\caC)_n$ and $\alpha\colon \widetilde{\caD \times \caD} \to \caD$ is the map in $\Wald^\gd_\infty$ adjoint to the identity of $\caD$ in $\Exact_\infty.$
\end{itemize}

The map $\beta\coloneqq \rho + \tau\colon S(\caC)_{n-2} \times \widetilde{\caC \times \caC} \to S(\caC)_n,$ given by  
$$ ((X_1 \to \dots \to X_{n-2}), (A,B)) \mapsto (A \to X_1 \oplus A \to \dots\to X_{n-2} \oplus A \to  X_{n-2} \oplus A \oplus B) $$ is a section of $\sigma\colon S(\caC)_n \to S(\caC)_{n-2} \times \widetilde{\caC \times \caC}.$

It also happens, as we show next, that $\phi(\sigma) $ is a section (and so an inverse) of $\phi(\beta)$.

The map $\beta \circ \sigma\colon S(\caC)_n \to S(\caC)_n$, given by sending 
$X\coloneqq(X_1 \to \dots \to X_n) \mapsto $$$ (X_1 \to (X_2/X_1) \oplus X_1 \to \dots \to (X_{n-1}/X_1) \oplus X_1 \to (X_{n-1}/X_1) \oplus X_1 \oplus (X_n/X_{n-1})) \simeq $$
$$ \delta(X_1) \oplus (0 \to X_2/X_1 \to \dots \to X_{n-1}/X_1 \to X_{n-1}/X_1) \oplus (0 \to \dots \to 0 \to X_n/X_{n-1})\simeq $$ 
$$ \delta(X_1) \oplus (\underbrace{(X_1 \to X_2 \to ... \to X_{n-1} \to X_{n-1} )}_{X'}/\delta(X_1)) \oplus (X / X') $$ 
lifts along the map $\xi$ via the map $$\kappa\colon \caD \to S(\caD)_3, \ X \mapsto ( \delta(X_1) \to X' \to X)$$ in $\Wald_\infty^\gd$.

The map $\kappa$ is a section of the map $$\ev_3\colon S(\caD)_3 \to \caD \simeq S(\caD)_1, \ (A \to B \to C) \mapsto C$$
induced by the map $[1] \to [3]$ in $\Delta^{hC_2}$ that sends $0\to 0$ and $1\to 3$. By \cref{dfghjkkkk}, for $\caD= S(\caC)_n$ the real functor $\phi$ identifies 
$ \xi $ with $\ev_3.$

So we obtain $$ \id \simeq \phi(\ev_3) \circ \phi(\kappa) \simeq \phi(\xi) \circ \phi(\kappa) \simeq \phi(\beta) \circ \phi(\sigma)$$
and therefore conclude that $\phi(\sigma) $ is a section (and so an inverse) of $\phi(\beta)$.

We are finally ready to prove (2). We do so by induction over $\frac{n}{2} \geq 0$ for $n$ even, respectively $\frac{n-1}{2} \geq 0$ for $n$ odd. For $n=0$ and $n=1$ there is nothing to show. 

Assume that (2) holds for $\frac{n}{2} \geq 0$ for $n$ even and $\frac{n-1}{2} \geq 0$ for $n$ odd. The canonical maps $$\gamma_{n+2}\colon S(\caC)_{n+2} \to (\widetilde{\caC \times \caC})^{\frac{n+2}{2}}$$
for $n $ even, and $$\gamma_{n+2}\colon S(\caC)_{n+2}  \to  (\widetilde{\caC \times \caC})^{\frac{n+1}{2}} \times \caC $$ for $n $ odd, factor canonically as
$$ S(\caC)_{n+2} \xrightarrow{\sigma} S(\caC)_{n} \times \widetilde{\caC \times \caC} \xrightarrow{\gamma_n \times \widetilde{\caC \times \caC} } (\widetilde{\caC \times \caC})^{\frac{n+2}{2}} $$
for $n $ even, and $$ S(\caC)_{n+2} \xrightarrow{\sigma} S(\caC)_{n} \times \widetilde{\caC \times \caC} \xrightarrow{ \gamma_n \times \widetilde{\caC \times \caC} }  (\widetilde{\caC \times \caC})^{\frac{n+1}{2}} \times \caC $$ for $n $ odd.
This concludes the inductive step \textemdash since then (2) also holds for $\frac{n}{2}+1 $ for $n$ even and $\frac{n-1}{2} +1$ for $n$ odd\textemdash, and with it the proof.
\end{proof}

\cref{dfghjhgk} allows us to show that the
$\caS^{1,1}$-construction localizes to the real suspension $\Sigma^{1,1}$ after inverting the morphisms $\gamma_3\colon S(\caD)_3\to\widetilde{\caD\times \caD}\times \caD$
for all $\caD \in \Wald_\infty^{\gd}.$ To make this precise, we need to introduce the following localization of $\D\Wald^\gd_\infty$.

\begin{definition}
The real $\infty$-category of fissile Waldhausen $\infty$-categories with genuine duality is the full subcategory 
$$\D\Wald_\infty^{\gd,\fiss} \subset \D\Wald_\infty^{\gd}$$ 
spanned by the real functors $(\Wald_\infty^{\gd})^\op \to \widehat{\Spc}^{C_2}$
that invert the morphisms  $\gamma_3\colon S(\caD)_3 \to \widetilde{\caD \times \caD} \times \caD$ that sends $(A \to B \to C) \mapsto (A, C/B, B/A)$ in $\Wald^\gd_\infty $ for $D \in \Wald^\gd_\infty .$
\end{definition}

\begin{remark}
The real $\infty$-category $\D\Wald_\infty^{\gd}$ is presentable in a larger universe. Observe that $\D\Wald_\infty^{\gd,\fiss}$ is the real accessible localization of $\D\Wald_\infty^{\gd}$ with respect to the morphisms $\gamma_3\colon S(\caD)_3 \to \widetilde{\caD \times \caD} \times \caD$ for $\caD \in \Wald^\gd_\infty.$ Hence $\D\Wald_\infty^{\gd,\fiss}$ is a reflective full subcategory of $\D\Wald_\infty^{\gd}$ and is in particular presentable in a larger universe, as well as genuine preadditive.
\end{remark}

\begin{remark}
The localized $\infty$-category $\D\Wald_\infty^{\gd,\fiss}$ is closed under sifted colimits in $\D\Wald_\infty^{\gd}$, since sifted colimits commute with finite products and cotensors with $C_2.$ Being $\D\Wald_\infty^{\gd,\fiss}$ preadditive, it is also closed under finite sums and so closed in $\D\Wald_\infty^{\gd}$ under all colimits.
\end{remark}

\begin{notation}
We denote by $L\colon\D\Wald_\infty^{\gd} \to \D\Wald_\infty^{\gd,\fiss}$ the real left adjoint of the embedding $\D\Wald_\infty^{\gd,\fiss}\subset \D\Wald_\infty^{\gd}.$
\end{notation}

\begin{corollary}\label{cor:S_suspension_commute}
There exists a commutative square
\[
\begin{tikzcd}
\D\Wald_\infty^\gd \ar[r, "S"] \ar[d, "L"]      & \rs\D\Wald_\infty^{\gd} \ar[d, "L"] \\
\D\Wald_\infty^{\gd,\fiss} \ar[r, "B"]      &\rs\D\Wald_\infty^{\gd,\fiss}.
\end{tikzcd}
\]
\end{corollary}

\begin{proof}
The real left adjoint $L$ of the embedding $\D\Wald_\infty^{\gd,\fiss}\subset \D\Wald_\infty^{\gd}$
is the universal semiadditive theory. So by \cref{dfghjhgk}
the map $L \circ \theta\colon L \circ B \to L \circ S $ is an equivalence. Moreover, as $L$ is left adjoint, we have $ B \circ L \simeq L \circ B $ so that we obtain a canonical equivalence $ B \circ L \simeq L \circ S.$
\end{proof}

\begin{corollary}\label{comm_sq_BL_LS}
There exists a commutative square

\[
\begin{tikzcd}
\D\Wald_\infty^\gd \ar[r, "\caS^{1,1}"] \ar[d, "L"]        & \D\Wald_\infty^{\gd} \ar[d, "L"] \\
\D\Wald_\infty^{\gd,\fiss} \ar[r, "\Sigma^{1,1}"]       & \D\Wald_\infty^{\gd,\fiss}.
\end{tikzcd}
\]

\end{corollary}

\begin{proof}
After composing the canonical equivalence $ B\circ L \simeq L \circ S$ of \cref{cor:S_suspension_commute} with geometric realization, we obtain a canonical equivalence $$\Sigma^{1,1} \circ L \simeq |-| \circ B \circ L \simeq |-| \circ L \circ S \simeq L \circ S^{1,1}. $$
\end{proof}

\begin{remark}
In particular, the real functor $\caS^{1,1}\colon \D\Wald_\infty^\gd \to \D\Wald_\infty^\gd$ preserves local equivalences, which we view as a weak form of additivity.
\end{remark}

\subsection{Some theory about left actions of a Segal space}\label{subsec:left_actions}

It is known that given a connected space $X$ and a map of spaces $Y \to X$ there is a canonical left action of the space $\Omega(X)$ on the fiber of $Y \to X$, since any loop in $X$ gives rise to an autoequivalence of the fiber. Moreover the map $Y \to X$ can be recovered from the left action of $\Omega(X)$. This is, the functor $\Spc_{/X} \to \Spc$ given by taking the fiber lifts to an equivalence
$$ \Spc_{/X} \simeq \LMod_{\Omega(X)}(\Spc), $$
where the right hand side denotes the $\infty$-category of spaces equipped with a left action of the loop space $\Omega(X).$

In this section we provide a similar account for Segal spaces instead of $\infty$-categories, and where spaces are replaced by $\infty$-pregroupoids in the sense of \cref{def:pregroupoid}; see \cref{prop:left_action_equiv_Segal_spaces}. This account will serve us to prove an identification criterion for when a functor $F\colon\caC\to\caD$ is genuine excisive, which in turn allows us to show in \cref{subsec:additivization} that the additivization of a theory is, as one would expect, an additive theory.

In the example above we considered an action of the loop space $\Omega(X)$, which is a grouplike $A_\infty$-space, on a space. More generally, one can consider actions of an $A_\infty$-space or a many object version of
$A_\infty$-space, which is a Segal space. We will encode such an action by a right fibration of simplicial spaces.

Given a right fibration of simplicial spaces $B\to A$, where $A$ is a Segal space, for every $n \in \Delta$ there is a canonical equivalence $ B_n \simeq A_n \times_{A_0} B_0$ induced by the map $\{n\} \subset [n]$ and so we obtain a map $\mu_n\colon A_n \times_{A_0} B_0 \simeq B_n \to B_0$ induced by the map $\{0\} \subset [n]$.

Let us now try to get a better feeling of what the action map $\mu_n\colon A_n \times_{A_0} B_0 \simeq B_n \to B_0$ does. Consider $s\colon V \to U$ a morphism of $A$, i.e.\ an object of the space $A_1$ with source $V $ and target $U$ in $A_0$. 
The map $\mu_1\colon A_1 \times_{A_0} B_0 \simeq B_1 \to B_0 $ gives rise to a map $ A_1 \times_{A_0} B_0 \to A_1 \times_{A_0} B_0 $ over $A_1$ that induces on the fibers over $s\in A_1$ a map $s^\ast\colon (B_0)_U \to (B_0)_V$. Given now two composable morphisms $t\colon W \to V, s\colon V \to U$ of $A$ with composition $ s \circ t$, i.e.\ the image under the map $ A_1 \times_{A_0} A_1 \simeq A_2 \to A_1$ induced by the map in $\Delta$ that sends 0 to 0 and 1 to 2, we have $$ (s \circ t)^\ast \simeq t^\ast \circ s^\ast\colon (B_0)_U \to (B_0)_V \to (B_0)_W.$$ Moreover for every object $U$ of $A_0$ the identity $\id_U$ of $U$, i.e. the image of $U$  under the map $A_0 \to A_1$ induced by the map $[1] \to [0]$ in $\Delta$, yields an equivalence $(\id_U)^\ast\colon (B_0)_U \to (B_0)_U.$ Hence for every invertible morphism $s\colon V \to U$ of $A$ the map $s^\ast\colon (B_0)_U \to (B_0)_V $ is an equivalence.

We will be interested in the case that $A$ is as in the following definition. 

\begin{definition}\label{def:pregroupoid}
We call a simplicial space $X$ an $\infty$-pregroupoid if for every $[m],[n],[l] \in \Delta$ such that $[n] = \{ 0 <...< m =0' < 1' < \cdots <l'\}= [m] \ast [l-1]$
the following square is a pullback square:
\[
\begin{tikzcd}
X_n \ar[r] \ar[d]       & X_m \ar[d] \\
X_{l} \ar[r]            & X_0.
\end{tikzcd}
\]
\end{definition}

\begin{notation}For an $\infty$-pregroupoid $A$ we write $\LMod_{A}(\Spc) \subset \caP(\Delta)_{/A} $ for the full subcategory spanned by the right fibrations of simplicial spaces with target $A$.
\end{notation}

Note that for $A \simeq \ast$ the final simplicial space, we have $\LMod_{A}(\Spc) \simeq \Spc \subset \caP(\Delta) $ since a map of simplicial spaces $B \to \ast$ is a 
right fibration of simplicial spaces if and only if $B$ is constant. Observe, too, that a pullback of a right fibration of simplicial spaces $B \to A$ along any map of simplicial spaces $\tau\colon A' \to A$ is a right fibration of simplicial spaces so that $\tau$ yields a change of base functor $\LMod_{A}(\Spc) \to \LMod_{A'}(\Spc)$. Especially the map $A \to \ast$ gives rise to a functor
$\mathrm{triv}_{A}\colon \Spc \simeq \LMod_{\ast}(\Spc) \to \LMod_{A}(\Spc)$ that sends
a space $X$ to $A \times \delta(X)$ with $\delta $ the diagonal embedding $\Spc \subset \caP(\Delta).$

\begin{proposition}\label{prop:left_action_equiv_Segal_spaces} Let $A$ be an $\infty$-pregroupoid. There is a a canonical equivalence
\[
\begin{tikzcd}
\Spc_{/{|A|}}\ar[r, "\theta", "\simeq"']     &\LMod_{A}(\Spc)
\end{tikzcd}
\]
natural in $A$ with respect to pullback.
\end{proposition}

\begin{proof}
The simplicial object $A$ is classified by a left fibration $\mathfrak{A} \to \Delta^\op$.The realization of $A$ is the $\infty$-category arising from $\mathfrak{A} $ after inverting all equivalences. Thus the canonical functor $\mathfrak{A} \to |A|$ yields an embedding 
$\Fun(|A|, \Spc) \subset \Fun(\mathfrak{A}, \Spc) $ with essential image those functors that invert every morphism.

Given a small $\infty$-category $\caD$, denote by $\mathrm{LFib}_{\caD} \subset \Cat_{\infty / \caD}$ the full subcategory spanned by the left fibrations with target $\caD.$ Then for every left fibration $\caE \to \caD$ the forgetful functor 
$(\mathrm{LFib}_{\caD})_{/\caE} \to \mathrm{LFib}_{\caE}$ is an equivalence.

So we obtain an equivalence 
$$ \Fun(\mathfrak{A}, \Spc) \simeq \mathrm{LFib}_{ \mathfrak{A} } \simeq (\mathrm{LFib}_{\Delta^\op})_{/ \mathfrak{A} } \simeq \caP(\Delta)_{/A}. $$

We will show next that a functor $\mathfrak{A} \to \Spc$ inverts every morphism if and only if its corresponding map of simplicial spaces $B \to A$ is a right fibration. And therefore conclude that such equivalence restricts to $$\theta\colon \Spc_{/{|A|}} \simeq \Fun(|A|, \Spc) \simeq \LMod_{A_1}(\Spc).$$

Consider a functor $h\colon\mathfrak{A} \to \Spc$, and its corresponding map of simplicial spaces $\xi\colon B \to A$. The map $\xi$ is classified by a map of left fibrations $\mathfrak{B} \to \mathfrak{A}$ over $\Delta^\op$, that is automatically itself a left fibration classifying the functor $h\colon \mathfrak{A} \to \Spc.$ 

Then, for any map $f\colon[n] \to [m]$ in $\Delta^\op$ we can construct a commutative square 
\begin{equation}\label{fkjhgp}
\begin{tikzcd}
\mathfrak{B}_n \ar[r] \ar[d]        & \mathfrak{B}_m \ar[d] \\
\mathfrak{A}_n \ar[r, "f_\ast"]      & \mathfrak{A}_m.
\end{tikzcd}
\end{equation}

The proof consists on showing that, on one hand, the functor $h\colon\mathfrak{A} \to \Spc$ inverts every morphism lying over $f$ if and only if the square (\ref{fkjhgp}) is a pullback square; and in the other hand, that the square (\ref{fkjhgp}) is a pullback square for every morphism $f\colon [n] \to [m]$ in $\Delta^\op$ if and only if $B\to A$ is a right fibration.

Since the functor $\mathfrak{A}\to \Delta^\op$ is a left fibration, every morphism lying over $f\colon [n] \to [m]$ in $\Delta^\op$ is of the form $X \to f_\ast(X)$ in $\mathfrak{A}$. Given such a morphism $X \to f_\ast(X)$,
the square (\ref{fkjhgp}) gives a map between the fibers of the vertical maps over $X$ and $f_\ast(X)$, which is equivalent to the induced map   $$\mathfrak{B}_X \simeq (\mathfrak{B}_n)_X \to \mathfrak{B}_{f_\ast(X)} \simeq (\mathfrak{B}_m)_{f_\ast(X)}$$ on fibers of the left fibration
$\mathfrak{B} \to \mathfrak{A}.$

For the second part, recall that a map $\xi\colon B \to A$ is a right fibration of simplicial spaces if square (\ref{fkjhgp}) is a pullback square for every map $\alpha_n\colon [n] \to [0]$ in $\Delta^\op$ corresponding to the map $\{n\} \subset [n]$ in $\Delta.$ Then if if $h$ inverts all morphisms, in particular $\xi\colon B \to A$ is a right fibration of simplicial spaces.

Let us assume now that $\xi\colon B \to A$ is a right fibration of simplicial spaces. Given a map $g\colon [n] \to [m]$ in $\Delta^\op$ 
by the pasting law square (\ref{fkjhgp}) is a pullback square for $g$ if it is a pullback square for the composition $[n] \xrightarrow{g} [m] \xrightarrow{\alpha_m} [0] $ in $\Delta^\op$.
So it is enough to check that square (\ref{fkjhgp}) is a pullback square for every map $[n] \to [0] $ in $\Delta^\op$. Every map $[n] \to [0] $ in $\Delta^\op$
factors as $[n] \xrightarrow{\beta} [1] \to [0] $ for some map $\beta$ in $\Delta^\op$ corresponding to a map in $\Delta$ that sends $1$ to $n.$ As $\alpha_1 \circ \beta = \alpha_n$ in $\Delta^\op,$ by the pasting law square 
(\ref{fkjhgp}) is a pullback square for $\beta.$ Consequently it is enough to show that square (\ref{fkjhgp}) is a pullback square for the map $\nu\colon [1] \to [0]$ in $\Delta^\op$ corresponding to the map $\{0\} \subset [1].$

In this case we need to check that for every object $X$ of $\mathfrak{A}_1$ encoding a morphism $s$ of $A$ with source $V $ and target $U$, the map $\nu\colon [1] \to [0]$ in $\Delta^\op$ induces an equivalence on fibers $$s^\ast\colon \mathfrak{B}_X \simeq (\mathfrak{B}_0)_U  \to \mathfrak{B}_{f_\ast(X)} \simeq (\mathfrak{B}_0)_V.$$ Since $A$ is an $\infty$-pregroupoid, $s$ is invertible. So $s^\ast$ is an equivalence.
\end{proof}

\begin{remark}
The equivalence $\theta$ is inverse to the functor $\LMod_A(\Spc) \to \Sp_{|A|}$ sending a right fibration $B \to A$ to $|B| \to |A|$.
\end{remark}

\section{\texorpdfstring{The connective real $K$-theory spectrum}{The connective real K-theory spectrum}}\label{sec:spectrum}

In this section we define the connective genuine $C_2$-spectrum of real $K$-theory. We begin by defining the real $K$-theory genuine $C_2$-space associated to a Waldhausen $\infty$-category with genuine duality in \cref{subsec:KR_space}; we also derive from it the corresponding hermitian $K$-theory. \cref{subsec:enriched_spectra} is devoted to introduce enriched spectra. It will be especially relevant for us one case of the $\Spc^{C_2}$-enriched variant, which we call genuine spectra (\cref{def:genuine_spectra}) and its connective variant (\cref{def:connective_genuine_spectra}). The most relevant result therein being an enriched version of the fact that excisive functors lift to spectra; see \cref{spec}. 

Finally, in \cref{subsec:KR_spectrum} we promote the real $K$-theory spaces defined in \cref{subsec:KR_space} to genuine spectra (see \cref{fghjkmno}). For this, we use that the real $K$-theory functor $\KR\colon \D\Wald_\infty^\gd \to \Spc^{C_2}$ is an additive theory and in consequence a genuine excisive functor, which implies that it lifts to monoid objects in $\Spc^{C_2}$ whose underlying monoid in $\Spc$ is grouplike. We strengthen the latter by showing in \cref{fhjkkkl} that it actually lifts to group objects in $\Spc^{C_2}$. The resulting lift $\D\Wald_\infty^\gd \to \Grp(\Spc^{C_2} )$ is genuine excisive, and thus by our recount of enriched spectra, it lifts to $\Spc^{C_2}$-enriched $S^{1,1}$-spectra in $\Grp_{E_\infty}(\Spc^{C_2})$. Then we identify $S^{1,1}$-spectra in $\Grp_{E_\infty}(\Spc^{C_2})$ with connective genuine $C_2$-spectra.

\subsection{Real $K$-theory genuine $C_2$-space}\label{subsec:KR_space}
In this subsection we define the real $K$-theory genuine $C_2$-space associated to a Waldhausen $\infty$-category with genuine duality, and describe its component parts in \cref{description_KR_u_C2}. In particular, we define and obtain a description of the corresponding hermitian $K$-theory space.

We recall \cref{def:KR_space}, where we first introduced the real $K$-theory functor:
The real $K$-theory functor
$$\KR\colon\D\Wald_\infty^\gd\to \Spc^{C_2}_\ast$$
is defined as the additivization $\KR\coloneqq\add(\iota)$ of the functor $\iota\colon\D\Wald_\infty^\gd\to\Spc^{C_2}_\ast$, which sends an $\infty$-category with genuine duality $\caC$ to its maximal subspace in $\caC$ equipped with the restricted genuine duality.

\begin{remark}
By \cref{fghkkkkl}, the real functor $\KR$ preserves sifted colimits and thus is uniquely determined by its restriction to $\Wald_\infty^\gd$.
\end{remark}

For any $\caC \in \Wald_\infty^\gd$ the real
K-theory space admits the following description, that the reader will find familiar
$$\KR(\caC) =\add(\iota)(\caC) = (\Omega^{1,1}\circ\iota\circ\caS^{1,1})(\caC) \simeq \Omega^{1,1}(|\iota \circ S(\caC)|).$$

In the following result, we study each of its components: its $C_2$-fixed points, and its underlying $\infty$-category. For this, we consider the functor $S(C)^\simeq\colon\Delta^\op\to  \Spc$ given by the composite 
\[
\begin{tikzcd}
\Delta^\op\ar[r, "S(\caC)"]     &\Exact_\infty\ar[r, "(-)^\simeq"]        &\Spc,
\end{tikzcd}
\]
where we also use $S(\caC)$ for the underlying functor of the real functor $S(\caC)\colon\underline{\Delta}^{\op}\to\Wald^\gd_\infty$ and $(-)^\simeq$ takes the maximal subspace; as well as the functor $S(\caC)^h\colon\Delta^\op\to  \Spc$, defined as the composite 
\[
\begin{tikzcd}
\Delta^\op\ar[r, "e"]       &(\Delta^{\op})^{hC_2}\ar[r, "S(\caC)"]        &\D\Wald^\gd_\infty\ar[r, "\iota"]        &\Spc^{C_2}\ar[r, "(-)^{C_2}"]&\Spc.
\end{tikzcd}
\]

\begin{proposition}\label{description_KR_u_C2}
Let $\caC$ be a Waldhausen $\infty$-category with genuine duality. Then there exist canonical equivalences
$$ \KR(\caC)^u \simeq \Omega(|S(\caC)^\simeq|)\hspace{1em}\text{and}\hspace{1em}\KR^{C_2}(\caC) \simeq \fib(|S(\caC)^h| \to |S(\caC)^\simeq|).$$
\end{proposition}

\begin{proof}For the left most equivalence, we observe that
$$ \KR(\caC)^u = \Omega^{1,1}(|\iota \circ S(\caC)|)^u \simeq\Omega(\vert \iota\circ S(\caC)\vert^u) \simeq\Omega(|S(\caC)^\simeq|).$$

For the right most equivalence, by \cref{lem:description_fixed_points_Omega11}
we have 
$$\KR(\caC)^{C_2}= \Omega^{1,1}(|\iota \circ S(\caC)|)^{C_2} \simeq \fib(|S(\caC)^h| \to |S(\caC)^\simeq|).$$

The statement then follows from the canonical equivalence 
$$|\iota \circ S(\caC)|^{C_2} \simeq |\iota \circ S(\caC) \circ e|^{C_2} \simeq |S(\caC)^h|,$$
where on the most left term we are considering real geometric realization \textemdash since in this case $S(\caC)$ refers to the real functor $S(\caC)\colon\underline{\Delta}^\op\to\Wald^\gd_\infty$\textemdash while in the other terms we are taking usual geometric realization, and $S(\caC)$ refers to the underlying functor $(\Delta^\op)^{hC_2}\to\Wald^\gd_\infty$.
\end{proof}

As in Hesselholt and Madsen's work \cite{hesselholt-madsen}, we define hermitian $K$-theory from the real $K$-theory space by by taking the corresponding space of $C_2$-fixed points. Therefore, \cref{description_KR_u_C2} provides an explicit description of the hermitian $K$-theory space.

\begin{definition}The hermitian $K$-theory functor, $\KH\colon\D\Wald^\gd_\infty\to\Spc_\ast$ is defined as
$$\KH\coloneqq(-)^{C_2}\circ\KR.$$
\end{definition}

With this terminology, we deduce from \cref{description_KR_u_C2} a comparison with the Grothendieck-Witt space of a stable $\infty$-category with duality defined by the second author in \cite[Definition 3.4]{spitzweck.gw}.

\begin{corollary}Let $\caC$ be a stable $\infty$-category with duality. Then 
$\KH(\caC)$, when $\caC$ is considered with the standard genuine duality, and the Grothendieck-Witt space of $\caC$ as defined in \cite{spitzweck.gw} coincide.
\end{corollary}

\subsection{Enriched spectra}\label{subsec:enriched_spectra}
We dedicate this subsection to introduce generalities of genuine $C_2$-spectra, \cref{def:genuine_spectra}. Throughout the rest of the subsection we let $\caV$ be a presentable symmetric monoidal $\infty$-category.

\subsubsection{Excisive functors} In this subsection we introduce enriched excisive functors, which not only will play a crucial role in the definition of enriched spectra but also allow for an interpretation of real $K$-theory as a first Goodwillie derivative in \cref{prop:Kth_Goodwillie}.
 
\begin{definition}
Let $\caC,\caD \in \Cat_\infty^\caV$ and $K \in \caV$. Let $\caC$ be such that admits tensors with $K$ and $\caD$ such that it admits cotensors with $K$. We say that a $\caV$-enriched functor $F \colon \caC \to \caD$ is $K$-excisive if for every $X \in \caC$ the natural map $F(X) \to F(K \otimes X)^K$ is an equivalence.
\end{definition}

We denote by 
$$\Exc_K^\caV(\caC,\caD) \subset \Fun_\caV(\caC,\caD)$$ 
the full subcategory spanned by the $K$-excisive functors.

\begin{remark}
Since $\caD$ admits cotensors with $K$ and in $ \Fun_\caV(\caC,\caD)$ those are computed objectwise, also $ \Fun_\caV(\caC,\caD)$ admits cotensors with $K$. Moreover for every $K$-excisive $\caV$-enriched functor $F\colon \caC \to \caD$, its cotensor with $K$ is $K$-excisive. Then $\Exc_K^\caV(\caC,\caD)$ admits cotensors with $K$, too.
\end{remark}

\begin{proposition}\label{prop:Kexcisive_approximation}
Let $\caC,\caD \in \Cat_\infty^\caV$ and $K \in \caV$. Let $\caC$ be such that admits tensors with $K$ and $\caD$ such that it admits cotensors with $K$.
The embedding
$$\Exc_K^\caV(\caC,\caD) \subset \Fun_\caV(\caC,\caD)$$
admits a left adjoint $P$, which we call $K$-excisive approximation,
and which sends $F$ to the colimit of the sequential diagram
$$ F \to F(K \otimes (-))^K \to (F(K \otimes K \otimes (-))^K)^K \to\cdots  $$
\end{proposition}

\begin{proof}
Given any $\caV$-enriched functor $F\colon \caC\to \caD$, consider $P(F) $ to be the colimit of the diagram in the statement. This assembles into a functor $P\colon\Fun_\caV(\caC,\caD) \to \Fun_\caV(\caC,\caD).$ 

Now, we want to show that the $\caV$-enriched functor $P(F)$ is actually $K$-excisive. Indeed, the composite
$$\colim(F \to F(K \otimes (-))^K \to (F(K \otimes K \otimes (-))^K)^K \to\cdots) \longrightarrow $$$$ \colim(F(K \otimes (-)) \to F(K \otimes K \otimes (-))^K \to\cdots)^K $$$$ \simeq \colim(F(K \otimes (-))^K \to (F(K \otimes K \otimes (-))^K)^K \to\cdots)$$ is the canonical equivalence.
Thus $P$ defines a functor $\Fun_\caV(\caC,\caD) \to \Exc_K^\caV(\caC,\caD).$

Finally, we need to prove that $P$ is indeed left adjoint to the embedding. Let $\eta\colon F \to P(F)$ be the canonical $\caV$-enriched natural transformation. Then $\eta$ defines a natural transformation $\id \to P$ of functors $\Fun_\caV(\caC,\caD) \to \Fun_\caV(\caC,\caD).$ The $\caV$-enriched natural transformation
$\eta\colon F \to P(F)$ is an equivalence if $F$ is $K$-excisive
since then the diagram
$$F \to F(K \otimes (-))^K \to (F(K \otimes K \otimes (-))^K)^K \to\cdots $$ becomes constant. Therefore to prove the statement it suffices to prove that the $\caV$-enriched natural transformation $\eta\colon F \to P(F)$ is a $P$-equivalence.
Now, the induced $\caV$-enriched natural transformation
$P(\eta)\colon P(F) \to P(P(F))$ identifies with the
$\caV$-enriched natural transformation $\eta_{P(F)}\colon P(F) \to P(P(F)),$
which is an equivalence by what we have proven.
\end{proof}

The theory of enriched Goodwillie calculus has not yet been developed, and it is not a straightforward generalization of the classical one for levels higher than one. However, it is to be expected that any such theory has the excisive approximation as above as the first Goodwillie derivative. With this in mind, and inspired by Barwick's description of algebraic $K$-theory as a first Goodwillie derivative in  \cite[Theorem 7.4]{barwick.wald} we prove the following result.

\begin{theorem}\label{prop:Kth_Goodwillie}
The real $K$-theory functor $\KR\colon \D\Wald_\infty^\gd \to \Spc_*^{C_2}$ 
factors as $$ \D\Wald_\infty^\gd \xrightarrow{L} \D\Wald_\infty^{\gd,\fiss} \xrightarrow{P(\iota)} \Spc_*^{C_2},$$
where $\iota$ is the real functor $\D\Wald_\infty^{\gd,\fiss} \to \Spc_*^{C_2}$ that sends a small Waldhausen $\infty$-category with genuine duality to its space of objects.
\end{theorem}

\begin{proof}
As seen in \cref{prop:Kexcisive_approximation}, the genuine excisive approximation of $\iota$ is computed
as $$ \colim_{n\geq 0} (\Omega^{1,1})^n \circ \iota \circ (\Sigma^{1,1})^n.$$

Therefore by \cref{comm_sq_BL_LS} there is a canonical equivalence
$$ P(\iota)\circ L \simeq \colim_{n\geq 0} (\Omega^{1,1})^n \circ \iota \circ (\Sigma^{1,1})^n \circ L $$
$$ \simeq \colim_{n\geq 0} (\Omega^{1,1})^n \circ \iota \circ L \circ (\mathcal{S}^{1,1})^n $$
$$ \simeq \colim_{n\geq 0} (\Omega^{1,1})^n \circ \iota \circ(\mathcal{S}^{1,1})^n. $$

The last equivalence holds because the functor
$$ \mathcal{S}^{1,1} \colon \D\Wald_\infty^\gd \to \D\Wald_\infty^\gd $$
and therefore also the iterate $$(\mathcal{S}^{1,1})^n \colon \D\Wald_\infty^\gd \to \D\Wald_\infty^\gd $$
land in the full subcategory $ \D\Wald_\infty^{\gd,\fiss}.$
Indeed, this is equivalent to say that for every $\caC \in \Wald^\gd_\infty$
the functor $ \mathcal{S}^{1,1}(\caC) = |\map_{\Wald^\gd_\infty}(-,S(\caC)) |$ is local with respect to $\gamma_3$, which follows from the real addivity theorem \cref{thm:add}.
\end{proof}

\subsubsection{Enriched spectra}
Denote by $\caV_K \subset \caV$ the smallest symmetric monoidal full subcategory containing $K.$ Set $$ \Pre\Sp^\caV_K(\caD) \coloneqq \Fun_\caV(\caV_K, \caD)$$

\begin{definition}\label{def:V_enriched_K_spectra}Let $K\in\caV$ and $\caD$ a $\caV$-category that admits cotensors with $K$. We define the $\caV$-enriched $\infty$-category of $\caV$-enriched $K$-spectra in $\caD$ as 
$$\Sp^\caV_K(\caD)\coloneqq \Exc^\caV_K(\caV_K, \caD).$$
\end{definition}

\begin{example}\label{ex:enriched_spectra}
Let $\caV$ be the $\infty$-category of pointed spaces with the smash product
and $K= S^1$. Let also $\caC, \caD$ be pointed $\infty$-categories such that $\caC$ admits pushouts and $\caD$ admits pullbacks. Then we know that $\caC$ and $\caD$ are canonically $\Spc_\ast$-enriched $\infty$-categories such that $\caC$ admits tensors with $S^1$ (given by suspension) and $\caD$ admits cotensors with $S^1$ (given by loops).

In this case a reduced functor $F\colon \caC \to \caD$ is $S^1$-excisive if and only if it sends suspensions to loops, which by \cite[Proposition 1.4.2.13]{lurie.higheralgebra} is equivalent to it sending pushout squares to pullback squares. We conclude then that a reduced functor $F\colon \caC \to \caD$ is $S^1$-excisive (with the new definition) if and only if it is excisive (in the usual sense). Moreover, we have a canonical equivalence $\Sp^\caV_K(\caD) \simeq \Sp(\caD),$ where the right hand side denotes spectra objects in $\caD$.
\end{example}

\begin{remark}
Note that if $\caC$ is small and $\caD$ admits small colimits and tensors, the $\caV$-enriched $\infty$-category $\Exc_K^\caV(\caC, \caD) \subset \Fun_\caV(\caC, \caD) $ is the $\caV$-enriched localization with respect to the set of morphisms
$$ \{ \mathrm{lan}_{K \otimes Y}(K \otimes Z) \to \mathrm{lan}_{Y}(Z) \mid Y \in \caC, \ Z \in \caD \} .$$

In particular, $ \Sp^\caV_K \subset \Pre \Sp^\caV_K$ is a $\caV$-enriched localization.
\end{remark}

Let $\caD$ be a $\caV$-category and $K \in \caV.$
There is a forgetful $\caV$-functor $\Sp^\caV_K(\caD) \to \caD$
given by the composition
$\Sp^\caV_K(\caD) \subset \Fun^\caV(\caV_K,\caD) \to \caD$,
where the last $\caV$-functor evaluates at $1 \in \caV_K.$

\begin{lemma}\label{D_presentable_then_Sp(D_to_D_admits_left_adj}
If $\caD$ is presentable, the forgetful $\caV$-functor $\Sp^\caV_K(\caD) \to \caD$ admits a left adjoint.
\end{lemma}
\begin{proof}
We know that the embedding $\Sp^\caV_K(\caD) \subset \Fun^\caV(\caV_K,\caD)$
admits a left adjoint and by \cref{Venr:ev_admits_adjoint}
the evaluation $\caV$-functor $\Fun^\caV(\caV_K,\caD) \to \caD$
admits a left adjoint as well.
\end{proof}

For cases that are of special interest for us, we will lighten the notation to refer to spectra.
\begin{notation}When we consider $\caD=\caV$, we write $ \Sp^\caV_K$ for $\Sp^\caV_K(\caD)$ and
$ \Pre \Sp^\caV_K $ for $ \Pre \Sp^\caV_K(\caD).$
Moreover, we will write $\Sp_{K}(\caD)$ for $\Sp_{K}^{\Spc^{C_2}_\ast}(\caD)$.
\end{notation}

As in the non-enriched case, we can consider the functor 
$$\Omega_K^\infty \colon\Sp^\caV_K(\caD) \subset \Fun_\caV(\caV_K, \caD) \to \caD$$ the $\caV$-enriched functor that evaluates at the tensor unit of $\caV$. To identify when such a functor is an equivalence, we introduce the notion of stability.

\begin{definition}
We say that a $\caV$-enriched $\infty$-category $\caD$ is $K$-stable
if $\caD$ admits cotensors with $K$ and the cotensor functor $(-)^K\colon \caD \to \caD $ is an equivalence of $\caV$-enriched $\infty$-categories.
\end{definition}

\begin{remark}\label{rmk:K_stable_admits_tensors}
A $K$-stable $\caV$-category $\caD$ admits tensors with $K$. 
\end{remark}

\begin{definition}
Let $\caV=\Spc^{C_2}_\ast$ be endowed with the smash product and set $K=S^{1,1}$. We say that $D$ is genuine stable if it is $S^{1,1}$-stable.
\end{definition}

\begin{example}In the scenario of \cref{ex:enriched_spectra}, 
if $\caD$ is in addition admits also pushouts, we find by taking $F$ to be the identity of $\caD$, $\caD^\op$, respectively, that $\caD$ is stable if and only if $D$ is
$S^1$-stable.
\end{example}

\begin{example}\label{ex:gen_exc_S11_exc}
Let $\caV= \Spc^{C_2}_\ast$ be endowed with the smash product and set $K=S^{1,1}$. Consider $C$ to be a pointed real $\infty$-category that admits tensors with $S^{1,1}$ as $\Spc^{C_2}_\ast$-enriched $\infty$-category, and $\caD$ a pointed real $\infty$-category that admits cotensors with $S^{1,1}$ as $\Spc^{C_2}_\ast$-enriched $\infty$-category. A reduced real functor $F\colon \caC \to \caD$ is $S^{1,1}$-excisive if and only if $F$ is genuine excisive.
\end{example}

Now we dedicate ourselves to prove some structural results about $K$-excisive functors, that will allow us to identify when the functor $\Omega_K^\infty \colon\Sp^\caV_K(\caD) \to \caD$ is an equivalence

\begin{lemma}\label{mfhe3hto}
Let $\caC,\caD \in \Cat_\infty^\caV$ and $K \in \caV$. Suppose that the tensor with $K$ exists in $\caC$ and the cotensor with $K$ exists in $\caD$. Then the $\caV$-enriched $\infty$-category $\Exc_K^\caV(\caC,\caD)$ is $K$-stable. In particular, $\Sp^\caV_K(\caD) $ is $K$-stable.
\end{lemma}

\begin{proof}
By definition, we need to show that the $\caV$-enriched functor $$(-)^K\colon \Exc_K^\caV(\caC,\caD) \to \Exc_K^\caV(\caC,\caD)$$ is an equivalence. Since the cotensor with $K$ is computed objectwise, $ (-)^K $ is equivalent to the $\caV$-enriched functor $$(-)^K_\ast\colon \Exc_K^\caV(\caC,\caD) \to \Exc_K^\caV(\caC,\caD)$$ induced by $(-)^K\colon \caD \to \caD$.

The $\caV$-enriched functor $K \otimes - \colon \caC \to \caC$ induces a $\caV$-enriched functor $$(K \otimes -)^\ast\colon \Exc_K^\caV(\caC,\caD) \to \Exc_K^\caV(\caC,\caD)$$ that is inverse to the $\caV$-enriched functor $(-)^K_\ast$. Indeed, there is a canonical equivalence
$$ (K \otimes -)^\ast \circ (-)^K_\ast \simeq (-)^K_\ast \circ (K \otimes -)^\ast$$
as well as a canonical map $ \id \to (-)^K_\ast \circ (K \otimes -)^\ast$
that is an equivalence.
\end{proof}

\begin{proposition}\label{fghjkl}
Let $\caD$ be a $\caV$-enriched $\infty$-category that admits cotensors with $K.$ Then $\caD$ is $K$-stable if and only if the $\caV$-enriched functor $ \Omega_K^\infty\colon \Sp^\caV_K(\caD) \to \caD  $ is an equivalence.
\end{proposition}

\begin{proof}
If we assume that $\Omega_K^\infty\colon \Sp^\caV_K(\caD) \to \caD$, by \cref{mfhe3hto} we conclude that $\caD$ is $K$-stable.	

Assume now that $\caD$ is $K$-stable. Then $\caD$ must admit tensors with $K$,
and we know that $K \otimes -$ is the inverse of $(-)^K\colon \caD \to \caD.$ 
But we know that a $\caV$-enriched functor $F\colon  \caV_K \to \caD$ is $K$-excisive if and only if it preserves tensors with $K$ and in this case $F$ is determined by its value on the tensorunit of $\caV.$
\end{proof}

In the next result, we identify conditions under which we can lift an excisive functor to spectra. 
\begin{proposition}\label{spec}
Let $\caC,\caD \in \Cat_\infty^\caV$ and $K \in \caV$. Suppose that the tensor with $K$ exists in $\caC$ and that the cotensor with $K$ exists in $\caD$. Then the $\caV$-enriched functor $\Omega_K^\infty\colon \Sp^\caV_K(\caD) \to \caD$ induces an equivalence $$\Exc_K^\caV(C,\Sp_K^\caV(\caD)) \to \Exc_K^\caV(C,\caD).$$

\end{proposition}

\begin{proof}
This follows from \cref{mfhe3hto,fghjkl}, together with the fact that the canonical equivalence
$\Fun_\caV(\caC, \Fun_\caV(\caV_K, \caD)) \simeq \Fun_\caV(\caV_K, \Fun_\caV(\caC, \caD))$ restricts to an equivalence
$$\Exc_K^\caV(\caC, \Sp_K^\caV(\caD)) \simeq \Sp^\caV_K(\Exc_K^\caV(\caC,\caD)).$$
\end{proof}

We come back now to the case when stability conditions do not necessarily hold, and we will prove some general facts about enriched spectra. 

\begin{corollary}\label{lem:Sp_compatibility}Let $\caV$ be a symmetric monoidal $\infty$-category, $\caD$ a $\caV$-enriched $\infty$-category, and $K,T\in\caV$. Then there is a canonical $\caV$-enriched equivalence
$$\Sp_K^\caV(\Sp_{T}(\caD)) \simeq \Sp_{K \otimes T}^\caV(\caD).$$
\end{corollary}
\begin{proof}
It follows directly by noticing that the forgetful functors $$\Sp^\caV_{K \otimes T}(\Sp^\caV_K(\Sp^\caV_{T}(\caD))) \to \Sp^\caV_{K \otimes T}(\caD),$$
$$\Sp_{K \otimes T}(\Sp^\caV_K(\Sp_{T}(\caD))) \to \Sp^\caV_K(\Sp^\caV_{T}(\caD))$$
are equivalences.
\end{proof}

\begin{proposition}
Let $\caD \in \Cat_\infty^\caV$ and $K \in \caV$ be such that $\caD$ admits cotensors with $K$. There is a canonical equivalence of $\caV$-categories
$$\Sp_K^\caV(\caD) \simeq \lim(\cdots\caD \xrightarrow{(-)^K} \caD \xrightarrow{(-)^K} \caD),$$
where the limit is taken in the $\infty$-category $\Cat_\infty^\caV.$
\end{proposition}

\begin{proof}
We start with proving that the $\caV$-category 
$ \caB\coloneqq \lim(\cdots\caD \xrightarrow{(-)^K} \caD \xrightarrow{(-)^K} \caD)$
is $K$-stable.
The $\caV$-functor $(-)^K\colon\caD \to \caD$ gives rise to a 
$\caV$-functor $\xi\colon\caB \to\caB$, which is the cotensor with $K$ of $\caB$.
Let $\pr_i\colon \caB \to \caD $ denote the $i$-th component.
The $\caV$-functor $\pr_i$ factors as $$\xi\colon\caB \xrightarrow{(-)^K} \caD) \to \caB \xrightarrow{(-)^K} \caD) \xrightarrow{\pr_{i+1}} \caD. $$
So $\xi$ is an equivalence: an inverse is given by the $\caV$-functor $\rho\colon \caB \to \caB$, which is determined by asking that the $\caV$-functor $\pr_{i+1}$ factors as $\rho\colon \caB \to \caB \xrightarrow{\pr_{i}} \caD. $

By \cref{mfhe3hto} the $\caV$-category $\Sp_K^\caV(\caD)$ is $K$-stable. 
Consequently it is enough to see that for any 
$K$-stable $\caV$-category $\caC $ 
there is a canonical equivalence
$$\Exc_K^\caV(\caC,\Sp_K^\caV(\caD)) \simeq \Exc_K^\caV(\caC,\caB).$$

Since $\Fun^\caV(-,-)$ is an internal hom, it commutes with limits in the second variable and therefore the canonical $\caV$-functor
$$\Fun^\caV(\caC,\caB) \to \lim (\cdots\Fun^\caV(\caC,\caD) \xrightarrow{(-)^K}\Fun^\caV(C,\caD) \xrightarrow{(-)^K}\Fun^\caV(C,\caD)) $$
is an equivalence, which restricts to an equivalence
$$\Exc_K^\caV(\caC,\caB) \to  \lim (\cdots\Exc_K^\caV(\caC,\caD) \xrightarrow{(-)^K}\Exc_K^\caV(\caC,\caD) \xrightarrow{(-)^K}\Exc_K^\caV(\caC,\caD)). $$
By \cref{mfhe3hto} the $\caV$-category $\Exc_K^\caV(\caC,\caD)$ is $K$-stable.
Thus the $\caV$-functor $$\pr_0\colon  \lim (\cdots\Exc_K^\caV(\caC,\caD) \xrightarrow{(-)^K}\Exc_K^\caV(\caC,\caD)\xrightarrow{(-)^K}\Exc_K^\caV(\caC,\caD)) \to \Exc_K^\caV(\caC,\caD)$$
is an equivalence.
By \cref{spec} the canonical $\caV$-functor
$$\Exc_K^\caV(\caC,\Sp_K^\caV(\caD)) \to \Exc_K^\caV(\caC,\caD)$$
is an equivalence.
The desired equivalence is the composition
$$ \Exc_K^\caV(\caC,\Sp_K^\caV(\caD)) \simeq \Exc_K^\caV(\caC,\caD) \simeq\Exc_K^\caV(\caC,\caB). $$
\end{proof}

The unique left adjoint symmetric monoidal functor $\Spc \to \caV$ yields a left adjoint symmetric monoidal functor $\tau\colon \Spc_\ast \to \caV_\ast,$ where $\caV_\ast$ carries the smash product.

\begin{notation}\label{not:Sp_V_C_canonical_enrichment}
Let $\caD$ be a $\caV$-category, which admits a zero object, that we see as a $\caV_\ast$-category via \cref{Venr:pointed_vs_not_pointed}. Assume that $\caD$ admits cotensors with $\tau(S^1)$ as a $\caV_\ast$-category.

We write 
$$\Sp^\caV(\caD) \coloneqq \lim(\cdots\caD \xrightarrow{(-)^{\tau(S^1)}} \caD \xrightarrow{(-)^{\tau(S^1)}} \caD),$$
where the limit is taken in $\Cat_\infty^{\caV_\ast}.$
\end{notation}

\begin{proposition}\label{prop:Sp(C)_admits_real_structure}
Let $\caD$ be a $\caV$-category, which admits a zero object, that we see as a $\caV_\ast$-category via \cref{Venr:pointed_vs_not_pointed}. Assume that $\caD$ admits cotensors with $\tau(S^1)$ as a $\caV_\ast$-category. Then the $\caV$-category $\Sp^\caV(\caD)$ has underlying $\infty$-category the usual $\infty$-category of spectra on the underlying $\infty$-category of $\caC$. 
\end{proposition}

\begin{proof}
First note that the $\caV_\ast$-functor $ (-)^{\tau(S^1)} \colon \caD \to \caD$ induces on underlying $\Spc_\ast$-enriched $\infty$-categories the functor $\Omega= (-)^{S^1}.$
So we conclude the proof by using that the forgetful functor $\Cat_\infty^{\caV_\ast} \to \Cat_\infty^{\Spc_\ast}$ preserves small limits, as it is a right adjoint. 
\end{proof}

Finally, we define a very special case for us of enriched spectra: that of genuine $C_2$-spectra.

\begin{definition}\label{def:genuine_spectra}
For any real $\infty$-category $\caD$ we define the real $\infty$-category of genuine $C_2$-spectra in $\caD$ as
$$\Sp^{C_2}(\caD)\coloneqq \Sp^{\Spc_\ast^{C_2}}_{S^{2,1}}(\caD),$$
where $S^{2,1}\coloneqq S^{1,1} \wedge S^{1,0}$ and $S^{1,0}$ is the sphere $S^{1}$ with trivial $C_2$-action and genuine refinement
given by the diagonal map $S^1 \to (S^1)^{hC_2}=\Fun(BC_2,S^1).$
\end{definition}

\begin{notation}When we consider $\caD=\Spc^{C_2}_\ast$, we simply write $\Sp^{C_2}$ for $\Sp^{C_2}(\Spc_\ast^{C_2})$, the real $\infty$-category of genuine $C_2$-spectra.
\end{notation}

The $\infty$-category of genuine $C_2$-spectra admits the following alternative description that we give in \cref{char_genuine_spectra}, for which we present some results.

Finally, we are able to describe genuine $C_2$-spectra as stated below.

\begin{proposition}\label{char_genuine_spectra}
There is a canonical equivalence of real $\infty$-categories as follows
$$\Sp^{C_2}(\caD) \simeq \Sp_{S^{1,1}}(\Sp(\caD)).$$
\end{proposition}

\begin{proof}
This follows from \cref{prop:Sp(C)_admits_real_structure} by applying \cref{lem:Sp_compatibility}.
\end{proof}

\subsection{Connective real $K$-theory genuine $C_2$-spectrum}\label{subsec:KR_spectrum}

Finally, we lift the real $K$-theory genuine $C_2$-space to a connective genuine $C_2$-spectrum. We achieve this in two steps. Firstly, we functorially refine the $K$-theory genuine $C_2$-space to a (non-genuine) connective $C_2$-spectrum in \cref{fhjkkkl}; secondly, we refine such connective $C_2$-spectrum to a connective \emph{genuine} $C_2$-spectrum (\cref{def:connective_genuine_spectra}) in \cref{fghjkmno}. For the latter, it will essential to know that that real $K$-theory functor is genuine excisive, which is implied by \cref{additivization_is_additive}.

For this section, it will be useful to have at hand a canonical real functor $\Grp_{E_\infty}(\caD) \to \Sp(\caD)$ for any real $\infty$-category $\caD$. We construct such a functor in what follows.

\begin{proposition}
Let $\caD$ be a presentable $\caV$-category. There is a canonical $\caV$-functor $\Grp^\caV_{E_\infty}(\caD) \to \Sp^\caV(\caD)$
that is left adjoint in the enriched sense to a lift $ \Sp^\caV(\caD) \to\Grp^\caV_{E_\infty}(\caD) $ of the forgetful functor
$ \Sp^\caV(\caD) \to \caD$.
\end{proposition}
\begin{proof}

By \cref{D_presentable_then_Sp(D_to_D_admits_left_adj} there is a pair of right and left adjoints $R\colon\caD \rightleftarrows \Sp^\caV(\caD)\colon L$ that induces, by \cref{Venr:Grp_is_reflective}, a pair of right and left adjoints $\Grp^\caV_{E_\infty}(\caD) \rightleftarrows \Grp^\caV_{E_\infty}(\Sp^\caV(\caD))$. Now, \cref{prop:Sp(C)_admits_real_structure,Venr:V_pre_add_when_underlying_is_V_pre_add} imply that the $\caV$-category $\Sp^\caV(\caD)$ is $\caV$-additive, and therefore \cref{Venr:forgetful_Mon_Grp_are_equiv} says that the forgetful functor $\Grp^\caV_{E_\infty}(\Sp^\caV(\caD)) \to \Sp^\caV(\caD)$
is an equivalence.

This concludes the proof, since we can consider the composition 
$$\Grp^\caV_{E_\infty}(\caD) \rightleftarrows \Grp^\caV_{E_\infty}(\Sp^\caV(\caD)) \simeq \Sp^\caV(\caD).$$
\end{proof}

We now proceed with the meaning of connective in this context.

\begin{definition}\label{def:connective_genuine_spectra}
Let $\caD$ be a real $\infty$-category. We call a genuine $C_2$-spectrum in $\caD$, say $X$, connective if seen as an $S^{1,1}$-spectrum in $\Sp(\caD)$, each level $X_n\coloneqq X((S^{1,1})^{\wedge n})$ \textemdash which is a spectrum in $\caD$\textemdash belongs to the essential image of the functor $\Grp_{E_\infty}(\caD) \to \Sp(\caD)$.
\end{definition}
We will denove by $\Sp^{C_2}_{\geq0} (\caD) \subset  \Sp^{C_2}(\caD) $  the full subcategory spanned by the connective genuine $C_2$-spectra in $\caD$.

Let us go into a short digression in order to justify \cref{def:connective_genuine_spectra}. Recall that in the classical setting, the forgetful functor $\Sp \to \Spc$ lifts to an equivalence $\Sp_{\geq 0} \simeq \Grp_{E_\infty}(\Spc)$, between connective spectra and grouplike $E_\infty$-spaces. Using that there exists an equivalence $\Spc^{C_2} \simeq \Fun(BC_2^{\triangleleft},\Spc) $, we see that it is also true that the forgetful functor $\Sp(\Spc^{C_2}) \to \Spc^{C_2}$ lifts to an equivalence $\Sp(\Spc^{C_2})_{\geq 0} \simeq \Grp_{E_\infty}(\Spc^{C_2})$,
where $\Sp(\Spc^{C_2})_{\geq 0}$ is the full subcategory spanned by those spectra in $\Spc^{C_2}$ which induce on $C_2$-fixed points and on the underlying space a connective spectrum.

With that in mind, we observe that \cref{def:connective_genuine_spectra} for $\caD=\Spc^{C_2}$ calls a genuine $C_2$-spectrum $X$ connective if seen as a $S^{1,1}$-spectrum in $\Sp(\Spc^{C_2})$ each level $X_n$, which is a spectrum in $\Spc^{C_2}$, induces on $C_2$-fixed points and on the underlying space a connective spectrum.

Moreover, this situation emerges for $\D\Wald_\infty^\gd$ as well. 

\begin{lemma}\label{rmk:embedding_grplike_objects_in_Sp}
Let $\caC$ be a small real $\infty$-cateogry, then
the canonical real functor below is fully faithful. $$\Grp_{E_\infty}(\D\caC) \to \Sp(\D\caC)$$ 
\end{lemma}

\begin{proof}
By definition the $\infty$-category $\D\Wald_\infty^\gd $ is a real localization of $\Fun_{ \Spc^{C_2}} (\caC^\op, \Spc^{C_2})$. Thus it suffices to prove the result for
$\caD= \Fun_{ \Spc^{C_2}} (\caC^\op, \Spc^{C_2}).$
In this case the real functor $\Grp_{E_\infty}(\caD) \to \Sp(\caD)$ factors as
\[
\begin{tikzcd}[column sep=small]
\Grp_{E_\infty}(\caD)\ar[dr, "\simeq"]      &       &       &\Sp(\caD).\\
\           &\Fun_{ \Spc^{C_2}} (\caC^\op,\Grp_{E_\infty}( \Spc^{C_2}))\ar[r]      &\Fun_{ \Spc^{C_2}} (\caC^\op,\Sp( \Spc^{C_2}))\ar[ur, "\simeq"]       &\ 
\end{tikzcd}
\]
We conclude then that it is fully faithful.
\end{proof}

After enlarging the universe, we can apply \cref{rmk:embedding_grplike_objects_in_Sp} to $\D\Wald^\gd_\infty$.

To achieve the goal of lifting the real $K$-theory space to a genuine $C_2$-spectrum, we present in what follows an intermediate result about lifting the real functor $$\Omega^{1,1}\circ\caS^{1,1}\colon \D\Wald_\infty^\gd\to \D\Wald_\infty^\gd$$ to grouplike $E_\infty$-objects. It will be useful to have in mind the fact that given a real functor $F\colon\caC\to\caD$ that preserves finite products, where $\caC$ is preadditive and $\caD$ has final products (both conditions in the enriched sense), there is a unique lift $F\colon\caC\to\Mon_{E_\infty}(\caD)$ (see \cref{Venr:lifting_to_monoids}).

\begin{proposition}\label{fhjkkkl}
The unique lift of the real functor $\Omega^{1,1} \circ \caS^{1,1}\colon \D\Wald_\infty^\gd \to \D\Wald_\infty^\gd$ to $E_\infty$-monoids in $\D\Wald_\infty^\gd$ factors as 
\[\begin{tikzcd}[column sep=tiny]
\D\Wald_\infty^\gd\ar[rr]\ar[dr]       &           &\Mon_{E_\infty}(\D\Wald_\infty^\gd).\\
\                               &\Grp_{E_\infty}(\D\Wald_\infty^\gd)\ar[ru, hook]
\end{tikzcd}
\]
\end{proposition}

\begin{proof}
We begin the proof by showing that the real functor $\Omega^{1,1}\circ\caS^{1,1}\colon\D\Wald^\gd_\infty\to\D\Wald_\infty^\gd$ lifting to a real functor whose target is $\Grp_{E_\infty}(\D\Wald_\infty^\gd)$, is equivalent to the real functor 
\begin{equation}\label{functor_thats_equiv_to_lift}
\Wald_\infty^\gd \xrightarrow{S} \rs\Wald_\infty^\gd \xrightarrow{\iota}  \rs\Spc_\ast \xrightarrow{|-| } \Spc^{C_2}_\ast \xrightarrow{\Omega^{1,1}}  \Spc^{C_2}_\ast 
\end{equation}
lifting to commutative group objects in $\Spc_\ast^{C_2}$, where $\iota\colon\rs\Wald_\infty^\gd\to\rs\Spc_\ast$ takes the underlying space in each level.

Since commutative group objects are closed under sifted colimits
and by \cref{fghkkkkl} the real functor $ \Omega^{1,1} \circ \caS^{1,1}\colon  \D \Wald_\infty^\gd \to \D \Wald_\infty^\gd$ preserves sifted colimits, it is enough to check that the restriction
$$ \Wald_\infty^\gd \subset \D \Wald_\infty^\gd \xrightarrow{\Omega^{1,1} \circ \caS^{1,1}} \D\Wald_\infty^\gd $$ lifts to grouplike $E_\infty$-objects or, equivalently, that for every $\caD\in \Wald_\infty^\gd$ the real functor
$$\psi\colon \Wald_\infty^\gd \subset \D \Wald_\infty^\gd \xrightarrow{\Omega^{1,1} \circ \caS^{1,1}} \D\Wald_\infty^\gd \xrightarrow{\map_{\D\Wald_\infty^\gd}(\caD,-)} \Spc_\ast^{C_2}$$
lifts to grouplike $E_\infty$-objects in $\Spc_\ast^{C_2}$. 
The real functor $\psi$ can be written as 
$$ \Wald_\infty^\gd \xrightarrow{S} \rs\Wald_\infty^\gd \xrightarrow{\map_{\Wald_\infty^\gd}(\caD,-)_\ast} \rs\Spc_\ast \xrightarrow{|-| } \Spc^{C_2}_\ast \xrightarrow{\Omega^{1,1}}  \Spc^{C_2}_\ast $$
and so as 
$$ \Wald_\infty^\gd \xrightarrow{\Hom_{\Wald_\infty^\gd}(\caD,-)} \Wald_\infty^\gd \xrightarrow{S} \rs\Wald_\infty^\gd \xrightarrow{\iota}  \rs\Spc_\ast \xrightarrow{|-| } \Spc^{C_2}_\ast \xrightarrow{\Omega^{1,1}}  \Spc^{C_2}_\ast $$
by \cref{SSSSSS}. We conclude, as announced, that it is enough to check that the real functor (\ref{functor_thats_equiv_to_lift}) lifts to grouplike $E_\infty$-objects in $\Spc_\ast^{C_2}$.

As usual, we will verify that both the underlying and the $C_2$-fixed point versions do so. This is, it is enough to show that both functors below lift to grouplike $E_\infty$-spaces.
$$\alpha\colon \Wald_\infty^\gd \xrightarrow{S} \rs\Wald_\infty^\gd \xrightarrow{\iota}  \rs\Spc_\ast \xrightarrow{|-| } \Spc^{C_2}_\ast \xrightarrow{\Omega^{1,1}}  \Spc^{C_2}_\ast \xrightarrow{(-)^u} \Spc_\ast,$$
$$\beta\colon \Wald_\infty^\gd \xrightarrow{S} \rs\Wald_\infty^\gd \xrightarrow{\iota}  \rs\Spc_\ast \xrightarrow{|-| } \Spc^{C_2}_\ast \xrightarrow{\Omega^{1,1}}  \Spc^{C_2}_\ast \xrightarrow{(-)^{C_2}} \Spc_\ast$$
The first of these functors, $\alpha\colon\Wald_\infty^\gd\to \Spc_\ast$, factors as 
$$ \Wald_\infty^\gd \xrightarrow{\mathrm{forget}} \Exact_\infty \xrightarrow{S} \s\Exact_\infty \xrightarrow{(-)^\simeq_\ast } \s\Spc_\ast \xrightarrow{|-|} \Spc_\ast \xrightarrow{\Omega} \Spc_\ast$$ 
and so lifts to grouplike $E_\infty$-spaces since $\Omega$ does. Observe that even the functor 
$$\varphi\colon \Exact_\infty \xrightarrow{S} \s\Exact_\infty \xrightarrow{(-)^\simeq_\ast } \s\Spc_\ast \xrightarrow{|-|} \Spc_\ast$$ lifts to grouplike $E_\infty$-spaces. Indeed, the composition $ \Exact_\infty \xrightarrow{\varphi} \Spc_\ast \xrightarrow{\pi_0} \Set_\ast$ being equivalent to 
$$ \Exact_\infty \xrightarrow{S} \s\Exact_\infty \xrightarrow{ \pi_0 \circ (-)^\simeq } \s\Set_\ast \xrightarrow{|-|} \Set_\ast$$ trivially lifts to abelian groups, since the latter is constant with value the contractible space as for every $\caC \in \Exact_\infty$ the simplicial set $\pi_0(S(\caC)^\simeq)$ is in degree zero the set with one element. This remark, although extra at this point, will be of use later in the proof.

So it remains to check that $\beta\colon\Wald^\gd_\infty\to\Spc_\ast$ lifts to grouplike $E_\infty$-spaces too. Again, we will use a conveniently chosen factorization, although this time the conclusion will not follow as directly. It follows from \cref{lem:description_fixed_points_Omega11}
that the functor $\beta$ factors as 
$$ \Wald_\infty^\gd \xrightarrow{S} \rs\Wald_\infty^\gd \xrightarrow{\iota} \rs\Spc_\ast \xrightarrow{|-| } \Spc^{C_2}_\ast \xrightarrow{\rho} \Fun(\Delta^1, \Spc_\ast) \xrightarrow{\fib} \Spc_\ast, $$
where $\rho $ sends a pointed genuine $C_2$-space $X$ to the map $X^{C_2} \to X^u$ of pointed spaces and $\fib$ takes the fiber of a map of pointed spaces. Since the functor $\fib\colon\Fun(\Delta^1,\Spc_\ast)$ preserves products, it is enough to see that the functor 
$$\zeta\colon \Wald_\infty^\gd \xrightarrow{S} \rs\Wald_\infty^\gd \xrightarrow{\iota} \rs\Spc_\ast \xrightarrow{|-| } \Spc^{C_2}_\ast \xrightarrow{\rho} \Fun(\Delta^1, \Spc_\ast) $$ lifts to grouplike $E_\infty$-spaces. The composition $\Wald_\infty^\gd \xrightarrow{\zeta} \Fun(\Delta^1, \Spc_\ast) \xrightarrow{\ev_1} \Spc_\ast $
is equivalent to the functor $ \Wald_\infty^\gd \xrightarrow{\mathrm{forget}} \Exact_\infty \xrightarrow{\varphi} \Spc_\ast$ that lifts to grouplike $E_\infty$-spaces as we have seen before. The composition $\Wald_\infty^\gd \xrightarrow{\zeta} \Fun(\Delta^1, \Spc_\ast) \xrightarrow{\ev_0} \Spc_\ast $, on the other hand, factors as
$$ \Wald_\infty^\gd \xrightarrow{S} \rs\Wald_\infty^\gd \xrightarrow{\iota} \rs\Spc_\ast \xrightarrow{(-)^{C_2}} \s\Spc_\ast \xrightarrow{|-| } \Spc_\ast,$$
where the functor $(-)^{C_2}\colon \rs\Spc_\ast \xrightarrow{} \s\Spc_\ast$
sends a real simplicial space $X$ to the composition $ \Delta^\op \xrightarrow{e} (\Delta^{hC_2})^\op \xrightarrow{X} \Spc^{C_2} \xrightarrow{(-)^{C_2}} \Spc$
with $e\colon [n] \mapsto [n] \ast [n]^\op$ the edgewise subdivision.

So it is enough to check that the composition 
$$ \Wald_\infty^\gd \xrightarrow{S} \rs\Wald_\infty^\gd \xrightarrow{\iota} \rs\Spc_\ast \xrightarrow{(-)^{C_2}} \s\Spc_\ast \xrightarrow{|-| } \Spc_\ast \xrightarrow{\pi_0} \Set_\ast$$ being equivalent to the 
composition 
$$ \Wald_\infty^\gd \xrightarrow{S} \rs\Wald_\infty^\gd \xrightarrow{\iota} \rs\Spc_\ast \xrightarrow{(-)^{C_2}} \s\Spc_\ast \xrightarrow{\pi_0 } \sSet_\ast \xrightarrow{|-|} \Set_\ast$$
lifts to abelian groups.

We conclude the proof by noticing that for every $\caC \in \Wald_\infty^\gd$ and $Y\coloneqq (\iota \circ S(\caC))^{C_2} \in \s \Spc_\ast$, the geometric realization of the simplicial commutative monoid $\pi_0 \circ Y$ is an abelian group by \cref{lem:things_are_abelian_grps}.
\end{proof}

The reader may be interested in looking at \cite[Proposition 5.8]{spitzweck.gw}, where a similar proof is given in the context of stable $\infty$-categories.

\begin{lemma}\label{lem:things_are_abelian_grps}
Let $\caC$ be in $\Wald^\gd_\infty$ and consider the pointed simplicial space $Y\coloneqq (\iota \circ S(\caC))^{C_2}$, where the functor $(-)^{C_2}\colon \rs\Spc_\ast \xrightarrow{} \s\Spc_\ast$
sends a real simplicial space $X$ to the composition $$ \Delta^\op \xrightarrow{e} (\Delta^{hC_2})^\op \xrightarrow{X} \Spc^{C_2} \xrightarrow{(-)^{C_2}} \Spc$$
with $e\colon [n] \mapsto [n] \ast [n]^\op$ the edgewise subdivision.
 Then the geometric realization of the simplicial commutative monoid $\pi_0 \circ Y$ is an abelian group.
\end{lemma}
\begin{proof}

We consider the commutative square
\[
\begin{tikzcd}
\Wald_\infty^\gd\ar[r, "\iota"]\ar[d, "\caH"']      &\Spc_\ast^{C_2}\ar[d, "(-)^{C_2}"]\\
\Cat_{\infty\ast}\ar[r, "(-)^{\simeq}"]      &\Spc_\ast.
\end{tikzcd}
\]

The cartesian structure on $\caC$ gives rise to a symmetric monoidal structure on $\caH(\caC)$ that restricts to $\caH^\mathrm{nd}(\caC)$ and yields the $\mathbb{E}_\infty$-space structure on $(\caH^\mathrm{nd}(\caC))^\simeq \simeq \iota(\caC)^{C_2}.$ The tensor product
of any $\alpha, \beta \in \caH(\caC)$ lying over $A,B \in C$ is the sum in the $E_\infty$-space $\caH(\caC)_{A \oplus B}$ of the pullbacks of $\alpha \in \H(\caC)_A $ along the first projection $A \oplus B \to A$ and $ \beta \in \H(\caC)_B $ along the second projection $A \oplus B \to B$.

We want to see that the geometric realization of the simplicial commutative monoid $\pi_0 \circ Y$ is an abelian group. The geometric realization of $ \pi_0 \circ Y$ is the coequalizer $K$ of commutative monoids (which is computed in sets) of the source and target maps $ \pi_0(Y_1) \to \pi_0(Y_0)$.

Let $\varphi \in Y_0 \simeq \caH^\mathrm{nd}(\caC)^\simeq $ lying over $X \in \caC$.
The additive inverse $-\varphi$ of $\varphi$ in $\caH(\caC)_X \simeq \map_\caC(X, X^\dual)^{hC_2}$ belongs to $ \caH^\mathrm{nd}(\caC)^\simeq \simeq Y_0. $

We claim that $\varphi \in Y_0$ and $- \varphi \in Y_0$ become mutually inverse in the coequalizer $K$. Before diving in showing this, we recall that an object of $S(\caC)_2$
is a cofiber sequence $T\coloneqq(Z \to Z' \to Z'')$ in $\caC$, and that a non-degenerated object in $S(\caC)_2$ consists of both an object $T$ of $S(\caC)_2$ and a map $T \to T^\dual= (Z''^\dual \to Z'^\dual \to Z^\dual)$
of cofiber sequences in $\caC$, plus some coherence data.

We leave as an exercise to the reader to show that there exists a $\psi\in \H^\mathrm{nd}(S(\caC)_2)$
with underlying cofiber sequence
$$ X\xrightarrow{\delta} X \oplus X \xrightarrow{\varphi +(-\varphi) } X^\dual$$ in $\caC$, where $\delta$ is the diagonal map, and map of cofiber sequences
$$Z\coloneqq (X \xrightarrow{\delta} X \oplus X \xrightarrow{\varphi +(-\varphi) } X^\dual)
\ \to \ Z^\dual\coloneqq (X \xrightarrow{(\varphi, -\varphi)} X^\dual \oplus X^\dual \xrightarrow{ +} X^\dual)$$ in $\caC$ given by the identity on the extremes
and the map $\varphi \oplus (- \varphi)$ in the middle. Moreover, one can check that the map $[1] \to [2]$ in $\Delta^{hC_2}$ sending $0\mapsto 0$ and $1\mapsto 2$
yields a functor $$\caH^\mathrm{nd}(S(\caC)_2) \to \caH^\mathrm{nd}(S(\caC)_1) \simeq \caH^\mathrm{nd}(\caC)$$
that sends $\psi$ to the tensor product of $\varphi$ and $-\varphi$
in $\caH^\mathrm{nd}(\caC)$ or, in other words, their sum in the $E_\infty$-space
$Y_0 \simeq \caH^\mathrm{nd}(\caC)^\simeq.$

Finally, the map $[3] \mapsto [2]$ in $\Delta^{hC_2}$ that sends $0 \mapsto 0$, $1,2 \mapsto 1$ and $3 \mapsto 2$ gives rise to a map $ \H^\mathrm{nd}(S(\caC)_2)^\simeq \to  \H^\mathrm{nd}(S(\caC)_3)^\simeq \simeq Y_1$ that sends $\psi$ to an object $\psi'.$ The source and target maps $ Y_1 \to Y_0$ send $\psi'$, respectively, to the sum of $\varphi$ and $-\varphi$ in $Y_0$ and $0$.
\end{proof}

\begin{lemma}\label{lem:to_lift_KR}
Let $X$  be in $\Sp_{S^{1,1}}(\Grp_{E_\infty}( \D\Wald_\infty^\gd) ) \subset \Pre\Sp_{S^{1,1}}(\Sp( \D\Wald_\infty^\gd) ).$ Assume that for every $Z \in \Wald_\infty^\gd$ the image of $X$ under the real functor $$ \Sp_{S^{1,1}}(\Grp_{E_\infty}( \D\Wald_\infty^\gd) ) \xrightarrow{\map_{\D\Wald_\infty^\gd}(Z,-)} \Sp_{S^{1,1}}(\Grp_{E_\infty}(\Spc^{C_2})) \xrightarrow{(-)^u} $$$$ \Sp(\Grp_{E_\infty}(\Spc)) \simeq \Sp,$$ which we denote by $Y$, has the property that for every $n \geq 1$ the space $Y_n$ is connected. Then $X$ is a connective genuine $C_2$-spectrum in $\D\Wald_\infty^\gd$.
\end{lemma}

\begin{proof}
By \cref{char_genuine_spectra} we know that 
$$\Sp_{S^{1,1}}(\Sp( \D\Wald_\infty^\gd) ) \simeq \Sp_{S^{2,1}}( \D\Wald_\infty^\gd).$$
The proof will consist on showing that under the hypotheses, $X$ is an object in the left hand side $\infty$-category.

We want to see that for every $n \geq 0$ the canonical map
$X_n \to \Omega^{1,1}(X_{n+1})$ in $\Sp( \D\Wald_\infty^\gd) $ is an equivalence, where $\Omega^{1,1}$ is formed in $\Sp(\D\Wald_\infty^\gd)$ via the embedding 
$$\Grp_{E_\infty}(\D\Wald_\infty^\gd)\hookrightarrow \Sp(\D\Wald_\infty^\gd)$$ 
of \cref{rmk:embedding_grplike_objects_in_Sp}. This is equivalent to ask that $\Omega^{1,1}(X_{n+1}) \in \Sp( \D\Wald_\infty^\gd)$ belong to $\Grp_{E_\infty}( \D\Wald_\infty^\gd).$
To show this, it suffices to check that for every $Z \in \Wald_\infty^\gd$ the
$C_2$-spectrum 
$$\D\Wald_\infty^\gd(Z,\Omega^{1,1}(X_{n+1})) \simeq \Omega^{1,1}(\D\Wald_\infty^\gd(Z,X_{n+1})) \in \Sp(\Spc^{C_2})$$ belongs to $\Grp_{E_\infty}(\Spc^{C_2})$. This, in turn, is equivalent to prove that the spectra 
$$ \Omega^{1,1}(\D\Wald_\infty^\gd(Z,X_{n+1}))^u \simeq \Omega(\D\Wald_\infty^\gd(Z,X_{n+1})^u), $$$$ \Omega^{1,1}(\D\Wald_\infty^\gd(Z,X_{n+1}))^{C_2} \simeq \fib(\D\Wald_\infty^\gd(Z,X_{n+1})^{C_2} \to \D\Wald_\infty^\gd(Z,X_{n+1})^u)$$ belong to $\Grp_{E_\infty}(\Spc)$ \textemdash or, equivalently, that are connective.

By assumption we know that $\D\Wald_\infty^\gd(Z,X_{n+1})^u$ is connected, and since  the fiber of a map of spectra from a connective spectrum to a $0$-connective spectrum is again connective, we conclude the proof.
\end{proof}

\begin{theorem}\label{thm:lifting_Omega11_S11}
The real functor $\Omega^{1,1} \circ \caS^{1,1}\colon \D\Wald_\infty^\gd \to \D\Wald_\infty^\gd $ lifts to a real functor $$\D\Wald_\infty^\gd \to \Sp^{C_2}_{\geq0}(\D\Wald_\infty^\gd).$$
\end{theorem}

\begin{proof}
The functor $\Omega^{1,1}\circ\caS^{1,1}$ is an additive theory by \cref{dfghjlk}, and therefore genuine excisive (see \cref{ex:gen_exc_S11_exc}). It follows from \cref{fhjkkkl} that the functor $$\Omega^{1,1} \circ \caS^{1,1}\colon \D\Wald_\infty^\gd \to \D\Wald_\infty^\gd$$ uniquely lifts to an additive theory $ \D\Wald_\infty^\gd \to \Grp_{E_\infty}(\D\Wald_\infty^\gd) $ and then by \cref{spec} uniquely lifts to a genuine excisive functor
$$\psi\colon \D\Wald_\infty^\gd \to \Sp_{S^{1,1}}(\Grp_{E_\infty}(\D\Wald_\infty^\gd)).$$

We check now that, for every $\caC\in\D\Wald^\gd_\infty$, its image $\psi(\caC)$ verifies the hypotheses of \cref{lem:to_lift_KR}. For every $\caC\in \D\Wald_\infty^\gd$ the real functor
$$\D \Wald_\infty^\gd \xrightarrow{\Omega^{1,1} \circ \caS^{1,1}} \D\Wald_\infty^\gd \xrightarrow{\map_{\D\Wald_\infty^\gd}(\caC,-)} \Spc_\ast^{C_2}$$
can be written as 
$$ \D\Wald_\infty^\gd \xrightarrow{\Hom_{\D\Wald_\infty^\gd}(\caC,-)} \D\Wald_\infty^\gd \xrightarrow{S} \rs\D\Wald_\infty^\gd \xrightarrow{\iota}  \rs\Spc_\ast \xrightarrow{|-| } \Spc^{C_2}_\ast \xrightarrow{\Omega^{1,1}}  \Spc^{C_2}_\ast $$
by \cref{SSSSSS}. 
Hence the functor
\begin{equation}\label{first_functor_to_lift}
\D \Wald_\infty^\gd \xrightarrow{\Omega^{1,1} \circ \caS^{1,1}} \D\Wald_\infty^\gd \xrightarrow{\map_{\D\Wald_\infty^\gd}(\caC,-)} \Spc_\ast^{C_2} \xrightarrow{u} \Spc_\ast
\end{equation}
can be written as 
\begin{equation}\label{the_one_to_show_it_is_connective}
\D\Wald_\infty^\gd \xrightarrow{\Hom_{\D\Wald_\infty^\gd}(\caC,-)} \D\Wald_\infty^\gd \xrightarrow{u} \D\Wald_\infty \xrightarrow{S} \s\D\Wald_\infty \xrightarrow{(-)^\simeq}  \s\Spc_\ast \xrightarrow{|-| } \Spc_\ast \xrightarrow{\Omega}  \Spc_\ast.
\end{equation}
The last part of (\ref{the_one_to_show_it_is_connective}), the functor 
$$\D\Wald_\infty \xrightarrow{S} \s\D\Wald_\infty \xrightarrow{(-)^\simeq}  \s\Spc_\ast \xrightarrow{|-| } \Spc_\ast \xrightarrow{\Omega}  \Spc_\ast,$$
is excisive and so uniquely lifts to the $K$-theory $K \colon \D\Wald_\infty \to \Sp,$ which is known to assigne a connective spectrum.

Now, the composition
\begin{equation}\label{the_one}
\D\Wald_\infty^\gd \xrightarrow{\psi} \Sp_{S^{1,1}}(\Grp_{E_\infty}(\D\Wald_\infty^\gd))  \xrightarrow{\map_{\D\Wald_\infty^\gd}(D,-)} \Sp_{S^{1,1}}(\Grp_{E_\infty}(\Spc^{C_2})) \xrightarrow{u} \Sp(\Grp_{E_\infty}(\Spc))
\end{equation}
uniquely lifts the functor (\ref{first_functor_to_lift}), whereas the  composition 
$$\D\Wald_\infty^\gd \xrightarrow{\Hom_{\D\Wald_\infty^\gd}(\caC,-)} \D\Wald_\infty^\gd \xrightarrow{u} \D\Wald_\infty \xrightarrow{\alpha}\Sp$$ lifts the expression (\ref{the_one_to_show_it_is_connective}). Since both are liftings of the same functor, just different expressions of it, and $K$ takes values in connective spectra, we conclude that the composition (\ref{the_one}) also takes values in connective spectra.

Therefore, we conclude that the functor $\psi$ is the lifting of $\Omega^{1,1}\circ\caS^{1,1}$ that we sought by applying \cref{lem:to_lift_KR}.
\end{proof}

\begin{corollary}\label{fghjkmno}
The real functor $\mathrm{KR}\colon \D\Wald_\infty^\gd \to \Spc^{C_2}$ lifts to a real functor $$\mathrm{KR}\colon \D\Wald_\infty^\gd \to \Sp_{\geq0}^{C_2}.$$
\end{corollary}

\begin{proof}
Observe first that since $\iota\colon\D\Wald_\infty^\gd\to\Spc_\ast^{C_2}$ preserves cotensors, it commutes with the functor $\Omega^{1,1}$, and therefore we can write $\KR\coloneqq\add(\iota)=\caS^{1,1}\circ\Omega^{1,1}\circ\iota$. On one hand, the functor $\iota$ induces a functor between connective genuine $C_2$-spectra
$$\Sp^{C_2}_{\geq0}(\D\Wald_\infty^\gd) \to  \Sp^{C_2}_{\geq0}(\Spc^{C_2}_\ast),$$
for it preserves finite limits and cotensors. On the other hand, by \cref{thm:lifting_Omega11_S11}, the composite $\Omega^{1,1} \circ \caS^{1,1}$ lifts to a real functor 
$$\D\Wald_\infty^\gd \to \Sp^{C_2}_{\geq0}(\D\Wald_\infty^\gd).$$

The composition of these functors is the lift sought, $\mathrm{KR}\colon \D\Wald_\infty^\gd \to \Spc^{C_2}.$
\end{proof}

\part{Appendices}\label{part:app}
\appendix

\section{Enriched infinity category theory}\label{Appx:enriched}
In this appendix we aim to present references, or proofs when that is more convenient, to results of enriched $\infty$-category theory that we use thorough the paper. This is by no means an exhaustive nor didactic introduction to the topic, but a guide to find such results in the literature. 

\subsection{The $\infty$-category of $\mathcal{V}$-enriched $\infty$-categories}

We begin with a brief presentation of the definition of enriched $\infty$-categories given by the first author in \cite{heine-enriched1, heine-enriched2}, which is very similar to that presented by Gepner and Haugseng in \cite{GEPNER2015575} but it contains a proof of an enriched Yoneda lemma as does Hinich's work on the topic \cite{hinich2020yoneda}.

We must warn the reader of a slight change in nomenclature with respect to \cite{GEPNER2015575}. We are calling $\caV$-enriched $\infty$-precategories, or $\caV$-precategories, to what Gepner and Haugseng call $\caV$-enriched $\infty$-categories (or categorical algebras in $\caV$); in the same line, we call $\caV$-enriched $\infty$-categories (or simply $\caV$-categories) to what is \emph{complete} $\caV$-enriched $\infty$-categories in \cite{GEPNER2015575}. We justify this choice in \cref{justification:shift_of_nomenclature}.

We start by defining a many object version of the final non-symmetric $\infty$-operad $\Ass\coloneqq\Delta^\op$, where we use the non-identity involution $(-)^\op$ on $\Ass=\Delta^\op.$ See \cite[Definition 4.1.1]{GEPNER2015575} for a similar description \textemdash we provide in \cref{def_of_Venriched_are_the_same} a clarification on how these relate when defining $\caV$-enriched $\infty$-precategories. We will use ``(generalized) $\infty$-operad'' for non-symmetric (generalized) $\infty$-operad, and say explicitly when it is symmetric. 

Before proceeding, we direct the reader to \cite[Definition 2.4.1]{GEPNER2015575} for the definition of generalized $\infty$-operad, and \cite[Definition 2.4.2]{GEPNER2015575} for the morphisms between them. We also recall that given two generalized $\infty$-operads $\mathcal{O}$ and $\mathcal{O'}$, we denote by $\Alg_{\mathcal{O}}(\mathcal{O'})$ the full subcategory of $\Fun_{\Delta^\op}(\mathcal{O}, \mathcal{O'})$ spanned by the maps of generalized $\infty$-operads.

\begin{construction}
For any space $X$ the forgetful functor $\nu\colon \Delta \to \Set$ that sends $[n] \mapsto \{0,\dots, n \}$ gives rise to a functor $$ \Ass \xrightarrow{(-)^\op} \Ass \xrightarrow{\nu^\op } \Set^\op  \xrightarrow{\Fun(-,X)} \Spc $$ classifying a left fibration 
$ \Ass_{X} \to \Ass$ that is a generalized non-symmetric $\infty$-operad.
If $X$ is contractible, the functor $\Ass_X \to \Ass$ is an equivalence.
\end{construction}

We are now ready to define $\caV$-enriched $\infty$-precategory.

\begin{definition}\label{GH}
Given a space $X$ and an $\infty$-operad $\caV^\otimes$, an $\infty$-precategory enriched in
$\caV$ (or $\caV$-precategory) with space of objects $X$ is a map $\Ass_{X} \to \caV^\otimes$ of generalized $\infty$-operads.
\end{definition}

The assignment $X \mapsto \Alg_{\Ass_X}(\caV)$ defines a contravariant functor from small spaces to large $\infty$-categories classifying a cartesian fibration $\psi\colon \PreCat_\infty^\caV  \to \Spc.$ The $\infty$-category $\PreCat_\infty^\caV$ is the $\infty$-categry of $\caV$-precategories.

Given a space $X,$ an $\infty$-operad $\caV^\otimes$ and a $\caV$-precategory 
$\caC\colon \Ass_{X} \to \caV^\otimes$ with space of objects $X$, we call the induced functor
$(\Ass_\X)_{[1]} \simeq X \times X \to \caV$ the underlying graph of $\caC.$

\begin{remark}\label{def_of_Venriched_are_the_same}
Definition \ref{GH} is a slight variant from the original definition \cite[Definition 2.4.4]{GEPNER2015575}. Gepner and Haugseng define homotopy-coherent enrichment in an $\infty$-operad $\caV^\otimes$ (under the name of categorical algebra in $\caV$ or $\caV$-enriched $\infty$-category) with space of objects $X$, as a map of generalized $\infty$-operads $\Ass_{X} \to (\caV^\otimes)^\rev$, where $(\caV^\otimes)^\rev$ is the reversed $\infty$-operad structure defined as the composition $\caV^\otimes \to \Ass \xrightarrow{(-)^\op} \Ass.$

We have decided to present it this way because of how composition works. Given objects $A, B, C \in X$ the composition maps of a $\caV$-precategory
$\caC$ with space of objects $X$ are multi-morphisms
$$\mathsf{Mor}_\caC(B, C), \mathsf{Mor}_\caC(A, B) \to  \mathsf{Mor}_\caC(A, C)$$ in our definition whereas they are maps
$$\mathsf{Mor}_\caC(A, B), \mathsf{Mor}_\caC(B, C) \to  \mathsf{Mor}_\caC(A, C)$$ in 
the definition of Gepner and Haugseng.
\end{remark}

\begin{remark}\label{rmk:base_change}
There is a functor $\Op_\infty \to \widehat{\Cat}_{\infty/\Spc}^\Cart$ that maps  $\caV \mapsto \PreCat _\infty^\caV$
given as the composition

\[
\begin{tikzcd}[row sep=tiny, column sep=small]
\Op_\infty\ar[r]        &\Fun(\Spc^\op, \Cat_\infty)\ar[r, "\simeq"]       &\widehat{\Cat}_{\infty/\Spc}^\Cart\\
\caV\ar[r, mapsto]         & (X\mapsto \Alg_{\Ass_X}(\caV))      &
\end{tikzcd}
\]
\end{remark}

\begin{remark}\label{adj_pass_to_precat}
An adjunction $F\colon \caV \rightleftarrows \mathcal{W} \colon G$ of $\infty$-operads canonically induces an adjunction $\PreCat_\infty^\caV \rightleftarrows \PreCat_\infty^\mathcal{W}.$
\end{remark}

\begin{construction}\label{const:op_V_precategory}
We will see in this construction that, given a non symmetric $\infty$-operad $\caV$, the $\infty$-category $ \PreCat_\infty^\caV$ carries a canonical involution denoted by $(-)^\op$ sending a $\caV$-enriched $\infty$-precategory to its opposite
$\caV$-enriched $\infty$-precategory $\caC^\op$ with $\map_{\caC^\op}(A, B) \simeq \map_{\caC}(B,A)$ for any $A,B \in \caC^\simeq.$

There is a canonical equivalence
$ \nu \circ (-)^\op \simeq \nu$ that yields for any space $X$ a canonical equivalence
$$\Ass_{X}^\rev = ((-)^\op)^\ast(\Ass_X)\simeq \Ass_X$$ over $\Ass$ that yields on the fiber over
$[1] \in \Ass$ the equivalence $X \times X \simeq X \times X $ switching the factors.

Then any $\caV$-enriched $\infty$-precategory $\Ass_X \to \caV^\otimes$ with space of objects $X$ corresponds to a $\caV^\rev$-enriched $\infty$-precategory $\Ass_X \simeq \Ass_X^\rev \to (\caV^\otimes)^\rev$ with with the same space of objects but reversed graph, which we call the opposite $\caV$-enriched $\infty$-precategory.
\end{construction}

We now embark to define $\caV$-categories, for which we continue following \cite{GEPNER2015575}.

When $\caV=\Spc$, we call $\caV$-enriched $\infty$-precategories simply $\infty$-precategories.
Now, consider $\ast$ to be the final $\infty$-precategory and $\mathcal{J} \to \ast$ a cartesian lift in $\PreCat_\infty^\Spc$ of the map of spaces $\{0,1\} \to [0]$. So $\mathcal{J}$ is the $\infty$-precategory with contractible mapping spaces and space of objects a two-element set.

Note that for any $\infty$-precategory $\caD$ with small space of objects $X$, the canonical map $ \PreCat_\infty(\ast, \caD) \to \Spc([0],X) \simeq X$ is an equivalence since the fiber over any $Z \in X$ is the contractible space $ \Alg(\Spc)(\ast, \caD(Z,Z)).$  

\begin{definition}\label{comp}
We call an $\infty$-precategory $\caD$ with small space of objects $X$ complete if the canonical map
$$\theta\colon X \simeq \PreCat^\Spc_\infty(\ast, \caD) \to \PreCat^\Spc_\infty(\mathcal{J}, \caD)$$ is an equivalence
\end{definition}

Given a $\caV$-category $\caC$ with space of objects $X$, we write $\caC^\simeq $ for $X$.

We think of $\PreCat^\Spc_\infty(\mathcal{J}, \caD) $ as the space of equivalences in $\caD$ and so as the genuine space of objects of $\caD$. Therefore $\theta$ may be thought of as a comparison map between the space of objects and the genuine space of objects, and a complete $\infty$-precategory one whose genuine space of objects and space of objects coincide.

We will reduce the definition of $\caV$-enriched $\infty$-category to a requirement on its underlying $\infty$-precategory. For this, we need to define the functor 
$$\PreCat^\caV_\infty\to\PreCat^\Spc_\infty$$
induced by the functor $\caV\to\Spc$ that sends $Z\mapsto \caV(1_\caV,Z)$.

\begin{definition}\label{def:enriched_V_category}
Let $\caV$ be a monoidal $\infty$-category. A $\caV$-enriched $\infty$-category is a $\caV$-enriched $\infty$-precategory such that its image under the functor $\PreCat^\caV_\infty\to\PreCat^\Spc_\infty$ is a complete $\infty$-precategory.
\end{definition}

We denote by $\Cat_\infty^\caV$ the full subcategory of $\PreCat_\infty^\caV$ spanned by the $\caV$-categories. We call a morphism in $ \Cat_\infty^\caV$ a $\caV$-enriched functor.

\begin{remark}\label{justification:shift_of_nomenclature}
Note that the extra requirement on a $\caV$-precategory to be a $\caV$-category solely relies on the underlying $\infty$-category with the canonical enrichment on $\Spc$, and it does not involve $\caV$ whatsoever. Moreover, there is an embedding $\Cat_\infty \subset \PreCat_\infty^\Spc$ with essential image the complete $\infty$-precategories. 
\end{remark}

\begin{definition}
We call a $\caV$-enriched functor $\caC \to \caD$ fully faithful, or embedding, if it is cartesian with respect to the canonical functor $\Cat^\caV_\infty\to\Spc$.
\end{definition}

\begin{remark}\label{rmk:Cat_V_reflective_subcat}
The full subcategory $\Cat_\infty^\caV \subset \PreCat_\infty^\caV$ is reflective, where a morphism $\caC \to \caD$ with $\caD$ complete is a local equivalence if and only if it is
cartesian over $\Spc$ and induces an essentially surjective map on spaces of objects.
\end{remark}

\begin{remark}\label{rmk:base_change_Vcat}
Importantly, \cref{rmk:Cat_V_reflective_subcat} pairs with \cref{rmk:base_change} to prove that we also count with base change for $\caV$-categories. 
\end{remark}

\begin{remark}\label{adj_pass_to_Vcat}
By \cref{adj_pass_to_precat} we know that an adjunction $F\colon \caV \rightleftarrows \mathcal{W} \colon G$ of $\infty$-operads induces an adjunction $\PreCat^\caV_\infty\rightleftarrows\PreCat^\caW_\infty$. The left adjoint preserves local equivalences for the localization of \cref{rmk:Cat_V_reflective_subcat} and so induces an adjunction $\Cat_\infty^\caV \rightleftarrows \Cat_\infty^\mathcal{W}$,
where the right adjoint is the restriction of the canonical functor $\PreCat^\caW_\infty\to\PreCat^\caV_\infty.$
\end{remark}

The following remark guarantees that given a $\caV$-category $\caC$ with space of objects $X$, for any embedding $Y\subset X$ one can consider the full $\caV$-subcategory of $\caC$ with space of objects $Y$.

\begin{remark}
Given a $\caV$-enriched $\infty$-category $\caC$ and an embedding $\iota\colon Y \subset \caC^\simeq$
the pullback $\iota^\ast(\caC) \in \PreCat_\infty^\caV$ is a $\caV$-enriched $\infty$-category. We call $\iota^\ast(\caC)$ the full $\caV$-subcategory spanned by the objects of $Y.$
\end{remark}

We culminate here the introduction of the elements the rest of the appendix will refer to. It is worth noticing that the $\infty$-categories of $\caV$-enriched categories obtained in all three models are equivalent.

\begin{remark}
The full subcategory inclusion $\Op_\infty \subset \Op_\infty^{\mathrm{gen}}$ of $\infty$-operads into generalized $\infty$-operads admits a left adjoint denoted by $\L$ and called operadic localization. So a $\caV$-precategory with space of objects $X$ is likewise a map of $\infty$-operads $\L(\Ass_X) \to \caV^\otimes.$

The abstract $\infty$-operad $\L(\Ass_X) $ can be explicitely constructed:
Gepner-Haugseng construct a presentation of $\L(\Ass_X) $ as a simplicial operad
(\cite[Definition 4.2.2 and Corollary 4.2.8]{GEPNER2015575})
that is identified by Macpherson \cite{macpherson2019operad} with a construction
of Hinich (\cite[3.2.11]{hinich2020yoneda}).
We have decided to work with $\Ass_X$ instead of $\L(\Ass_X) $ 
as $\Ass_X$ seems a simpler object, is much easier to define and has the  technically convenient advantage that the functor $\Ass_X \to \Ass$ is a left fibration.
\end{remark}
\subsection{Pointed versus unpointed}
Let $\caV$ be a (symmetric) monoidal $\infty$-category compatible with small colimits
that admits a final object, and $\caV_\ast \subset \Fun([1],\caV)$ the full subcategory spanned by the morphisms whose source is a final object. One can observe that $\caV_\ast \subset \Fun([1],\caV)$ is a (symmetric) monoidal localization, where the left adjoint sends a morphism $A \to B$ in $\caV$ to the pushout $B/A$ of $A\to B$ along $A \to \ast$.
Consequently $\caV_\ast$ carries a canonical (symmetric) monoidal structure such that the embedding $\caV_\ast \subset \Fun([1],\caV)$ is lax (symmetric) monoidal.

\begin{proposition}\label{Venr:pointed_vs_not_pointed}
Let $\caV$ be a monoidal $\infty$-category compatible with small colimits
that admits a final object.
Then there is a canonical embedding 
$$\Cat_\infty^{\caV_\ast} \subset (\Cat_\infty^{\caV})_\ast.$$
\end{proposition}

\begin{proof}
By \cref{Venr:Fun_enriched_is_Fun_to_enriched} there is a canonical embedding 
$$\Cat_\infty^{\Fun([1],\caV)} \subset \Fun([1], \Cat_\infty^{\caV})$$
over $$\Cat_\infty^{\Fun(\{0\},\caV)} \simeq \Fun(\{0\}, \Cat_\infty^{\caV}).$$
Thus the composition $$\Cat_\infty^{\caV_\ast} \subset\Cat_\infty^{\Fun([1],\caV)} \subset \Fun([1], \Cat_\infty^{\caV})$$
is induced by the lax monoidal functor
$\caV_\ast \to \ast \to \caV$, where the last functor selects the final object.
Since $\Cat_\infty^\ast$ is contractible, 
we conclude that the embedding above restricts to the claimed embedding
$$\Cat_\infty^{\caV_\ast} \subset  (\Cat_\infty^{\caV})_\ast.$$
\end{proof}

\subsection{Enriched functor $\infty$-categories}
Let us consider $\caV$ to be a presentable symmetric  monoidal $\infty$-category. In this case $\Cat^\caV_\infty$ is a presentable symmetric monoidal $\infty$-category as well (see \cite[Corollary 5.6.16]{GEPNER2015575}), and the assignment $\caV \mapsto \Cat^\caV $ defines an endofunctor of the very large $\infty$-category of presentable symmetric monoidal $\infty$-categories and (lax) symmetric monoidal functors. Especially $\Cat^\caV_\infty$ is a closed symmetric monoidal $\infty$-category, whose internal hom we denote by $ \Fun^\caV(-,-).$

Given $\caV$-categories $\caC$ and $\caD$, we will call a morphism in $\Fun^\caV(\caC,\caD)$ a $\caV$-transformation.

\subsection{Enriched slice $\infty$-categories}
 \begin{definition}
Let $\caC$ be a $\caV$-category and $Y $ an object of $\caC$.
Then we define the $\caV$-enriched slice category of $\caC$ over $Y$ as
$$\caC_{/Y}\coloneqq\{Y\} \times_{\caC^{\{1\}}} \caC^{[1]}.$$
\end{definition}

\begin{proposition}\label{Venr:u_of_slice category}
Let $\caV$ and $\mathcal{W}$ be presentable symmetric monoidal $\infty$-categories, and let $G\colon\caV\to\mathcal{W}$ be a functor that admits a symmetric monoidal left adjoint $F$, let $\caC$ be a $\caV$-category and $Y$ an object of $\caC$. Then there is an equivalence $$G_\ast(\caC_{/Y})\simeq G_\ast(\caC)_{/Y}.$$
\end{proposition}

\begin{proof}
The adjunction $F\colon\caW \rightleftarrows \caV\colon G$
induces an adjunction $F_\ast \colon \Cat_\infty^\caV \rightleftarrows \Cat_\infty^\caW \colon G_\ast$ between their respective $\infty$-categories of enriched $\infty$-categories.
Note that the symmetric monoidal functor $F$ is canonically a symmetric monoidal functor under $\Spc$, since this is the initial presentable symmetric monoidal $\infty$-category. Then $F_\ast$ is canonically a symmetric monoidal functor
under $\Cat_\infty$ and so a $\Cat_\infty$-linear functor.
This implies by adjointness that its right adjoint $G_\ast$ is a lax $\Cat_\infty$-linear functor and preserves cotensors with
$\Cat_\infty$.

So there is a canonical equivalence $$G_\ast(\caC_{/Y})= G_\ast(\{Y\} \times_{\caC^{\{1\}}} \caC^{[1]}) \simeq \{Y\} \times_{G_\ast(\caC^{\{1\}})} G_\ast(\caC^{[1]}) \simeq $$
$$\{Y\}\times_{G_\ast(\caC)^{\{1\}}} G_\ast(\caC)^{[1]}\eqqcolon G_\ast(\caC)_{/Y}.$$
\end{proof}

\subsection{Enriched $\infty$-presheaves}

We now consider $\caV$ to be a presentable monoidal $\infty$-category. In this section we stufy the $\caV$-enriched $\infty$-category of $\caV$-enriched presheaves following \cite{heine-enriched1}, which was also studied in \cite{hinich2020yoneda} and compared with other models by \cite{heine-enriched5}.

For this section we will need a special kind of colimits, that we now define. 

\begin{definition}
Let $\caC$ be a $\caV$-category. The tensor of an object $V \in \caV$ with an object $X \in \caC$
is an object $V \otimes X$ and a morphism
$V \to Mor_\caC(X,V \otimes X)$ such that the canonical morphism
$Mor_\caC(V \otimes X,Y) \to Mor_\caV(V,Mor_\caC(X,Y))$
in $\caV$ is an equivalence.
\end{definition}

\begin{theorem}\label{Venr:thm:defining_presheaves}
Let $\caV$ be a presentable monoidal $\infty$-category, $\caC$ a small $\caV$-category and $\caD$ a $\caV$-category with small colimits and tensors with $\caV$. Then there exists a $\caV$-category $\caP_\caV(\caC)$, which admits small colimits and tensors, together with a  $\caV$-functor $\iota\colon \caC \to \caP_\caV(\caC)$
that induces an equivalence of $\infty$-categories
$$\Fun_\caV^L(\caP_\caV(\caC), \caD) \simeq \Fun_\caV(\caC,\caD),$$
where the $L$ in the left hand term indicates that we consider $\caV$-functors that preserve small colimits and tensors with $\caV.$

\end{theorem}
\begin{proof}
This follows directly from \cite[Theorem 4.2]{heine-enriched1}. 
\end{proof}

\begin{definition}\label{Venr:def:presheaves}
We call $\caP_\caV(\caC)$ the $\infty$-category of $\caV$-enriched presheaves on $\caC$.
\end{definition}

\begin{definition}\label{Venr:def:Yoneda_embedding}
We call the universal $\caV$-functor $\iota\colon \caC \to \caP_\caV(\caC)$ the $\caV$-enriched Yoneda embedding.
\end{definition}

Note that by \cite[Proposition 4.15]{heine-enriched1} the $\caV$-enriched Yoneda embedding is indeed an embedding. 

\begin{remark}
When $\caV$ is a presentable symmetric monoidal $\infty$-category, there is a canonical equivalence $$\caP_\caV(\caC) \simeq \Fun^\caV(\caC^\op,\caV).$$
This follows from \cite[Theorem 5.20]{heine-enriched1} for the particular case of having a symmetric monoidal structure on $\caV$.
\end{remark}

The next is \cite[Corollary 4.85.]{heine-enriched4}:

\begin{proposition}\label{Venr:descr_right_adjoint}
Let $\caV$ be a presentable monoidal $\infty$-category, $\caC$ a small $\caV$-category and $\caD$ a $\caV$-category with small colimits and tensors with $\caV$. Given a $\caV$-functor $F\colon  \caC \to \caD$
its extension $\caP_\caV(\caC) \to \caD$ is left adjoint to
the composition of $\caV$-functors $\caD \to \caP_\caV(\caD) \xrightarrow{F^\ast}\caP_\caV(\caC)$.
\end{proposition}


For $\caD=\caP_\caV(\caC)$ and $F$ the Yoneda-embedding, it follows from \cref{Venr:descr_right_adjoint} a non functorial version of the enriched Yoneda lemma; for a complete version see, for example, and \cite[Lemma 5.2.4]{heine-enriched1}. See also \cite[Section 6.2.7]{hinich2020yoneda} for this result proven for the Yoneda embedding defined therein.

\begin{lemma}[Yoneda lemma]\label{Yoneda_lemma}
Let $\caV$ be a presentable monoidal $\infty$-category, $\caC$ a small $\caV$-category, $Z \in \caC$ and $H \in \caP_\caV(\caC)$.
There is a canonical equivalence in $\caV$ as below.
$$\mathsf{Mor}_{ \caP_\caV(\caC)}(y(Z), H) \simeq H(Z)$$
\end{lemma}

\begin{proposition}\label{Venr:Yoneda_embedding_preserves_limits_and_cotensors}
Let $\caC$ be a $\caV$-category. Then the Yoneda embedding $\caC \subset \caP_\caV(\caC)$ preserves small limits and cotensors with $\caV$.
\end{proposition}

\begin{proof}
This is \cite[Corollary 3.68.]{heine-enriched3}.
\end{proof}

Given a lax monoidal functor $F\colon\caV \to \caW$
and a small $\caV$-category $\caC$
there is a left adjoint $\caV$-functor (in the sense of \cref{Venr:def:left_adjoint})
$$\caP_\caV(\caC) \to F^\ast(\caP_\caW(F_\ast(\caC))) $$
extending the $\caV$-functor
$\caC \to F^\ast(\caP_\caW(F_\ast(\caC))).$

\begin{remark}For any small $\caV$-category $\caC$. The forgetful functor above commutes with the Yoneda embedding $\caC\subset\caP_\caV(\caC)$.
\end{remark}

\subsection{Enriched adjoint functors}
\begin{definition}\label{Venr:def:left_adjoint}
Given $\caV$-functors $F\colon\caC\to\caD$ and $G\colon\caD \to \caC$, we say that $F$ is left adjoint to $G$, or that $G$ is right adjoint to $F$, if there is a $\caV$-transformation $\varepsilon\colon FG \to \id$ such that for any $X \in \caC, Y \in \caD$ the morphism
$$\mathsf{Mor}_\caC(X,G(Y)) \to \mathsf{Mor}_\caD(F(X),FG(Y)) \to \mathsf{Mor}_\caD(F(X),Y) $$ is an equivalence in $\caV$.
\end{definition}

\begin{example}\label{Venr:ev_admits_adjoint}
Let $\caC$ and $\caD$ be $\caV$-categories, such that $\caD$ admits small colimits and tensors with $\caV$. Then the functor evaluation at any object $\ev\colon\Fun^\caV(\caC,\caD)$ admits a left adjoint by \cite[Proposition 4.41.]{heine-enriched4}.
\end{example}

\begin{remark}\label{Venr:F_equiv_iff_left_adj_and_underlying_is_equiv}
A $\caV$-functor $F\colon \caC \to \caD$ is an equivalence if and only if $F$ admits a right adjoint and the underlying functor of $F$ is an equivalence.
\end{remark}

The next is \cite[Remark 2.75.]{heine-enriched4}:

\begin{proposition}
\label{Venr:F_admits_right_adjoint}
A $\caV$-functor $F\colon\caC\to\caD$ admits a right adjoint
if and only if for any $Y \in \caD$ there is a $Z \in \caC$
and a morphism $F(Z) \to Y $ such that for any $X \in \caC$ the morphism
$$\mathsf{Mor}_\caC(X,Z) \to \mathsf{Mor}_\caD(F(X),F(Z)) \to \mathsf{Mor}_\caD(F(X),Y) $$ is an equivalence in $\caV$.
\end{proposition}

The following result allows us to upgrade adjunctions to enriched adjunctions. 

\begin{corollary}\label{Venr:upgrade_adjunctions}
Let $\caC, \caD$ be $\caV$-categories that admit tensors and $F\colon\caC\to\caD$ a $\caV$-enriched functor, whose underlying functor admits a right adjoint. If $F$ preserves tensors with $\caV$, then it admits a $\caV$-enriched right adjoint.
\end{corollary}

\begin{proposition}\label{Venr:left_adjs_of_embeddings}
Let $\caV$ be a monoidal $\infty$-category, $\caC$ and $\caD$ $\caV$-categories admitting cotensors with $\caV$, and a $\caV$-enriched embedding $\iota\colon\caD \to \caC$ that preserves cotensors. Consider also an essentially surjective $\caV$-functor $L\colon \caC \to \caD$ and a $\caV$-transformation and $\eta\colon \id \to \iota \circ L$. If the $\caV$-transformations $L \circ \iota$ and $\iota \circ L$ are equivalences, then $\eta$ exhibits $L$ as left to adjoint to $\iota$.
\end{proposition}
\begin{proof}
We need to see that for any $X \in \caC, Y \in \caD$ the canonical morphism
$$\mathsf{Mor}_\caD(L(X),Y) \to \mathsf{Mor}_\caC(X,\iota(Y))$$
in $\caV$ is an equivalence.
To check this, it is enough to see that for any $Z \in \caV$
the induced map
$\caV(Z, \mathsf{Mor}_\caD(L(X),Y)) \to \caV(Z,\mathsf{Mor}_\caC(X,\iota(Y)))$
is an equivalence.
Since $\caC, \caD$ admit cotensors and $\iota$ preserves them, this map is canonically equivalent to the map 
$$\caD(L(X),Y^Z) \to \caC(X,\iota(Y^Z)) \simeq \caC(X,\iota(Y)^Z).$$ So we can reduce to show that
the map $\caD(L(X),Y) \to \caC(X,\iota(Y))$ is an equivalence,
which follows from \cite[Proposition 5.2.7.4]{lurie.HTT}.
\end{proof}

\subsection{Weighted limits and colimits}

\begin{definition}\label{Venr:def:colimit}
Let $F\colon\mathcal{J} \to \caC$ and $H\colon\mathcal{J}^\op \to \caC$ be $\caV$-functors, consider also $X \in \caC$ and $\psi\colon H \to F^\ast(X)$ a morphism in $\caP_\caV(\caC).$ We say that $\psi$ exhibits $X$ as the $H$-weighted colimit of $F$
if for any $Y \in \caC$ the canonical morphism
$$\mathsf{Mor}_\caC(X,Y) \to \mathsf{Mor}_{\caP_\caV(\caJ)}(F^\ast(X),  F^\ast(Y)) \to \mathsf{Mor}_{\caP_\caV(\caJ)}(H, F^\ast(Y)) $$ in $\caV$ is an equivalence.
\end{definition}

\begin{remark}Dually we can define weighted limits.
\end{remark}

\begin{remark}\label{Venr:left_adj_preserve_tensors}
It is easy to check that left $\caV$-adjoints preserve weighted colimits. 
\end{remark}

Let us distinguish two important classes of weighted colimits that allow us to get a feeling of how this notion works: tensors and conical.

\begin{remark}
For $\caJ$ the $\caV$-category with connected space of objects, and the tensor unit of $V$ as endomorphisms of such object, we have that $\caV$-functors $H\colon \caJ \to \caV$ and $F\colon \caJ^\op \to \caC$ correspond to objects $V \in \caV$ and $X \in \caC$. We call the $H$-weighted colimit of $F$ the tensor of $V$ and $X$.

The tensor of an object $V \in \caV$ with an object $X \in \caC$
is an object $V \otimes X$ and a morphism
$V \to \mathsf{Mor}_\caC(X,V \otimes X)$ such that the canonical morphism
$\mathsf{Mor}_\caC(V \otimes X,Y) \to \mathsf{Mor}_\caV(V,\mathsf{Mor}_\caC(X,Y))$
is an equivalence.
\end{remark}

\begin{definition}\label{Venr:def:presentable}
We call a $\caV$-category presentable if it admits tensors and its underlying $\infty$-category is presentable.
\end{definition}

We now proceed to describe conical colimits. 

\begin{remark}
Let $K$ be an $\infty$-category and $\caJ\coloneqq K^\caV$, where the superscript $\caV$ indicates that we see $K$ as a $\caV$-category in the canonical way.
Then $\caV$-functors $H\colon \caJ \to \caV$ and $F\colon \caJ^\op \to \caC$ correspond to functors $H'\colon K \to \caV$ and $T\colon K^\op \to \caC$. 
For $H$ corresponding to the constant functor $H'\colon K \to \caV$ with value the tensor unit, we call the $H$-weighted colimit of $F$ the conical colimit of $F$.

So a conical colimit of a functor $F\colon K_\caV^\op \to \caC$
corresponding to a functor $T\colon K^\op \to \caC$
is an object $X$ of $\caC$ and a map $T \to \delta(X)$ in 
$\caP(K)$, where $\delta\colon \caC \to \Fun(K^\op,\caC) $ is the diagonal functor, such that for any $Y \in \caC$ the canonical morphism
$$ \mathsf{Mor}_\caC(X,Y) \to \lim(\mathsf{Mor}_\caC(-,Y) \circ T^\op) $$
is an equivalence.

\end{remark}

\subsection{$\caV$-linear categories and $\caV$-categories}
This part is included in \cite[Section 2]{heine-enriched1}.

Let 
$$(\Cat_\infty^\caV)^{\mathrm{cc}} \subset \Cat_\infty^\caV$$
be the subcategory with objects the $\caV$-categories that admit
small colimits and tensors with $\caV$ and with morphisms the $\caV$-functors preserving small colimits and tensors with $\caV$.

The next result follows from \cite[Proposition 5.4]{heine-enriched1}.
\begin{theorem}\label{Venr:equiv_left_modules_enriched_cats}
There is a canonical equivalence 
$$\LMod_\caV(\Cat_\infty) \simeq (\Cat_\infty^\caV)^{\mathrm{cc}},$$
where the left hand side is the $\infty$-category of left $\caV$-modules in $ \Cat_\infty $, that sends an $\infty$-category $\mathcal{M}$ left tensored over $\caV$ to a $\caV$-enriched $\infty$-category with space of objects $\mathcal{M}^\simeq.$
\end{theorem}

\subsection{Enriched higher algebra}
\begin{notation}
Let $\caV$ be a monoidal $\infty$-category compatible with small coimits.
Let $1\otimes (-)\colon \Spc \to \caV$ denote the unique left adjoint monoidal functor, which induces a functor
$$(-)^\caV\colon \Cat_\infty^\Spc \to \Cat_\infty^\caV.$$
\end{notation}

\begin{definition}\label{Venr:def:monoids}
Let $\caC$ be a $\caV$-category that admits finite products. We denote by $\Mon^\caV_{E_\infty}(\caC)$ the $\caV$-enriched $\infty$-category of $\caV$-enriched commutative monoids in $\caC$, that we define as the full subcategory  
$$\Mon_{E_\infty}^\caV(\caC) \subset \Fun^\caV(\Fin_\ast^\caV, \caC) $$ 
spanned by the $\caV$-functors 
$F\colon \Fin_\ast^\caV \to \caC$
(corresponding to functors $\Fin_\ast \to \caC$) such that for any $n \geq 0$ the family of n morphisms $\langle n \rangle \to \langle 1 \rangle$, whose fiber over $1 \in \langle 1 \rangle$ consists precisely of one element, 
yields an equivalence in $\caC$ as below.
$$ F(\langle n \rangle) \to  \prod_{i=1}^n F(\langle 1 \rangle).$$ 
\end{definition}

\begin{remark}\label{Venr:rmk:enriched_vs_unenriched_monoids}
Let $\caV$ be a monoidal $\infty$-category compatible with small colimits.
Then the canonical equivalence $\Fun^\caV(\Fin_\ast^\caV, \caC) \simeq \Fun^\caV(\Fin_\ast, \caC)$
restricts to an equivalence
$$\Mon_{E_\infty}^\caV(\caC) \to \Mon_{E_\infty}(\caC).$$
\end{remark}

We can consider the forgetful $\caV$-functor $\Mon_{E_\infty}^\caV(\caC) \subset \Fun^\caV(\Fin_\ast^\caV, \caC) \to \caC$, where the last $\caV$-functor evaluates at
$\langle 1 \rangle \in \Fin_\ast^\caV.$

\begin{remark}\label{Venr:Mon_closed_under_co_tensors}
Let $F \in \Mon_{E_\infty}^\caV(\caC) $ and $K \in \caV$, then the cotensor
$F^K \in \Fun^\caV(\Fin_\ast^\caV, \caC) $ belongs to $\Mon_{E_\infty}^\caV(\caC) $ since 
for any $n \geq 0$ the family of n morphisms $\langle n \rangle \to \langle 1 \rangle$, whose fiber over $1 \in \langle 1 \rangle$ consists precisely of one element, 
induces the morphism 
$$ F^K(\langle n \rangle) \to  \prod_{i=1}^n F^K(\langle 1 \rangle),$$
which is the result of cotensoring with $K$ the morphism below $$ F(\langle n \rangle) \to  \prod_{i=1}^n F(\langle 1 \rangle).$$

\end{remark}

\begin{remark}\label{Venr:Mon_is_reflective}
For $\caC$ a presentable $\caV$-category, the embedding 
$$\iota\colon\Mon_{E_\infty}^\caV(\caC) \hookrightarrow \Fun^\caV(\Fin_\ast^\caV, \caC)$$ 
admits a left adjoint in the enriched sense. Indeed, by the dual of \cref{Venr:upgrade_adjunctions}, it suffices to show that $\iota$ preserves cotensors with $\caV$ and its underlying functor admits a left adjoint. The former is the content of \cref{Venr:Mon_closed_under_co_tensors} whereas the latter, since $\caC$ is presentable, follows from \cref{Venr:rmk:enriched_vs_unenriched_monoids}.
\end{remark}

\begin{definition}\label{Venr:def:group_objects}
Let $\caC$ be a $\caV$-category that admits finite products. We denote by $\Grp^\caV_{E_\infty}(\caC)$ the $\caV$-category of $\caV$-enriched group objects in $\caC$, that we define as the full subcategory 
$$\Grp^\caV_{E_\infty}(\caC) \subset \Mon^\caV_{E_\infty}(\caC)$$
spanned by the $\caV$-enriched commutative monoids $F$ in $\caC$,
whose underlying monoid is grouplike, i.e.\ the shear map
$F(1) \times F(1) \to F(1) \times F(1)$ is an equivalence. 
\end{definition}

We can consider the forgetful $\caV$-functor $$\Grp_{E_\infty}^\caV(\caC) \subset \Fun^\caV(\Fin_\ast^\caV, \caC) \to \caC,$$ where the last $\caV$-functor evaluates at
$\langle 1 \rangle \in \Fin_\ast^\caV.$

\begin{remark}\label{Venr:Grp_is_reflective}
Similarly to \cref{Venr:Mon_is_reflective}, one can check that for $\caC$ a presentable $\caV$-category, the embedding 
$$\iota\colon\Grp_{E_\infty}^\caV(\caC) \hookrightarrow \Fun^\caV(\Fin_\ast^\caV, \caC)$$ admits a left adjoint in the enriched sense.
\end{remark}

\begin{definition}\label{Venr:def:V_pre_additive}
We say that a $\caV$-category $\caC$ is $\caV$-(pre)additive if it has finite products and finite coproducts, and its underlying $\infty$-category is (pre)additive.
\end{definition}

\begin{remark}\label{Venr:V_pre_add_when_underlying_is_V_pre_add}
Let $\caC$ be a $\caV$-category that admits cotensors, then
$\caC$ is $\caV$-(pre)additive if and only if the underlying $\infty$-category of $\caC$ is (pre)additive.
\end{remark}

\begin{theorem}\label{Venr:lifting_to_monoids}
Let $\caC$ and $\caD$ be $\caV$-enriched $\infty$-categories such that
$\caC$ is $\caV$-preadditive and $\caD$ has finite products. The canonical $\caV$-functor 
$$\theta\colon \Fun^{\caV,{\prod}}(\caC, \Mon^\caV_{E_\infty}(\caD)) \to \Fun^{\caV,{\prod}}(\caC, \caD),$$ is an equivalence, where the superscript $\prod$ indicates that we consider product preserving functors.
\end{theorem}

\begin{proof}
We first show that it suffices to check that $\theta$ induces a bijection between equivalence classes of objects in the underlying $\infty$-categories.
Indeed, by Yoneda $\theta$ is an equivalence if for any $\caV$-category $\caA$ the $\caV$-functor $\Fun^\caV(\caA, \theta)$ induces a bijection between equivalence classes of objects in the underlying $\infty$-categories.
The $\caV$-functor $\Fun^\caV(\caA, \theta)$ is canonically equivalent to the $\caV$-functor
$$\Fun^{\caV,{\prod}}(\caC, \Mon^\caV_{E_\infty}(\Fun^\caV(\caA, \caD))) \to  \Fun^{\caV,{\prod}}(\caC, \Fun^\caV(\caA, \caD)),$$
which is $\theta$ when $\caD$ is replaced by $\Fun^\caV(\caA, \caD).$

Since $\caC$ is $\caV$-preadditive, the forgetful functor
$\Mon_{E_\infty}(\caC) \to \caC$
is an equivalence.
Then the $\caV$-functor 
$$\Fun^{\caV,{\prod}}(\caC, \caD) \to \Fun^{\caV,{\prod}}(\Mon^\caV_{E_\infty}(\caC), \Mon^\caV_{E_\infty}(\caD)) \simeq \Fun^{\caV,{\prod}}(\caC, \Mon^\caV_{E_\infty}(\caD))$$
provides liftings of functors $\caC\to\caD$. It remains to show that these lifts are unique up to equivalence. Suppose we have two lifts $\alpha,\beta\colon\caC\to\Mon^\caV_{E_\infty}(\caD)$ of a $\caV$-functor $F\colon\caC\to\caD$.
Then the induced $\caV$-functors $\Mon^\caV_{E_\infty}(\alpha), \Mon^\caV_{E_\infty}(\beta)$ lift the $\caV$-functor $\Mon^\caV_{E_\infty}(F)$ along the forgetful functor $\Mon^\caV_{E_\infty}(\Mon^\caV_{E_\infty}(\caD)) \to \Mon^\caV_{E_\infty}(\caD),$ which is an equivalence since $\Mon^\caV_{E_\infty}(\caD)$ is preadditive.
Hence the $\caV$-functors $\Mon^\caV_{E_\infty}(\alpha), \Mon^\caV_{E_\infty}(\beta)$
are equivalent.
Since both forgetful functors $\Mon^\caV_{E_\infty}(\caC) \to \caC$ and $\Mon^\caV_{E_\infty}(\Mon^\caV_{E_\infty}(\caD)) \to \Mon^\caV_{E_\infty}(\caD)$
are equivalences, also $\alpha$ and $\beta$ are equivalent. 
\end{proof}

\begin{theorem}\label{Venr:lifting_to_grp}
For any $\caV$-category $\caD$ and any additive $\caV$-category $\caC$
the forgetful functor
$\Grp^\caV_{E_\infty}(\caD) \to \caD$ induces an equivalence
$$\Fun_\caV^{\prod}(\caC, \Grp^\caV_{E_\infty}(\caD)) \simeq \Fun_\caV^{\prod}(\caC, \caD).$$
\end{theorem}

\begin{proof}
The proof is analogous to the proof of \cref{Venr:lifting_to_monoids},
where we use that the $\caV$-functor $\Grp^\caV_{E_\infty}(\caC) \to \caC$ is an equivalence since $\caC$ is $\caV$-additive.
\end{proof}

\begin{proposition}\label{Venr:forgetful_Mon_Grp_are_equiv}
Let $\caD$ be a $\caV$-category that has colimits and admits tensors with $\caV$. Then 
\begin{enumerate}
\item If $\caD$ is $\caV$-preadditve, the $\caV$-functor $\Mon^\caV_{E_\infty}(\caD) \to \caD$ is an equivalence.
\item If $\caD$ is $\caV$-additive, the $\caV$-functor $\Grp^\caV_{E_\infty}(\caD) \to \caD$ is an equivalence.
\end{enumerate}
\end{proposition}

\begin{proof}
We only prove item 1, item 2 is similar. By \cref{Venr:rmk:enriched_vs_unenriched_monoids} the forgetful functor $\Mon^\caV_{E_\infty}(\caD) \to \caD$ induces on underlying $\infty$-categories the forgetful functor $\Mon_{E_\infty}(\caD) \to \caD$,
which is an equivalence by \cite[Proposition 2.3]{gepner-groth-nikolaus}. Hence by \cref{Venr:F_equiv_iff_left_adj_and_underlying_is_equiv} it will be enough to show that the $\caV$-functor
$\Mon^\caV_{E_\infty}(\caD) \to \caD$ admits a left adjoint. 
This follows from \cref{Venr:Mon_is_reflective} and \cref{Venr:ev_admits_adjoint}.
\end{proof}

\subsection{Other results we use}

\begin{remark}\label{gghhjjn}
For every small $\infty$-category $\caC$ and monoidal $\infty$-category $\caV$ there is a canonical equivalence of cartesian fibrations over $\Spc$

\[
\begin{tikzcd}[column sep=tiny]
\mathsf{PreCat}_\infty^{\Fun(\caC, \caV)}\ar[rr,"\simeq"]\ar[dr]       &           &(\mathsf{PreCat}^\caV_\infty)^\caC\ar[dl]\\
\         &\Spc         &
\end{tikzcd}
\]

natural in both $\caC$, with respect to restriction, and in $\caV$, with respect to pushforward \textemdash where $(-)^\caC$ denotes the cotensoring of cartesian fibrations over $\Spc$, i.e.\ the fiberwise functor category.
\end{remark}

Then we can consider a canonical functor $\theta\colon\mathsf{PreCat}_\infty^{\Fun(\caC, \caV)}\to \Fun(\caC, \mathsf{PreCat}^\caV_\infty)$ given by the composition
\begin{equation}\label{theta_in_Vapp}
    \mathsf{PreCat}_\infty^{\Fun(\caC, \caV)} \simeq (\mathsf{PreCat}^\caV_\infty)^\caC \to \Fun(\caC, \mathsf{PreCat}^\caV_\infty). 
\end{equation}

\begin{lemma}\label{fghjkfb}
Let $\mathcal{J}$ be an $\infty$-category and $\phi\colon\caC \to \caD$ a cartesian fibration. If $\caD$ admits all limits indexed by $\mathcal{J}$, then there is an adjunction

\[
\begin{tikzcd}[column sep = huge]
	{\caC^{\mathcal{J}}} & {\Fun(\mathcal{J}, \caC),}
	\arrow[""{name=0, anchor=east, inner sep=0}, "{L}", bend left=30, from=1-1, to=1-2]
	\arrow[""{name=1, anchor=center, inner sep=0}, "{R}", bend left=30, from=1-2, to=1-1]
	\arrow["\dashv"{anchor=center, rotate=-90}, draw=none, from=0, to=1]
\end{tikzcd}
\]

where $\caC^\mathcal{J}$ is the cotensor of the cartesian fibration $\phi$ with $\mathcal{J}$, such that $L$ is the projection and the counit $L(R(\beta)) \to \beta $
at a given functor $\beta\colon \mathcal{J} \to \caC$ is the cartesian lift of the natural transformation $f\colon \delta(\lim(\phi \circ \beta)) \to \phi \circ \beta$\textemdash we use $\delta(\lim(\alpha))$ for the constant functor on $\lim(\phi\circ\beta)$.
\end{lemma}

\begin{proof}
For every $(Y,\gamma) \in \caC^{\mathcal{J}}$ the canonical map below is an equivalence: 
$$ \map_{ \caC^{\mathcal{J}}}(  (Y,\gamma), (\lim(\alpha),f^\ast(\beta))) \simeq $$
$$\map_{  \caD }( Y, \lim(\alpha)) \times_{ \map_{  \Fun(\mathcal{J} , \caD) }( \delta(Y), \delta(\lim(\alpha))) } \map_{ \Fun(\mathcal{J} , \caC) }(\gamma, f^\ast(\beta)) $$
$$ \simeq \map_{  \caD }( Y, \lim(\alpha)) \times_{ \map_{  \Fun(\mathcal{J} , \caD) }( \delta(Y), \alpha) } \map_{ \Fun(\mathcal{J} , \caC) }(\gamma, \beta) $$ 

Since the map $\map_{  \caD }( Y, \lim(\alpha)) \to \map_{  \Fun(\mathcal{J} , \caD)}( \delta(Y), \alpha)$ is an equivalence, we obtain an equivalence 

$$\map_{ \caC^{\mathcal{J}}}(  (Y,\gamma), (\lim(\alpha),f^\ast(\beta)))\simeq \map_{ \Fun(\mathcal{J} , \caC)} ( \gamma, \beta ).$$
\end{proof}

We apply \cref{fghjkfb} to the cartesian fibration $\psi\colon \PreCat_\infty^\mathcal{V} \to \Spc$ that takes the corresponding space of objects; see \cite[Remark 3.17]{heine-enriched1}.

\begin{corollary}\label{adj_theta_rho}
Let $\caV$ be a monoidal $\infty$-category and $\caC$ an $\infty$-category. Then the functor $\theta$ of \cref{theta_in_Vapp} is part of an adjunction

\[
\begin{tikzpicture}
\node[](A) {$\mathsf{PreCat}_\infty^{\Fun(\caC, \caV)}$};
\node[right of=A,xshift=4cm](B){$\Fun(\caC,\mathsf{PreCat}^\caV_\infty)$};
\draw[->] ($(A.east)+(0,.25cm)$) to [bend left=25] node[above]{$\theta$} ($(B.west)+(0,.25cm)$);
\draw[->] ($(B.west)+(0,-.25cm)$) to [bend left=25] node[below]{$\rho$} ($(A.east)-(0,.25cm)$);
\node[] at ($(A.east)!0.5!(B.west)$) {$\bot$};
\end{tikzpicture}
\]
whose counit component $(\theta\circ\rho)(\beta) \to \beta $
at $\beta\colon \caC \to \mathsf{PreCat}^\caV_\infty$ lying over a functor $\alpha\colon\caC \to \Spc$ is the objectwise $\psi$-cartesian lift of the natural transformation $f\colon\delta(\lim(\alpha)) \to \alpha. $  
\end{corollary}

In the next two lemmas we describe the behavior of cartesian morphisms wtih the functors $\theta$ and $\rho$.

\begin{lemma}\label{cartesian_iff_objwise_cartesian_under_theta}
A morphism of $\mathsf{PreCat}_\infty^{\Fun(\caC, \caV)}$ is cartesian over $\Spc$
if and only if its image under the functor $\theta$ of \cref{theta_in_Vapp} is objectwise cartesian over $\Spc$.
\end{lemma}

\begin{proof}
It follows from \cref{gghhjjn} together with the fact that the functor $\theta$ is simply the projection.
\end{proof}

\begin{lemma}\label{rho_respects_cartesian}
If a morphism of $\Fun(\caC, \mathsf{PreCat}^\caV_\infty) $ is objectwise cartesian over $\Spc$, then its image under $\rho$ is cartesian over $\Spc$.
\end{lemma}

\begin{proof}
This follows from \cref{cartesian_iff_objwise_cartesian_under_theta} and the fact that the counit of the adjunction
$$\theta\colon \mathsf{PreCat}_\infty^{\Fun(\caC, \caV)} \rightleftarrows \Fun(\caC, \mathsf{PreCat}^\caV_\infty) \colon\rho $$ is objectwise cartesian over $\Spc.$
\end{proof}

\begin{lemma}\label{unit_theta_rho}
The unit of the adjunction $$\theta\colon \mathsf{PreCat}_\infty^{\Fun(\caC, \caV)} \rightleftarrows \Fun(\caC, \mathsf{PreCat}^\caV_\infty) \colon\rho $$ is cartesian over $\Spc.$
\end{lemma}

\begin{proof}
This follows from \cref{adj_theta_rho}, \cref{cartesian_iff_objwise_cartesian_under_theta}, and the triangle identities.
\end{proof}

As we see in what follows, the adjunction $\theta\dashv\rho$ induces a similar adjunction when considering the corresponding subcategories $\Cat_\infty^\mathcal{V}\subset \PreCat^\mathcal{V}_\infty$ spanned by the $\mathcal{V}$-categories. 

\begin{lemma}\label{Venr:lem:restriction_of_G}
Let $\caV$ be a monoidal $\infty$-category and $\caC$ an $\infty$-category.
The right adjoint $\rho\colon \Fun(\caC, \mathsf{PreCat}^\caV_\infty) \to \mathsf{PreCat}_\infty^{\Fun(\caC, \caV)}$ of \cref{adj_theta_rho}
restricts to a functor $$G\colon \Fun(\caC, \mathrm{Cat}^\caV_\infty) \to \Cat_\infty^{\Fun(\caC, \caV)}$$ that admits a left adjoint $F$, which is natural in $\caC$ with respect to restriction and natural in $\caV$ with respect to pushforward.
\end{lemma}

\begin{proof}
The functor 
$$\theta\colon \mathsf{PreCat}_\infty^{\Fun(\caC, \caV)} \simeq (\mathsf{PreCat}^\caV_\infty)^\caC \to \Fun(\caC, \mathsf{PreCat}^\caV_\infty) $$
preserves local equivalences, as for every $X \in \caC$ there is a commutative triangle
\[
\begin{tikzcd}[column sep=tiny]
\mathsf{PreCat}_\infty^{\Fun(\caC, \caV)}  \ar[rd, "(\ev_X)_\ast"']\ar[rr, "\theta"]         &     &\Fun(\caC,  \mathsf{PreCat}^\caV_\infty) \ar[ld, "\ev_X"]\\ 
&   \mathsf{PreCat}^\caV_\infty. 
\end{tikzcd}
\]

Then $\theta$ induces a functor $$F\colon\Cat_\infty^{\Fun(\caC, \caV)} \to \Fun(\caC, \mathrm{Cat}^\caV_\infty) $$
natural in $\caC$ with respect to restriction and in $\caV$ with respect to pushforward.

As $\theta$ preserves local equivalences, its right adjoint too admits a restriction to a functor
$$G\colon\Fun(\caC, \mathrm{Cat}^\caV_\infty) \to \Cat_\infty^{\Fun(\caC, \caV)},$$ right adjoint to $F.$
\end{proof}

\begin{remark}\label{equivalence from adjunction}
If $\caV$ admits finite products, $F$ preserves finite products. So if $\caV$ is a cartesian symmetric monoidal $\infty$-category, for every $X \in \Cat_\infty^{\Fun(\caC, \caV)}, \ Y \in \Fun(\caC, \Cat^\caV_\infty)$ there is a canonical equivalence
$$ G(\hom_{  \Fun(\caC, \Cat^\caV_\infty)}(F(X),Y)) \simeq \Fun^{\Fun(\caC, \caV)}(X, G(Y)). $$
\end{remark}

\begin{remark}\label{remmm}
Let $\caV$ be a monoidal $\infty$-category and $\caC$ an $\infty$-category. By adjointness if $\caV$ admits small limits, there is a commutative triangle
\[
\begin{tikzcd}[column sep=small]
\Fun(\caC, \Cat^\caV_\infty) \ar[rd, "\lim"'] \ar[rr, "G"]        &       &  \Cat_\infty^{\Fun(\caC, \caV)} \ar[ld, "\lim_\ast"] \\ 
&  \Cat^\caV_\infty.
\end{tikzcd}
\]
\end{remark}

\begin{lemma}\label{unit_restriction_theta_rho}
The unit of the adjunction 
$$F\colon \mathrm{Cat}_\infty^{\Fun(\caC, \caV)} \rightleftarrows \Fun(\caC, \mathrm{Cat}^\caV_\infty) \colon G $$ is cartesian over $\Spc$.
\end{lemma}

\begin{proof}
This follows from \cref{rho_respects_cartesian,unit_theta_rho} using that every local equivalence of $\mathsf{PreCat}^\caV_\infty$ is 
cartesian over $\Spc.$
\end{proof}

\begin{lemma}\label{Venr:Fun_enriched_is_Fun_to_enriched}
If $\caC$ admits an initial object, then the unit of the adjunction $$F\colon \mathrm{Cat}_\infty^{\Fun(\caC, \caV)} \rightleftarrows \Fun(\caC, \mathrm{Cat}^\caV_\infty) \colon G $$ is an equivalence.
\end{lemma}

\begin{proof}
If $\caC$ admits an initial object, evaluation at the initial object of $\caC$ computes the limit over $\caC$. So $F$ and $G$ canonically commute with evaluation at the initial object (\cref{remmm}). Thus the unit is especially essentially surjective and so by \cref{unit_restriction_theta_rho} an equivalence.
\end{proof}

The following results refer to the functors $F$ and $G$ of \cref{Venr:Fun_enriched_is_Fun_to_enriched}.

\begin{lemma}
For every $\beta\colon \caC \to \Cat^\caV_\infty$ lying over $\alpha\colon \caC \to \Spc $ the counit $F G(\beta) \to \beta$ is objectwise cartesian over $\Spc$.
\end{lemma}

\begin{corollary}
For every $\beta \colon\caC \to \Cat^\caV_\infty$ lying over $\alpha\colon \caC \to \Spc $ the counit $F G(\beta) \to \beta$
is an equivalence if and only if for every $Z \in \caC$
the canonical map $ \lim(\alpha) \to \alpha(Z)$ is essentially surjective.
\end{corollary}

\section{Real spine inclusions}\label{Appx:real_spine_inc}

In this appendix we show that the definition of real spine inclusions, \cref{def:real_spine_inc}, is the natural one to consider when we draw inspiration from the usual real spine inclusions. More precisely, we will show the following result.

\begin{proposition}[\cref{ex:real_spine_inc}]\label{prop:spine_inc_usual}
For every $n \geq 2$ the real spine inclusion $\Lambda^n \subset \Delta^n$ is the unique lift with target $\Delta^n$ of the canonical inclusion $$\Delta^1 \sqcup_{\Delta^0 } \Delta^1 \sqcup_{ \Delta^0 } \ldots \sqcup_{\Delta^0 } \Delta^1 \subset \Delta^n$$ of simplicial sets along the forgetful functor $\sSet^{hC_2} \to \sSet$. 

\end{proposition}

We defer the proof of \cref{prop:spine_inc_usual} to the end of the appendix, and recall now the definitions involved, as well as present the results the proof relies on. 

Unsurprisingly, we first need to bring to mind the main components for, and the definition itself of, real spine inclusions given in \cref{subsec:real_Segal_objects}. We know that for every $n \geq 0 $ the category $[n]$ has a unique strict duality that gives its nerve $\Delta^n \in \sSet$ the structure of a real simplicial set. With this structure in $\Delta^n$, we write below \cref{def:real_spine_inc}.

\begin{definition}[\cref{def:real_spine_inc}]\label{def:spine_inc_A2}
We define the spine inclusions $\Lambda^n \to \Delta^n$ for $n \geq 1$. For $n$ even, we set $\Lambda^n$ to be the coproduct of $[0]$ with $\frac{n}{2}$ copies of $\widetilde{[1] \coprod [1]}$, constructed by using alternatively the two maps $[0]\to[1]$. We have
$$\Lambda^n \coloneqq [0] \coprod_{\widetilde{[0] \coprod [0]}} \widetilde{[1] \coprod [1]} \coprod_{\widetilde{[0] \coprod [0]}}\dots \coprod_{\widetilde{[0] \coprod [0]}} \widetilde{[1] \coprod [1]}. $$
 
Given $\Lambda^n$, we define the $n-th$ spine inclusion $\Lambda^n\to\Delta^n$ by the universal property of the coproduct, where the compatible maps considered are the ones induced by considering for every $1 \leq i \leq \frac{n}{2}$, the maps $ [1] \to [n]$ in $\Delta$ sending $0 $ to $i-1$ and $1$ to $i$, that yield maps $\widetilde{[1] \coprod [1]} \to [n]$ in $\rsSpc$; and the map $[0] \to [n]$ in $\Delta^{hC_2}$ sending $0$ to $\frac{n}{2}$.

For $n \geq 3 $ odd, we set
$$\Lambda^n \coloneqq [1] \coprod_{\widetilde{[0] \coprod [0]}} \widetilde{[1] \coprod [1]} \coprod_{\widetilde{[0] \coprod [0]}}\dots \coprod_{\widetilde{[0] \coprod [0]}} \widetilde{[1] \coprod [1]}.$$

In this case, the spine inclusion $\Lambda^n\to \Delta^n$ is induced by considering, for every $1 \leq i \leq \frac{n-1}{2}$, the maps $ [1] \to [n]$ sending $0 $ to $i-1$ and $1$ to $i$ in $\Delta$; and the map $[1] \to [n]$ in $\Delta^{hC_2}$, sending $0$ to $ \frac{n-1}{2}$ and $1$ to $\frac{n-1}{2}+1$ in $\Delta^{hC_2}$.
\end{definition}

The following three simple results will make for a short proof of \cref{prop:spine_inc_usual}.

\begin{lemma}\label{lemA1}
Let $F,G\colon J \to \Spc$ be functors and
$\psi\colon F \to G$ a natural transformation, whose components are all embeddings. Then the induced commutative square
\[
\begin{tikzcd}
\lim_{j \in J} F \ar[r]\ar[d]  &\lim_{j \in J} G \ar[d]\\
\prod_{j \in J}F(j)\ar[r]&\prod_{j \in J}G(j)
\end{tikzcd}
\]
is a pullback square.
\end{lemma}

\begin{proof}
Let $\theta\colon X \to Y$ be the map of left fibrations over $J$ classifying the natural transformation $\psi\colon F \to G$. The commutative square in the statement
is equivalent to the following one
\[
\begin{tikzcd}
\Fun_J(J,X) \ar[r]\ar[d]  &\Fun_J(J,Y) \ar[d]\\
\prod_{j \in J}X_j\ar[r]&\prod_{j \in J}Y_j.
\end{tikzcd}
\]
Since all components of $\psi$ are embeddings, so is the map $\theta$ and therefore both horizontal maps in the last square are embeddings as well. Moreover every section
$\gamma$ of $ Y \to J$ factors through $X$
if and only if for every $j \in J$ the image
$\gamma(j) \in Y_j$ belongs to $X_j$.
Hence the last square is a pullback square.
\end{proof}

\begin{lemma}\label{lemA2}
Let $J$ be an $\infty$-category and $F\colon J \to \Cat_\infty$ a functor.
Any morphism $\phi\colon A \to B$ in $\lim(F)$ whose image in $\prod\coloneqq\prod_{j \in J} F(j)$ is a monomorphism, is cartesian with respect to the functor $\rho\colon\lim F \to \prod_{j \in J} F(j)$.
\end{lemma}

\begin{proof}
In other words, we want to show that the diagram below is a pullback square.
\[
\begin{tikzcd}
\lim_{j \in J} F(j)(C_j,A_j) \ar[r]\ar[d]    &\lim_{j \in J} F(j)(C_j,B_j)\ar[d]\\
\prod_{j \in J}F(j)(C_j,A_j)\ar[r]&\prod_{j \in J}F(j)(C_j,B_j)
\end{tikzcd}
\]
So the claim follows from \cref{lemA1}.
\end{proof}

\begin{corollary}\label{cor:forgetful_cartesian}
Let $\alpha\colon A\to B$ be a map of real simplicial spaces 
whose underlying map of simplicial spaces is an embedding in each simplicial degree.
Then $\alpha$ is cartesian with respect to the functor $\phi\colon \s\Spc^{hC_2} \to \s\Spc.$
\end{corollary}

\begin{proof}
The claim follows from \cref{lemA2} applied to $J=BC_2$ and $F\colon BC_2 \to \Cat_\infty$ the functor encoding the $C_2$-action on $\s\Spc$.
\end{proof}

We conclude the appendix as announced.

\begin{proof}[Proof of \cref{prop:spine_inc_usual}]
From the description in \cref{def:spine_inc_A2}, it can easily be checked that the inclusion $\Lambda^n \subset \Delta^n$ of real simplicial sets lifts the canonical inclusion 
$$\Delta^1 \sqcup_{\Delta^0 } \Delta^1 \sqcup_{ \Delta^0 } \ldots \sqcup_{\Delta^0 } \Delta^1 \subset \Delta^n.$$  
Since by \cref{cor:forgetful_cartesian} any lift of such map is cartesian with respect to the forgetful functor $\s\Spc^{hC_2} \to \s\Spc$, we conclude that this lift is unique. 
\end{proof}

\bibliographystyle{alpha}
\bibliography{add}

\end{document}